\documentclass{amsart}
\usepackage{amssymb,amsmath,graphicx,color}
\usepackage[all]{xy}
\usepackage{enumerate}
\usepackage{mathrsfs}

\theoremstyle{plain}
\newtheorem{theorem}{Theorem}[section]
\newtheorem{lemma}[theorem]{Lemma}
\newtheorem{corollary}[theorem]{Corollary}
\newtheorem{proposition}[theorem]{Proposition}

\theoremstyle{definition}
\newtheorem{definition}[theorem]{Definition}
\newtheorem{example}[theorem]{Example}

\theoremstyle{remark}
\newtheorem{remark}[theorem]{Remark}

\newtheorem{claim}{Claim}

\newcommand\Aut{\operatorname{Aut}}

\newcommand\Out{\operatorname{Out}}
\newcommand\Hom{\operatorname{Hom}}
\newcommand\End{\operatorname{End}}
\newcommand\Mat{\operatorname{Mat}}
\newcommand\GL{\operatorname{GL}}
\newcommand\IA{\operatorname{IA}}
\newcommand\Ind{\operatorname{Ind}}
\newcommand\Sym{\operatorname{Sym}}
\newcommand\sym{\operatorname{sym}}
\newcommand\alt{\operatorname{alt}}

\newcommand\id{\operatorname{id}}

\newcommand\opegr{\operatorname{gr}}
\newcommand\Span{\operatorname{Span}}
\newcommand\Rad{\operatorname{Rad}}
\newcommand\Soc{\operatorname{Soc}}
\newcommand\im{\operatorname{im}}
\newcommand\Der{\operatorname{Der}}

\newcommand\gr{\mathrm{gr}}

\newcommand\op{\mathrm{op}}
\newcommand\ab{\mathrm{ab}}
\newcommand\sgn{\mathrm{sgn}}

\newcommand\Z{\mathbb{Z}}
\newcommand\N{\mathbb{N}}

\newcommand\K{\Bbbk}
\newcommand\R{\mathbb{R}}
\newcommand\Q{\mathbb{Q}}
\newcommand\repS{\mathbb{S}} 

\newcommand\A{\mathbf{A}} 
\newcommand\prA{\wti{\mathbf{A}^{L}}}
\newcommand\F{\mathbf{F}} 
\newcommand\FAb{\mathbf{FAb}} 

\newcommand\fVect{\mathbf{fVect}}
\newcommand\gVect{\mathbf{gVect}}

\newcommand\catC{\mathcal{C}} 

\newcommand\jfA{\mathcal{A}}
\newcommand\jfE{\mathcal{E}}
\newcommand\Lie{\mathcal{L}}
\newcommand\g{\mathfrak{g}}

\newcommand\gpS{\mathfrak{S}} 

\newcommand\ti{\tilde}
\newcommand\wti{\widetilde}
\def\edge{{-}\kern-.2em{-}\kern-.2em{-}\kern-.2em{-}\kern-.2em{-}}

\newcommand\centre[1]{\begin{array}{c} #1 \end{array}}

\title[Actions of automorphism groups of free groups on Jacobi diagrams.II]{Actions of automorphism groups of free groups on spaces of Jacobi diagrams. II}
\author{Mai Katada}
\address{Department of Mathematics, Kyoto University, Kyoto 606-8502, Japan}
\email{katada.mai.36s@st.kyoto-u.ac.jp}
\date{June 12, 2021 (First version: May 19, 2021)}

\begin{document}

\begin{abstract}
  The automorphism group $\Aut(F_n)$ of the free group $F_n$ acts on a space $A_d(n)$ of Jacobi diagrams of degree $d$ on $n$ oriented arcs.
  We study the $\Aut(F_n)$-module structure of $A_d(n)$ by using two actions on the associated graded vector space of $A_d(n)$:
  an action of the general linear group $\GL(n,\Z)$
  and an action of the graded Lie algebra $\gr(\IA(n))$ of the IA-automorphism group $\IA(n)$ of $F_n$ associated with its lower central series.
  We extend the action of $\gr(\IA(n))$ to an action of the associated graded Lie algebra of the Andreadakis filtration of the endomorphism monoid of $F_n$. By using this action, we study the $\Aut(F_n)$-module structure of $A_d(n)$. We obtain an indecomposable decomposition of $A_d(n)$ as $\Aut(F_n)$-modules for $n\geq 2d$. Moreover, we obtain the radical filtration of $A_d(n)$ for $n\geq 2d$ and the socle of $A_3(n)$.
\end{abstract}

\maketitle
\setcounter{tocdepth}{1}
\tableofcontents

\section{Introduction} \label{sec:intro}
 \emph{Jacobi diagrams} are uni-trivalent graphs, which graphically encode the algebraic structures of Lie algebras and their representations.
 Jacobi diagrams were introduced for the \emph{Kontsevich integral}, which is a universal finite type link invariant and unifies all quantum link invariants \cite{BN_v, Kontsevich, Kassel, Ohtsuki}.
 The associated graded vector space of finite type link invariants is isomorphic to the space of \emph{weight systems}, which is the dual to the space of Jacobi diagrams.

 Let $\K$ be a field of characteristic $0$. We study the $\K$-vector space $A(n)$ of Jacobi diagrams on $n$-component oriented arcs, which is the target space of the Kontsevich integral for \emph{string links} \cite{HL,BN_sl} or \emph{bottom tangles} \cite{Habiro_b}.
 We consider the degree $d$ part $A_d(n)$ of $A(n)$, where the degree of a Jacobi diagram is determined by half the number of its vertices.

 We consider a filtration for $A_d(n)$ defined by the number of trivalent vertices.
 The associated graded vector space of $A_d(n)$ is identified via the PBW (Poincar\'{e}--Birkhoff--Witt) map \cite{BN_sl} with a graded vector space $B_d(n)$ of \emph{open Jacobi diagrams} of degree $d$ that are colored by elements of an $n$-dimensional $\K$-vector space.

 Habiro and Massuyeau \cite{HM_k} extended the Kontsevich integral to a functor $Z$ from the category of bottom tangles in handlebodies to the category $\A$ of \emph{Jacobi diagrams in handlebodies}.
 By using the functor $Z$, we obtain an action on $A_d(n)$ of the \emph{automorphism group} $\Aut(F_n)$ of the free group $F_n$ of rank $n$, which is induced from an action of handlebody groups on bottom tangles.
 In fact, the $\Aut(F_n)$-action on $A_d(n)$ induces an action on $A_d(n)$ of the \emph{outer automorphism group} $\Out(F_n)$ of $F_n$.

 In a previous paper \cite{Mai1}, we observed that the $\Aut(F_n)$-action on $A_d(n)$ induces two actions on $B_d(n)$:
 an action of the \emph{general linear group} $\GL(n;\Z)$
 and an action of the graded Lie algebra $\gr(\IA(n))$ of the \emph{IA-automorphism group} $\IA(n)$ of $F_n$ associated with the lower central series.
 We used these two actions on $B_d(n)$ to study the $\Aut(F_n)$-module structure of $A_d(n)$ for $d=2$.
 However, it is rather difficult to compute the $\gr(\IA(n))$-action on $B_d(n)$ directly for general $d$.

 The aim of the present paper is to study the $\Aut(F_n)$-module structure of $A_d(n)$ for general $d$ and especially $d=3$ in detail.
 We consider the \emph{Andreadakis filtration} $\jfE_{\ast}(n)$ of the endomorphism monoid $\End(F_n)$ of $F_n$.
 We extend the action of the graded Lie algebra $\gr(\IA(n))$ to
 an action of the associated graded Lie algebra $\gr(\jfE_{\ast}(n))$ of the Andreadakis filtration.
 On the other hand, we construct a graphical version of the $\gr(\jfE_{\ast}(n))$-action on $B_d(n)$.
 By using this graphical action, we study the $\Aut(F_n)$-module structure of $A_d(n)$.
 We obtain an indecomposable decomposition of $A_d(n)$ as $\Aut(F_n)$-modules for $n\geq 2d$. Moreover, we obtain the radical filtration of $A_d(n)$ for $n\geq 2d$ and the socle of $A_3(n)$.

 \subsection{Andreadakis filtration of $\End(F_n)$}
  Let $\Gamma_r:=\Gamma_r(F_n)$ denote the $r$-th term of the lower central series of the free group $F_n$.
  Let $\Lie_r(n) :=\Gamma_r/\Gamma_{r+1}$ for $r\geq 1$ and set $H:=\Lie_1(n)$.
  Note that $\Lie_r(n)$ is the degree $r$ part of the free Lie algebra $\Lie_{\ast}(n)$ on $H$.

  Let $\IA(n)$ denote the IA-automorphism group of $F_n$, which is the kernel of the canonical homomorphism $\Aut(F_n)\rightarrow \GL(n;\Z)$.

  The \emph{Andreadakis filtration} $\jfA_{\ast}(n)$ of $\Aut(F_n)$ \cite{Andreadakis, Satoh_survey}
  $$\Aut(F_n)=\jfA_0(n)\supset \jfA_{1}(n)=\IA(n)\supset \jfA_2(n)\supset \cdots$$
  is defined by
  $$\jfA_{r}(n)=\ker(\Aut(F_n)\rightarrow \Aut(F_n/\Gamma_{r+1})).$$
  For $r\geq 1$, we have an injective homomorphism
  $$\tau_r: \opegr^r(\jfA_{\ast}(n))\hookrightarrow \Hom(H,\Lie_{r+1}(n)),$$
  which is called the \emph{Johnson homomorphism}.
  By Andreadakis \cite{Andreadakis} and Kawazumi \cite{Kawazumi},
  we have $\opegr^1(\IA(n))\cong \opegr^1(\jfA_{\ast}(n))\cong \Hom(H,\Lie_{2}(n))$.

  We construct the Andreadakis filtration $\jfE_{\ast}(n)$ of $\End(F_n)$ in a similar way by
  $$\jfE_{r}(n)=\ker(\End(F_n)\rightarrow \End(F_n/\Gamma_{r+1})).$$
  We define an equivalence relation on the monoid $\jfE_{r}(n)$ and consider the quotient group $\opegr^r(\jfE_{\ast}(n))$, which includes $\opegr^r(\jfA_{\ast}(n))$ (see Section \ref{ss32}).
  We also construct the Johnson homomorphism
  $$
   \ti{\tau}_r:\opegr^r(\jfE_{\ast}(n))\xrightarrow{\cong}\Hom(H,\Lie_{r+1}(n))
  $$
  of $\End(F_n)$,
  which turns out to be an abelian group isomorphism (see Proposition \ref{p712}).

  The target group $\Hom(H,\Lie_{r+1}(n))\cong H^{\ast}\otimes \Lie_{r+1}(n)$ of the Johnson homomorphism is identified with the degree $r$ part $\Der_r(\Lie_{\ast}(n))$ of the derivation Lie algebra $\Der(\Lie_{\ast}(n))$ of the free Lie algebra $\Lie_{\ast}(n)$ and with the tree module $T_r(n)$, which we define in Section \ref{ss315}.
  From the above, we have abelian group isomorphisms
  $$
    \opegr^r(\jfE_{\ast}(n))\cong H^{\ast}\otimes \Lie_{r+1}(n)\cong \Der_r(\Lie_{\ast}(n)) \cong T_r(n).
  $$
  Thus, we have
  $$\opegr^1(\IA(n))\cong \opegr^1(\jfE_{\ast}(n))\cong H^{\ast}\otimes \Lie_{2}(n)\cong \Der_1(\Lie_{\ast}(n))\cong T_1(n).$$

  Moreover, we have isomorphisms of graded Lie algebras
  \begin{equation}\label{dereq}
    \gr(\jfE_{\ast}(n))=\bigoplus_{r\geq 1}\opegr^r(\jfE_{\ast}(n))\cong \Der(\Lie_{\ast}(n))\cong \bigoplus_{r\geq 1}T_r(n)
  \end{equation}
  (see Section \ref{ss34}).
  In what follows, we identify these three graded Lie algebras.

 \subsection{Actions of the derivation Lie algebra on $B_d(n)$}
  Let $A_d(n)$ be the $\K$-vector space spanned by Jacobi diagrams of degree $d$ on $n$ oriented arcs.
  We consider a filtration for $A_d(n)$
  $$A_d(n)=A_{d,0}(n)\supset A_{d,1}(n)\supset A_{d,2}(n)\supset \cdots,$$
  where $A_{d,k}(n)$ is the subspace of $A_d(n)$ spanned by Jacobi diagrams with at least $k$ trivalent vertices.
  The $\Aut(F_n)$-action on $A_d(n)$ that we considered in \cite{Mai1} can be naturally extended to an action of the endomorphism monoid $\End(F_n)$ of $F_n$ on $A_d(n)$. (See Section \ref{s4}.)

  Let $V_n$ be an $n$-dimensional $\K$-vector space, which will be identified with the first cohomology of a handlebody of genus $n$.
  The associated graded vector space of $A_d(n)$ is isomorphic via the PBW map \cite{BN_sl} to a graded vector space $B_d(n)=\bigoplus_{k\geq 0}B_{d,k}(n)$ of $V_n$-colored open Jacobi diagrams of degree $d$, where $B_{d,k}(n)$ is the subspace of $B_d(n)$ spanned by open Jacobi diagrams with exactly $k$ trivalent vertices.

  We defined in \cite{Mai1} a $\gr(\IA(n))$-action on $B_d(n)$ by using the bracket map
  $$
    [\cdot,\cdot]:B_{d,k}(n)\otimes_{\Z} \opegr^r(\IA(n))\rightarrow B_{d,k+r}(n).
  $$
  We extend the $\gr(\IA(n))$-action to an action of $\gr(\jfE_{\ast}(n))$  on $B_d(n)$.

  We define a $\K$-linear map
  $$
    [\cdot,\cdot]:B_{d,k}(n)\otimes_{\Z}\opegr^r(\jfE_{\ast}(n))\rightarrow B_{d,k+r}(n)
  $$
  by using the following theorem.
  \begin{theorem}[see Theorem \ref{th731}]
   For any $r\geq 1$, we have
   $$[A_{d,k}(n),\jfE_r(n)]\subset A_{d,k+r}(n).$$
  \end{theorem}
  To prove this theorem, we introduce a category $\A^{L}$, which includes as full subcategories the category $\A$ of Jacobi diagrams in handlebodies and the category isomorphic to the PROP for Casimir Lie algebras \cite{Hinich}. (see Section \ref{s4} and Appendix \ref{sA}).

  By using the bracket maps, we obtain $\K$-linear maps
  $$
    \ti{\beta}_{d,k}^r: \opegr^r(\jfE_{\ast}(n))\rightarrow\Hom(B_{d,k}(n),B_{d,k+r}(n)),
  $$
  which form an action of the graded Lie algebra $\gr(\jfE_{\ast}(n))$ on the graded vector space $B_d(n)$.

  We also define a $\K$-linear map
  $$
   c:B_{d,k}(n)\otimes_{\Z} T_r(n) \rightarrow B_{d,k+r}(n),
  $$
  which is an analogue of the contraction map for a vector space and its dual vector space (see Section \ref{s5}).
  By using the map $c$, we obtain $\K$-linear maps
  $$\gamma_{d,k}^r: T_r(n)\rightarrow \Hom(B_{d,k}(n),B_{d,k+r}(n)),$$
  which form an action of the graded Lie algebra $\bigoplus_{r\geq 1}T_r(n)$ on the graded vector space $B_d(n)$.

  Via the isomorphisms \eqref{dereq}, these two actions of the derivation Lie algebra $\Der(\Lie_{\ast}(n))$ on $B_d(n)$ coincide up to sign. (See Theorem \ref{th91}.)

  By using the linear map $c$ for computation, we obtain the surjectivity of the bracket map.
  \begin{proposition}[see Proposition \ref{pbracketsurj}]
    For $n\geq 2d-k$, the bracket map
    $$[\cdot,\cdot]:B_{d,k}(n)\otimes_{\Z} \opegr^1(\IA(n))\rightarrow B_{d,k+1}(n)$$
    is surjective.
  \end{proposition}

 \subsection{The $\GL(n;\Z)$-module structure of $B_d(n)$}
  The $\GL(n;\Z)$-action on $B_d(n)$ that is induced by the $\Aut(F_n)$-action on $A_d(n)$ naturally extends to a polynomial $\GL(V_n)$-action on $B_d(n)$ \cite{Mai1}.
  Therefore, the $\GL(V_n)$-module $B_d(n)$ can be decomposed into the direct sum of images of the Schur functors.
  In general, however, it remains open to obtain an irreducible decomposition of $B_d(n)$ as $\GL(V_n)$-modules.
  We can reduce this problem to the connected parts $B_{d,k}^c(n)\subset B_{d,k}(n)$ (see Theorem \ref{th521}).

  For a partition $\lambda\vdash N$, let $V_{\lambda}$ denote the image of $V_n$ under the Schur functor $\repS_{\lambda}$.
  By using the results by Bar-Natan \cite{BN}, we have isomorphisms of $\GL(V_n)$-modules
  $$B_3(n)=B_{3,0}(n)\oplus\cdots\oplus B_{3,4}(n),$$
  where
  \begin{gather*}
    \begin{split}
     B_{3,0}(n)&\cong V_{(6)}\oplus V_{(4,2)}\oplus V_{(2^3)},\\
     B_{3,1}(n)&\cong V_{(3,1^2)}\oplus V_{(2,1^3)},\\
     B_{3,2}(n)&\cong V_{(4)}\oplus V_{(3,1)}\oplus (V_{(2^2)})^{\oplus 2},\\
     B_{3,3}(n)&=B^c_{3,3}\cong V_{(1^3)},\\
     B_{3,4}(n)&=B^c_{3,4}\cong V_{(2)}
    \end{split}
  \end{gather*}
  (see Proposition \ref{p531} for the cases $d=3,4,5$).

  In general degrees, we obtain irreducible decompositions of $B_{d,k}(n)$ as $\GL(V_n)$-modules for $k=0,1$.
  \begin{proposition}[see Proposition \ref{pdecompgen}]
   For any $d\geq 1$, we have
   $$
     B_{d,0}(n)\cong \bigoplus_{\lambda\vdash d} V_{2\lambda},
   $$
   where $2\lambda=(2\lambda_1,\cdots,2\lambda_r)\vdash 2d$ for $\lambda=(\lambda_1,\cdots,\lambda_r)\vdash d$.
   For any $d\geq 2$, we have
   $$
     B_{d,1}(n)\cong \bigoplus_{\lambda\vdash 2d-1 \text{ with exactly $3$ odd parts}}V_{\lambda}.
   $$
  \end{proposition}

 \subsection{The $\Aut(F_n)$-module structure of $A_d(n)$}
  We consider the $\Aut(F_n)$-module structure of $A_d(n)$ and give an indecomposable decomposition of $A_d(n)$.
  We have $$A_0(n)=\K \quad(n\geq 0),\quad A_d(0)=0 \quad (d\geq 1)$$ and we studied the cases where $d=1, 2$ in \cite{Mai1}.
  Thus, we mainly consider the cases where $d\geq 3,n\geq 1$.

  In \cite{Mai1}, we constructed a functor $A_d:\F^{\op}\rightarrow \fVect$ from the opposite category $\F^{\op}$ of the category $\F$ of finitely generated free groups to the category $\fVect$ of filtered vector spaces, which includes the $\Aut(F_n)$-module structures of $A_d(n)$ for all $n\geq 0$ (see Section \ref{sspartI}).

  For $X\in A_d(2d)$, let
  $$A_d X:\F^{\op}\rightarrow \fVect$$
  denote the subfunctor of $A_d$ generated by $X$.
  That is, for any $n\in\N$, $A_d X(n)$ is the $\Aut(F_n)$-submodule of $A_d(n)$ defined by
  $$A_d X (n):=\Span_\K\{A_d(f)(X)\mid f\in\F^\op(2d,n)\}.$$

  Set
  $$P=\centre{
\begingroup%
  \makeatletter%
  \providecommand\color[2][]{%
    \errmessage{(Inkscape) Color is used for the text in Inkscape, but the package 'color.sty' is not loaded}%
    \renewcommand\color[2][]{}%
  }%
  \providecommand\transparent[1]{%
    \errmessage{(Inkscape) Transparency is used (non-zero) for the text in Inkscape, but the package 'transparent.sty' is not loaded}%
    \renewcommand\transparent[1]{}%
  }%
  \providecommand\rotatebox[2]{#2}%
  \newcommand*\fsize{\dimexpr\f@size pt\relax}%
  \newcommand*\lineheight[1]{\fontsize{\fsize}{#1\fsize}\selectfont}%
  \ifx\svgwidth\undefined%
    \setlength{\unitlength}{53.25118286bp}%
    \ifx\svgscale\undefined%
      \relax%
    \else%
      \setlength{\unitlength}{\unitlength * \real{\svgscale}}%
    \fi%
  \else%
    \setlength{\unitlength}{\svgwidth}%
  \fi%
  \global\let\svgwidth\undefined%
  \global\let\svgscale\undefined%
  \makeatother%
  \begin{picture}(1,0.51725529)%
    \lineheight{1}%
    \setlength\tabcolsep{0pt}%
    \put(0,0){\includegraphics[width=\unitlength,page=1]{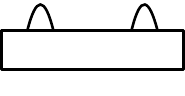}}%
    \put(0.2183161,0.21126291){\makebox(0,0)[lt]{\lineheight{1.45000005}\smash{\begin{tabular}[t]{l}$\sym_{2d}$\end{tabular}}}}%
    \put(0,0){\includegraphics[width=\unitlength,page=2]{P.pdf}}%
  \end{picture}%
\endgroup%
},\quad Q=\centre{
\begingroup%
  \makeatletter%
  \providecommand\color[2][]{%
    \errmessage{(Inkscape) Color is used for the text in Inkscape, but the package 'color.sty' is not loaded}%
    \renewcommand\color[2][]{}%
  }%
  \providecommand\transparent[1]{%
    \errmessage{(Inkscape) Transparency is used (non-zero) for the text in Inkscape, but the package 'transparent.sty' is not loaded}%
    \renewcommand\transparent[1]{}%
  }%
  \providecommand\rotatebox[2]{#2}%
  \newcommand*\fsize{\dimexpr\f@size pt\relax}%
  \newcommand*\lineheight[1]{\fontsize{\fsize}{#1\fsize}\selectfont}%
  \ifx\svgwidth\undefined%
    \setlength{\unitlength}{69.73846498bp}%
    \ifx\svgscale\undefined%
      \relax%
    \else%
      \setlength{\unitlength}{\unitlength * \real{\svgscale}}%
    \fi%
  \else%
    \setlength{\unitlength}{\svgwidth}%
  \fi%
  \global\let\svgwidth\undefined%
  \global\let\svgscale\undefined%
  \makeatother%
  \begin{picture}(1,0.39757915)%
    \lineheight{1}%
    \setlength\tabcolsep{0pt}%
    \put(0,0){\includegraphics[width=\unitlength,page=1]{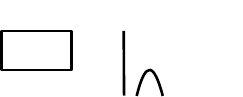}}%
    \put(0.03069267,0.1459089){\makebox(0,0)[lt]{\lineheight{1.45000005}\smash{\begin{tabular}[t]{l}$\alt_{2}$\end{tabular}}}}%
    \put(0,0){\includegraphics[width=\unitlength,page=2]{Q.pdf}}%
  \end{picture}%
\endgroup%
}\in A_d(2d).$$
  Then, we have the following direct decomposition of $A_d(n)$ as $\Aut(F_n)$-modules, which is indecomposable for $n\geq 2d$.
  \begin{theorem}[see Theorem \ref{decompositionofAd}] \label{introdirect}
   We have $A_d(n)=A_d P(n)\oplus A_d Q(n)$ for any $d,n\geq 1$.
   This direct decomposition is indecomposable for $n\geq 2d$.
  \end{theorem}
  In degree $1$, we have $A_1 Q(n)=0$ and $A_1(n)\cong \Sym^2(V_n)$ is simple for $n\geq 1$.
  In \cite{Mai1}, we obtained that the direct decomposition of $A_2(n)$ is indecomposable for $n\geq 3$ (see Theorem 6.9 in \cite{Mai1}).
  We improve Theorem \ref{introdirect} for $d=3,4$ (see Theorems \ref{decomp3} and \ref{decomp4}).

  In general degree $d$, we obtain the radical of $A_{d,k}(n)$ for any $k\geq 0$ if $n\geq 2d$.
  \begin{theorem}[see Theorem \ref{rad}]
   Let $n\geq 2d$.
   The filtration of $A_d(n)$ by the number of trivalent vertices coincides with the radical filtration of $A_d(n)$.
  \end{theorem}

  In degree $3$, we obtain the socle of $A_3(n)$ as well (see Proposition \ref{soc}).

 \subsection{Direct decomposition of the functor $A_d$}
  Lastly, we give an indecomposable decomposition of the functor $A_d$.

  By Theorem \ref{introdirect}, we obtain an indecomposable decomposition of the functor $A_d$.
  \begin{theorem}[see Theorem \ref{directdecomp}]
    We have an indecomposable decomposition
    \begin{equation}\label{introdirectdecomp}
      A_d=A_d P\oplus A_d Q
    \end{equation}
    in the functor category $\fVect^{\F^{\op}}$.
  \end{theorem}
  In degree $1$, we have $A_1 Q=0$ and $A_1=A_1 P$.
  In \cite{Mai1}, we obtained the direct decomposition \eqref{introdirectdecomp} of the functor $A_2$ and proved that \eqref{introdirectdecomp} is indecomposable (see Theorems 6.5 and 6.14 of \cite{Mai1}).

 \subsection{Organization of the paper}
  In Section \ref{sec:pre}, we recall the category $\A$ of Jacobi diagrams in handlebodies, N-series and graded Lie algebras, contents of the previous paper \cite{Mai1}, Hopf algebras and Lie algebras in a linear symmetric strict monoidal category.
  In Section \ref{s3}, we construct the Andreadakis filtration and the Johnson homomorphism of $\End(F_n)$.
  In Section \ref{s4}, we construct an action of the derivation Lie algebra $\Der(\Lie_{\ast}(n))$ on $B_d(n)$, which is defined by the bracket map. In preparation for the definition of the bracket map, we construct an extended category $\A^{L}$ of the category $\A$, which includes a Lie algebra structure.
  In Section \ref{s5}, we define a contraction map, which forms another action of $\Der(\Lie_{\ast}(n))$ on $B_d(n)$.
  In Section \ref{s6}, we prove that two actions of $\Der(\Lie_{\ast}(n))$ on $B_d(n)$ defined in Sections \ref{s4} and \ref{s5} coincide up to sign.
  In Section \ref{s7}, we compute the $\GL(n;\Z)$-module structure of $B_d(n)$.
  In Section \ref{s8}, we study the $\Aut(F_n)$-module structure of $A_d(n)$ by using the $\GL(n;\Z)$-module structure of $B_d(n)$ and the action of $\Der(\Lie_{\ast}(n))$ on $B_d(n)$.
  In Section \ref{s9}, we give an indecomposable decomposition of the functor $A_d$.
  In Appendix \ref{sA}, we study an expected presentation of the category $\A^{L}$.

 \subsection{Acknowledgments}
  The author would like to thank Kazuo Habiro for careful reading and valuable advice.

\section{Preliminaries}\label{sec:pre}
 In this section, we recall the contents of the previous paper \cite{Mai1} and definitions of the category $\A$ of Jacobi diagrams in handlebodies, Hopf algebras and Lie algebras in a symmetric strict monoidal category, and an action of an N-series on a filtered vector space and that of a graded Lie algebra on a graded vector space.

 In what follows, we work over a fixed field $\K$ of characteristic $0$.
 For a vector space $V$ and an abelian group $G$, we just write $V\otimes G$ instead of $V\otimes_{\Z} G$.
 For vector spaces $V$ and $W$, we also write $V\otimes W$ instead of $V\otimes_{\K} W$.

 For $n\geq 0$, let $[n]:=\{1,\cdots,n\}$.

 \subsection{The category $\A$ of Jacobi diagrams in handlebodies}
  Here, we briefly review the category $\A$ of Jacobi diagrams in handlebodies defined in \cite{HM_k}. We use the same notations as in \cite{Mai1}.

  For $n\geq 0$, let $X_n=\centre{
\begingroup%
  \makeatletter%
  \providecommand\color[2][]{%
    \errmessage{(Inkscape) Color is used for the text in Inkscape, but the package 'color.sty' is not loaded}%
    \renewcommand\color[2][]{}%
  }%
  \providecommand\transparent[1]{%
    \errmessage{(Inkscape) Transparency is used (non-zero) for the text in Inkscape, but the package 'transparent.sty' is not loaded}%
    \renewcommand\transparent[1]{}%
  }%
  \providecommand\rotatebox[2]{#2}%
  \newcommand*\fsize{\dimexpr\f@size pt\relax}%
  \newcommand*\lineheight[1]{\fontsize{\fsize}{#1\fsize}\selectfont}%
  \ifx\svgwidth\undefined%
    \setlength{\unitlength}{68.96467443bp}%
    \ifx\svgscale\undefined%
      \relax%
    \else%
      \setlength{\unitlength}{\unitlength * \real{\svgscale}}%
    \fi%
  \else%
    \setlength{\unitlength}{\svgwidth}%
  \fi%
  \global\let\svgwidth\undefined%
  \global\let\svgscale\undefined%
  \makeatother%
  \begin{picture}(1,0.35281587)%
    \lineheight{1}%
    \setlength\tabcolsep{0pt}%
    \put(0,0){\includegraphics[width=\unitlength,page=1]{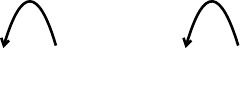}}%
    \put(0.09653969,0.00511188){\makebox(0,0)[lt]{\lineheight{1.45000005}\smash{\begin{tabular}[t]{l}$1$\end{tabular}}}}%
    \put(0.85441458,0.00521592){\makebox(0,0)[lt]{\lineheight{1.45000005}\smash{\begin{tabular}[t]{l}$n$\end{tabular}}}}%
    \put(0.5629133,0.22591732){\makebox(0,0)[lt]{\lineheight{1.45000005}\smash{\begin{tabular}[t]{l}$\cdots$\end{tabular}}}}%
    \put(0,0){\includegraphics[width=\unitlength,page=2]{X_n.pdf}}%
    \put(0.37031179,0.00700557){\makebox(0,0)[lt]{\lineheight{1.45000005}\smash{\begin{tabular}[t]{l}$2$\end{tabular}}}}%
  \end{picture}%
\endgroup%
}$ be the oriented $1$-manifold consisting of $n$ arc components.

\begin{figure}[h]
       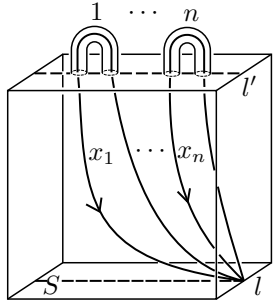
       \caption{The handlebody $U_n$.}
       \label{Fig1}
  \end{figure}

  Let $I=[-1,1]$.
  For $n\geq 0$, let $U_n \subset \R^3$ denote the handlebody of genus $n$ that is obtained from the cube $I^3$ by attaching $n$ handles on the top square $I^2 \times \{1\}$ as depicted in Figure \ref{Fig1}.
  We call $l := I \times \{0\} \times \{-1\}$ the \emph{bottom line} of $U_n$ and $l':= I\times \{0\} \times \{1\}$ the \emph{upper line} of $U_n$.
  We call $S := I^2 \times \{-1\}$ the \emph{bottom square} of $U_n$.
  For $i\in[n]$, let $x_i$ be a loop which goes through only the $i$-th handle of the handlebody $U_n$ just once and let $x_i$ denote its homotopy class as well. In what follows, for loops $\gamma_1$ and $\gamma_2$ with base points on $l$, let $\gamma_2\gamma_1$ denote the loop that goes through $\gamma_1$ first and then goes through $\gamma_2$. That is, we write a product of elements of the fundamental group of $U_n$ in the opposite order to the usual one.
  Let $H=H_1(U_n;\Z)$ and let $\bar{x}_i\in H$ be the homology class of $x_i$.
  We have $H=\bigoplus_{i=1}^n \Z\bar{x}_i$ and $\pi_1(U_n)=\langle x_1,\cdots,x_n\rangle$.
  Let
  $$V_n=H^1(U_n;\K)=\Hom(H,\K)$$
  and let $\{v_1,\cdots,v_n\}$ be the dual basis of $\{\bar{x}_1,\cdots,\bar{x}_n\}$.

  The objects in $\A$ are nonnegative integers.

  For $m,n\geq 0$, the hom-set $\A(m,n)$ is the $\K$-vector space spanned by \emph{$(m,n)$-Jacobi diagrams} modulo the STU relation.
  An $(m,n)$-Jacobi diagram is a Jacobi diagram on $X_n$ mapped into $U_m$ in such a way that the endpoints of $X_n$ are uniformly distributed on the bottom line $l$ of $U_m$ (see \cite{HM_k, Mai1} for further details).
  We usually depict $(m,n)$-Jacobi diagrams by drawing their images under the orthogonal projection of $\R^3$ onto $\R\times \{0\}\times \R$.

  The \emph{degree} of an $(m,n)$-Jacobi diagram is the degree of its Jacobi diagram.
  Let $\A_d(m,n)\subset \A(m,n)$ be the subspace spanned by $(m,n)$-Jacobi diagrams of degree $d$.
  We have $\A(m,n)=\bigoplus_{d\geq0} \A_d(m,n)$.

  The category $\A$ has a structure of a linear symmetric strict monoidal category.
  The tensor product on objects is addition. The monoidal unit is $0$. The tensor product on morphisms is juxtaposition followed by horizontal rescaling and relabelling of indices.
  The symmetry is determined by
  $$P_{1,1}=\scalebox{0.5}{$\centre{
\begingroup%
  \makeatletter%
  \providecommand\color[2][]{%
    \errmessage{(Inkscape) Color is used for the text in Inkscape, but the package 'color.sty' is not loaded}%
    \renewcommand\color[2][]{}%
  }%
  \providecommand\transparent[1]{%
    \errmessage{(Inkscape) Transparency is used (non-zero) for the text in Inkscape, but the package 'transparent.sty' is not loaded}%
    \renewcommand\transparent[1]{}%
  }%
  \providecommand\rotatebox[2]{#2}%
  \newcommand*\fsize{\dimexpr\f@size pt\relax}%
  \newcommand*\lineheight[1]{\fontsize{\fsize}{#1\fsize}\selectfont}%
  \ifx\svgwidth\undefined%
    \setlength{\unitlength}{68.24999903bp}%
    \ifx\svgscale\undefined%
      \relax%
    \else%
      \setlength{\unitlength}{\unitlength * \real{\svgscale}}%
    \fi%
  \else%
    \setlength{\unitlength}{\svgwidth}%
  \fi%
  \global\let\svgwidth\undefined%
  \global\let\svgscale\undefined%
  \makeatother%
  \begin{picture}(1,1.27812969)%
    \lineheight{1}%
    \setlength\tabcolsep{0pt}%
    \put(0,0){\includegraphics[width=\unitlength,page=1]{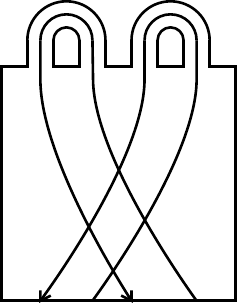}}%
  \end{picture}%
\endgroup%
}$}: 2\rightarrow 2.$$

 \subsection{N-series and graded Lie algebras}
  Here, we briefly review the definition of an action of an N-series on a filtered vector space and the induced action of the graded Lie algebra on the graded vector space (see \cite{Mai1} for details).

  An \emph{N-series} $K_{\ast}=(K_n)_{n\geq 1}$ of a group $K$ is a descending series
  $$
    K=K_1\supset K_2\supset \cdots
  $$
  such that $[K_n,K_m]\subset K_{n+m}$ for all $n,m\geq 1$.

  A \emph{morphism} $f: G_{\ast}\rightarrow K_{\ast}$ between N-series is a group homomorphism $f: G_1 \rightarrow K_1$ such that we have  $f(G_n)\subset K_n$ for all $n\geq 1.$

  For a filtered vector space $W_{\ast}$, set
  $$
   \Aut_n(W_{\ast}):=\{\phi\in\Aut_{\fVect}(W_{\ast})\mid [\phi,w]\in W_{k+n}\text{ for all }w\in W_k, k\geq 0\}\quad (n\geq 1),
  $$
  where $[\phi,w]:=\phi(w)-w$ for $w\in W_k$.
  We can easily check that $\Aut_{\ast}(W_{\ast}):=(\Aut_n(W_{\ast}))_{n\geq 1}$ is an N-series.

  \begin{definition}(Action of N-series on filtered vector spaces)
   Let $K_{\ast}$ be an N-series and $W_{\ast}$ be a filtered vector space.
   An \emph{action} of $K_{\ast}$ on $W_{\ast}$ is a morphism $f:K_{\ast}\rightarrow \Aut_{\ast}(W_{\ast})$ between N-series.
  \end{definition}

  For an N-series $K_{\ast}$, we have a graded Lie algebra $\gr(K_{\ast})=\bigoplus_{n\geq 1} K_n/K_{n+1}$, where the Lie bracket is defined by the commutator.

  For a graded vector space $W=\bigoplus_{k\geq 0}W_k$, set
  $$
  \End_n(W):=\{\phi\in\End(W)\mid \phi(W_k)\subset W_{k+n}\text{ for }k\geq 0\}\quad(n\geq 1).
  $$
  We can check that $\End_{+}(W)=\bigoplus_{n\geq 1}\End_n(W)$ is a graded Lie algebra, where the Lie bracket is defined by
  $$[f,g]:=f\circ g-g\circ f \quad \text{for} \quad f\in\End_k(W),g\in\End_l(W)\; (k,l\geq 1).$$

  \begin{definition}(Action of graded Lie algebras on graded vector spaces)
   Let $L_{+}=\bigoplus_{n\geq 1}L_n$ be a graded Lie algebra and $W=\bigoplus_{k\geq 0}W_k$ be a graded vector space.
   An \emph{action} of $L_{+}$ on $W$ is a morphism $f:L_{+}\rightarrow \End_{+}(W)$ between graded Lie algebras.
  \end{definition}

  \begin{proposition}
   An action of an N-series $K_{\ast}$ on a filtered vector space $W_{\ast}$ induces an action of the graded Lie algebra $\gr(K_{\ast})$ on the graded vector space $\gr(W_{\ast})$, which is a morphism
   $$\rho_{+}:\bigoplus_{n\geq 1}\opegr^n(K_{\ast}) \rightarrow\bigoplus_{n\geq 1}\End_n(\gr(W_{\ast}))$$ defined by $\rho_{+}(gK_{n+1})([v]_{W_{k+1}})=[[g,v]]_{W_{k+n+1}}$ for $gK_{n+1}\in \opegr^n(K_{\ast})$, $[v]_{W_{k+1}}\in \opegr^k(W_{\ast})$.
  \end{proposition}
  The proof can be seen in Proposition 5.14 of \cite{Mai1}.

 \subsection{Contents of the previous paper}\label{sspartI}
  Here, we briefly review the notations and contents of the previous paper \cite{Mai1}.
  Let $\Aut(F_n)$ denote the automorphism group of the free group $F_n$ of rank $n$ and $\GL(n;\Z)$ the general linear group of degree $n$.
  Let $\IA(n)$ denote the IA-automorphism group of $F_n$, that is the kernel of the canonical surjection
  $$\Aut(F_n)\rightarrow \Aut(H_1(F_n;\Z))\cong \GL(n;\Z).$$
  Let $\Gamma_{\ast}(\IA(n))=(\Gamma_{r}(\IA(n)))_{r\geq 1}$ denote the lower central series of $\IA(n)$, and
  $\gr(\IA(n))=\bigoplus_{r\geq 1}\opegr^r(\IA(n))$ the associated graded Lie algebra, where $\opegr^r(\IA(n))=\Gamma_{r}(\IA(n))/\Gamma_{r+1}(\IA(n))$.

  Let $A_d(n)=\A_d(0,n)$ denote the $\K$-vector space of Jacobi diagrams of degree $d$ on $X_n$.
  We consider a filtration for $A_d(n)$
  $$
    A_d(n)=A_{d,0}(n)\supset A_{d,1}(n)\supset\cdots\supset A_{d,2d-2}(n)\supset A_{d,2d-1}(n)=0,
  $$
  such that $A_{d,k}(n)\subset A_d(n)$ is the subspace spanned by Jacobi diagrams with at least $k$ trivalent vertices. Hence, $A_d(n)$ is a filtered vector space.

  Let $\F$ denote the category of finitely generated free groups, and $\fVect$ the category of filtered vector spaces over $\K$.

  We have a $\K$-vector space isomorphism
  $$Z:\K\F^{\op}(m,n)\xrightarrow{\cong} \A_0(m,n)$$
  from the hom-set $\K\F^{\op}(m,n)$ of the $\K$-linearization of the opposite category of $\F$ to the degree $0$ part of the hom-set $\A(m,n)$ \cite{HM_k}.
  We define a functor
  $$A_d:\F^{\op}\rightarrow \fVect$$
  by $A_d(n)=\A_d(0,n)$ for an object $n\in \N$ and
  $A_d(f)=Z(f)_{\ast}$ for a morphism $f\in \F^{\op}(m,n)$, where $Z(f)_{\ast}$ denotes the post-composition with $Z(f)$.
  By restricting this functor to the automorphism group, we obtain an action of the opposite group $\Aut(F_n)^{\op}$ of $\Aut(F_n)$ on $A_d(n)$ for each $n\geq 0$. We consider this action as a right action of $\Aut(F_n)$ on $A_d(n)$.
  The $\Aut(F_n)$-action on $A_d(n)$ induces an action on $A_d(n)$ of the outer automorphism group $\Out(F_n)$ of $F_n$ (see Theorem 5.1 in \cite{Mai1}).

  On the other hand, the associated graded vector space $\gr(A_d(n))$ of $A_d(n)$ is identified via the PBW map \cite{BN_v, BN_sl}
  \begin{equation}\label{PBW}
    \theta_{d,n}:\gr(A_d(n))\xrightarrow{\cong} B_d(n)
  \end{equation}
  with the graded $\K$-vector space $B_d(n)=\bigoplus_{k\geq 0}B_{d,k}(n)= \bigoplus_{k=0}^{2d-2}B_{d,k}(n)$ of $V_n$-colored open Jacobi diagrams of degree $d$, where the grading is determined by the number of trivalent vertices.
  Note that we have $\theta_{d,n}=\bigoplus_{k}\theta_{d,n,k}$, where $$\theta_{d,n,k}:\opegr^k(A_d(n))\xrightarrow{\cong} B_{d,k}(n).$$

  Let $\FAb$ denote the category of finitely generated free abelian groups, and $\gVect$ the category of graded vector spaces over $\K$.

  We define a functor
  $$B_d: \FAb^\op\rightarrow \gVect$$
  by sending an object $n\in\N$ to the graded vector space $B_d(n)$
  and a morphism $f\in\FAb^\op(m,n)=\Mat(m,n;\Z)$ to
  $B_d(f)$, which is a right action on each coloring, where we consider an element of $V_n$ as a $1\times n$ matrix.
  By restricting this functor to the automorphism group, we obtain an action of the opposite group $\GL(n;\Z)^{\op}$ of $\GL(n;\Z)$ on $B_d(n)$ for each $n\geq 0$. We consider this action as a right action of $\GL(n;\Z)$ on $B_d(n)$. Note that the $\GL(n;\Z)$-action on $B_d(n)$ naturally extends to a $\GL(V_n)$-action on $B_d(n)$.

  \begin{proposition}[see Proposition 3.2 of \cite{Mai1}]
   For $d\geq 0$, the PBW maps \eqref{PBW} give a natural isomorphism
   $$\theta_d:\gr\circ A_d\overset{\cong}\Rightarrow B_d\circ\ab^\op,$$
   where $\ab^{\op}$ denotes the opposite functor of the abelianization functor and $\gr$ denote the functor that sends a filtered vector space to its associated graded vector space.
  \end{proposition}

  By this proposition, it turns out that the $\Aut(F_n)$-action on $A_d(n)$ induces two actions on $B_d(n)$, which form an action of an extended graded Lie algebra on a graded vector space (see Theorem 5.15 in \cite{Mai1}).
  One of them is the $\GL(n;\Z)$-action and the other of them is an action of the graded Lie algebra $\gr(\IA(n))$ on the graded vector space $B_d(n)$, which consists of $\GL(n;\Z)$-module homomorphisms
  \begin{equation}\label{bracket}
    [\cdot,\cdot]: B_{d,k}(n)\otimes \opegr^r(\IA(n)) \rightarrow B_{d,k+r}(n)
  \end{equation}
  for $k\geq 0, r\geq 1$ (see Proposition 5.9 and Corollary 5.15 of \cite{Mai1}).
  By using these two actions on $B_d(n)$, we obtained an indecomposable decomposition of $A_2(n)$ as $\Aut(F_n)$-modules (see Theorem 6.9 of \cite{Mai1}).

 \subsection{Hopf algebra in a symmetric strict monoidal category}\label{ssHopf}
  We review the definition of a Hopf algebra in a symmetric strict monoidal category.
  Let $\catC=(\catC,\otimes,I,P)$ be a symmetric strict monoidal category.
  A \emph{Hopf algebra} in $\catC$ is an object $H$ in $\catC$ equipped with morphisms
  $$
    \mu:H\otimes H\rightarrow H, \quad \eta:I\rightarrow H, \quad \Delta:H\rightarrow H\otimes H, \quad \epsilon:H\rightarrow I,\quad S:H\rightarrow H,
  $$
  called the \emph{multiplication, unit, comultiplication, counit}, and \emph{antipode}, respectively, satisfying
  \begin{enumerate}
   \item $\mu(\mu\otimes \id_H)=\mu(\id_H \otimes\mu),\quad \mu(\eta\otimes \id_H)=\id_H=\mu(\id_H\otimes \eta),$
   \item $(\Delta \otimes \id_H)\Delta=(\id_H\otimes \Delta)\Delta,\quad (\epsilon\otimes \id_H)\Delta=\id_H=(\id_H\otimes \epsilon)\Delta,$
   \item $\epsilon\eta=\id_I, \quad \epsilon\mu=\epsilon\otimes\epsilon,\quad \Delta\eta=\eta\otimes\eta,$
   \item $\Delta\mu=(\mu\otimes\mu)(\id_H\otimes P_{H,H}\otimes \id_H)(\Delta\otimes\Delta),$
   \item $\mu(\id_H\otimes S)\Delta=\mu(S\otimes \id_H)\Delta=\eta\epsilon.$
  \end{enumerate}
  A Hopf algebra $H$ is said to be
  \emph{cocommutative} if
  $P_{H,H}\Delta=\Delta.$

  Define $\mu_n : H^{\otimes{n}} \otimes H^{\otimes{n}} \rightarrow H^{\otimes{n}}$ and $\Delta_m : H^{\otimes{m}} \rightarrow H^{\otimes{m}} \otimes H^{\otimes{m}}$ inductively
  by
  $$\mu_0=\id_I, \quad \mu_{n+1}=(\mu_n \otimes \mu)(\id_{H^{\otimes{n}}} \otimes P_{H,H^{\otimes{n}}} \otimes \id_H)$$ for $n\geq 0$ and by
  $$
  \Delta_0=\id_I,\quad \Delta_{m+1}=(\id_{H^{\otimes{m}}} \otimes P_{H^{\otimes{m}},H} \otimes \id_H)(\Delta_m \otimes \Delta)
  $$ for $m\geq 0$.

  For morphisms $f,\ f': H^{\otimes{m}} \rightarrow H^{\otimes{n}}$, $m,n\geq 0$, the \emph{convolution} $f\ast f'$ of $f$ and $f'$ is defined by
  $$
    f\ast f':= \mu_n (f\otimes f')\Delta_m.
  $$


  The category $\A$ has a cocommutative Hopf algebra with the object $1$, where
  \begin{gather*}
    \mu=\centre{
\begingroup%
  \makeatletter%
  \providecommand\color[2][]{%
    \errmessage{(Inkscape) Color is used for the text in Inkscape, but the package 'color.sty' is not loaded}%
    \renewcommand\color[2][]{}%
  }%
  \providecommand\transparent[1]{%
    \errmessage{(Inkscape) Transparency is used (non-zero) for the text in Inkscape, but the package 'transparent.sty' is not loaded}%
    \renewcommand\transparent[1]{}%
  }%
  \providecommand\rotatebox[2]{#2}%
  \newcommand*\fsize{\dimexpr\f@size pt\relax}%
  \newcommand*\lineheight[1]{\fontsize{\fsize}{#1\fsize}\selectfont}%
  \ifx\svgwidth\undefined%
    \setlength{\unitlength}{30.7499995bp}%
    \ifx\svgscale\undefined%
      \relax%
    \else%
      \setlength{\unitlength}{\unitlength * \real{\svgscale}}%
    \fi%
  \else%
    \setlength{\unitlength}{\svgwidth}%
  \fi%
  \global\let\svgwidth\undefined%
  \global\let\svgscale\undefined%
  \makeatother%
  \begin{picture}(1,1.30247283)%
    \lineheight{1}%
    \setlength\tabcolsep{0pt}%
    \put(0,0){\includegraphics[width=\unitlength,page=1]{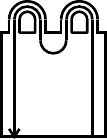}}%
  \end{picture}%
\endgroup%
},\;
    \eta=\centre{
\begingroup%
  \makeatletter%
  \providecommand\color[2][]{%
    \errmessage{(Inkscape) Color is used for the text in Inkscape, but the package 'color.sty' is not loaded}%
    \renewcommand\color[2][]{}%
  }%
  \providecommand\transparent[1]{%
    \errmessage{(Inkscape) Transparency is used (non-zero) for the text in Inkscape, but the package 'transparent.sty' is not loaded}%
    \renewcommand\transparent[1]{}%
  }%
  \providecommand\rotatebox[2]{#2}%
  \newcommand*\fsize{\dimexpr\f@size pt\relax}%
  \newcommand*\lineheight[1]{\fontsize{\fsize}{#1\fsize}\selectfont}%
  \ifx\svgwidth\undefined%
    \setlength{\unitlength}{30.7499995bp}%
    \ifx\svgscale\undefined%
      \relax%
    \else%
      \setlength{\unitlength}{\unitlength * \real{\svgscale}}%
    \fi%
  \else%
    \setlength{\unitlength}{\svgwidth}%
  \fi%
  \global\let\svgwidth\undefined%
  \global\let\svgscale\undefined%
  \makeatother%
  \begin{picture}(1,1.01165838)%
    \lineheight{1}%
    \setlength\tabcolsep{0pt}%
    \put(0,0){\includegraphics[width=\unitlength,page=1]{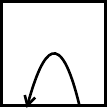}}%
  \end{picture}%
\endgroup%
},\;
    \Delta=\centre{
\begingroup%
  \makeatletter%
  \providecommand\color[2][]{%
    \errmessage{(Inkscape) Color is used for the text in Inkscape, but the package 'color.sty' is not loaded}%
    \renewcommand\color[2][]{}%
  }%
  \providecommand\transparent[1]{%
    \errmessage{(Inkscape) Transparency is used (non-zero) for the text in Inkscape, but the package 'transparent.sty' is not loaded}%
    \renewcommand\transparent[1]{}%
  }%
  \providecommand\rotatebox[2]{#2}%
  \newcommand*\fsize{\dimexpr\f@size pt\relax}%
  \newcommand*\lineheight[1]{\fontsize{\fsize}{#1\fsize}\selectfont}%
  \ifx\svgwidth\undefined%
    \setlength{\unitlength}{30.74999894bp}%
    \ifx\svgscale\undefined%
      \relax%
    \else%
      \setlength{\unitlength}{\unitlength * \real{\svgscale}}%
    \fi%
  \else%
    \setlength{\unitlength}{\svgwidth}%
  \fi%
  \global\let\svgwidth\undefined%
  \global\let\svgscale\undefined%
  \makeatother%
  \begin{picture}(1,1.45409838)%
    \lineheight{1}%
    \setlength\tabcolsep{0pt}%
    \put(0,0){\includegraphics[width=\unitlength,page=1]{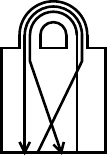}}%
  \end{picture}%
\endgroup%
},\;
    \epsilon=\centre{
\begingroup%
  \makeatletter%
  \providecommand\color[2][]{%
    \errmessage{(Inkscape) Color is used for the text in Inkscape, but the package 'color.sty' is not loaded}%
    \renewcommand\color[2][]{}%
  }%
  \providecommand\transparent[1]{%
    \errmessage{(Inkscape) Transparency is used (non-zero) for the text in Inkscape, but the package 'transparent.sty' is not loaded}%
    \renewcommand\transparent[1]{}%
  }%
  \providecommand\rotatebox[2]{#2}%
  \newcommand*\fsize{\dimexpr\f@size pt\relax}%
  \newcommand*\lineheight[1]{\fontsize{\fsize}{#1\fsize}\selectfont}%
  \ifx\svgwidth\undefined%
    \setlength{\unitlength}{30.74999894bp}%
    \ifx\svgscale\undefined%
      \relax%
    \else%
      \setlength{\unitlength}{\unitlength * \real{\svgscale}}%
    \fi%
  \else%
    \setlength{\unitlength}{\svgwidth}%
  \fi%
  \global\let\svgwidth\undefined%
  \global\let\svgscale\undefined%
  \makeatother%
  \begin{picture}(1,1.43902441)%
    \lineheight{1}%
    \setlength\tabcolsep{0pt}%
    \put(0,0){\includegraphics[width=\unitlength,page=1]{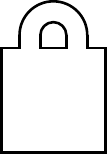}}%
  \end{picture}%
\endgroup%
},\;
    S=\centre{
\begingroup%
  \makeatletter%
  \providecommand\color[2][]{%
    \errmessage{(Inkscape) Color is used for the text in Inkscape, but the package 'color.sty' is not loaded}%
    \renewcommand\color[2][]{}%
  }%
  \providecommand\transparent[1]{%
    \errmessage{(Inkscape) Transparency is used (non-zero) for the text in Inkscape, but the package 'transparent.sty' is not loaded}%
    \renewcommand\transparent[1]{}%
  }%
  \providecommand\rotatebox[2]{#2}%
  \newcommand*\fsize{\dimexpr\f@size pt\relax}%
  \newcommand*\lineheight[1]{\fontsize{\fsize}{#1\fsize}\selectfont}%
  \ifx\svgwidth\undefined%
    \setlength{\unitlength}{30.74999894bp}%
    \ifx\svgscale\undefined%
      \relax%
    \else%
      \setlength{\unitlength}{\unitlength * \real{\svgscale}}%
    \fi%
  \else%
    \setlength{\unitlength}{\svgwidth}%
  \fi%
  \global\let\svgwidth\undefined%
  \global\let\svgscale\undefined%
  \makeatother%
  \begin{picture}(1,1.45121951)%
    \lineheight{1}%
    \setlength\tabcolsep{0pt}%
    \put(0,0){\includegraphics[width=\unitlength,page=1]{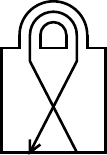}}%
  \end{picture}%
\endgroup%
}.
  \end{gather*}

 \subsection{Lie algebra in a linear symmetric strict monoidal category}\label{ssLie}
  We review the definition of a Lie algebra in a linear symmetric strict monoidal category.
  Let $\catC=(\catC,\otimes,I,P)$ be a linear symmetric strict monoidal category.
  A \emph{Lie algebra} in $\catC$ is an object $L$ in $\catC$ equipped with a morphism
  $$[\cdot,\cdot]:L\otimes L\rightarrow L$$
  satisfying
  \begin{enumerate}
    \item $[\cdot,\cdot](\id_{L\otimes L}+P_{L,L})=0$,
    \item $[\cdot,\cdot](\id_L\otimes[\cdot,\cdot])(\id_{L^{\otimes 3}}+\sigma+\sigma^2)=0,$ where $\sigma=(1,2,3):L^{\otimes 3}\rightarrow L^{\otimes 3}$.
  \end{enumerate}


\section{Andreadakis filtration $\jfE_{\ast}(n)$ of $\End(F_n)$}\label{s3}
 We briefly review the Andreadakis filtration and the Johnson homomorphism of $\Aut(F_n)$. See \cite{Satoh_survey} for further details.
 Then we consider its extension to the endomorphism monoid $\End(F_n)$ of $F_n$.

 \subsection{Andreadakis filtration $\jfA_{\ast}(n)$ of $\Aut(F_n)$}\label{ss31}
  In what follows, we consider the left action of $\Aut(F_n)$ on $F_n$.
  Let $\Gamma_r:=\Gamma_r(F_n)$ denote the $r$-th term of the lower central series of the free group $F_n$ of rank $n$.
  Let $\Lie_r(n) :=\Gamma_r/\Gamma_{r+1}$ for $r\geq 1$.
  Note that $H=\Lie_1(n)$ and that $\Lie_r(n)$ is the degree $r$ part of the free Lie algebra $\Lie_{\ast}(n)$ on $H$.

  For $r\geq 0$, the left action of $\Aut(F_n)$ on each nilpotent quotient $F_n/\Gamma_{r+1}$ induces a group homomorphism
  $$\Aut(F_n)\rightarrow \Aut(F_n/\Gamma_{r+1}).$$
  Set $$\jfA_{r}(n) := \ker(\Aut(F_n)\rightarrow \Aut(F_n/\Gamma_{r+1})) \lhd \Aut(F_n).$$
  Then we have a filtration, which is called the \emph{Andreadakis filtration} of $\Aut(F_n)$: $$\Aut(F_n)=\jfA_0(n)\supset\jfA_1(n)=\IA(n)\supset\jfA_2(n)\supset\cdots.$$
  For $r\geq 1$, the \emph{Johnson homomorphism} $$\tau_r:\gr^r(\jfA_{\ast}(n))\hookrightarrow \Hom(H,\Lie_{r+1}(n))$$
  is the injective homomorphism induced by the group homomorphism
  $$\tau'_r:\jfA_r(n)\rightarrow\Hom(H,\Lie_{r+1}(n))$$
  defined by
  $$\tau'_r(f)(x\Gamma_2):=f(x)x^{-1}\Gamma_{r+2}\quad\text{for}\ f\in\jfA_r(n), x\in F_n.$$

 \subsection{The target group of the Johnson homomorphism}\label{ss315}
  The target group $\Hom(H,\Lie_{r+1}(n))\cong H^{\ast}\otimes \Lie_{r+1}(n)$ of the Johnson homomorphism is identified with the degree $r$ part $\Der_r(\Lie_{\ast}(n))$ of the derivation Lie algebra $\Der(\Lie_{\ast}(n))$ of the free Lie algebra $\Lie_{\ast}(n)$ and with the tree module $T_r(n)$ via abelian group isomorphisms
  \begin{equation}\label{equiv}
    H^{\ast}\otimes \Lie_{r+1}(n)\cong \Der_r(\Lie_{\ast}(n)) \cong T_r(n).
  \end{equation}
  Here, we briefly review the derivation Lie algebra and the tree module. (See \cite{Satoh_survey} for details.)

  A \emph{derivation} $f$ of $\Lie_{\ast}(n)$ is a $\Z$-linear map
  $f:\Lie_{\ast}(n)\rightarrow \Lie_{\ast}(n)$ such that $f([a,b])=[f(a),b]+[a,f(b)]$ for any $a,b\in \Lie_{\ast}(n)$.
  The derivation Lie algebra $\Der(\Lie_{\ast}(n))$ of the Lie algebra $\Lie_{\ast}(n)$ is the set of all derivations of $\Lie_{\ast}(n)$.
  The degree $r$ part $\Der_r(\Lie_{\ast}(n))$ of the derivation Lie algebra is defined to be
  $$\Der_r(\Lie_{\ast}(n))=\{f\in \Der(\Lie_{\ast}(n))\mid f(a)\in \Lie_{r+1}(n) \text{ for any }a\in H\}.$$
  Then we have $\Der(\Lie_{\ast}(n))=\bigoplus_{r\geq 0}\Der_r(\Lie_{\ast}(n))$ and abelian group isomorphisms
  $$\Der_r(\Lie_{\ast}(n))\cong \Hom(H,\Lie_{r+1}(n))\cong H^{\ast}\otimes \Lie_{r+1}(n).$$

  We call a connected Jacobi diagram with no cycle a \emph{trivalent tree}.
  For $r\geq 0$, a trivalent tree is called a \emph{rooted trivalent tree of degree $r$} if it has one univalent vertex (called the \emph{root}) that is colored by an element of $H^{\ast}$ and $r+1$ univalent vertices (called \emph{leaves}) that are colored by elements of $H$.
  Let $T_r(n)$ denote the $\Z$-module spanned by rooted trivalent trees of degree $r$ modulo the AS, IHX and multilinearity relations.
  We have an abelian group isomorphism $$\Phi:H^{\ast}\otimes\Lie_{r+1}(n)\xrightarrow{\cong} T_{r}(n)$$
  defined by
  $$\Phi(v_i\otimes [\bar{x}_{i_1},\cdots,[\bar{x}_{i_r},\bar{x}_{i_{r+1}}]\cdots])=\centre{
\begingroup%
  \makeatletter%
  \providecommand\color[2][]{%
    \errmessage{(Inkscape) Color is used for the text in Inkscape, but the package 'color.sty' is not loaded}%
    \renewcommand\color[2][]{}%
  }%
  \providecommand\transparent[1]{%
    \errmessage{(Inkscape) Transparency is used (non-zero) for the text in Inkscape, but the package 'transparent.sty' is not loaded}%
    \renewcommand\transparent[1]{}%
  }%
  \providecommand\rotatebox[2]{#2}%
  \newcommand*\fsize{\dimexpr\f@size pt\relax}%
  \newcommand*\lineheight[1]{\fontsize{\fsize}{#1\fsize}\selectfont}%
  \ifx\svgwidth\undefined%
    \setlength{\unitlength}{64.90570833bp}%
    \ifx\svgscale\undefined%
      \relax%
    \else%
      \setlength{\unitlength}{\unitlength * \real{\svgscale}}%
    \fi%
  \else%
    \setlength{\unitlength}{\svgwidth}%
  \fi%
  \global\let\svgwidth\undefined%
  \global\let\svgscale\undefined%
  \makeatother%
  \begin{picture}(1,0.7267743)%
    \lineheight{1}%
    \setlength\tabcolsep{0pt}%
    \put(0,0){\includegraphics[width=\unitlength,page=1]{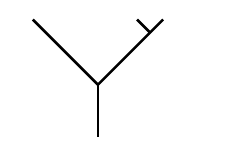}}%
    \put(-0.00233211,0.703062){\makebox(0,0)[lt]{\lineheight{1.45000005}\smash{\begin{tabular}[t]{l}$\bar{x}_{i_1}$\end{tabular}}}}%
    \put(0.35449094,0.01086311){\makebox(0,0)[lt]{\lineheight{1.45000005}\smash{\begin{tabular}[t]{l}$v_i$\end{tabular}}}}%
    \put(0.46271858,0.6989744){\makebox(0,0)[lt]{\lineheight{1.45000005}\smash{\begin{tabular}[t]{l}$\bar{x}_{i_r}$\end{tabular}}}}%
    \put(0,0){\includegraphics[width=\unitlength,page=2]{rootedtree1.pdf}}%
    \put(0.65665316,0.70061975){\makebox(0,0)[lt]{\lineheight{1.45000005}\smash{\begin{tabular}[t]{l}$\bar{x}_{i_{r+1}}$\end{tabular}}}}%
  \end{picture}%
\endgroup%
}$$
  for $v_i\in H^{\ast}$,  $[\bar{x}_{i_1},\cdots,[\bar{x}_{i_r},\bar{x}_{i_{r+1}}]\cdots]\in \Lie_{r+1}(n)$.

 \subsection{Andreadakis filtration $\jfE_{\ast}(n)$ of $\End(F_n)$}\label{ss32}
  We extend the above construction to the endomorphism monoid $\End(F_n)$ of $F_n$.
  For $r\geq 0$, consider the canonical map
  $$\rho_r:\End(F_n)\rightarrow \End(F_n/\Gamma_{r+1})$$
  and set $\jfE_r(n):=\ker(\rho_r)$.
  Then we have a filtration of monoids
  $$\End(F_n)=\jfE_0(n)\supset\jfE_1(n)\supset\cdots$$
  and we call $\jfE_{\ast}(n)=(\jfE_{r}(n))_{r\geq 0}$ the \emph{Andreadakis filtration} of $\End(F_n)$.

  For $f\in\End(F_n)$ and $x,y\in F_n$, set
  $$[f,x]:=f(x)x^{-1},\quad {}^y x=yxy^{-1},$$
  and for a subset $T\subset F_n$,
  set $$[f,T]=\{[f,x]\in F_n \mid x\in T\}.$$
  We can easily check the following lemma.
  \begin{lemma}
    $$f\in \jfE_{r}(n)\quad \Leftrightarrow\quad [f, F_n]\subset \Gamma_{r+1}\quad \Leftrightarrow\quad [f,x_i]\in \Gamma_{r+1}\; (\text{for any }i\in [n]).$$
  \end{lemma}
  For subsets $S\subset \End(F_n)$ and $T\subset F_n$,
  let $[S,T]$ denote the subgroup of $F_n$ generated by the elements $[f,x]$ for $f\in S, x\in T$.
  \begin{lemma}\label{l711}
   We have
   $$
     [\jfE_r(n),\Gamma_k]\subset \Gamma_{k+r}
   $$
   for $r\geq 0, k\geq 1$.
  \end{lemma}
  \begin{proof}
   It is well known that $[\jfA_r(n),\Gamma_k]\subset \Gamma_{k+r}$ by Andreadakis \cite{Andreadakis}. The same proof can be applied to $\jfE_r(n)$. We use induction on $k$.
   When $k=1$, we have $[\jfE_r(n),F_n]\subset\Gamma_{r+1}$ by the definition of $\jfE_r(n)$.
   Suppose that $[\jfE_r(n),\Gamma_{k-1}]\subset\Gamma_{k-1+r}$.
   We will show that $[\jfE_r(n),\Gamma_{k}]\subset\Gamma_{k+r}$.
   Let $f\in\jfE_r(n)$.
   Recall that $\Gamma_k$ is generated by the commutator $[x,y]$ with $x\in \Gamma_{k-1},y\in F_n$.
   We can check that for $x\in\Gamma_{k-1},y\in F_n$, we have
   $$
     [f,[x,y]]={}^{[f,y]}([[f,y]^{-1},f(x)]\cdot
     [[f,x],{}^{x}y]\cdot[[x,y],[f,y]^{-1}])\in \Gamma_{k+r}.
   $$
   For $z,w\in\Gamma_k$, we have
   \begin{gather*}
     [f,zw]=[f,z]\cdot{}^{z}[f,w]\equiv [f,z][f,w] \quad(\bmod\:\Gamma_{k+r+1})
   \end{gather*}
   and by letting $w=z^{-1}$, we have
   \begin{gather*}
     [f,z^{-1}]\equiv [f,z]^{-1}\quad(\bmod\:\Gamma_{k+r+1}).
   \end{gather*}
   Therefore, we have $[f,z]\in\Gamma_{k+r}$ for any $z\in\Gamma_k$.
  \end{proof}

  Define a map
  $$\sigma:\End(F_n)\rightarrow \End(F_n)$$
  by
  $\sigma(f)=\ti{f}$ for $f\in \End(F_n)$,
  where
  $$\ti{f}(x_i)=[f,x_i]^{-1}x_i=x_i f(x_i)^{-1} x_i$$
  for $i\in [n]$.

  \begin{lemma}\label{l712}
   We have
   \begin{gather}\label{l7121}
     \sigma^2=\id_{\End(F_n)}
   \end{gather}
   \begin{gather}\label{l7122}
     f\in \jfE_r(n)\quad\Rightarrow\quad \sigma(f)\in \jfE_r(n)
   \end{gather}
   \begin{gather}\label{l7123}
     f\in \jfE_r(n)\quad\Rightarrow\quad f\sigma(f),\: \sigma(f) f\in \jfE_{2r}(n).
   \end{gather}
  \end{lemma}
  \begin{proof}
    We have (\ref{l7121}) since for any $f\in \End(F_n)$ and $i\in [n]$,
    we have
    $$\sigma^2(f)(x_i)=x_i\ti{f}(x_i)^{-1}x_i=x_ix_i^{-1}f(x_i)x_i^{-1}x_i=f(x_i).$$

    We have (\ref{l7122}) since for any $f\in \jfE_r(n)$ and $i\in[n]$,
    we have
    $$[\ti{f}, x_i]=[f,x_i]^{-1}\in \Gamma_{r+1}.$$

    We prove (\ref{l7123}). Let $f\in \jfE_r(n)$.
    We have
    $$[f\ti{f},x_i]=f([\ti{f},x_i])[f,x_i]=f([f,x_i]^{-1})[f,x_i]=[f,[f,x_i]^{-1}]\in \Gamma_{2r+1}$$
    for any $i\in [n]$.
    Thus, we have
    \begin{equation}\label{l71231}
      f\ti{f}\in\jfE_{2r}(n).
    \end{equation}
    By (\ref{l7122}), we have $\ti{f}\in \jfE_r(n)$, and by (\ref{l7121}) and (\ref{l71231}),
    $$\ti{f}f=\ti{f}\ti{\ti{f}}\in\jfE_{2r}(n).$$
  \end{proof}
  For $N\geq r\geq 0$, we define an equivalence relation $\sim_{N}$ on the monoid $\jfE_r(n)$ by
  $$f\sim_{N}g \quad \overset{\text{def}}\Leftrightarrow \quad
     [f,x]\equiv[g,x]\; (\bmod\: \Gamma_{N+1}) \quad \text{ for any }x\in F_n
  $$
  for $f,g\in\jfE_r(n)$.
  Thus we have
  $$f\sim_{N} \id_{F_n} \quad \Leftrightarrow \quad [f,x]\in \Gamma_{N+1} \quad \text{ for any }x\in F_n \quad \Leftrightarrow \quad f\in \jfE_N(n).$$

  \begin{lemma}\label{lrlinv}
   Let $r\geq 1$. For $f\in \jfE_r(n)$, define $f^{R}_{N}$ and $f^{L}_{N}$ for $N\geq r+1$ inductively by
   $$
     f^{R}_{N}=
     \begin{cases}
       \ti{f} & (N=r+1)\\
       f^{R}_{N-1}\wti{f f^{R}_{N-1}} & (N\geq r+2),
     \end{cases}
   $$
   $$
     f^{L}_{N}=
     \begin{cases}
       \ti{f} & (N=r+1)\\
       \wti{f^{L}_{N-1} f}f^{L}_{N-1} & (N\geq r+2).
     \end{cases}
   $$
   Then we have
   \begin{gather*}
     f^{R}_{N}\in \jfE_r(n),\quad f f^{R}_{N}\in \jfE_N(n),\quad f^{R}_{N}\sim_{N-1} f^{R}_{N-1},\\
     f^{L}_{N}\in \jfE_r(n),\quad f^{L}_{N} f\in \jfE_N(n),\quad f^{L}_{N}\sim_{N-1} f^{L}_{N-1}.
   \end{gather*}
  \end{lemma}
  \begin{proof}
   We use induction on $N\geq r+1$.
   When $N=r+1$, by Lemma \ref{l712}, we have $\ti{f}\in\jfE_r(n)$ and $f \ti{f}\in\jfE_{2r}(n)\subset\jfE_{r+1}(n)$.
   Suppose that $f^{R}_{N-1}\in \jfE_r(n)$ satisfies $f f^{R}_{N-1}\in\jfE_{N-1}(n)$.
   By Lemma \ref{l712}, we have
   $\wti{f f^{R}_{N-1}}\in \jfE_{N-1}(n)$ and $f f^{R}_{N-1} \wti{f f^{R}_{N-1}}\in \jfE_{2N-2}(n)\subset \jfE_{N}(n)$.
   Then we have $f^{R}_N= f^{R}_{N-1}\wti{f f^{R}_{N-1}}\in \jfE_{r}(n)$ and $f f^{R}_N\in \jfE_{N}(n)$.
   Since $\wti{f f^{R}_{N-1}}\in \jfE_{N-1}(n)$, we have $f^{R}_{N}\sim_{N-1} f^{R}_{N-1}$.
   The case for $f^{L}_N$ is similar.
  \end{proof}

  \begin{proposition}\label{p711}
   For $N\geq 1$, we have a filtration of groups
   $$\jfE_1(n)/\sim_{N} \;\supset \jfE_2(n)/\sim_{N} \;\supset \cdots \supset \jfE_{N-1}(n)/\sim_{N} \;\supset \jfE_N(n)/\sim_{N}\;=1.$$
   Moreover, this is an N-series.
  \end{proposition}
  \begin{proof}
   Firstly, we show that $\jfE_r(n)/\sim_{N}$ is a group for each $r\geq 1$.
   For $f,f',g\in \jfE_r(n)$ such that $f\sim_N f'$, we can easily check that $fg\sim_N f'g$ and $gf\sim_N gf'$.
   Thus, the composition makes the set $\jfE_r(n)/\sim_{N}$ a monoid.
   For $[f]\in \jfE_r(n)/\sim_{N}$, by Lemma \ref{lrlinv}, it follows that $[f][f^{R}_N]=[f^{L}_N][f]=1 \in\jfE_r(n)/\sim_{N}$.
   Since $\jfE_r(n)/\sim_{N}$ is a monoid, we have $[f^{R}_N]=[f^{L}_N]$ and this is the inverse of $[f]$. Therefore, $\jfE_r(n)/\sim_{N}$ is a group for each $r\geq 1$.

   Since $\jfE_r(n)\supset \jfE_{r+1}(n)$, we have $\jfE_r(n)/\sim_{N} \; \supset \jfE_{r+1}(n)/\sim_{N}$.
   Secondly, we show that the descending series is an N-series.
   It suffices to show that for $f\in \jfE_r(n), g\in \jfE_s(n)$, we have
   $$
   [[f],[g]]=[f][g][f]^{-1}[g]^{-1}=[fgf^{R}_{N}g^{R}_{N}]\in\jfE_{r+s}(n)/\sim_{N}.
   $$
   Note that, by Lemma \ref{lrlinv}, we can take $f^{R}_{N}, g^{R}_{N}\in \jfE_r(n)$ such that $ff^{R}_{N}, gg^{R}_{N}\in \jfE_{N}(n)\cap\jfE_{r+s}(n)$.
   By commutator calculus, for $x\in F_n$, we have
   $$
   [fg,x]=[f,[g,x]]\;[g,x]\;[f,x]\equiv [g,x]\;[f,x] \quad (\bmod\: \Gamma_{r+s+1}),
   $$
   $$
   [g,[g^{R}_{N},x]\;[f^{R}_{N},x]]=[g,[g^{R}_{N},x]]\; {}^{[g^{R}_{N},x]}[g,[f^{R}_{N},x]]\equiv [g,[g^{R}_{N},x]] \quad (\bmod\: \Gamma_{r+s+1}).
   $$
   Similarly, we have
   $$[f^{R}_{N} g^{R}_{N},x]\equiv [g^{R}_{N},x]\;[f^{R}_{N},x]\quad (\bmod\: \Gamma_{r+s+1}),$$
   $$[f,[g^{R}_{N},x]\;[f^{R}_{N},x]]\equiv [f,[f^{R}_{N},x]]\quad (\bmod\: \Gamma_{r+s+1}).
   $$
   Thus, we have
   \begin{gather*}
    \begin{split}
      [fg,[f^{R}_{N} g^{R}_{N},x]]&\equiv
      [g,[f^{R}_{N} g^{R}_{N},x]]\;[f,[f^{R}_{N} g^{R}_{N},x]]\\
      &\equiv
      [g,[g^{R}_{N},x]\;[f^{R}_{N},x]]\;[f,[g^{R}_{N},x]\;[f^{R}_{N},x]]\\
      &\equiv
      [g,[g^{R}_{N},x]]\;[f,[f^{R}_{N},x]] \quad (\bmod\: \Gamma_{r+s+1}).
    \end{split}
   \end{gather*}
   Therefore, we have
   \begin{gather*}
    \begin{split}
      [f g f^{R}_{N} g^{R}_{N},x]&=
      [fg,[f^{R}_{N} g^{R}_{N},x]]\;[f^{R}_{N} g^{R}_{N},x]\;[fg,x]\\
      &\equiv [g,[g^{R}_{N},x]]\;[f,[f^{R}_{N},x]]\;[g^{R}_{N},x]\;[f^{R}_{N},x]\;[g,x]\;[f,x]\\
      &\equiv [g,[g^{R}_{N},x]]\;[g^{R}_{N},x]\;[g,x]\;[f,[f^{R}_{N},x]]\;[f^{R}_{N},x]\;[f,x]\\
      &=[gg^{R}_{N},x]\;[ff^{R}_{N},x]\\
      &\equiv 1 \quad (\bmod\: \Gamma_{r+s+1})
    \end{split}
   \end{gather*}
   and the proof is complete.
  \end{proof}

  For $N\geq r\geq 1$, we have a canonical projection
  $$p_{N+1}:\jfE_r(n)/\sim_{N+1}\twoheadrightarrow\jfE_r(n)/\sim_{N}.$$
  Let $\hat{\jfE}_r(n)$ denote the projective limit
  $\underset{N}{\varprojlim} (\jfE_r(n)/\sim_N)$ and  $$\pi_{N}:\hat{\jfE}_r(n)\twoheadrightarrow\jfE_r(n)/\sim_N$$
  denote the projection.
  By Proposition \ref{p711}, we have a descending series of groups
  $$\hat{\jfE}_1(n)\supset \hat{\jfE}_2(n)\supset \cdots$$
  satisfying
  $$\bigcap_{r\geq 1} \hat{\jfE}_r(n)=\{\id\}.$$

  \begin{proposition}\label{p713}
   The descending series $\hat{\jfE}_{\ast}(n):=(\hat{\jfE}_r(n))_{r\geq 1}$ is an N-series.
  \end{proposition}
  \begin{proof}
   By Proposition \ref{p711}, we have $[\jfE_{r}(n)/\sim_{N},\jfE_{s}(n)/\sim_{N}]\subset \jfE_{r+s}(n)/\sim_{N}$ for each $N>r,s$. By taking the projective limits, we have $[\hat{\jfE}_{r}(n),\hat{\jfE}_{s}(n)]\subset\hat{\jfE}_{r+s}(n)$.
  \end{proof}

  We have a graded Lie algebra $\gr(\hat{\jfE}_{\ast}(n))$ associated to the N-series $\hat{\jfE}_{\ast}(n)$.
  Let $\opegr^r(\jfE_{\ast}(n)):=\jfE_r(n)/\sim_{r+1}$ for $r\geq 1$
  and $\gr(\jfE_{\ast}(n)):=\bigoplus_{r\geq 1}\opegr^r(\jfE_{\ast}(n))$.

  \begin{proposition}\label{p714}
    We have a group isomorphism
    $$\bar{\pi}_{r+1}:\opegr^r(\hat{\jfE}_{\ast}(n))\xrightarrow{\cong}\opegr^r(\jfE_{\ast}(n))$$ induced by the projection $\pi_{r+1}:\hat{\jfE}_{r}(n)\rightarrow\opegr^r(\jfE_{\ast}(n))$.
    Therefore, $\gr(\jfE_{\ast}(n))$ is a graded Lie algebra.
  \end{proposition}
  \begin{proof}
   The projection $\pi_{r+1}$ induces $\bar{\pi}_{r+1}$ since for $f\in \hat{\jfE}_{r+1}(n)$, we have $\pi_{r+1}(f)\in \jfE_{r+1}(n)/\sim_{r+1}=1$.

   We will check that $\bar{\pi}_{r+1}$ is surjective.
   For any $f\in\jfE_r(n)$, let $\Phi(f)\in\hat{\jfE}_r(n)$ satisfy $\pi_N(\Phi(f))=[f]\in\jfE_r(n)/\sim_{N}$ for each $N>r$.
   We have $\bar{\pi}_{r+1}([\Phi(f)])=\pi_{r+1}(\Phi(f))=[f]\in\jfE_r(n)/\sim_{r+1}$. Therefore, $\bar{\pi}_{r+1}$ is surjective.

   Finally, we show that $\bar{\pi}_{r+1}$ is injective.
   Let $f\in\hat{\jfE}_r(n)$ satisfy $\bar{\pi}_{r+1}([f])=1\in \jfE_r(n)/\sim_{r+1}$ and $\pi_{N}(f)=[f_N]\in \jfE_r(n)/\sim_{N}$ for $f_N\in\jfE_r(n)$.
   Then, we have $f_{r+1}\in\jfE_{r+1}(n)$ and $f_N\sim_{r+1}f_{r+1}$ for any $N>r$.
   Therefore, we have $\pi_N(f)=[f_N]\in\jfE_{r+1}(n)/\sim_{N}$ for each $N>r$ and thus $[f]=1\in \opegr^r(\hat{\jfE}_{\ast}(n))$.
   The proof is complete.
  \end{proof}

 \subsection{Johnson homomorphism of $\End(F_n)$}\label{ss33}
  For $r\geq 1$, by using Lemma \ref{l711}, we can define a monoid homomorphism
  $$\ti{\tau}'_r:\jfE_r(n)\rightarrow\Hom(H, \Lie_{r+1}(n))$$
  by $\ti{\tau}'_r(f)(x\Gamma_2):=[f,x]\Gamma_{r+2}$ for $f\in\jfE_r(n), x\in F_n$.
  It is easily checked that the monoid homomorphism $\ti{\tau}'_r$ induces an injective group homomorphism
  $$
   \ti{\tau}_r:\opegr^r(\jfE_{\ast}(n))\hookrightarrow\Hom(H,\Lie_{r+1}(n)).
  $$
  We call it the \emph{$r$-th Johnson homomorphism} of $\End(F_n)$.

  \begin{proposition} \label{p712}
    The map $\ti{\tau}_r:\opegr^r(\jfE_{\ast}(n))\hookrightarrow\Hom(H, \Lie_{r+1}(n))$ is an abelian group isomorphism.
  \end{proposition}
  \begin{proof}
   It suffices to show that $\ti{\tau}_r$ is surjective.
   For any $\varphi\in\Hom(H,\Lie_{r+1}(n))$, we fix a representative of $\varphi(x_i\Gamma_2)\in \Lie_{r+1}(n)$ and write it $\varphi(x_i)\in \Gamma_{r+1},$ for $i\in[n]$.
   Define $\psi\in\End(F_n)$ by
   $$\psi(x_i)= \varphi(x_i)x_i \text{ for }i\in[n].$$
   It turns out that $[\psi,x]\Gamma_{r+2}=\varphi(x\Gamma_{2})\in \Lie_{r+1}(n)$ for any $x\in F_n$ by induction on the word length of $x\in F_n$.
   Therefore, we have $\ti{\tau}_r(\psi)=\varphi$ and thus the map $\ti{\tau}_r$ is surjective.
  \end{proof}
  Then we obtain the following commutative diagram:
  \begin{gather*}
   \xymatrix{
        \opegr^r(\jfA_{\ast}(n))\ar@{^{(}->}[d]_{\mathrm{inclusion}}\ar[rrd]^{\tau_r}
        \\
        \opegr^r(\jfE_{\ast}(n))\ar[rr]^{\cong}_{\ti{\tau}_r}
        &&
       \Hom(H, \Lie_{r+1}(n)).
      }
  \end{gather*}

  \begin{remark}\label{r711}
   It is well known that the Andreadakis filtration $\jfA_{\ast}(n)$ of $\Aut(F_n)$ includes the lower central series of $\IA(n)$:
   $$
     \Gamma_r(\IA(n))\subset\jfA_{r}(n).
   $$
   We have $\jfA_1(n)=\IA(n)$ by definition.
   Andreadakis \cite{Andreadakis} conjectured that
   \begin{equation}\label{Aconj}
     \jfA_r(n)=\Gamma_r(\IA(n))
   \end{equation} for all $r\geq 2,n\geq 2$.
   Andreadakis \cite{Andreadakis} ($n=3$) and Kawazumi \cite{Kawazumi} (for any $n$) showed that \eqref{Aconj} holds for $r=2$.
   Moreover, Andreadakis \cite{Andreadakis} showed that the first Johnson homomorphism $\tau_1$ of $\Aut(F_n)$ is an isomorphism.
   Therefore, we have abelian group isomorphisms
   \begin{equation}\label{bracketeq}
     \opegr^1(\IA(n))\cong \Hom(H,\Lie_2(n))\cong \opegr^1(\jfE_{\ast}(n)).
   \end{equation}
   Recently, Satoh \cite{Satoh_third} showed that \eqref{Aconj} holds for $r=3$.
   On the other hand, Bartholdi \cite{Bartholdi_e} showed that $$(\jfA_{5}(3)/ \Gamma_5(\IA(3)))\otimes \Q\cong \Q^{\oplus 3},$$ which is a counterexample of the Andreadakis conjecture.
   Now, the Andreadakis conjecture remains open for $n\gg r$.
  \end{remark}

 \subsection{The derivation Lie algebra}\label{ss34}
  By \eqref{equiv} and Proposition \ref{p712}, we have abelian group isomorphisms
  $$\opegr^r(\jfE_{\ast}(n))\cong H^{\ast}\otimes \Lie_{r+1}(n)\cong \Der_r(\Lie_{\ast}(n)) \cong T_r(n).$$
  We write $\ti{\tau}_r:\opegr^r(\jfE_{\ast}(n))\xrightarrow{\cong} \Der_r(\Lie_{\ast}(n))$ as well.
  \begin{proposition}
    The abelian group isomorphism
    $$\ti{\tau}=\bigoplus_{r\geq 1} \ti{\tau}_r:\gr(\jfE_{\ast}(n))\xrightarrow{\cong} \Der(\Lie_{\ast}(n))$$ is an isomorphism of graded Lie algebras.
  \end{proposition}
  \begin{proof}
    We only need to check that the Lie bracket of $\gr(\jfE_{\ast}(n))$ is sent to the Lie bracket of $\Der(\Lie_{\ast}(n))$.
    For $f\in\hat{\jfE}_r(n),g\in\hat{\jfE}_s(n)$ and $x\in F_n$, we have
    \begin{gather*}
      \begin{split}
          [\ti{\tau}_r([f]),\ti{\tau}_s([g])](x\Gamma_2)&=\ti{\tau}_r([f])\ti{\tau}_s([g])(x\Gamma_2)-\ti{\tau}_s([g])\ti{\tau}_r([f])(x\Gamma_2)\\ &=[f,[g,x]]\;[g,[f,x]]^{-1}=[[f,g],x]\in \Lie_{r+s+1}(n).
      \end{split}
    \end{gather*}
    On the other hand, we have $$\ti{\tau}_{r+s}([[f,g]])(x\Gamma_2)=[[f,g],x]\in \Lie_{r+s+1}(n).$$
    Therefore, $\ti{\tau}$ is an isomorphism of graded Lie algebras.
  \end{proof}

  \begin{remark}
    The tree module $\bigoplus_{r\geq 1}T_r(n)$ also has a graded Lie algebra structure which is induced by the Lie algebra structure of $\Der(\Lie_{\ast}(n))$.
    The Lie bracket $$[\cdot,\cdot]:T_r(n)\times T_s(n)\rightarrow T_{r+s}(n)$$
    is defined by the difference between two linear sums obtained by contracting the root of one of the trees and the leaves of the other tree
    $$\left[\scalebox{0.8}{$\centre{
\begingroup%
  \makeatletter%
  \providecommand\color[2][]{%
    \errmessage{(Inkscape) Color is used for the text in Inkscape, but the package 'color.sty' is not loaded}%
    \renewcommand\color[2][]{}%
  }%
  \providecommand\transparent[1]{%
    \errmessage{(Inkscape) Transparency is used (non-zero) for the text in Inkscape, but the package 'transparent.sty' is not loaded}%
    \renewcommand\transparent[1]{}%
  }%
  \providecommand\rotatebox[2]{#2}%
  \newcommand*\fsize{\dimexpr\f@size pt\relax}%
  \newcommand*\lineheight[1]{\fontsize{\fsize}{#1\fsize}\selectfont}%
  \ifx\svgwidth\undefined%
    \setlength{\unitlength}{56.12587567bp}%
    \ifx\svgscale\undefined%
      \relax%
    \else%
      \setlength{\unitlength}{\unitlength * \real{\svgscale}}%
    \fi%
  \else%
    \setlength{\unitlength}{\svgwidth}%
  \fi%
  \global\let\svgwidth\undefined%
  \global\let\svgscale\undefined%
  \makeatother%
  \begin{picture}(1,0.84365607)%
    \lineheight{1}%
    \setlength\tabcolsep{0pt}%
    \put(0,0){\includegraphics[width=\unitlength,page=1]{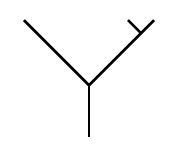}}%
    \put(-0.00269693,0.81304277){\makebox(0,0)[lt]{\lineheight{1.45000005}\smash{\begin{tabular}[t]{l}$x_{i_1}$\end{tabular}}}}%
    \put(0.36382187,0.01256243){\makebox(0,0)[lt]{\lineheight{1.45000005}\smash{\begin{tabular}[t]{l}$v_i$\end{tabular}}}}%
    \put(0.48496904,0.81433171){\makebox(0,0)[lt]{\lineheight{1.45000005}\smash{\begin{tabular}[t]{l}$x_{i_r}$\end{tabular}}}}%
    \put(0,0){\includegraphics[width=\unitlength,page=2]{treelier.pdf}}%
    \put(0.72327826,0.81623445){\makebox(0,0)[lt]{\lineheight{1.45000005}\smash{\begin{tabular}[t]{l}$x_{i_{r+1}}$\end{tabular}}}}%
  \end{picture}%
\endgroup%
}$},\scalebox{0.8}{$\centre{
\begingroup%
  \makeatletter%
  \providecommand\color[2][]{%
    \errmessage{(Inkscape) Color is used for the text in Inkscape, but the package 'color.sty' is not loaded}%
    \renewcommand\color[2][]{}%
  }%
  \providecommand\transparent[1]{%
    \errmessage{(Inkscape) Transparency is used (non-zero) for the text in Inkscape, but the package 'transparent.sty' is not loaded}%
    \renewcommand\transparent[1]{}%
  }%
  \providecommand\rotatebox[2]{#2}%
  \newcommand*\fsize{\dimexpr\f@size pt\relax}%
  \newcommand*\lineheight[1]{\fontsize{\fsize}{#1\fsize}\selectfont}%
  \ifx\svgwidth\undefined%
    \setlength{\unitlength}{61.62742285bp}%
    \ifx\svgscale\undefined%
      \relax%
    \else%
      \setlength{\unitlength}{\unitlength * \real{\svgscale}}%
    \fi%
  \else%
    \setlength{\unitlength}{\svgwidth}%
  \fi%
  \global\let\svgwidth\undefined%
  \global\let\svgscale\undefined%
  \makeatother%
  \begin{picture}(1,0.77105277)%
    \lineheight{1}%
    \setlength\tabcolsep{0pt}%
    \put(0,0){\includegraphics[width=\unitlength,page=1]{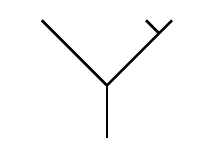}}%
    \put(-0.00245617,0.7460791){\makebox(0,0)[lt]{\lineheight{1.45000005}\smash{\begin{tabular}[t]{l}$x_{j_1}$\end{tabular}}}}%
    \put(0.41535324,0.0134059){\makebox(0,0)[lt]{\lineheight{1.45000005}\smash{\begin{tabular}[t]{l}$v_j$\end{tabular}}}}%
    \put(0.52568547,0.74360036){\makebox(0,0)[lt]{\lineheight{1.45000005}\smash{\begin{tabular}[t]{l}$x_{j_s}$\end{tabular}}}}%
    \put(0,0){\includegraphics[width=\unitlength,page=2]{treelies.pdf}}%
    \put(0.74272057,0.74533324){\makebox(0,0)[lt]{\lineheight{1.45000005}\smash{\begin{tabular}[t]{l}$x_{j_{s+1}}$\end{tabular}}}}%
  \end{picture}%
\endgroup%
}$}\right]
    =\sum_{l=1}^{s+1}\langle v_i,x_{j_l} \rangle \scalebox{0.95}{$\centre{
\begingroup%
  \makeatletter%
  \providecommand\color[2][]{%
    \errmessage{(Inkscape) Color is used for the text in Inkscape, but the package 'color.sty' is not loaded}%
    \renewcommand\color[2][]{}%
  }%
  \providecommand\transparent[1]{%
    \errmessage{(Inkscape) Transparency is used (non-zero) for the text in Inkscape, but the package 'transparent.sty' is not loaded}%
    \renewcommand\transparent[1]{}%
  }%
  \providecommand\rotatebox[2]{#2}%
  \newcommand*\fsize{\dimexpr\f@size pt\relax}%
  \newcommand*\lineheight[1]{\fontsize{\fsize}{#1\fsize}\selectfont}%
  \ifx\svgwidth\undefined%
    \setlength{\unitlength}{73.64846901bp}%
    \ifx\svgscale\undefined%
      \relax%
    \else%
      \setlength{\unitlength}{\unitlength * \real{\svgscale}}%
    \fi%
  \else%
    \setlength{\unitlength}{\svgwidth}%
  \fi%
  \global\let\svgwidth\undefined%
  \global\let\svgscale\undefined%
  \makeatother%
  \begin{picture}(1,0.94130387)%
    \lineheight{1}%
    \setlength\tabcolsep{0pt}%
    \put(0,0){\includegraphics[width=\unitlength,page=1]{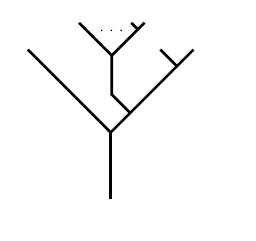}}%
    \put(-0.00092819,0.81500389){\makebox(0,0)[lt]{\lineheight{1.45000005}\smash{\begin{tabular}[t]{l}\tiny$x_{j_1}$\end{tabular}}}}%
    \put(0.3620858,0.01121777){\makebox(0,0)[lt]{\lineheight{1.45000005}\smash{\begin{tabular}[t]{l}$v_j$\end{tabular}}}}%
    \put(0.57817395,0.79790341){\makebox(0,0)[lt]{\lineheight{1.45000005}\smash{\begin{tabular}[t]{l}\tiny$x_{j_s}$\end{tabular}}}}%
    \put(0,0){\includegraphics[width=\unitlength,page=2]{treelieimr.pdf}}%
    \put(0.72090798,0.79711567){\makebox(0,0)[lt]{\lineheight{1.45000005}\smash{\begin{tabular}[t]{l}\tiny$x_{j_{s+1}}$\end{tabular}}}}%
    \put(0,0){\includegraphics[width=\unitlength,page=3]{treelieimr.pdf}}%
    \put(0.22141141,0.93085515){\makebox(0,0)[lt]{\lineheight{1.45000005}\smash{\begin{tabular}[t]{l}\tiny{$x_{i_1}$}\end{tabular}}}}%
    \put(0.41353197,0.9302944){\makebox(0,0)[lt]{\lineheight{1.45000005}\smash{\begin{tabular}[t]{l}\tiny{$x_{i_r}$}\end{tabular}}}}%
    \put(0.56210663,0.93007141){\makebox(0,0)[lt]{\lineheight{1.45000005}\smash{\begin{tabular}[t]{l}\tiny{$x_{i_{r+1}}$}\end{tabular}}}}%
    \put(0.55108431,0.45810217){\makebox(0,0)[lt]{\lineheight{1.45000005}\smash{\begin{tabular}[t]{l}\tiny$l$\end{tabular}}}}%
    \put(0.45641378,0.3845538){\makebox(0,0)[lt]{\lineheight{1.45000005}\smash{\begin{tabular}[t]{l}\tiny$1$\end{tabular}}}}%
    \put(0.72608893,0.64065021){\makebox(0,0)[lt]{\lineheight{1.45000005}\smash{\begin{tabular}[t]{l}\tiny$s$\end{tabular}}}}%
  \end{picture}%
\endgroup%
}$}
    -\sum_{l=1}^{r+1} \langle v_j,x_{i_l} \rangle \scalebox{0.95}{$\centre{
\begingroup%
  \makeatletter%
  \providecommand\color[2][]{%
    \errmessage{(Inkscape) Color is used for the text in Inkscape, but the package 'color.sty' is not loaded}%
    \renewcommand\color[2][]{}%
  }%
  \providecommand\transparent[1]{%
    \errmessage{(Inkscape) Transparency is used (non-zero) for the text in Inkscape, but the package 'transparent.sty' is not loaded}%
    \renewcommand\transparent[1]{}%
  }%
  \providecommand\rotatebox[2]{#2}%
  \newcommand*\fsize{\dimexpr\f@size pt\relax}%
  \newcommand*\lineheight[1]{\fontsize{\fsize}{#1\fsize}\selectfont}%
  \ifx\svgwidth\undefined%
    \setlength{\unitlength}{73.32424933bp}%
    \ifx\svgscale\undefined%
      \relax%
    \else%
      \setlength{\unitlength}{\unitlength * \real{\svgscale}}%
    \fi%
  \else%
    \setlength{\unitlength}{\svgwidth}%
  \fi%
  \global\let\svgwidth\undefined%
  \global\let\svgscale\undefined%
  \makeatother%
  \begin{picture}(1,0.94381457)%
    \lineheight{1}%
    \setlength\tabcolsep{0pt}%
    \put(0,0){\includegraphics[width=\unitlength,page=1]{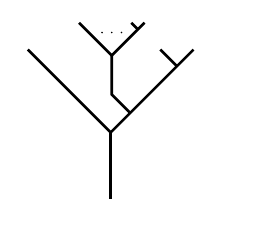}}%
    \put(-0.00093229,0.81695614){\makebox(0,0)[lt]{\lineheight{1.45000005}\smash{\begin{tabular}[t]{l}\tiny$x_{i_1}$\end{tabular}}}}%
    \put(0.36368684,0.00961589){\makebox(0,0)[lt]{\lineheight{1.45000005}\smash{\begin{tabular}[t]{l}$v_i$\end{tabular}}}}%
    \put(0.58073047,0.79978004){\makebox(0,0)[lt]{\lineheight{1.45000005}\smash{\begin{tabular}[t]{l}\tiny$x_{i_r}$\end{tabular}}}}%
    \put(0,0){\includegraphics[width=\unitlength,page=2]{treelieims.pdf}}%
    \put(0.72409563,0.79898882){\makebox(0,0)[lt]{\lineheight{1.45000005}\smash{\begin{tabular}[t]{l}\tiny$x_{i_{r+1}}$\end{tabular}}}}%
    \put(0,0){\includegraphics[width=\unitlength,page=3]{treelieims.pdf}}%
    \put(0.22239043,0.93331966){\makebox(0,0)[lt]{\lineheight{1.45000005}\smash{\begin{tabular}[t]{l}\tiny{$x_{j_1}$}\end{tabular}}}}%
    \put(0.41536049,0.93275642){\makebox(0,0)[lt]{\lineheight{1.45000005}\smash{\begin{tabular}[t]{l}\tiny{$x_{j_s}$}\end{tabular}}}}%
    \put(0.56459211,0.93253245){\makebox(0,0)[lt]{\lineheight{1.45000005}\smash{\begin{tabular}[t]{l}\tiny{$x_{j_{s+1}}$}\end{tabular}}}}%
    \put(0.55205694,0.45272065){\makebox(0,0)[lt]{\lineheight{1.45000005}\smash{\begin{tabular}[t]{l}\tiny$l$\end{tabular}}}}%
    \put(0.45696781,0.37884707){\makebox(0,0)[lt]{\lineheight{1.45000005}\smash{\begin{tabular}[t]{l}\tiny$1$\end{tabular}}}}%
    \put(0.72783539,0.63607587){\makebox(0,0)[lt]{\lineheight{1.45000005}\smash{\begin{tabular}[t]{l}\tiny$r$\end{tabular}}}}%
  \end{picture}%
\endgroup%
}$}.
    $$
  \end{remark}


\section{Action of $\gr(\jfE_{\ast}(n))$ on $B_{d}(n)$}\label{s4}
 We defined the bracket maps \eqref{bracket} in \cite{Mai1}.
 In this section, we extend them to linear maps
 $$[\cdot,\cdot]:B_{d,k}(n) \otimes \opegr^r(\jfE_{\ast}(n))\rightarrow B_{d,k+r}(n).$$

 In Section \ref{ss41}, we state Theorem \ref{th731}, which we use to obtain the extended bracket map.
 In Section \ref{ss42}, we extend the category $\A$ to a category $\A^{L}$, which includes a Lie algebra structure besides the Hopf algebra structure in $\A$.
 In Section \ref{ss425}, we observe some relations for morphisms of $\A^{L}$. By using these relations, we prove Theorem \ref{th731} in Section \ref{ss43}.

 \subsection{Bracket map $[\cdot,\cdot]: B_{d,k}(n)\otimes \opegr^r(\jfE_{\ast}(n))\rightarrow B_{d,k+r}(n)$}\label{ss41}
  We have a right $\End(F_n)$-action on $A_d(n)$ by letting
  $$u\cdot g:=A_d(g)(u)$$
  for $u\in A_d(n), g\in \End(F_n)$.
  We define
  \begin{equation}\label{eqbraE}
    [\cdot,\cdot]:A_d(n)\times\End(F_n)\rightarrow A_d(n)
  \end{equation}
  by $[u,g]:=u\cdot g -u$ for $u\in A_d(n),g\in \End(F_n)$, which we call the \emph{bracket map}.

  \begin{theorem}\label{th731}
   The N-series $\hat{\jfE}_{\ast}(n)$ acts on the right on the filtered vector space $A_d(n)$.
   That is, we have
   $$[A_{d,k}(n),\jfE_r(n)]\subset A_{d,k+r}(n)$$
   for any $r\geq 1$.
  \end{theorem}
  Note that we have $[A_{d,k}(n),\Gamma_{r}(\IA(n))]\subset A_{d,k+r}(n)$ (see Lemma 5.7 in \cite{Mai1}).
  We will prove Theorem \ref{th731} in Section \ref{ss43}.

  By using Theorem \ref{th731}, we can extend the bracket map $$[\cdot,\cdot]:B_{d,k}(n)\otimes\opegr^r(\IA(n))\rightarrow B_{d,k+r}(n)$$
  to $\opegr^r(\jfE_{\ast}(n))$.
  \begin{corollary}
   Let $r\geq 1$.
   The bracket map (\ref{eqbraE}) induces a $\K$-linear map
   \begin{equation*}
     [\cdot,\cdot]:B_{d,k}(n)\otimes\opegr^r(\jfE_{\ast}(n))\rightarrow B_{d,k+r}(n).
   \end{equation*}
  \end{corollary}

  We also extend the $\GL(n;\Z)$-module map
  $$\beta_{d,k}^r:\opegr^r(\IA(n))\rightarrow\Hom(B_{d,k}(n),B_{d,k+r}(n))$$ defined by $\beta_{d,k}^r(g)(u)=[u,g]$ for $g\in\opegr^r(\IA(n)),u\in B_{d,k}(n)$
  to a group homomorphism
  $$
    \ti{\beta}_{d,k}^r:\opegr^r(\jfE_{\ast}(n))\rightarrow\Hom(B_{d,k}(n),B_{d,k+r}(n)).
  $$
  \begin{remark}
   The right action of the N-series $\hat{\jfE}_{\ast}(n)$ on $A_d(n)$ induces an action of the graded Lie algebra $\gr(\jfE_{\ast}(n))$ on the graded vector space $B_{d}(n)$:
   $$\gr(\jfE_{\ast}(n))\xrightarrow{\cong}\gr(\hat{\jfE}_{\ast}(n))\rightarrow \bigoplus_{r\geq 1}\End_r(B_{d}(n)),$$
   which is given by the group homomorphisms $\ti{\beta}_{d,k}^r$.
   This induced action can be regarded as an action of the derivation Lie algebra $\Der(\Lie_{\ast}(n))$ on the graded vector space $B_d(n)$ by the identification in Section \ref{ss34}.
  \end{remark}

 \subsection{The category $\A^{L}$ of extended Jacobi diagrams in handlebodies}\label{ss42}
  In order to prove Theorem \ref{th731}, we extend the category $\A$ of Jacobi diagrams in handlebodies to another category $\A^{L}$.
  In Appendix \ref{sA}, we consider an expected presentation of the category $\A^{L}$.

  Construct the category $\A^{L}$ as follows.
  The set of objects of $\A^{L}$ is the free monoid generated by two objects $H$ and $L$, where multiplication is denoted by $\otimes$.
  The category $\A^{L}$ includes the category $\A$ as a full subcategory with the free monoid generated by $H$ as the set of objects.
  (On the other hand, the full subcategory with the free monoid generated by $L$ is isomorphic to a category in \cite{Hinich}. See Remark \ref{casimirlie}.)
  In the category $\A^{L}$, we consider diagrams that are obtained from Jacobi diagrams in handlebodies by attaching univalent vertices of the Jacobi diagrams to the bottom line $l$ and the upper line $l'$.
  \begin{example}\label{e722}
    Here is a morphism in $\A^{L}(H\otimes L\otimes H\otimes L\otimes H,H\otimes L^{\otimes 2}\otimes H)$:
    $$
\begingroup%
  \makeatletter%
  \providecommand\color[2][]{%
    \errmessage{(Inkscape) Color is used for the text in Inkscape, but the package 'color.sty' is not loaded}%
    \renewcommand\color[2][]{}%
  }%
  \providecommand\transparent[1]{%
    \errmessage{(Inkscape) Transparency is used (non-zero) for the text in Inkscape, but the package 'transparent.sty' is not loaded}%
    \renewcommand\transparent[1]{}%
  }%
  \providecommand\rotatebox[2]{#2}%
  \newcommand*\fsize{\dimexpr\f@size pt\relax}%
  \newcommand*\lineheight[1]{\fontsize{\fsize}{#1\fsize}\selectfont}%
  \ifx\svgwidth\undefined%
    \setlength{\unitlength}{135.74999889bp}%
    \ifx\svgscale\undefined%
      \relax%
    \else%
      \setlength{\unitlength}{\unitlength * \real{\svgscale}}%
    \fi%
  \else%
    \setlength{\unitlength}{\svgwidth}%
  \fi%
  \global\let\svgwidth\undefined%
  \global\let\svgscale\undefined%
  \makeatother%
  \begin{picture}(1,0.80236043)%
    \lineheight{1}%
    \setlength\tabcolsep{0pt}%
    \put(0,0){\includegraphics[width=\unitlength,page=1]{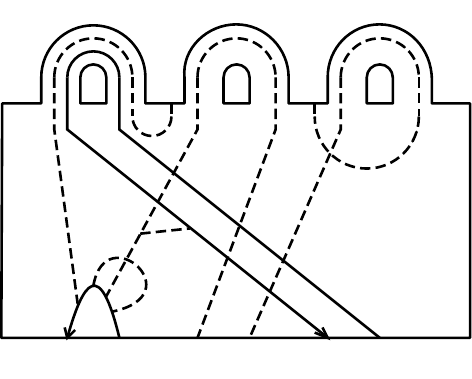}}%
    \put(0.16851353,0.77639492){\makebox(0,0)[lt]{\lineheight{1.45000005}\smash{\begin{tabular}[t]{l}$H$\end{tabular}}}}%
    \put(0.47238099,0.7781124){\makebox(0,0)[lt]{\lineheight{1.45000005}\smash{\begin{tabular}[t]{l}$H$\end{tabular}}}}%
    \put(0.77453083,0.77968549){\makebox(0,0)[lt]{\lineheight{1.45000005}\smash{\begin{tabular}[t]{l}$H$\end{tabular}}}}%
    \put(0.16493411,0.00792232){\makebox(0,0)[lt]{\lineheight{1.45000005}\smash{\begin{tabular}[t]{l}$H$\end{tabular}}}}%
    \put(0.72271724,0.00620486){\makebox(0,0)[lt]{\lineheight{1.45000005}\smash{\begin{tabular}[t]{l}$H$\end{tabular}}}}%
    \put(0.32278366,0.6354006){\makebox(0,0)[lt]{\lineheight{1.45000005}\smash{\begin{tabular}[t]{l}$L$\end{tabular}}}}%
    \put(0.62665111,0.6354006){\makebox(0,0)[lt]{\lineheight{1.45000005}\smash{\begin{tabular}[t]{l}$L$\end{tabular}}}}%
    \put(0.37459732,0.00691141){\makebox(0,0)[lt]{\lineheight{1.45000005}\smash{\begin{tabular}[t]{l}$L$\end{tabular}}}}%
    \put(0.50398694,0.00519394){\makebox(0,0)[lt]{\lineheight{1.45000005}\smash{\begin{tabular}[t]{l}$L$\end{tabular}}}}%
  \end{picture}%
\endgroup%
$$
  \end{example}
  As depicted in Figure \ref{source}, the objects $H$ and $L$ in the source of a morphism of $\A^{L}$ correspond to a handle of the handlebody and a univalent vertex attached to the upper line $l'$, respectively.
  \begin{figure}[h]
\begingroup%
  \makeatletter%
  \providecommand\color[2][]{%
    \errmessage{(Inkscape) Color is used for the text in Inkscape, but the package 'color.sty' is not loaded}%
    \renewcommand\color[2][]{}%
  }%
  \providecommand\transparent[1]{%
    \errmessage{(Inkscape) Transparency is used (non-zero) for the text in Inkscape, but the package 'transparent.sty' is not loaded}%
    \renewcommand\transparent[1]{}%
  }%
  \providecommand\rotatebox[2]{#2}%
  \newcommand*\fsize{\dimexpr\f@size pt\relax}%
  \newcommand*\lineheight[1]{\fontsize{\fsize}{#1\fsize}\selectfont}%
  \ifx\svgwidth\undefined%
    \setlength{\unitlength}{53.85826772bp}%
    \ifx\svgscale\undefined%
      \relax%
    \else%
      \setlength{\unitlength}{\unitlength * \real{\svgscale}}%
    \fi%
  \else%
    \setlength{\unitlength}{\svgwidth}%
  \fi%
  \global\let\svgwidth\undefined%
  \global\let\svgscale\undefined%
  \makeatother%
  \begin{picture}(1,0.78947368)%
    \lineheight{1}%
    \setlength\tabcolsep{0pt}%
    \put(0,0){\includegraphics[width=\unitlength,page=1]{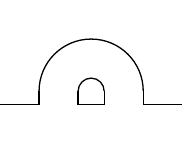}}%
    \put(0.41160563,0.66237604){\makebox(0,0)[lt]{\lineheight{1.45000005}\smash{\begin{tabular}[t]{l}$H$\end{tabular}}}}%
  \end{picture}%
\endgroup%
\;,\quad 
\begingroup%
  \makeatletter%
  \providecommand\color[2][]{%
    \errmessage{(Inkscape) Color is used for the text in Inkscape, but the package 'color.sty' is not loaded}%
    \renewcommand\color[2][]{}%
  }%
  \providecommand\transparent[1]{%
    \errmessage{(Inkscape) Transparency is used (non-zero) for the text in Inkscape, but the package 'transparent.sty' is not loaded}%
    \renewcommand\transparent[1]{}%
  }%
  \providecommand\rotatebox[2]{#2}%
  \newcommand*\fsize{\dimexpr\f@size pt\relax}%
  \newcommand*\lineheight[1]{\fontsize{\fsize}{#1\fsize}\selectfont}%
  \ifx\svgwidth\undefined%
    \setlength{\unitlength}{44.9999983bp}%
    \ifx\svgscale\undefined%
      \relax%
    \else%
      \setlength{\unitlength}{\unitlength * \real{\svgscale}}%
    \fi%
  \else%
    \setlength{\unitlength}{\svgwidth}%
  \fi%
  \global\let\svgwidth\undefined%
  \global\let\svgscale\undefined%
  \makeatother%
  \begin{picture}(1,0.83307344)%
    \lineheight{1}%
    \setlength\tabcolsep{0pt}%
    \put(0,0){\includegraphics[width=\unitlength,page=1]{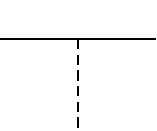}}%
    \put(0.42487081,0.76467071){\makebox(0,0)[lt]{\lineheight{1.45000005}\smash{\begin{tabular}[t]{l}$L$\end{tabular}}}}%
  \end{picture}%
\endgroup%

     \caption{Source of a morphism}
     \label{source}
  \end{figure}
  As depicted in Figure \ref{target}, the objects $H$ and $L$ in the target of a morphism of $\A^{L}$ correspond to an arc component mapped into the handlebody and a univalent vertex attached to the bottom line $l$, respectively.
  \begin{figure}[ht]
\begingroup%
  \makeatletter%
  \providecommand\color[2][]{%
    \errmessage{(Inkscape) Color is used for the text in Inkscape, but the package 'color.sty' is not loaded}%
    \renewcommand\color[2][]{}%
  }%
  \providecommand\transparent[1]{%
    \errmessage{(Inkscape) Transparency is used (non-zero) for the text in Inkscape, but the package 'transparent.sty' is not loaded}%
    \renewcommand\transparent[1]{}%
  }%
  \providecommand\rotatebox[2]{#2}%
  \newcommand*\fsize{\dimexpr\f@size pt\relax}%
  \newcommand*\lineheight[1]{\fontsize{\fsize}{#1\fsize}\selectfont}%
  \ifx\svgwidth\undefined%
    \setlength{\unitlength}{44.9999983bp}%
    \ifx\svgscale\undefined%
      \relax%
    \else%
      \setlength{\unitlength}{\unitlength * \real{\svgscale}}%
    \fi%
  \else%
    \setlength{\unitlength}{\svgwidth}%
  \fi%
  \global\let\svgwidth\undefined%
  \global\let\svgscale\undefined%
  \makeatother%
  \begin{picture}(1,0.94138869)%
    \lineheight{1}%
    \setlength\tabcolsep{0pt}%
    \put(0,0){\includegraphics[width=\unitlength,page=1]{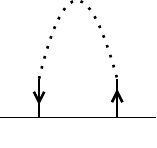}}%
    \put(0.41394335,0.01566839){\makebox(0,0)[lt]{\lineheight{1.45000005}\smash{\begin{tabular}[t]{l}$H$\end{tabular}}}}%
  \end{picture}%
\endgroup%
\;,\quad 
\begingroup%
  \makeatletter%
  \providecommand\color[2][]{%
    \errmessage{(Inkscape) Color is used for the text in Inkscape, but the package 'color.sty' is not loaded}%
    \renewcommand\color[2][]{}%
  }%
  \providecommand\transparent[1]{%
    \errmessage{(Inkscape) Transparency is used (non-zero) for the text in Inkscape, but the package 'transparent.sty' is not loaded}%
    \renewcommand\transparent[1]{}%
  }%
  \providecommand\rotatebox[2]{#2}%
  \newcommand*\fsize{\dimexpr\f@size pt\relax}%
  \newcommand*\lineheight[1]{\fontsize{\fsize}{#1\fsize}\selectfont}%
  \ifx\svgwidth\undefined%
    \setlength{\unitlength}{45bp}%
    \ifx\svgscale\undefined%
      \relax%
    \else%
      \setlength{\unitlength}{\unitlength * \real{\svgscale}}%
    \fi%
  \else%
    \setlength{\unitlength}{\svgwidth}%
  \fi%
  \global\let\svgwidth\undefined%
  \global\let\svgscale\undefined%
  \makeatother%
  \begin{picture}(1,0.77032093)%
    \lineheight{1}%
    \setlength\tabcolsep{0pt}%
    \put(0,0){\includegraphics[width=\unitlength,page=1]{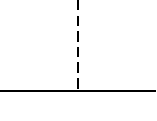}}%
    \put(0.42193663,0.01566839){\makebox(0,0)[lt]{\lineheight{1.45000005}\smash{\begin{tabular}[t]{l}$L$\end{tabular}}}}%
  \end{picture}%
\endgroup%

     \caption{Target of a morphism}
     \label{target}
  \end{figure}

  In the category $\A^{L}$, the object $H$ is considered as a Hopf algebra and $L$ is considered as a Lie algebra. See Section \ref{ss425} and Appendix \ref{sA}.

  To define morphisms of the category $\A^{L}$ precisely, we give the following definition.
  \begin{definition}
    For a finite set $T$, an \emph{$(X_m,T)$-diagram} is a quadruple $(D,V,f,g)$, where
    \begin{itemize}
      \item $D$ is a vertex-oriented uni-trivalent graph such that each connected component has at least one univalent vertex,
      \item $V$ is a subset of $\partial D=\{\text{univalent vertices of }D\}$,
      \item $f$ is an embedding of $V$ into the interior of $X_m$,
      \item $g$ is a bijection from $T$ to $\partial D\setminus V$.
    \end{itemize}
  \end{definition}
  Note that an $(X_m,\emptyset)$-diagram is a Jacobi diagram on $X_m$.

  For an object $w=H^{\otimes m_1}\otimes L^{\otimes n_1}\otimes\cdots\otimes H^{\otimes m_r}\otimes L^{\otimes n_r}\in\A^{L}$,
  let $m:=\sum_{i=1}^r m_i$ and $n:=\sum_{i=1}^r n_i$.
  For $p\geq 0$, let $[p]^{+}:=\{1^{+},\cdots,p^{+}\}$ and $[p]^{-}:=\{1^{-},\cdots,p^{-}\}$ be two copies of $[p]$.
  \begin{definition}
    For objects $w=H^{\otimes m_1}\otimes L^{\otimes n_1}\otimes\cdots\otimes H^{\otimes m_r}\otimes L^{\otimes n_r}\in \A^{L}$ and
    $w'=H^{\otimes m'_1}\otimes L^{\otimes n'_1}\otimes\cdots\otimes H^{\otimes m'_s}\otimes L^{\otimes n'_s}\in \A^{L}$,
    a \emph{$(w,w')$-diagram} consists of
    \begin{itemize}
      \item an $(X_{m'},[n]^{+}\sqcup [n']^{-})$-diagram $(D,V,f,g)$ such that each connected component of $D$ has at least one univalent vertex in $V\cup g([n']^{-})$
      \item a map $\varphi:X_{m'}\cup D\rightarrow U_m$ such that
       \begin{enumerate}
         \item the pair (the empty set $\emptyset$, the restriction $\varphi\mid _{X_{m'}}$) is an $(m,m')$-Jacobi diagram;
         that is, $\varphi$ maps $X_m'$ into $U_m$ in such a way that endpoints of $X_m'$ are arranged in the bottom line $l$ from left to right,
         \item $g([n]^{+})$ is mapped into $l'$ so that the corresponding object in $\A^{L}$ with respect to Figure \ref{source} will be $w$ when we look at the top line $l'$ from left to right,
         \item $g([n']^{-})$ is mapped into $l$ so that the corresponding object in $\A^{L}$ with respect to Figure \ref{target} will be $w'$ when we look at the bottom line $l$ from left to right.
       \end{enumerate}
    \end{itemize}
  \end{definition}
  We identify two $(w,w')$-diagrams if they are homotopic in $U_m$ relative to the endpoints of $X_m'\cup D$.
  In what follows, we simply write $D$ for a $(w,w')$-diagram.
  For objects $w$ and $w'$, the hom-set $\A^{L}(w,w')$ is the $\K$-vector space spanned by $(w,w')$-diagrams modulo the STU, AS and IHX relations.

  The composition of $\A^{L}$ is defined in a similar way to that of the category $\A$. We can define a square diagram for an $(w,w')$-diagram similarly.
  Let $D$ be a diagram in $\A^{L}(w, w')$ and $D'$ a diagram in $\A^{L}(w', w'')$.
  Deform $D'$ to have only the parallel copies of the handle cores in each handle.
  Then the composition $D'\circ D$ is a diagram obtained by stacking the cabling of $D$ on top of the square presentation of $D'$.
  \begin{example}\label{e723}
     For $D=\centre{
\begingroup%
  \makeatletter%
  \providecommand\color[2][]{%
    \errmessage{(Inkscape) Color is used for the text in Inkscape, but the package 'color.sty' is not loaded}%
    \renewcommand\color[2][]{}%
  }%
  \providecommand\transparent[1]{%
    \errmessage{(Inkscape) Transparency is used (non-zero) for the text in Inkscape, but the package 'transparent.sty' is not loaded}%
    \renewcommand\transparent[1]{}%
  }%
  \providecommand\rotatebox[2]{#2}%
  \newcommand*\fsize{\dimexpr\f@size pt\relax}%
  \newcommand*\lineheight[1]{\fontsize{\fsize}{#1\fsize}\selectfont}%
  \ifx\svgwidth\undefined%
    \setlength{\unitlength}{76.25811469bp}%
    \ifx\svgscale\undefined%
      \relax%
    \else%
      \setlength{\unitlength}{\unitlength * \real{\svgscale}}%
    \fi%
  \else%
    \setlength{\unitlength}{\svgwidth}%
  \fi%
  \global\let\svgwidth\undefined%
  \global\let\svgscale\undefined%
  \makeatother%
  \begin{picture}(1,0.67661581)%
    \lineheight{1}%
    \setlength\tabcolsep{0pt}%
    \put(0,0){\includegraphics[width=\unitlength,page=1]{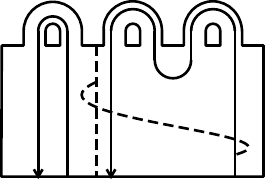}}%
  \end{picture}%
\endgroup%
}$ and $D'=\centre{
\begingroup%
  \makeatletter%
  \providecommand\color[2][]{%
    \errmessage{(Inkscape) Color is used for the text in Inkscape, but the package 'color.sty' is not loaded}%
    \renewcommand\color[2][]{}%
  }%
  \providecommand\transparent[1]{%
    \errmessage{(Inkscape) Transparency is used (non-zero) for the text in Inkscape, but the package 'transparent.sty' is not loaded}%
    \renewcommand\transparent[1]{}%
  }%
  \providecommand\rotatebox[2]{#2}%
  \newcommand*\fsize{\dimexpr\f@size pt\relax}%
  \newcommand*\lineheight[1]{\fontsize{\fsize}{#1\fsize}\selectfont}%
  \ifx\svgwidth\undefined%
    \setlength{\unitlength}{64.20720327bp}%
    \ifx\svgscale\undefined%
      \relax%
    \else%
      \setlength{\unitlength}{\unitlength * \real{\svgscale}}%
    \fi%
  \else%
    \setlength{\unitlength}{\svgwidth}%
  \fi%
  \global\let\svgwidth\undefined%
  \global\let\svgscale\undefined%
  \makeatother%
  \begin{picture}(1,0.86996657)%
    \lineheight{1}%
    \setlength\tabcolsep{0pt}%
    \put(0,0){\includegraphics[width=\unitlength,page=1]{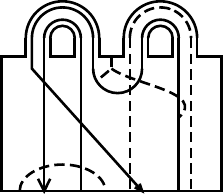}}%
  \end{picture}%
\endgroup%
}$,
     the composition $D'\circ D$ is $\centre{
\begingroup%
  \makeatletter%
  \providecommand\color[2][]{%
    \errmessage{(Inkscape) Color is used for the text in Inkscape, but the package 'color.sty' is not loaded}%
    \renewcommand\color[2][]{}%
  }%
  \providecommand\transparent[1]{%
    \errmessage{(Inkscape) Transparency is used (non-zero) for the text in Inkscape, but the package 'transparent.sty' is not loaded}%
    \renewcommand\transparent[1]{}%
  }%
  \providecommand\rotatebox[2]{#2}%
  \newcommand*\fsize{\dimexpr\f@size pt\relax}%
  \newcommand*\lineheight[1]{\fontsize{\fsize}{#1\fsize}\selectfont}%
  \ifx\svgwidth\undefined%
    \setlength{\unitlength}{103.00305813bp}%
    \ifx\svgscale\undefined%
      \relax%
    \else%
      \setlength{\unitlength}{\unitlength * \real{\svgscale}}%
    \fi%
  \else%
    \setlength{\unitlength}{\svgwidth}%
  \fi%
  \global\let\svgwidth\undefined%
  \global\let\svgscale\undefined%
  \makeatother%
  \begin{picture}(1,0.92273876)%
    \lineheight{1}%
    \setlength\tabcolsep{0pt}%
    \put(0,0){\includegraphics[width=\unitlength,page=1]{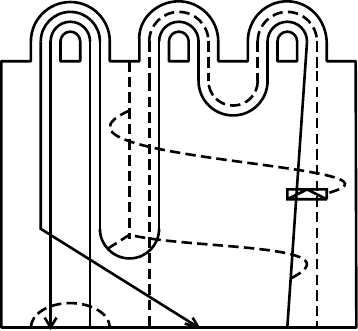}}%
  \end{picture}%
\endgroup%
}$, where the box notation represents a linear sum of Jacobi diagrams. (See \cite{HM_k} and \cite{Mai1} for the definition of the box notation.)
  \end{example}

  The identity morphism $\id_{H^{\otimes m_1}\otimes L^{\otimes n_1}\otimes \cdots \otimes H^{\otimes m_r}\otimes L^{\otimes n_r}}$ is the following diagram:
  \begin{figure}[ht]
\begingroup%
  \makeatletter%
  \providecommand\color[2][]{%
    \errmessage{(Inkscape) Color is used for the text in Inkscape, but the package 'color.sty' is not loaded}%
    \renewcommand\color[2][]{}%
  }%
  \providecommand\transparent[1]{%
    \errmessage{(Inkscape) Transparency is used (non-zero) for the text in Inkscape, but the package 'transparent.sty' is not loaded}%
    \renewcommand\transparent[1]{}%
  }%
  \providecommand\rotatebox[2]{#2}%
  \newcommand*\fsize{\dimexpr\f@size pt\relax}%
  \newcommand*\lineheight[1]{\fontsize{\fsize}{#1\fsize}\selectfont}%
  \ifx\svgwidth\undefined%
    \setlength{\unitlength}{168.37795708bp}%
    \ifx\svgscale\undefined%
      \relax%
    \else%
      \setlength{\unitlength}{\unitlength * \real{\svgscale}}%
    \fi%
  \else%
    \setlength{\unitlength}{\svgwidth}%
  \fi%
  \global\let\svgwidth\undefined%
  \global\let\svgscale\undefined%
  \makeatother%
  \begin{picture}(1,0.35353534)%
    \lineheight{1}%
    \setlength\tabcolsep{0pt}%
    \put(0,0){\includegraphics[width=\unitlength,page=1]{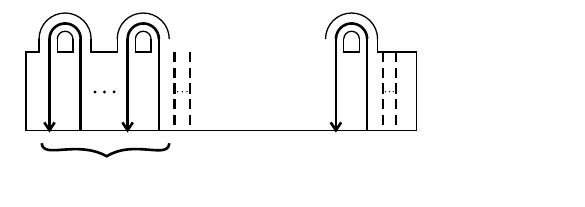}}%
    \put(0.14056891,0.04447137){\makebox(0,0)[lt]{\lineheight{1.45000005}\smash{\begin{tabular}[t]{l}$m_1$\end{tabular}}}}%
    \put(0,0){\includegraphics[width=\unitlength,page=2]{id.pdf}}%
    \put(0.27970622,0.0456565){\makebox(0,0)[lt]{\lineheight{1.45000005}\smash{\begin{tabular}[t]{l}$n_1$\end{tabular}}}}%
    \put(0.63317602,0.04639921){\makebox(0,0)[lt]{\lineheight{1.45000005}\smash{\begin{tabular}[t]{l}$n_r$\end{tabular}}}}%
    \put(0,0){\includegraphics[width=\unitlength,page=3]{id.pdf}}%
    \put(0.4907727,0.04689493){\makebox(0,0)[lt]{\lineheight{1.45000005}\smash{\begin{tabular}[t]{l}$m_r$\end{tabular}}}}%
  \end{picture}%
\endgroup%
.
  \end{figure}

  We can naturally extend the linear symmetric strict monoidal structure of $\A$ to the category $\A^{L}$, where the tensor product is defined to be the juxtaposition of the handlebodies.

  Note that the symmetries in $\A^{L}$ are determined by
  \begin{gather*}
   \begin{split}
     P_{H,H}&=\centre{
\begingroup%
  \makeatletter%
  \providecommand\color[2][]{%
    \errmessage{(Inkscape) Color is used for the text in Inkscape, but the package 'color.sty' is not loaded}%
    \renewcommand\color[2][]{}%
  }%
  \providecommand\transparent[1]{%
    \errmessage{(Inkscape) Transparency is used (non-zero) for the text in Inkscape, but the package 'transparent.sty' is not loaded}%
    \renewcommand\transparent[1]{}%
  }%
  \providecommand\rotatebox[2]{#2}%
  \newcommand*\fsize{\dimexpr\f@size pt\relax}%
  \newcommand*\lineheight[1]{\fontsize{\fsize}{#1\fsize}\selectfont}%
  \ifx\svgwidth\undefined%
    \setlength{\unitlength}{38.82448626bp}%
    \ifx\svgscale\undefined%
      \relax%
    \else%
      \setlength{\unitlength}{\unitlength * \real{\svgscale}}%
    \fi%
  \else%
    \setlength{\unitlength}{\svgwidth}%
  \fi%
  \global\let\svgwidth\undefined%
  \global\let\svgscale\undefined%
  \makeatother%
  \begin{picture}(1,1.27838214)%
    \lineheight{1}%
    \setlength\tabcolsep{0pt}%
    \put(0,0){\includegraphics[width=\unitlength,page=1]{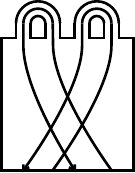}}%
  \end{picture}%
\endgroup%
}:H\otimes H\rightarrow H\otimes H,\quad
     P_{H,L}=\centre{
\begingroup%
  \makeatletter%
  \providecommand\color[2][]{%
    \errmessage{(Inkscape) Color is used for the text in Inkscape, but the package 'color.sty' is not loaded}%
    \renewcommand\color[2][]{}%
  }%
  \providecommand\transparent[1]{%
    \errmessage{(Inkscape) Transparency is used (non-zero) for the text in Inkscape, but the package 'transparent.sty' is not loaded}%
    \renewcommand\transparent[1]{}%
  }%
  \providecommand\rotatebox[2]{#2}%
  \newcommand*\fsize{\dimexpr\f@size pt\relax}%
  \newcommand*\lineheight[1]{\fontsize{\fsize}{#1\fsize}\selectfont}%
  \ifx\svgwidth\undefined%
    \setlength{\unitlength}{38.82448626bp}%
    \ifx\svgscale\undefined%
      \relax%
    \else%
      \setlength{\unitlength}{\unitlength * \real{\svgscale}}%
    \fi%
  \else%
    \setlength{\unitlength}{\svgwidth}%
  \fi%
  \global\let\svgwidth\undefined%
  \global\let\svgscale\undefined%
  \makeatother%
  \begin{picture}(1,1.27838214)%
    \lineheight{1}%
    \setlength\tabcolsep{0pt}%
    \put(0,0){\includegraphics[width=\unitlength,page=1]{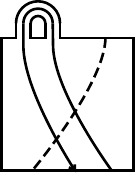}}%
  \end{picture}%
\endgroup%
}:H\otimes L\rightarrow L\otimes H, \\
     P_{L,H}&=\centre{
\begingroup%
  \makeatletter%
  \providecommand\color[2][]{%
    \errmessage{(Inkscape) Color is used for the text in Inkscape, but the package 'color.sty' is not loaded}%
    \renewcommand\color[2][]{}%
  }%
  \providecommand\transparent[1]{%
    \errmessage{(Inkscape) Transparency is used (non-zero) for the text in Inkscape, but the package 'transparent.sty' is not loaded}%
    \renewcommand\transparent[1]{}%
  }%
  \providecommand\rotatebox[2]{#2}%
  \newcommand*\fsize{\dimexpr\f@size pt\relax}%
  \newcommand*\lineheight[1]{\fontsize{\fsize}{#1\fsize}\selectfont}%
  \ifx\svgwidth\undefined%
    \setlength{\unitlength}{38.82448626bp}%
    \ifx\svgscale\undefined%
      \relax%
    \else%
      \setlength{\unitlength}{\unitlength * \real{\svgscale}}%
    \fi%
  \else%
    \setlength{\unitlength}{\svgwidth}%
  \fi%
  \global\let\svgwidth\undefined%
  \global\let\svgscale\undefined%
  \makeatother%
  \begin{picture}(1,1.27838213)%
    \lineheight{1}%
    \setlength\tabcolsep{0pt}%
    \put(0,0){\includegraphics[width=\unitlength,page=1]{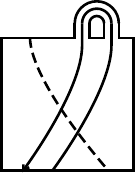}}%
  \end{picture}%
\endgroup%
}:L\otimes H\rightarrow H\otimes L,\quad
     P_{L,L}=\centre{
\begingroup%
  \makeatletter%
  \providecommand\color[2][]{%
    \errmessage{(Inkscape) Color is used for the text in Inkscape, but the package 'color.sty' is not loaded}%
    \renewcommand\color[2][]{}%
  }%
  \providecommand\transparent[1]{%
    \errmessage{(Inkscape) Transparency is used (non-zero) for the text in Inkscape, but the package 'transparent.sty' is not loaded}%
    \renewcommand\transparent[1]{}%
  }%
  \providecommand\rotatebox[2]{#2}%
  \newcommand*\fsize{\dimexpr\f@size pt\relax}%
  \newcommand*\lineheight[1]{\fontsize{\fsize}{#1\fsize}\selectfont}%
  \ifx\svgwidth\undefined%
    \setlength{\unitlength}{38.82448626bp}%
    \ifx\svgscale\undefined%
      \relax%
    \else%
      \setlength{\unitlength}{\unitlength * \real{\svgscale}}%
    \fi%
  \else%
    \setlength{\unitlength}{\svgwidth}%
  \fi%
  \global\let\svgwidth\undefined%
  \global\let\svgscale\undefined%
  \makeatother%
  \begin{picture}(1,0.99999996)%
    \lineheight{1}%
    \setlength\tabcolsep{0pt}%
    \put(0,0){\includegraphics[width=\unitlength,page=1]{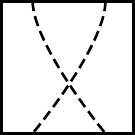}}%
  \end{picture}%
\endgroup%
}:L\otimes L\rightarrow L\otimes L.
   \end{split}
  \end{gather*}

  The degree of a $(w,w')$-diagram is defined by
  $$\frac{1}{2}\#\{\text{vertices}\}-\#\{\text{univalent vertices attached to the upper line }l'\}.$$
  Let $\A^{L}_{d}(w,w')\subset \A^{L}(w,w')$ be the subspace spanned by $(w,w')$-diagrams of degree $d$.
  We have
  $\A^{L}(w,w')=\bigoplus_{d\geq 0}\A^{L}_{d}(w,w')$.
  Since we have
  $$\A^{L}_{d'}(w',w'')\circ \A^{L}_{d}(w,w')\subset \A^{L}_{d+d'}(w,w'')$$ and
  $$\A^{L}_{d'}(w,w')\otimes \A^{L}_{d}(z,z')\subset \A^{L}_{d+d'}(w\otimes z,w'\otimes z')$$
  for any $w,w',w'',z,z'\in \A^{L}$,
  this grading is an $\N$-grading on $\A^{L}$.
  Note that we have $\A_d(m,n)= \A^{L}_d(H^{\otimes m},H^{\otimes n})$ for $m,n\geq 0$.

 \subsection{Relations for morphisms in $\A^{L}$}\label{ss425}
  Here, we observe some relations for morphisms of $\A^{L}$, which we use in the proof of Theorem \ref{th731}.

  The cocommutative Hopf algebra $(H,\mu,\eta,\Delta,\epsilon,S)$ in $\A$ naturally induces a cocommutative Hopf algebra in $\A^{L}$, such that
  \begin{gather*}
    \mu=\centre{
\begingroup%
  \makeatletter%
  \providecommand\color[2][]{%
    \errmessage{(Inkscape) Color is used for the text in Inkscape, but the package 'color.sty' is not loaded}%
    \renewcommand\color[2][]{}%
  }%
  \providecommand\transparent[1]{%
    \errmessage{(Inkscape) Transparency is used (non-zero) for the text in Inkscape, but the package 'transparent.sty' is not loaded}%
    \renewcommand\transparent[1]{}%
  }%
  \providecommand\rotatebox[2]{#2}%
  \newcommand*\fsize{\dimexpr\f@size pt\relax}%
  \newcommand*\lineheight[1]{\fontsize{\fsize}{#1\fsize}\selectfont}%
  \ifx\svgwidth\undefined%
    \setlength{\unitlength}{15.53033032bp}%
    \ifx\svgscale\undefined%
      \relax%
    \else%
      \setlength{\unitlength}{\unitlength * \real{\svgscale}}%
    \fi%
  \else%
    \setlength{\unitlength}{\svgwidth}%
  \fi%
  \global\let\svgwidth\undefined%
  \global\let\svgscale\undefined%
  \makeatother%
  \begin{picture}(1,1.22438894)%
    \lineheight{1}%
    \setlength\tabcolsep{0pt}%
    \put(0,0){\includegraphics[width=\unitlength,page=1]{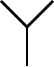}}%
  \end{picture}%
\endgroup%
}:=\centre{},\;
    \eta=\centre{
\begingroup%
  \makeatletter%
  \providecommand\color[2][]{%
    \errmessage{(Inkscape) Color is used for the text in Inkscape, but the package 'color.sty' is not loaded}%
    \renewcommand\color[2][]{}%
  }%
  \providecommand\transparent[1]{%
    \errmessage{(Inkscape) Transparency is used (non-zero) for the text in Inkscape, but the package 'transparent.sty' is not loaded}%
    \renewcommand\transparent[1]{}%
  }%
  \providecommand\rotatebox[2]{#2}%
  \newcommand*\fsize{\dimexpr\f@size pt\relax}%
  \newcommand*\lineheight[1]{\fontsize{\fsize}{#1\fsize}\selectfont}%
  \ifx\svgwidth\undefined%
    \setlength{\unitlength}{1.85746565bp}%
    \ifx\svgscale\undefined%
      \relax%
    \else%
      \setlength{\unitlength}{\unitlength * \real{\svgscale}}%
    \fi%
  \else%
    \setlength{\unitlength}{\svgwidth}%
  \fi%
  \global\let\svgwidth\undefined%
  \global\let\svgscale\undefined%
  \makeatother%
  \begin{picture}(1,10.74463112)%
    \lineheight{1}%
    \setlength\tabcolsep{0pt}%
    \put(0,0){\includegraphics[width=\unitlength,page=1]{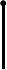}}%
  \end{picture}%
\endgroup%
}:=\centre{},\;
    \Delta=\centre{
\begingroup%
  \makeatletter%
  \providecommand\color[2][]{%
    \errmessage{(Inkscape) Color is used for the text in Inkscape, but the package 'color.sty' is not loaded}%
    \renewcommand\color[2][]{}%
  }%
  \providecommand\transparent[1]{%
    \errmessage{(Inkscape) Transparency is used (non-zero) for the text in Inkscape, but the package 'transparent.sty' is not loaded}%
    \renewcommand\transparent[1]{}%
  }%
  \providecommand\rotatebox[2]{#2}%
  \newcommand*\fsize{\dimexpr\f@size pt\relax}%
  \newcommand*\lineheight[1]{\fontsize{\fsize}{#1\fsize}\selectfont}%
  \ifx\svgwidth\undefined%
    \setlength{\unitlength}{15.62403785bp}%
    \ifx\svgscale\undefined%
      \relax%
    \else%
      \setlength{\unitlength}{\unitlength * \real{\svgscale}}%
    \fi%
  \else%
    \setlength{\unitlength}{\svgwidth}%
  \fi%
  \global\let\svgwidth\undefined%
  \global\let\svgscale\undefined%
  \makeatother%
  \begin{picture}(1,1.21338752)%
    \lineheight{1}%
    \setlength\tabcolsep{0pt}%
    \put(0,0){\includegraphics[width=\unitlength,page=1]{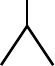}}%
  \end{picture}%
\endgroup%
}=\centre{},\\
    \epsilon=\centre{
\begingroup%
  \makeatletter%
  \providecommand\color[2][]{%
    \errmessage{(Inkscape) Color is used for the text in Inkscape, but the package 'color.sty' is not loaded}%
    \renewcommand\color[2][]{}%
  }%
  \providecommand\transparent[1]{%
    \errmessage{(Inkscape) Transparency is used (non-zero) for the text in Inkscape, but the package 'transparent.sty' is not loaded}%
    \renewcommand\transparent[1]{}%
  }%
  \providecommand\rotatebox[2]{#2}%
  \newcommand*\fsize{\dimexpr\f@size pt\relax}%
  \newcommand*\lineheight[1]{\fontsize{\fsize}{#1\fsize}\selectfont}%
  \ifx\svgwidth\undefined%
    \setlength{\unitlength}{2.05290117bp}%
    \ifx\svgscale\undefined%
      \relax%
    \else%
      \setlength{\unitlength}{\unitlength * \real{\svgscale}}%
    \fi%
  \else%
    \setlength{\unitlength}{\svgwidth}%
  \fi%
  \global\let\svgwidth\undefined%
  \global\let\svgscale\undefined%
  \makeatother%
  \begin{picture}(1,9.66835542)%
    \lineheight{1}%
    \setlength\tabcolsep{0pt}%
    \put(0,0){\includegraphics[width=\unitlength,page=1]{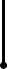}}%
  \end{picture}%
\endgroup%
}:=\centre{},\;
    S=\centre{
\begingroup%
  \makeatletter%
  \providecommand\color[2][]{%
    \errmessage{(Inkscape) Color is used for the text in Inkscape, but the package 'color.sty' is not loaded}%
    \renewcommand\color[2][]{}%
  }%
  \providecommand\transparent[1]{%
    \errmessage{(Inkscape) Transparency is used (non-zero) for the text in Inkscape, but the package 'transparent.sty' is not loaded}%
    \renewcommand\transparent[1]{}%
  }%
  \providecommand\rotatebox[2]{#2}%
  \newcommand*\fsize{\dimexpr\f@size pt\relax}%
  \newcommand*\lineheight[1]{\fontsize{\fsize}{#1\fsize}\selectfont}%
  \ifx\svgwidth\undefined%
    \setlength{\unitlength}{3.7511809bp}%
    \ifx\svgscale\undefined%
      \relax%
    \else%
      \setlength{\unitlength}{\unitlength * \real{\svgscale}}%
    \fi%
  \else%
    \setlength{\unitlength}{\svgwidth}%
  \fi%
  \global\let\svgwidth\undefined%
  \global\let\svgscale\undefined%
  \makeatother%
  \begin{picture}(1,4.99842601)%
    \lineheight{1}%
    \setlength\tabcolsep{0pt}%
    \put(0,0){\includegraphics[width=\unitlength,page=1]{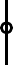}}%
  \end{picture}%
\endgroup%
}:=\centre{}.
  \end{gather*}
  Additionally, the triple $(L,[\cdot,\cdot],c_{L})$ is a Lie algebra with a symmetric invariant $2$-tensor in $\A^{L}$ (see Appendix \ref{ssA2}), where
  \begin{gather*}
    [\cdot,\cdot]=\centre{
\begingroup%
  \makeatletter%
  \providecommand\color[2][]{%
    \errmessage{(Inkscape) Color is used for the text in Inkscape, but the package 'color.sty' is not loaded}%
    \renewcommand\color[2][]{}%
  }%
  \providecommand\transparent[1]{%
    \errmessage{(Inkscape) Transparency is used (non-zero) for the text in Inkscape, but the package 'transparent.sty' is not loaded}%
    \renewcommand\transparent[1]{}%
  }%
  \providecommand\rotatebox[2]{#2}%
  \newcommand*\fsize{\dimexpr\f@size pt\relax}%
  \newcommand*\lineheight[1]{\fontsize{\fsize}{#1\fsize}\selectfont}%
  \ifx\svgwidth\undefined%
    \setlength{\unitlength}{15.37955329bp}%
    \ifx\svgscale\undefined%
      \relax%
    \else%
      \setlength{\unitlength}{\unitlength * \real{\svgscale}}%
    \fi%
  \else%
    \setlength{\unitlength}{\svgwidth}%
  \fi%
  \global\let\svgwidth\undefined%
  \global\let\svgscale\undefined%
  \makeatother%
  \begin{picture}(1,1.23639253)%
    \lineheight{1}%
    \setlength\tabcolsep{0pt}%
    \put(0,0){\includegraphics[width=\unitlength,page=1]{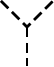}}%
  \end{picture}%
\endgroup%
}:=\centre{
\begingroup%
  \makeatletter%
  \providecommand\color[2][]{%
    \errmessage{(Inkscape) Color is used for the text in Inkscape, but the package 'color.sty' is not loaded}%
    \renewcommand\color[2][]{}%
  }%
  \providecommand\transparent[1]{%
    \errmessage{(Inkscape) Transparency is used (non-zero) for the text in Inkscape, but the package 'transparent.sty' is not loaded}%
    \renewcommand\transparent[1]{}%
  }%
  \providecommand\rotatebox[2]{#2}%
  \newcommand*\fsize{\dimexpr\f@size pt\relax}%
  \newcommand*\lineheight[1]{\fontsize{\fsize}{#1\fsize}\selectfont}%
  \ifx\svgwidth\undefined%
    \setlength{\unitlength}{30.75117737bp}%
    \ifx\svgscale\undefined%
      \relax%
    \else%
      \setlength{\unitlength}{\unitlength * \real{\svgscale}}%
    \fi%
  \else%
    \setlength{\unitlength}{\svgwidth}%
  \fi%
  \global\let\svgwidth\undefined%
  \global\let\svgscale\undefined%
  \makeatother%
  \begin{picture}(1,1)%
    \lineheight{1}%
    \setlength\tabcolsep{0pt}%
    \put(0,0){\includegraphics[width=\unitlength,page=1]{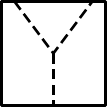}}%
  \end{picture}%
\endgroup%
}:L^{\otimes 2}\rightarrow L,\;
    c_{L}=\centre{
\begingroup%
  \makeatletter%
  \providecommand\color[2][]{%
    \errmessage{(Inkscape) Color is used for the text in Inkscape, but the package 'color.sty' is not loaded}%
    \renewcommand\color[2][]{}%
  }%
  \providecommand\transparent[1]{%
    \errmessage{(Inkscape) Transparency is used (non-zero) for the text in Inkscape, but the package 'transparent.sty' is not loaded}%
    \renewcommand\transparent[1]{}%
  }%
  \providecommand\rotatebox[2]{#2}%
  \newcommand*\fsize{\dimexpr\f@size pt\relax}%
  \newcommand*\lineheight[1]{\fontsize{\fsize}{#1\fsize}\selectfont}%
  \ifx\svgwidth\undefined%
    \setlength{\unitlength}{15.70958866bp}%
    \ifx\svgscale\undefined%
      \relax%
    \else%
      \setlength{\unitlength}{\unitlength * \real{\svgscale}}%
    \fi%
  \else%
    \setlength{\unitlength}{\svgwidth}%
  \fi%
  \global\let\svgwidth\undefined%
  \global\let\svgscale\undefined%
  \makeatother%
  \begin{picture}(1,0.74797506)%
    \lineheight{1}%
    \setlength\tabcolsep{0pt}%
    \put(0,0){\includegraphics[width=\unitlength,page=1]{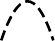}}%
  \end{picture}%
\endgroup%
}=:\centre{
\begingroup%
  \makeatletter%
  \providecommand\color[2][]{%
    \errmessage{(Inkscape) Color is used for the text in Inkscape, but the package 'color.sty' is not loaded}%
    \renewcommand\color[2][]{}%
  }%
  \providecommand\transparent[1]{%
    \errmessage{(Inkscape) Transparency is used (non-zero) for the text in Inkscape, but the package 'transparent.sty' is not loaded}%
    \renewcommand\transparent[1]{}%
  }%
  \providecommand\rotatebox[2]{#2}%
  \newcommand*\fsize{\dimexpr\f@size pt\relax}%
  \newcommand*\lineheight[1]{\fontsize{\fsize}{#1\fsize}\selectfont}%
  \ifx\svgwidth\undefined%
    \setlength{\unitlength}{30.75117737bp}%
    \ifx\svgscale\undefined%
      \relax%
    \else%
      \setlength{\unitlength}{\unitlength * \real{\svgscale}}%
    \fi%
  \else%
    \setlength{\unitlength}{\svgwidth}%
  \fi%
  \global\let\svgwidth\undefined%
  \global\let\svgscale\undefined%
  \makeatother%
  \begin{picture}(1,1)%
    \lineheight{1}%
    \setlength\tabcolsep{0pt}%
    \put(0,0){\includegraphics[width=\unitlength,page=1]{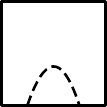}}%
  \end{picture}%
\endgroup%
}:I\rightarrow L^{\otimes 2}.
  \end{gather*}
  Moreover, $\A^{L}$ has two morphisms
  \begin{gather*}
    i=\centre{
\begingroup%
  \makeatletter%
  \providecommand\color[2][]{%
    \errmessage{(Inkscape) Color is used for the text in Inkscape, but the package 'color.sty' is not loaded}%
    \renewcommand\color[2][]{}%
  }%
  \providecommand\transparent[1]{%
    \errmessage{(Inkscape) Transparency is used (non-zero) for the text in Inkscape, but the package 'transparent.sty' is not loaded}%
    \renewcommand\transparent[1]{}%
  }%
  \providecommand\rotatebox[2]{#2}%
  \newcommand*\fsize{\dimexpr\f@size pt\relax}%
  \newcommand*\lineheight[1]{\fontsize{\fsize}{#1\fsize}\selectfont}%
  \ifx\svgwidth\undefined%
    \setlength{\unitlength}{10.83417053bp}%
    \ifx\svgscale\undefined%
      \relax%
    \else%
      \setlength{\unitlength}{\unitlength * \real{\svgscale}}%
    \fi%
  \else%
    \setlength{\unitlength}{\svgwidth}%
  \fi%
  \global\let\svgwidth\undefined%
  \global\let\svgscale\undefined%
  \makeatother%
  \begin{picture}(1,2.76901669)%
    \lineheight{1}%
    \setlength\tabcolsep{0pt}%
    \put(0,0){\includegraphics[width=\unitlength,page=1]{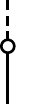}}%
    \put(0.45412913,1.44197913){\makebox(0,0)[lt]{\lineheight{1.45000005}\smash{\begin{tabular}[t]{l}$i$\end{tabular}}}}%
  \end{picture}%
\endgroup%
}:=\centre{
\begingroup%
  \makeatletter%
  \providecommand\color[2][]{%
    \errmessage{(Inkscape) Color is used for the text in Inkscape, but the package 'color.sty' is not loaded}%
    \renewcommand\color[2][]{}%
  }%
  \providecommand\transparent[1]{%
    \errmessage{(Inkscape) Transparency is used (non-zero) for the text in Inkscape, but the package 'transparent.sty' is not loaded}%
    \renewcommand\transparent[1]{}%
  }%
  \providecommand\rotatebox[2]{#2}%
  \newcommand*\fsize{\dimexpr\f@size pt\relax}%
  \newcommand*\lineheight[1]{\fontsize{\fsize}{#1\fsize}\selectfont}%
  \ifx\svgwidth\undefined%
    \setlength{\unitlength}{30.75117737bp}%
    \ifx\svgscale\undefined%
      \relax%
    \else%
      \setlength{\unitlength}{\unitlength * \real{\svgscale}}%
    \fi%
  \else%
    \setlength{\unitlength}{\svgwidth}%
  \fi%
  \global\let\svgwidth\undefined%
  \global\let\svgscale\undefined%
  \makeatother%
  \begin{picture}(1,1.00915176)%
    \lineheight{1}%
    \setlength\tabcolsep{0pt}%
    \put(0,0){\includegraphics[width=\unitlength,page=1]{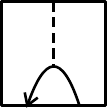}}%
  \end{picture}%
\endgroup%
}:L\rightarrow H,\;
    ad_{L}=\centre{
\begingroup%
  \makeatletter%
  \providecommand\color[2][]{%
    \errmessage{(Inkscape) Color is used for the text in Inkscape, but the package 'color.sty' is not loaded}%
    \renewcommand\color[2][]{}%
  }%
  \providecommand\transparent[1]{%
    \errmessage{(Inkscape) Transparency is used (non-zero) for the text in Inkscape, but the package 'transparent.sty' is not loaded}%
    \renewcommand\transparent[1]{}%
  }%
  \providecommand\rotatebox[2]{#2}%
  \newcommand*\fsize{\dimexpr\f@size pt\relax}%
  \newcommand*\lineheight[1]{\fontsize{\fsize}{#1\fsize}\selectfont}%
  \ifx\svgwidth\undefined%
    \setlength{\unitlength}{22.86421406bp}%
    \ifx\svgscale\undefined%
      \relax%
    \else%
      \setlength{\unitlength}{\unitlength * \real{\svgscale}}%
    \fi%
  \else%
    \setlength{\unitlength}{\svgwidth}%
  \fi%
  \global\let\svgwidth\undefined%
  \global\let\svgscale\undefined%
  \makeatother%
  \begin{picture}(1,1.14808236)%
    \lineheight{1}%
    \setlength\tabcolsep{0pt}%
    \put(0,0){\includegraphics[width=\unitlength,page=1]{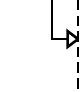}}%
    \put(-0.01324052,0.47725804){\makebox(0,0)[lt]{\lineheight{1.45000005}\smash{\begin{tabular}[t]{l}$ad_L$\end{tabular}}}}%
  \end{picture}%
\endgroup%
}:=\centre{
\begingroup%
  \makeatletter%
  \providecommand\color[2][]{%
    \errmessage{(Inkscape) Color is used for the text in Inkscape, but the package 'color.sty' is not loaded}%
    \renewcommand\color[2][]{}%
  }%
  \providecommand\transparent[1]{%
    \errmessage{(Inkscape) Transparency is used (non-zero) for the text in Inkscape, but the package 'transparent.sty' is not loaded}%
    \renewcommand\transparent[1]{}%
  }%
  \providecommand\rotatebox[2]{#2}%
  \newcommand*\fsize{\dimexpr\f@size pt\relax}%
  \newcommand*\lineheight[1]{\fontsize{\fsize}{#1\fsize}\selectfont}%
  \ifx\svgwidth\undefined%
    \setlength{\unitlength}{30.74999895bp}%
    \ifx\svgscale\undefined%
      \relax%
    \else%
      \setlength{\unitlength}{\unitlength * \real{\svgscale}}%
    \fi%
  \else%
    \setlength{\unitlength}{\svgwidth}%
  \fi%
  \global\let\svgwidth\undefined%
  \global\let\svgscale\undefined%
  \makeatother%
  \begin{picture}(1,1.41463416)%
    \lineheight{1}%
    \setlength\tabcolsep{0pt}%
    \put(0,0){\includegraphics[width=\unitlength,page=1]{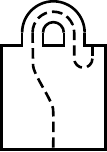}}%
  \end{picture}%
\endgroup%
}:H\otimes L\rightarrow L.
  \end{gather*}
  The degree of the morphism $c_{L}$ is $1$ and that of the others of the above morphisms is $0$.

  The iterated multiplications
  $$\mu^{[q]}=\centre{
\begingroup%
  \makeatletter%
  \providecommand\color[2][]{%
    \errmessage{(Inkscape) Color is used for the text in Inkscape, but the package 'color.sty' is not loaded}%
    \renewcommand\color[2][]{}%
  }%
  \providecommand\transparent[1]{%
    \errmessage{(Inkscape) Transparency is used (non-zero) for the text in Inkscape, but the package 'transparent.sty' is not loaded}%
    \renewcommand\transparent[1]{}%
  }%
  \providecommand\rotatebox[2]{#2}%
  \newcommand*\fsize{\dimexpr\f@size pt\relax}%
  \newcommand*\lineheight[1]{\fontsize{\fsize}{#1\fsize}\selectfont}%
  \ifx\svgwidth\undefined%
    \setlength{\unitlength}{18.53032939bp}%
    \ifx\svgscale\undefined%
      \relax%
    \else%
      \setlength{\unitlength}{\unitlength * \real{\svgscale}}%
    \fi%
  \else%
    \setlength{\unitlength}{\svgwidth}%
  \fi%
  \global\let\svgwidth\undefined%
  \global\let\svgscale\undefined%
  \makeatother%
  \begin{picture}(1,1.10711279)%
    \lineheight{1}%
    \setlength\tabcolsep{0pt}%
    \put(0,0){\includegraphics[width=\unitlength,page=1]{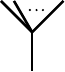}}%
  \end{picture}%
\endgroup%
}:H^{\otimes q}\rightarrow H$$
  and the iterated comultiplications
  $$\Delta^{[q]}=\centre{
\begingroup%
  \makeatletter%
  \providecommand\color[2][]{%
    \errmessage{(Inkscape) Color is used for the text in Inkscape, but the package 'color.sty' is not loaded}%
    \renewcommand\color[2][]{}%
  }%
  \providecommand\transparent[1]{%
    \errmessage{(Inkscape) Transparency is used (non-zero) for the text in Inkscape, but the package 'transparent.sty' is not loaded}%
    \renewcommand\transparent[1]{}%
  }%
  \providecommand\rotatebox[2]{#2}%
  \newcommand*\fsize{\dimexpr\f@size pt\relax}%
  \newcommand*\lineheight[1]{\fontsize{\fsize}{#1\fsize}\selectfont}%
  \ifx\svgwidth\undefined%
    \setlength{\unitlength}{18.61272512bp}%
    \ifx\svgscale\undefined%
      \relax%
    \else%
      \setlength{\unitlength}{\unitlength * \real{\svgscale}}%
    \fi%
  \else%
    \setlength{\unitlength}{\svgwidth}%
  \fi%
  \global\let\svgwidth\undefined%
  \global\let\svgscale\undefined%
  \makeatother%
  \begin{picture}(1,1.09958412)%
    \lineheight{1}%
    \setlength\tabcolsep{0pt}%
    \put(0,0){\includegraphics[width=\unitlength,page=1]{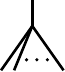}}%
  \end{picture}%
\endgroup%
}:H\rightarrow H^{\otimes q}$$
  for $q\geq 0$ are inductively defined by
  $$\mu^{[0]}=\eta,\quad\mu^{[1]}=\id_H, \quad \mu^{[q+1]}=\mu(\mu^{[q]}\otimes\id_H) \quad(q\geq 1),$$
  $$\Delta^{[0]}=\epsilon,\quad\Delta^{[1]}=\id_H, \quad \Delta^{[q+1]}=(\Delta^{[q]}\otimes\id_H)\Delta \quad(q\geq 1).$$
  Let
  $$
   ad_{H}=\centre{
\begingroup%
  \makeatletter%
  \providecommand\color[2][]{%
    \errmessage{(Inkscape) Color is used for the text in Inkscape, but the package 'color.sty' is not loaded}%
    \renewcommand\color[2][]{}%
  }%
  \providecommand\transparent[1]{%
    \errmessage{(Inkscape) Transparency is used (non-zero) for the text in Inkscape, but the package 'transparent.sty' is not loaded}%
    \renewcommand\transparent[1]{}%
  }%
  \providecommand\rotatebox[2]{#2}%
  \newcommand*\fsize{\dimexpr\f@size pt\relax}%
  \newcommand*\lineheight[1]{\fontsize{\fsize}{#1\fsize}\selectfont}%
  \ifx\svgwidth\undefined%
    \setlength{\unitlength}{26.01602248bp}%
    \ifx\svgscale\undefined%
      \relax%
    \else%
      \setlength{\unitlength}{\unitlength * \real{\svgscale}}%
    \fi%
  \else%
    \setlength{\unitlength}{\svgwidth}%
  \fi%
  \global\let\svgwidth\undefined%
  \global\let\svgscale\undefined%
  \makeatother%
  \begin{picture}(1,1.00899359)%
    \lineheight{1}%
    \setlength\tabcolsep{0pt}%
    \put(0,0){\includegraphics[width=\unitlength,page=1]{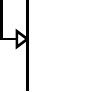}}%
    \put(0.41322261,0.52178294){\makebox(0,0)[lt]{\lineheight{1.45000005}\smash{\begin{tabular}[t]{l}$ad_H$\end{tabular}}}}%
  \end{picture}%
\endgroup%
}:= \centre{
\begingroup%
  \makeatletter%
  \providecommand\color[2][]{%
    \errmessage{(Inkscape) Color is used for the text in Inkscape, but the package 'color.sty' is not loaded}%
    \renewcommand\color[2][]{}%
  }%
  \providecommand\transparent[1]{%
    \errmessage{(Inkscape) Transparency is used (non-zero) for the text in Inkscape, but the package 'transparent.sty' is not loaded}%
    \renewcommand\transparent[1]{}%
  }%
  \providecommand\rotatebox[2]{#2}%
  \newcommand*\fsize{\dimexpr\f@size pt\relax}%
  \newcommand*\lineheight[1]{\fontsize{\fsize}{#1\fsize}\selectfont}%
  \ifx\svgwidth\undefined%
    \setlength{\unitlength}{15.67500013bp}%
    \ifx\svgscale\undefined%
      \relax%
    \else%
      \setlength{\unitlength}{\unitlength * \real{\svgscale}}%
    \fi%
  \else%
    \setlength{\unitlength}{\svgwidth}%
  \fi%
  \global\let\svgwidth\undefined%
  \global\let\svgscale\undefined%
  \makeatother%
  \begin{picture}(1,1.68899512)%
    \lineheight{1}%
    \setlength\tabcolsep{0pt}%
    \put(0,0){\includegraphics[width=\unitlength,page=1]{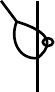}}%
  \end{picture}%
\endgroup%
}=\centre{
\begingroup%
  \makeatletter%
  \providecommand\color[2][]{%
    \errmessage{(Inkscape) Color is used for the text in Inkscape, but the package 'color.sty' is not loaded}%
    \renewcommand\color[2][]{}%
  }%
  \providecommand\transparent[1]{%
    \errmessage{(Inkscape) Transparency is used (non-zero) for the text in Inkscape, but the package 'transparent.sty' is not loaded}%
    \renewcommand\transparent[1]{}%
  }%
  \providecommand\rotatebox[2]{#2}%
  \newcommand*\fsize{\dimexpr\f@size pt\relax}%
  \newcommand*\lineheight[1]{\fontsize{\fsize}{#1\fsize}\selectfont}%
  \ifx\svgwidth\undefined%
    \setlength{\unitlength}{44.2499995bp}%
    \ifx\svgscale\undefined%
      \relax%
    \else%
      \setlength{\unitlength}{\unitlength * \real{\svgscale}}%
    \fi%
  \else%
    \setlength{\unitlength}{\svgwidth}%
  \fi%
  \global\let\svgwidth\undefined%
  \global\let\svgscale\undefined%
  \makeatother%
  \begin{picture}(1,0.95596904)%
    \lineheight{1}%
    \setlength\tabcolsep{0pt}%
    \put(0,0){\includegraphics[width=\unitlength,page=1]{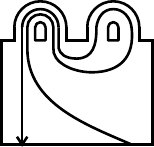}}%
  \end{picture}%
\endgroup%
},
  $$
  which denotes the \emph{adjoint action}, and
  \begin{equation}\label{comm}
    comm=\centre{
\begingroup%
  \makeatletter%
  \providecommand\color[2][]{%
    \errmessage{(Inkscape) Color is used for the text in Inkscape, but the package 'color.sty' is not loaded}%
    \renewcommand\color[2][]{}%
  }%
  \providecommand\transparent[1]{%
    \errmessage{(Inkscape) Transparency is used (non-zero) for the text in Inkscape, but the package 'transparent.sty' is not loaded}%
    \renewcommand\transparent[1]{}%
  }%
  \providecommand\rotatebox[2]{#2}%
  \newcommand*\fsize{\dimexpr\f@size pt\relax}%
  \newcommand*\lineheight[1]{\fontsize{\fsize}{#1\fsize}\selectfont}%
  \ifx\svgwidth\undefined%
    \setlength{\unitlength}{30.53033043bp}%
    \ifx\svgscale\undefined%
      \relax%
    \else%
      \setlength{\unitlength}{\unitlength * \real{\svgscale}}%
    \fi%
  \else%
    \setlength{\unitlength}{\svgwidth}%
  \fi%
  \global\let\svgwidth\undefined%
  \global\let\svgscale\undefined%
  \makeatother%
  \begin{picture}(1,1.11414333)%
    \lineheight{1}%
    \setlength\tabcolsep{0pt}%
    \put(0,0){\includegraphics[width=\unitlength,page=1]{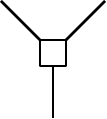}}%
  \end{picture}%
\endgroup%
}:= \centre{
\begingroup%
  \makeatletter%
  \providecommand\color[2][]{%
    \errmessage{(Inkscape) Color is used for the text in Inkscape, but the package 'color.sty' is not loaded}%
    \renewcommand\color[2][]{}%
  }%
  \providecommand\transparent[1]{%
    \errmessage{(Inkscape) Transparency is used (non-zero) for the text in Inkscape, but the package 'transparent.sty' is not loaded}%
    \renewcommand\transparent[1]{}%
  }%
  \providecommand\rotatebox[2]{#2}%
  \newcommand*\fsize{\dimexpr\f@size pt\relax}%
  \newcommand*\lineheight[1]{\fontsize{\fsize}{#1\fsize}\selectfont}%
  \ifx\svgwidth\undefined%
    \setlength{\unitlength}{33.51819455bp}%
    \ifx\svgscale\undefined%
      \relax%
    \else%
      \setlength{\unitlength}{\unitlength * \real{\svgscale}}%
    \fi%
  \else%
    \setlength{\unitlength}{\svgwidth}%
  \fi%
  \global\let\svgwidth\undefined%
  \global\let\svgscale\undefined%
  \makeatother%
  \begin{picture}(1,0.92566937)%
    \lineheight{1}%
    \setlength\tabcolsep{0pt}%
    \put(0,0){\includegraphics[width=\unitlength,page=1]{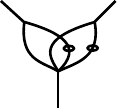}}%
  \end{picture}%
\endgroup%
}=\centre{
\begingroup%
  \makeatletter%
  \providecommand\color[2][]{%
    \errmessage{(Inkscape) Color is used for the text in Inkscape, but the package 'color.sty' is not loaded}%
    \renewcommand\color[2][]{}%
  }%
  \providecommand\transparent[1]{%
    \errmessage{(Inkscape) Transparency is used (non-zero) for the text in Inkscape, but the package 'transparent.sty' is not loaded}%
    \renewcommand\transparent[1]{}%
  }%
  \providecommand\rotatebox[2]{#2}%
  \newcommand*\fsize{\dimexpr\f@size pt\relax}%
  \newcommand*\lineheight[1]{\fontsize{\fsize}{#1\fsize}\selectfont}%
  \ifx\svgwidth\undefined%
    \setlength{\unitlength}{44.2499995bp}%
    \ifx\svgscale\undefined%
      \relax%
    \else%
      \setlength{\unitlength}{\unitlength * \real{\svgscale}}%
    \fi%
  \else%
    \setlength{\unitlength}{\svgwidth}%
  \fi%
  \global\let\svgwidth\undefined%
  \global\let\svgscale\undefined%
  \makeatother%
  \begin{picture}(1,0.95596901)%
    \lineheight{1}%
    \setlength\tabcolsep{0pt}%
    \put(0,0){\includegraphics[width=\unitlength,page=1]{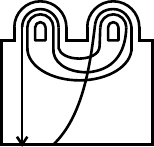}}%
  \end{picture}%
\endgroup%
},
  \end{equation}
  which denotes the \emph{commutator}.

  \begin{lemma}\label{l723}
   We have
    \begin{enumerate}
      \item $S\circ i =-i$\label{l7231}
      \item $\Delta \circ i=i\otimes \eta +\eta\otimes i$\label{l7232}
      \item $\epsilon\circ i=0$\label{l7233}
      \item $ad_{H}(i\otimes i)=-i\circ[\cdot,\cdot]$.\label{l7234}
    \end{enumerate}
   \end{lemma}
  \begin{proof}
   They can be checked by diagrammatic computation.
  \end{proof}

  Let $\g$ be a Lie algebra and $U=U(\g)$ be the universal enveloping algebra.
  We have a filtration $F_{\ast}(U)$ of $U$ induced by the usual filtration of the tensor algebra $T(\g)$ of $\g$.
  Since $U$ has a cocommutative Hopf algebra structure, we can define the commutator operator
  $$comm:U^{\otimes 2}\rightarrow U$$
  in a similar way as \eqref{comm}.
  For $x_1,\cdots,x_m,y_1,\cdots,y_n\in \g$, we have
  $$comm(x_1\cdots x_m,y_1\cdots y_n)\in F_{\min(m,n)}(U).$$
  The following lemma is a diagrammatic version of this fact.

  \begin{lemma}\label{l721}
     \begin{enumerate}
       \item Let $m,n\geq 1$.
        We have
        $$\centre{
\begingroup%
  \makeatletter%
  \providecommand\color[2][]{%
    \errmessage{(Inkscape) Color is used for the text in Inkscape, but the package 'color.sty' is not loaded}%
    \renewcommand\color[2][]{}%
  }%
  \providecommand\transparent[1]{%
    \errmessage{(Inkscape) Transparency is used (non-zero) for the text in Inkscape, but the package 'transparent.sty' is not loaded}%
    \renewcommand\transparent[1]{}%
  }%
  \providecommand\rotatebox[2]{#2}%
  \newcommand*\fsize{\dimexpr\f@size pt\relax}%
  \newcommand*\lineheight[1]{\fontsize{\fsize}{#1\fsize}\selectfont}%
  \ifx\svgwidth\undefined%
    \setlength{\unitlength}{100.80512746bp}%
    \ifx\svgscale\undefined%
      \relax%
    \else%
      \setlength{\unitlength}{\unitlength * \real{\svgscale}}%
    \fi%
  \else%
    \setlength{\unitlength}{\svgwidth}%
  \fi%
  \global\let\svgwidth\undefined%
  \global\let\svgscale\undefined%
  \makeatother%
  \begin{picture}(1,0.96191636)%
    \lineheight{1}%
    \setlength\tabcolsep{0pt}%
    \put(0,0){\includegraphics[width=\unitlength,page=1]{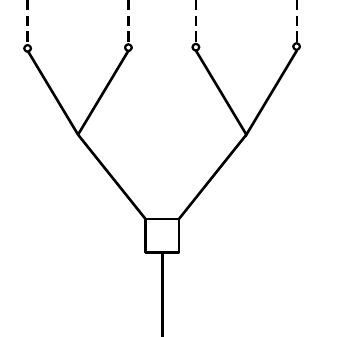}}%
    \put(-0.00087675,0.81286176){\makebox(0,0)[lt]{\lineheight{1.45000005}\smash{\begin{tabular}[t]{l}\small$i$\end{tabular}}}}%
    \put(0.88251602,0.81550962){\makebox(0,0)[lt]{\lineheight{1.45000005}\smash{\begin{tabular}[t]{l}\small$i$\end{tabular}}}}%
    \put(0.05888971,0.52237511){\makebox(0,0)[lt]{\lineheight{1.45000005}\smash{\begin{tabular}[t]{l}\small$\mu^{[m]}$\end{tabular}}}}%
    \put(0.71636931,0.53428534){\makebox(0,0)[lt]{\lineheight{1.45000005}\smash{\begin{tabular}[t]{l}\small$\mu^{[n]}$\end{tabular}}}}%
    \put(0,0){\includegraphics[width=\unitlength,page=2]{l72111.pdf}}%
  \end{picture}%
\endgroup%
}
        =\sum_{\alpha}c_\alpha\centre{
\begingroup%
  \makeatletter%
  \providecommand\color[2][]{%
    \errmessage{(Inkscape) Color is used for the text in Inkscape, but the package 'color.sty' is not loaded}%
    \renewcommand\color[2][]{}%
  }%
  \providecommand\transparent[1]{%
    \errmessage{(Inkscape) Transparency is used (non-zero) for the text in Inkscape, but the package 'transparent.sty' is not loaded}%
    \renewcommand\transparent[1]{}%
  }%
  \providecommand\rotatebox[2]{#2}%
  \newcommand*\fsize{\dimexpr\f@size pt\relax}%
  \newcommand*\lineheight[1]{\fontsize{\fsize}{#1\fsize}\selectfont}%
  \ifx\svgwidth\undefined%
    \setlength{\unitlength}{61.82644933bp}%
    \ifx\svgscale\undefined%
      \relax%
    \else%
      \setlength{\unitlength}{\unitlength * \real{\svgscale}}%
    \fi%
  \else%
    \setlength{\unitlength}{\svgwidth}%
  \fi%
  \global\let\svgwidth\undefined%
  \global\let\svgscale\undefined%
  \makeatother%
  \begin{picture}(1,1.60702524)%
    \lineheight{1}%
    \setlength\tabcolsep{0pt}%
    \put(0,0){\includegraphics[width=\unitlength,page=1]{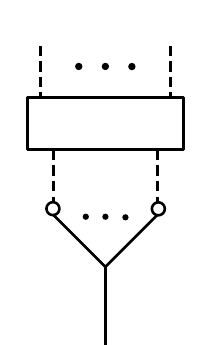}}%
    \put(0.37213376,0.95269146){\makebox(0,0)[lt]{\lineheight{1.45000005}\smash{\begin{tabular}[t]{l}$D_\alpha$\end{tabular}}}}%
    \put(0.0565834,0.59866999){\makebox(0,0)[lt]{\lineheight{1.45000005}\smash{\begin{tabular}[t]{l}$i$\end{tabular}}}}%
    \put(0.83964053,0.58506888){\makebox(0,0)[lt]{\lineheight{1.45000005}\smash{\begin{tabular}[t]{l}$i$\end{tabular}}}}%
    \put(0.55394326,0.32503001){\makebox(0,0)[lt]{\lineheight{1.45000005}\smash{\begin{tabular}[t]{l}$\mu^{[n]}$\end{tabular}}}}%
    \put(0.37317834,1.55723874){\makebox(0,0)[lt]{\lineheight{1.45000005}\smash{\begin{tabular}[t]{l}$m+n$\end{tabular}}}}%
    \put(0,0){\includegraphics[width=\unitlength,page=2]{l72112.pdf}}%
  \end{picture}%
\endgroup%
}
        =\sum_{\beta}c_\beta\centre{
\begingroup%
  \makeatletter%
  \providecommand\color[2][]{%
    \errmessage{(Inkscape) Color is used for the text in Inkscape, but the package 'color.sty' is not loaded}%
    \renewcommand\color[2][]{}%
  }%
  \providecommand\transparent[1]{%
    \errmessage{(Inkscape) Transparency is used (non-zero) for the text in Inkscape, but the package 'transparent.sty' is not loaded}%
    \renewcommand\transparent[1]{}%
  }%
  \providecommand\rotatebox[2]{#2}%
  \newcommand*\fsize{\dimexpr\f@size pt\relax}%
  \newcommand*\lineheight[1]{\fontsize{\fsize}{#1\fsize}\selectfont}%
  \ifx\svgwidth\undefined%
    \setlength{\unitlength}{63.18582361bp}%
    \ifx\svgscale\undefined%
      \relax%
    \else%
      \setlength{\unitlength}{\unitlength * \real{\svgscale}}%
    \fi%
  \else%
    \setlength{\unitlength}{\svgwidth}%
  \fi%
  \global\let\svgwidth\undefined%
  \global\let\svgscale\undefined%
  \makeatother%
  \begin{picture}(1,1.57245184)%
    \lineheight{1}%
    \setlength\tabcolsep{0pt}%
    \put(0,0){\includegraphics[width=\unitlength,page=1]{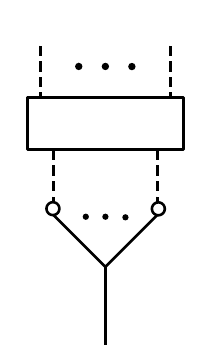}}%
    \put(0.36412771,0.9243274){\makebox(0,0)[lt]{\lineheight{1.45000005}\smash{\begin{tabular}[t]{l}$D_\beta$\end{tabular}}}}%
    \put(0.05536607,0.58579025){\makebox(0,0)[lt]{\lineheight{1.45000005}\smash{\begin{tabular}[t]{l}$i$\end{tabular}}}}%
    \put(0.82157657,0.57248176){\makebox(0,0)[lt]{\lineheight{1.45000005}\smash{\begin{tabular}[t]{l}$i$\end{tabular}}}}%
    \put(0.54202577,0.31803734){\makebox(0,0)[lt]{\lineheight{1.45000005}\smash{\begin{tabular}[t]{l}$\mu^{[m]}$\end{tabular}}}}%
    \put(0.36514982,1.52373644){\makebox(0,0)[lt]{\lineheight{1.45000005}\smash{\begin{tabular}[t]{l}$m+n$\end{tabular}}}}%
    \put(0,0){\includegraphics[width=\unitlength,page=2]{l72113.pdf}}%
  \end{picture}%
\endgroup%
},$$
        where $c_\alpha,c_\beta\in \Z$, and where $D_\alpha$ (resp. $D_\beta$) is a union of trees with $m$ (resp. $n$) trivalent vertices.
        Moreover, for $m=n=1$, we have
        $$\centre{
\begingroup%
  \makeatletter%
  \providecommand\color[2][]{%
    \errmessage{(Inkscape) Color is used for the text in Inkscape, but the package 'color.sty' is not loaded}%
    \renewcommand\color[2][]{}%
  }%
  \providecommand\transparent[1]{%
    \errmessage{(Inkscape) Transparency is used (non-zero) for the text in Inkscape, but the package 'transparent.sty' is not loaded}%
    \renewcommand\transparent[1]{}%
  }%
  \providecommand\rotatebox[2]{#2}%
  \newcommand*\fsize{\dimexpr\f@size pt\relax}%
  \newcommand*\lineheight[1]{\fontsize{\fsize}{#1\fsize}\selectfont}%
  \ifx\svgwidth\undefined%
    \setlength{\unitlength}{43.82791337bp}%
    \ifx\svgscale\undefined%
      \relax%
    \else%
      \setlength{\unitlength}{\unitlength * \real{\svgscale}}%
    \fi%
  \else%
    \setlength{\unitlength}{\svgwidth}%
  \fi%
  \global\let\svgwidth\undefined%
  \global\let\svgscale\undefined%
  \makeatother%
  \begin{picture}(1,1.19786668)%
    \lineheight{1}%
    \setlength\tabcolsep{0pt}%
    \put(0,0){\includegraphics[width=\unitlength,page=1]{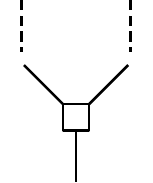}}%
    \put(-0.00345367,0.8008173){\makebox(0,0)[lt]{\lineheight{1.45000005}\smash{\begin{tabular}[t]{l}$i$\end{tabular}}}}%
    \put(0.93253083,0.80935697){\makebox(0,0)[lt]{\lineheight{1.45000005}\smash{\begin{tabular}[t]{l}$i$\end{tabular}}}}%
    \put(0,0){\includegraphics[width=\unitlength,page=2]{commi.pdf}}%
  \end{picture}%
\endgroup%
}=-\centre{
\begingroup%
  \makeatletter%
  \providecommand\color[2][]{%
    \errmessage{(Inkscape) Color is used for the text in Inkscape, but the package 'color.sty' is not loaded}%
    \renewcommand\color[2][]{}%
  }%
  \providecommand\transparent[1]{%
    \errmessage{(Inkscape) Transparency is used (non-zero) for the text in Inkscape, but the package 'transparent.sty' is not loaded}%
    \renewcommand\transparent[1]{}%
  }%
  \providecommand\rotatebox[2]{#2}%
  \newcommand*\fsize{\dimexpr\f@size pt\relax}%
  \newcommand*\lineheight[1]{\fontsize{\fsize}{#1\fsize}\selectfont}%
  \ifx\svgwidth\undefined%
    \setlength{\unitlength}{24.5865748bp}%
    \ifx\svgscale\undefined%
      \relax%
    \else%
      \setlength{\unitlength}{\unitlength * \real{\svgscale}}%
    \fi%
  \else%
    \setlength{\unitlength}{\svgwidth}%
  \fi%
  \global\let\svgwidth\undefined%
  \global\let\svgscale\undefined%
  \makeatother%
  \begin{picture}(1,1.8398101)%
    \lineheight{1}%
    \setlength\tabcolsep{0pt}%
    \put(0,0){\includegraphics[width=\unitlength,page=1]{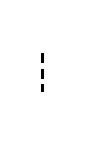}}%
    \put(0.20642933,0.63647302){\makebox(0,0)[lt]{\lineheight{1.45000005}\smash{\begin{tabular}[t]{l}$i$\end{tabular}}}}%
    \put(0,0){\includegraphics[width=\unitlength,page=2]{commib.pdf}}%
  \end{picture}%
\endgroup%
}.$$ \label{l7211}
       \item Let $m\geq 1$.
         We have $$\centre{
\begingroup%
  \makeatletter%
  \providecommand\color[2][]{%
    \errmessage{(Inkscape) Color is used for the text in Inkscape, but the package 'color.sty' is not loaded}%
    \renewcommand\color[2][]{}%
  }%
  \providecommand\transparent[1]{%
    \errmessage{(Inkscape) Transparency is used (non-zero) for the text in Inkscape, but the package 'transparent.sty' is not loaded}%
    \renewcommand\transparent[1]{}%
  }%
  \providecommand\rotatebox[2]{#2}%
  \newcommand*\fsize{\dimexpr\f@size pt\relax}%
  \newcommand*\lineheight[1]{\fontsize{\fsize}{#1\fsize}\selectfont}%
  \ifx\svgwidth\undefined%
    \setlength{\unitlength}{73.07242729bp}%
    \ifx\svgscale\undefined%
      \relax%
    \else%
      \setlength{\unitlength}{\unitlength * \real{\svgscale}}%
    \fi%
  \else%
    \setlength{\unitlength}{\svgwidth}%
  \fi%
  \global\let\svgwidth\undefined%
  \global\let\svgscale\undefined%
  \makeatother%
  \begin{picture}(1,1.14954435)%
    \lineheight{1}%
    \setlength\tabcolsep{0pt}%
    \put(0,0){\includegraphics[width=\unitlength,page=1]{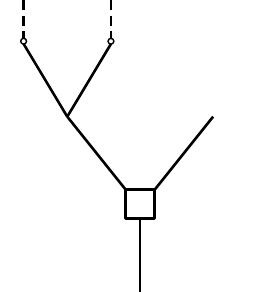}}%
    \put(-0.00232004,0.97141575){\makebox(0,0)[lt]{\lineheight{1.45000005}\smash{\begin{tabular}[t]{l}$i$\end{tabular}}}}%
    \put(0.48744234,0.9730179){\makebox(0,0)[lt]{\lineheight{1.45000005}\smash{\begin{tabular}[t]{l}$i$\end{tabular}}}}%
    \put(0.00702642,0.60439555){\makebox(0,0)[lt]{\lineheight{1.45000005}\smash{\begin{tabular}[t]{l}$\mu^{[m]}$\end{tabular}}}}%
    \put(0,0){\includegraphics[width=\unitlength,page=2]{l72121.pdf}}%
    \put(0.85767531,0.68404806){\makebox(0,0)[lt]{\lineheight{1.45000005}\smash{\begin{tabular}[t]{l}$\eta$\end{tabular}}}}%
    \put(0,0){\includegraphics[width=\unitlength,page=3]{l72121.pdf}}%
  \end{picture}%
\endgroup%
}=\centre{
\begingroup%
  \makeatletter%
  \providecommand\color[2][]{%
    \errmessage{(Inkscape) Color is used for the text in Inkscape, but the package 'color.sty' is not loaded}%
    \renewcommand\color[2][]{}%
  }%
  \providecommand\transparent[1]{%
    \errmessage{(Inkscape) Transparency is used (non-zero) for the text in Inkscape, but the package 'transparent.sty' is not loaded}%
    \renewcommand\transparent[1]{}%
  }%
  \providecommand\rotatebox[2]{#2}%
  \newcommand*\fsize{\dimexpr\f@size pt\relax}%
  \newcommand*\lineheight[1]{\fontsize{\fsize}{#1\fsize}\selectfont}%
  \ifx\svgwidth\undefined%
    \setlength{\unitlength}{68.29244762bp}%
    \ifx\svgscale\undefined%
      \relax%
    \else%
      \setlength{\unitlength}{\unitlength * \real{\svgscale}}%
    \fi%
  \else%
    \setlength{\unitlength}{\svgwidth}%
  \fi%
  \global\let\svgwidth\undefined%
  \global\let\svgscale\undefined%
  \makeatother%
  \begin{picture}(1,1.14214676)%
    \lineheight{1}%
    \setlength\tabcolsep{0pt}%
    \put(0,0){\includegraphics[width=\unitlength,page=1]{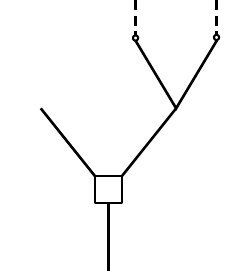}}%
    \put(0.49878437,0.9717469){\makebox(0,0)[lt]{\lineheight{1.45000005}\smash{\begin{tabular}[t]{l}$i$\end{tabular}}}}%
    \put(0.95496847,0.96830831){\makebox(0,0)[lt]{\lineheight{1.45000005}\smash{\begin{tabular}[t]{l}$i$\end{tabular}}}}%
    \put(0.77265102,0.62909059){\makebox(0,0)[lt]{\lineheight{1.45000005}\smash{\begin{tabular}[t]{l}$\mu^{[m]}$\end{tabular}}}}%
    \put(0,0){\includegraphics[width=\unitlength,page=2]{l72122.pdf}}%
    \put(-0.00354633,0.6752927){\makebox(0,0)[lt]{\lineheight{1.45000005}\smash{\begin{tabular}[t]{l}$\eta$\end{tabular}}}}%
    \put(0,0){\includegraphics[width=\unitlength,page=3]{l72122.pdf}}%
  \end{picture}%
\endgroup%
}=0.$$\label{l7212}
       \item
         We have $$\centre{
\begingroup%
  \makeatletter%
  \providecommand\color[2][]{%
    \errmessage{(Inkscape) Color is used for the text in Inkscape, but the package 'color.sty' is not loaded}%
    \renewcommand\color[2][]{}%
  }%
  \providecommand\transparent[1]{%
    \errmessage{(Inkscape) Transparency is used (non-zero) for the text in Inkscape, but the package 'transparent.sty' is not loaded}%
    \renewcommand\transparent[1]{}%
  }%
  \providecommand\rotatebox[2]{#2}%
  \newcommand*\fsize{\dimexpr\f@size pt\relax}%
  \newcommand*\lineheight[1]{\fontsize{\fsize}{#1\fsize}\selectfont}%
  \ifx\svgwidth\undefined%
    \setlength{\unitlength}{64.24983129bp}%
    \ifx\svgscale\undefined%
      \relax%
    \else%
      \setlength{\unitlength}{\unitlength * \real{\svgscale}}%
    \fi%
  \else%
    \setlength{\unitlength}{\svgwidth}%
  \fi%
  \global\let\svgwidth\undefined%
  \global\let\svgscale\undefined%
  \makeatother%
  \begin{picture}(1,0.8046459)%
    \lineheight{1}%
    \setlength\tabcolsep{0pt}%
    \put(0,0){\includegraphics[width=\unitlength,page=1]{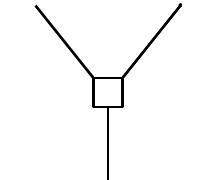}}%
    \put(0.83910243,0.76654876){\makebox(0,0)[lt]{\lineheight{1.45000005}\smash{\begin{tabular}[t]{l}$\eta$\end{tabular}}}}%
    \put(0,0){\includegraphics[width=\unitlength,page=2]{l72131.pdf}}%
    \put(-0.00374686,0.76639914){\makebox(0,0)[lt]{\lineheight{1.45000005}\smash{\begin{tabular}[t]{l}$\eta$\end{tabular}}}}%
  \end{picture}%
\endgroup%
}=\quad\centre{}.$$\label{l7213}
     \end{enumerate}
  \end{lemma}
  For example, we have
  \begin{gather*}
    \centre{
\begingroup%
  \makeatletter%
  \providecommand\color[2][]{%
    \errmessage{(Inkscape) Color is used for the text in Inkscape, but the package 'color.sty' is not loaded}%
    \renewcommand\color[2][]{}%
  }%
  \providecommand\transparent[1]{%
    \errmessage{(Inkscape) Transparency is used (non-zero) for the text in Inkscape, but the package 'transparent.sty' is not loaded}%
    \renewcommand\transparent[1]{}%
  }%
  \providecommand\rotatebox[2]{#2}%
  \newcommand*\fsize{\dimexpr\f@size pt\relax}%
  \newcommand*\lineheight[1]{\fontsize{\fsize}{#1\fsize}\selectfont}%
  \ifx\svgwidth\undefined%
    \setlength{\unitlength}{49.69710912bp}%
    \ifx\svgscale\undefined%
      \relax%
    \else%
      \setlength{\unitlength}{\unitlength * \real{\svgscale}}%
    \fi%
  \else%
    \setlength{\unitlength}{\svgwidth}%
  \fi%
  \global\let\svgwidth\undefined%
  \global\let\svgscale\undefined%
  \makeatother%
  \begin{picture}(1,1.14694804)%
    \lineheight{1}%
    \setlength\tabcolsep{0pt}%
    \put(0,0){\includegraphics[width=\unitlength,page=1]{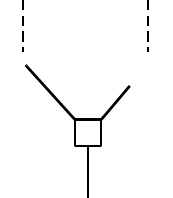}}%
    \put(-0.0030458,0.79055211){\makebox(0,0)[lt]{\lineheight{1.45000005}\smash{\begin{tabular}[t]{l}$i$\end{tabular}}}}%
    \put(0.94049889,0.7947275){\makebox(0,0)[lt]{\lineheight{1.45000005}\smash{\begin{tabular}[t]{l}$i$\end{tabular}}}}%
    \put(0,0){\includegraphics[width=\unitlength,page=2]{comm12.pdf}}%
    \put(0.50441104,0.79056381){\makebox(0,0)[lt]{\lineheight{1.45000005}\smash{\begin{tabular}[t]{l}$i$\end{tabular}}}}%
  \end{picture}%
\endgroup%
}=-\centre{
\begingroup%
  \makeatletter%
  \providecommand\color[2][]{%
    \errmessage{(Inkscape) Color is used for the text in Inkscape, but the package 'color.sty' is not loaded}%
    \renewcommand\color[2][]{}%
  }%
  \providecommand\transparent[1]{%
    \errmessage{(Inkscape) Transparency is used (non-zero) for the text in Inkscape, but the package 'transparent.sty' is not loaded}%
    \renewcommand\transparent[1]{}%
  }%
  \providecommand\rotatebox[2]{#2}%
  \newcommand*\fsize{\dimexpr\f@size pt\relax}%
  \newcommand*\lineheight[1]{\fontsize{\fsize}{#1\fsize}\selectfont}%
  \ifx\svgwidth\undefined%
    \setlength{\unitlength}{26.40102604bp}%
    \ifx\svgscale\undefined%
      \relax%
    \else%
      \setlength{\unitlength}{\unitlength * \real{\svgscale}}%
    \fi%
  \else%
    \setlength{\unitlength}{\svgwidth}%
  \fi%
  \global\let\svgwidth\undefined%
  \global\let\svgscale\undefined%
  \makeatother%
  \begin{picture}(1,1.20003351)%
    \lineheight{1}%
    \setlength\tabcolsep{0pt}%
    \put(0.57651564,0.37520897){\makebox(0,0)[lt]{\lineheight{1.45000005}\smash{\begin{tabular}[t]{l}$i$\end{tabular}}}}%
    \put(0,0){\includegraphics[width=\unitlength,page=1]{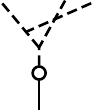}}%
  \end{picture}%
\endgroup%
},\quad
    \centre{
\begingroup%
  \makeatletter%
  \providecommand\color[2][]{%
    \errmessage{(Inkscape) Color is used for the text in Inkscape, but the package 'color.sty' is not loaded}%
    \renewcommand\color[2][]{}%
  }%
  \providecommand\transparent[1]{%
    \errmessage{(Inkscape) Transparency is used (non-zero) for the text in Inkscape, but the package 'transparent.sty' is not loaded}%
    \renewcommand\transparent[1]{}%
  }%
  \providecommand\rotatebox[2]{#2}%
  \newcommand*\fsize{\dimexpr\f@size pt\relax}%
  \newcommand*\lineheight[1]{\fontsize{\fsize}{#1\fsize}\selectfont}%
  \ifx\svgwidth\undefined%
    \setlength{\unitlength}{49.69710912bp}%
    \ifx\svgscale\undefined%
      \relax%
    \else%
      \setlength{\unitlength}{\unitlength * \real{\svgscale}}%
    \fi%
  \else%
    \setlength{\unitlength}{\svgwidth}%
  \fi%
  \global\let\svgwidth\undefined%
  \global\let\svgscale\undefined%
  \makeatother%
  \begin{picture}(1,1.14694805)%
    \lineheight{1}%
    \setlength\tabcolsep{0pt}%
    \put(0,0){\includegraphics[width=\unitlength,page=1]{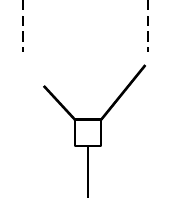}}%
    \put(-0.0030458,0.79055211){\makebox(0,0)[lt]{\lineheight{1.45000005}\smash{\begin{tabular}[t]{l}$i$\end{tabular}}}}%
    \put(0.94049889,0.7947275){\makebox(0,0)[lt]{\lineheight{1.45000005}\smash{\begin{tabular}[t]{l}$i$\end{tabular}}}}%
    \put(0,0){\includegraphics[width=\unitlength,page=2]{comm21.pdf}}%
    \put(0.45612836,0.80049104){\makebox(0,0)[lt]{\lineheight{1.45000005}\smash{\begin{tabular}[t]{l}$i$\end{tabular}}}}%
  \end{picture}%
\endgroup%
}=\centre{
\begingroup%
  \makeatletter%
  \providecommand\color[2][]{%
    \errmessage{(Inkscape) Color is used for the text in Inkscape, but the package 'color.sty' is not loaded}%
    \renewcommand\color[2][]{}%
  }%
  \providecommand\transparent[1]{%
    \errmessage{(Inkscape) Transparency is used (non-zero) for the text in Inkscape, but the package 'transparent.sty' is not loaded}%
    \renewcommand\transparent[1]{}%
  }%
  \providecommand\rotatebox[2]{#2}%
  \newcommand*\fsize{\dimexpr\f@size pt\relax}%
  \newcommand*\lineheight[1]{\fontsize{\fsize}{#1\fsize}\selectfont}%
  \ifx\svgwidth\undefined%
    \setlength{\unitlength}{23.09999922bp}%
    \ifx\svgscale\undefined%
      \relax%
    \else%
      \setlength{\unitlength}{\unitlength * \real{\svgscale}}%
    \fi%
  \else%
    \setlength{\unitlength}{\svgwidth}%
  \fi%
  \global\let\svgwidth\undefined%
  \global\let\svgscale\undefined%
  \makeatother%
  \begin{picture}(1,1.40584424)%
    \lineheight{1}%
    \setlength\tabcolsep{0pt}%
    \put(0.18753696,0.4176924){\makebox(0,0)[lt]{\lineheight{1.45000005}\smash{\begin{tabular}[t]{l}$i$\end{tabular}}}}%
    \put(0,0){\includegraphics[width=\unitlength,page=1]{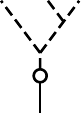}}%
  \end{picture}%
\endgroup%
},\\
    \centre{
\begingroup%
  \makeatletter%
  \providecommand\color[2][]{%
    \errmessage{(Inkscape) Color is used for the text in Inkscape, but the package 'color.sty' is not loaded}%
    \renewcommand\color[2][]{}%
  }%
  \providecommand\transparent[1]{%
    \errmessage{(Inkscape) Transparency is used (non-zero) for the text in Inkscape, but the package 'transparent.sty' is not loaded}%
    \renewcommand\transparent[1]{}%
  }%
  \providecommand\rotatebox[2]{#2}%
  \newcommand*\fsize{\dimexpr\f@size pt\relax}%
  \newcommand*\lineheight[1]{\fontsize{\fsize}{#1\fsize}\selectfont}%
  \ifx\svgwidth\undefined%
    \setlength{\unitlength}{49.69710912bp}%
    \ifx\svgscale\undefined%
      \relax%
    \else%
      \setlength{\unitlength}{\unitlength * \real{\svgscale}}%
    \fi%
  \else%
    \setlength{\unitlength}{\svgwidth}%
  \fi%
  \global\let\svgwidth\undefined%
  \global\let\svgscale\undefined%
  \makeatother%
  \begin{picture}(1,1.14694805)%
    \lineheight{1}%
    \setlength\tabcolsep{0pt}%
    \put(0,0){\includegraphics[width=\unitlength,page=1]{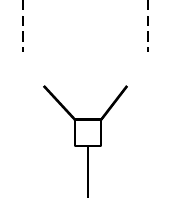}}%
    \put(-0.0030458,0.79055211){\makebox(0,0)[lt]{\lineheight{1.45000005}\smash{\begin{tabular}[t]{l}$i$\end{tabular}}}}%
    \put(0.94049889,0.7947275){\makebox(0,0)[lt]{\lineheight{1.45000005}\smash{\begin{tabular}[t]{l}$i$\end{tabular}}}}%
    \put(0,0){\includegraphics[width=\unitlength,page=2]{comm22.pdf}}%
  \end{picture}%
\endgroup%
}=\centre{
\begingroup%
  \makeatletter%
  \providecommand\color[2][]{%
    \errmessage{(Inkscape) Color is used for the text in Inkscape, but the package 'color.sty' is not loaded}%
    \renewcommand\color[2][]{}%
  }%
  \providecommand\transparent[1]{%
    \errmessage{(Inkscape) Transparency is used (non-zero) for the text in Inkscape, but the package 'transparent.sty' is not loaded}%
    \renewcommand\transparent[1]{}%
  }%
  \providecommand\rotatebox[2]{#2}%
  \newcommand*\fsize{\dimexpr\f@size pt\relax}%
  \newcommand*\lineheight[1]{\fontsize{\fsize}{#1\fsize}\selectfont}%
  \ifx\svgwidth\undefined%
    \setlength{\unitlength}{29.99498035bp}%
    \ifx\svgscale\undefined%
      \relax%
    \else%
      \setlength{\unitlength}{\unitlength * \real{\svgscale}}%
    \fi%
  \else%
    \setlength{\unitlength}{\svgwidth}%
  \fi%
  \global\let\svgwidth\undefined%
  \global\let\svgscale\undefined%
  \makeatother%
  \begin{picture}(1,1.15019249)%
    \lineheight{1}%
    \setlength\tabcolsep{0pt}%
    \put(0.14692804,0.47170206){\makebox(0,0)[lt]{\lineheight{1.45000005}\smash{\begin{tabular}[t]{l}$i$\end{tabular}}}}%
    \put(0,0){\includegraphics[width=\unitlength,page=1]{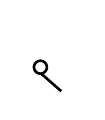}}%
    \put(0.90141574,0.47170206){\makebox(0,0)[lt]{\lineheight{1.45000005}\smash{\begin{tabular}[t]{l}$i$\end{tabular}}}}%
    \put(0,0){\includegraphics[width=\unitlength,page=2]{comm22r1.pdf}}%
  \end{picture}%
\endgroup%
}-\centre{
\begingroup%
  \makeatletter%
  \providecommand\color[2][]{%
    \errmessage{(Inkscape) Color is used for the text in Inkscape, but the package 'color.sty' is not loaded}%
    \renewcommand\color[2][]{}%
  }%
  \providecommand\transparent[1]{%
    \errmessage{(Inkscape) Transparency is used (non-zero) for the text in Inkscape, but the package 'transparent.sty' is not loaded}%
    \renewcommand\transparent[1]{}%
  }%
  \providecommand\rotatebox[2]{#2}%
  \newcommand*\fsize{\dimexpr\f@size pt\relax}%
  \newcommand*\lineheight[1]{\fontsize{\fsize}{#1\fsize}\selectfont}%
  \ifx\svgwidth\undefined%
    \setlength{\unitlength}{28.4949794bp}%
    \ifx\svgscale\undefined%
      \relax%
    \else%
      \setlength{\unitlength}{\unitlength * \real{\svgscale}}%
    \fi%
  \else%
    \setlength{\unitlength}{\svgwidth}%
  \fi%
  \global\let\svgwidth\undefined%
  \global\let\svgscale\undefined%
  \makeatother%
  \begin{picture}(1,1.21073964)%
    \lineheight{1}%
    \setlength\tabcolsep{0pt}%
    \put(0.10202157,0.49653288){\makebox(0,0)[lt]{\lineheight{1.45000005}\smash{\begin{tabular}[t]{l}$i$\end{tabular}}}}%
    \put(0,0){\includegraphics[width=\unitlength,page=1]{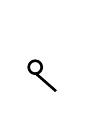}}%
    \put(0.89622618,0.49653288){\makebox(0,0)[lt]{\lineheight{1.45000005}\smash{\begin{tabular}[t]{l}$i$\end{tabular}}}}%
    \put(0,0){\includegraphics[width=\unitlength,page=2]{comm22r2.pdf}}%
  \end{picture}%
\endgroup%
}-\centre{
\begingroup%
  \makeatletter%
  \providecommand\color[2][]{%
    \errmessage{(Inkscape) Color is used for the text in Inkscape, but the package 'color.sty' is not loaded}%
    \renewcommand\color[2][]{}%
  }%
  \providecommand\transparent[1]{%
    \errmessage{(Inkscape) Transparency is used (non-zero) for the text in Inkscape, but the package 'transparent.sty' is not loaded}%
    \renewcommand\transparent[1]{}%
  }%
  \providecommand\rotatebox[2]{#2}%
  \newcommand*\fsize{\dimexpr\f@size pt\relax}%
  \newcommand*\lineheight[1]{\fontsize{\fsize}{#1\fsize}\selectfont}%
  \ifx\svgwidth\undefined%
    \setlength{\unitlength}{28.4949794bp}%
    \ifx\svgscale\undefined%
      \relax%
    \else%
      \setlength{\unitlength}{\unitlength * \real{\svgscale}}%
    \fi%
  \else%
    \setlength{\unitlength}{\svgwidth}%
  \fi%
  \global\let\svgwidth\undefined%
  \global\let\svgscale\undefined%
  \makeatother%
  \begin{picture}(1,1.26338048)%
    \lineheight{1}%
    \setlength\tabcolsep{0pt}%
    \put(0.10202157,0.49653288){\makebox(0,0)[lt]{\lineheight{1.45000005}\smash{\begin{tabular}[t]{l}$i$\end{tabular}}}}%
    \put(0,0){\includegraphics[width=\unitlength,page=1]{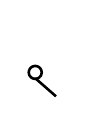}}%
    \put(0.89622618,0.49653288){\makebox(0,0)[lt]{\lineheight{1.45000005}\smash{\begin{tabular}[t]{l}$i$\end{tabular}}}}%
    \put(0,0){\includegraphics[width=\unitlength,page=2]{comm22r3.pdf}}%
  \end{picture}%
\endgroup%
}.
  \end{gather*}

  \begin{proof}[Proof of Lemma \ref{l721}]
    By using Lemma \ref{l723} (\ref{l7232}) and $\centre{
\begingroup%
  \makeatletter%
  \providecommand\color[2][]{%
    \errmessage{(Inkscape) Color is used for the text in Inkscape, but the package 'color.sty' is not loaded}%
    \renewcommand\color[2][]{}%
  }%
  \providecommand\transparent[1]{%
    \errmessage{(Inkscape) Transparency is used (non-zero) for the text in Inkscape, but the package 'transparent.sty' is not loaded}%
    \renewcommand\transparent[1]{}%
  }%
  \providecommand\rotatebox[2]{#2}%
  \newcommand*\fsize{\dimexpr\f@size pt\relax}%
  \newcommand*\lineheight[1]{\fontsize{\fsize}{#1\fsize}\selectfont}%
  \ifx\svgwidth\undefined%
    \setlength{\unitlength}{42.75574994bp}%
    \ifx\svgscale\undefined%
      \relax%
    \else%
      \setlength{\unitlength}{\unitlength * \real{\svgscale}}%
    \fi%
  \else%
    \setlength{\unitlength}{\svgwidth}%
  \fi%
  \global\let\svgwidth\undefined%
  \global\let\svgscale\undefined%
  \makeatother%
  \begin{picture}(1,0.776825)%
    \lineheight{1}%
    \setlength\tabcolsep{0pt}%
    \put(0,0){\includegraphics[width=\unitlength,page=1]{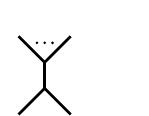}}%
    \put(0.3549849,0.31177703){\makebox(0,0)[lt]{\lineheight{1.45000005}\smash{\begin{tabular}[t]{l}$\mu^{[n]}$\end{tabular}}}}%
    \put(-0.00708055,0.12682643){\makebox(0,0)[lt]{\lineheight{1.45000005}\smash{\begin{tabular}[t]{l}$\Delta$\end{tabular}}}}%
    \put(0,0){\includegraphics[width=\unitlength,page=2]{deltamu.pdf}}%
    \put(0.21332821,0.70483182){\makebox(0,0)[lt]{\lineheight{1.45000005}\smash{\begin{tabular}[t]{l}$n$\end{tabular}}}}%
  \end{picture}%
\endgroup%
}=\centre{
\begingroup%
  \makeatletter%
  \providecommand\color[2][]{%
    \errmessage{(Inkscape) Color is used for the text in Inkscape, but the package 'color.sty' is not loaded}%
    \renewcommand\color[2][]{}%
  }%
  \providecommand\transparent[1]{%
    \errmessage{(Inkscape) Transparency is used (non-zero) for the text in Inkscape, but the package 'transparent.sty' is not loaded}%
    \renewcommand\transparent[1]{}%
  }%
  \providecommand\rotatebox[2]{#2}%
  \newcommand*\fsize{\dimexpr\f@size pt\relax}%
  \newcommand*\lineheight[1]{\fontsize{\fsize}{#1\fsize}\selectfont}%
  \ifx\svgwidth\undefined%
    \setlength{\unitlength}{29.19317395bp}%
    \ifx\svgscale\undefined%
      \relax%
    \else%
      \setlength{\unitlength}{\unitlength * \real{\svgscale}}%
    \fi%
  \else%
    \setlength{\unitlength}{\svgwidth}%
  \fi%
  \global\let\svgwidth\undefined%
  \global\let\svgscale\undefined%
  \makeatother%
  \begin{picture}(1,1.12863955)%
    \lineheight{1}%
    \setlength\tabcolsep{0pt}%
    \put(0,0){\includegraphics[width=\unitlength,page=1]{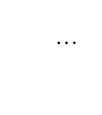}}%
    \put(0.0401983,0.08423301){\makebox(0,0)[lt]{\lineheight{1.45000005}\smash{\begin{tabular}[t]{l}$\mu$\end{tabular}}}}%
    \put(-0.01037003,0.5399708){\makebox(0,0)[lt]{\lineheight{1.45000005}\smash{\begin{tabular}[t]{l}$\Delta$\end{tabular}}}}%
    \put(0,0){\includegraphics[width=\unitlength,page=2]{mudelta.pdf}}%
    \put(0.52759836,1.02319975){\makebox(0,0)[lt]{\lineheight{1.45000005}\smash{\begin{tabular}[t]{l}$n$\end{tabular}}}}%
    \put(0,0){\includegraphics[width=\unitlength,page=3]{mudelta.pdf}}%
  \end{picture}%
\endgroup%
}$,
    we have
    $$
      \centre{}=\sum_{\sigma:(p,q)-\text{shuffle},p+q=n}
      \centre{
\begingroup%
  \makeatletter%
  \providecommand\color[2][]{%
    \errmessage{(Inkscape) Color is used for the text in Inkscape, but the package 'color.sty' is not loaded}%
    \renewcommand\color[2][]{}%
  }%
  \providecommand\transparent[1]{%
    \errmessage{(Inkscape) Transparency is used (non-zero) for the text in Inkscape, but the package 'transparent.sty' is not loaded}%
    \renewcommand\transparent[1]{}%
  }%
  \providecommand\rotatebox[2]{#2}%
  \newcommand*\fsize{\dimexpr\f@size pt\relax}%
  \newcommand*\lineheight[1]{\fontsize{\fsize}{#1\fsize}\selectfont}%
  \ifx\svgwidth\undefined%
    \setlength{\unitlength}{108.3096549bp}%
    \ifx\svgscale\undefined%
      \relax%
    \else%
      \setlength{\unitlength}{\unitlength * \real{\svgscale}}%
    \fi%
  \else%
    \setlength{\unitlength}{\svgwidth}%
  \fi%
  \global\let\svgwidth\undefined%
  \global\let\svgscale\undefined%
  \makeatother%
  \begin{picture}(1,0.97120876)%
    \lineheight{1}%
    \setlength\tabcolsep{0pt}%
    \put(0,0){\includegraphics[width=\unitlength,page=1]{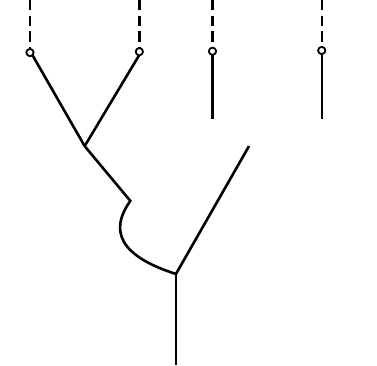}}%
    \put(-0.00088522,0.82071427){\makebox(0,0)[lt]{\lineheight{1.45000005}\smash{\begin{tabular}[t]{l}\small$i$\end{tabular}}}}%
    \put(0.88138105,0.82821777){\makebox(0,0)[lt]{\lineheight{1.45000005}\smash{\begin{tabular}[t]{l}\small$i$\end{tabular}}}}%
    \put(0.06612448,0.52546528){\makebox(0,0)[lt]{\lineheight{1.45000005}\smash{\begin{tabular}[t]{l}\small$\mu^{[m]}$\end{tabular}}}}%
    \put(0,0){\includegraphics[width=\unitlength,page=2]{l72114.pdf}}%
    \put(0.4178295,0.4323894){\makebox(0,0)[lt]{\lineheight{1.45000005}\smash{\begin{tabular}[t]{l}\tiny$S$\end{tabular}}}}%
    \put(0,0){\includegraphics[width=\unitlength,page=3]{l72114.pdf}}%
    \put(0.66008477,0.27297741){\makebox(0,0)[lt]{\lineheight{1.45000005}\smash{\begin{tabular}[t]{l}\tiny$S$\end{tabular}}}}%
    \put(0,0){\includegraphics[width=\unitlength,page=4]{l72114.pdf}}%
    \put(0.66963628,0.6040813){\makebox(0,0)[lt]{\lineheight{1.45000005}\smash{\begin{tabular}[t]{l}\small$\sigma$\end{tabular}}}}%
    \put(0,0){\includegraphics[width=\unitlength,page=5]{l72114.pdf}}%
    \put(0.57897636,0.54427937){\makebox(0,0)[lt]{\lineheight{1.45000005}\smash{\begin{tabular}[t]{l}\tiny$p$\end{tabular}}}}%
    \put(0.77378967,0.53991878){\makebox(0,0)[lt]{\lineheight{1.45000005}\smash{\begin{tabular}[t]{l}\tiny$q$\end{tabular}}}}%
  \end{picture}%
\endgroup%
}.
    $$
    By Lemma \ref{l723} (\ref{l7231}), it suffices to consider $D=\centre{
\begingroup%
  \makeatletter%
  \providecommand\color[2][]{%
    \errmessage{(Inkscape) Color is used for the text in Inkscape, but the package 'color.sty' is not loaded}%
    \renewcommand\color[2][]{}%
  }%
  \providecommand\transparent[1]{%
    \errmessage{(Inkscape) Transparency is used (non-zero) for the text in Inkscape, but the package 'transparent.sty' is not loaded}%
    \renewcommand\transparent[1]{}%
  }%
  \providecommand\rotatebox[2]{#2}%
  \newcommand*\fsize{\dimexpr\f@size pt\relax}%
  \newcommand*\lineheight[1]{\fontsize{\fsize}{#1\fsize}\selectfont}%
  \ifx\svgwidth\undefined%
    \setlength{\unitlength}{98.35891973bp}%
    \ifx\svgscale\undefined%
      \relax%
    \else%
      \setlength{\unitlength}{\unitlength * \real{\svgscale}}%
    \fi%
  \else%
    \setlength{\unitlength}{\svgwidth}%
  \fi%
  \global\let\svgwidth\undefined%
  \global\let\svgscale\undefined%
  \makeatother%
  \begin{picture}(1,1.13158016)%
    \lineheight{1}%
    \setlength\tabcolsep{0pt}%
    \put(0,0){\includegraphics[width=\unitlength,page=1]{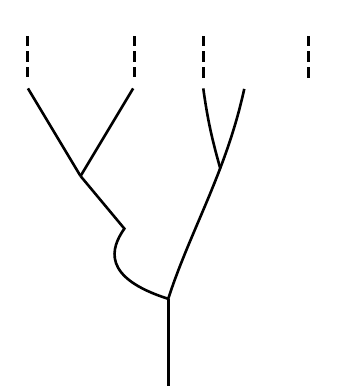}}%
    \put(-0.00307785,0.8676644){\makebox(0,0)[lt]{\lineheight{1.45000005}\smash{\begin{tabular}[t]{l}$i$\end{tabular}}}}%
    \put(0.93987268,0.87049075){\makebox(0,0)[lt]{\lineheight{1.45000005}\smash{\begin{tabular}[t]{l}$i$\end{tabular}}}}%
    \put(0.26696152,0.608588){\makebox(0,0)[lt]{\lineheight{1.45000005}\smash{\begin{tabular}[t]{l}$\mu^{[m]}$\end{tabular}}}}%
    \put(0,0){\includegraphics[width=\unitlength,page=2]{l72119.pdf}}%
    \put(0.41858766,0.46492852){\makebox(0,0)[lt]{\lineheight{1.45000005}\smash{\begin{tabular}[t]{l}$S$\end{tabular}}}}%
    \put(0,0){\includegraphics[width=\unitlength,page=3]{l72119.pdf}}%
    \put(0.62923734,1.10023132){\makebox(0,0)[lt]{\lineheight{1.45000005}\smash{\begin{tabular}[t]{l}$p$\end{tabular}}}}%
    \put(0.84925235,1.10028536){\makebox(0,0)[lt]{\lineheight{1.45000005}\smash{\begin{tabular}[t]{l}$q$\end{tabular}}}}%
    \put(0,0){\includegraphics[width=\unitlength,page=4]{l72119.pdf}}%
    \put(0.18762471,1.09738795){\makebox(0,0)[lt]{\lineheight{1.45000005}\smash{\begin{tabular}[t]{l}$m$\end{tabular}}}}%
    \put(0,0){\includegraphics[width=\unitlength,page=5]{l72119.pdf}}%
  \end{picture}%
\endgroup%
}$.
    By Lemma \ref{l723} (\ref{l7233}), we have  $\centre{
\begingroup%
  \makeatletter%
  \providecommand\color[2][]{%
    \errmessage{(Inkscape) Color is used for the text in Inkscape, but the package 'color.sty' is not loaded}%
    \renewcommand\color[2][]{}%
  }%
  \providecommand\transparent[1]{%
    \errmessage{(Inkscape) Transparency is used (non-zero) for the text in Inkscape, but the package 'transparent.sty' is not loaded}%
    \renewcommand\transparent[1]{}%
  }%
  \providecommand\rotatebox[2]{#2}%
  \newcommand*\fsize{\dimexpr\f@size pt\relax}%
  \newcommand*\lineheight[1]{\fontsize{\fsize}{#1\fsize}\selectfont}%
  \ifx\svgwidth\undefined%
    \setlength{\unitlength}{31.23534272bp}%
    \ifx\svgscale\undefined%
      \relax%
    \else%
      \setlength{\unitlength}{\unitlength * \real{\svgscale}}%
    \fi%
  \else%
    \setlength{\unitlength}{\svgwidth}%
  \fi%
  \global\let\svgwidth\undefined%
  \global\let\svgscale\undefined%
  \makeatother%
  \begin{picture}(1,1.36864198)%
    \lineheight{1}%
    \setlength\tabcolsep{0pt}%
    \put(0,0){\includegraphics[width=\unitlength,page=1]{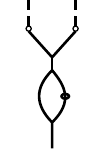}}%
    \put(-0.00969204,1.05518709){\makebox(0,0)[lt]{\lineheight{1.45000005}\smash{\begin{tabular}[t]{l}$i$\end{tabular}}}}%
    \put(0.78504094,1.04301078){\makebox(0,0)[lt]{\lineheight{1.45000005}\smash{\begin{tabular}[t]{l}$i$\end{tabular}}}}%
    \put(0.58205429,0.81166039){\makebox(0,0)[lt]{\lineheight{1.45000005}\smash{\begin{tabular}[t]{l}$\mu$\end{tabular}}}}%
    \put(0.08167744,0.64195396){\makebox(0,0)[lt]{\lineheight{1.45000005}\smash{\begin{tabular}[t]{l}$\Delta$\end{tabular}}}}%
    \put(0.66424454,0.43760817){\makebox(0,0)[lt]{\lineheight{1.45000005}\smash{\begin{tabular}[t]{l}$S$\end{tabular}}}}%
    \put(0.5441882,0.1718058){\makebox(0,0)[lt]{\lineheight{1.45000005}\smash{\begin{tabular}[t]{l}$\mu$\end{tabular}}}}%
    \put(0,0){\includegraphics[width=\unitlength,page=2]{l72116.pdf}}%
  \end{picture}%
\endgroup%
}=\centre{
\begingroup%
  \makeatletter%
  \providecommand\color[2][]{%
    \errmessage{(Inkscape) Color is used for the text in Inkscape, but the package 'color.sty' is not loaded}%
    \renewcommand\color[2][]{}%
  }%
  \providecommand\transparent[1]{%
    \errmessage{(Inkscape) Transparency is used (non-zero) for the text in Inkscape, but the package 'transparent.sty' is not loaded}%
    \renewcommand\transparent[1]{}%
  }%
  \providecommand\rotatebox[2]{#2}%
  \newcommand*\fsize{\dimexpr\f@size pt\relax}%
  \newcommand*\lineheight[1]{\fontsize{\fsize}{#1\fsize}\selectfont}%
  \ifx\svgwidth\undefined%
    \setlength{\unitlength}{31.23534272bp}%
    \ifx\svgscale\undefined%
      \relax%
    \else%
      \setlength{\unitlength}{\unitlength * \real{\svgscale}}%
    \fi%
  \else%
    \setlength{\unitlength}{\svgwidth}%
  \fi%
  \global\let\svgwidth\undefined%
  \global\let\svgscale\undefined%
  \makeatother%
  \begin{picture}(1,1.36864198)%
    \lineheight{1}%
    \setlength\tabcolsep{0pt}%
    \put(0,0){\includegraphics[width=\unitlength,page=1]{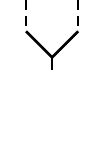}}%
    \put(-0.00969204,1.05518709){\makebox(0,0)[lt]{\lineheight{1.45000005}\smash{\begin{tabular}[t]{l}$i$\end{tabular}}}}%
    \put(0.78504094,1.04301078){\makebox(0,0)[lt]{\lineheight{1.45000005}\smash{\begin{tabular}[t]{l}$i$\end{tabular}}}}%
    \put(0.58205429,0.81166039){\makebox(0,0)[lt]{\lineheight{1.45000005}\smash{\begin{tabular}[t]{l}$\mu$\end{tabular}}}}%
    \put(0,0){\includegraphics[width=\unitlength,page=2]{l72117.pdf}}%
    \put(0.13511833,0.65843746){\makebox(0,0)[lt]{\lineheight{1.45000005}\smash{\begin{tabular}[t]{l}$\epsilon$\end{tabular}}}}%
    \put(0.15313485,0.43160933){\makebox(0,0)[lt]{\lineheight{1.45000005}\smash{\begin{tabular}[t]{l}$\eta$\end{tabular}}}}%
    \put(0,0){\includegraphics[width=\unitlength,page=3]{l72117.pdf}}%
  \end{picture}%
\endgroup%
}=\centre{
\begingroup%
  \makeatletter%
  \providecommand\color[2][]{%
    \errmessage{(Inkscape) Color is used for the text in Inkscape, but the package 'color.sty' is not loaded}%
    \renewcommand\color[2][]{}%
  }%
  \providecommand\transparent[1]{%
    \errmessage{(Inkscape) Transparency is used (non-zero) for the text in Inkscape, but the package 'transparent.sty' is not loaded}%
    \renewcommand\transparent[1]{}%
  }%
  \providecommand\rotatebox[2]{#2}%
  \newcommand*\fsize{\dimexpr\f@size pt\relax}%
  \newcommand*\lineheight[1]{\fontsize{\fsize}{#1\fsize}\selectfont}%
  \ifx\svgwidth\undefined%
    \setlength{\unitlength}{45.71011522bp}%
    \ifx\svgscale\undefined%
      \relax%
    \else%
      \setlength{\unitlength}{\unitlength * \real{\svgscale}}%
    \fi%
  \else%
    \setlength{\unitlength}{\svgwidth}%
  \fi%
  \global\let\svgwidth\undefined%
  \global\let\svgscale\undefined%
  \makeatother%
  \begin{picture}(1,0.9352416)%
    \lineheight{1}%
    \setlength\tabcolsep{0pt}%
    \put(0,0){\includegraphics[width=\unitlength,page=1]{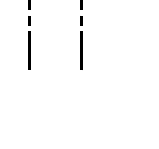}}%
    \put(0.01337448,0.72104676){\makebox(0,0)[lt]{\lineheight{1.45000005}\smash{\begin{tabular}[t]{l}$i$\end{tabular}}}}%
    \put(0.55644371,0.71272625){\makebox(0,0)[lt]{\lineheight{1.45000005}\smash{\begin{tabular}[t]{l}$i$\end{tabular}}}}%
    \put(0,0){\includegraphics[width=\unitlength,page=2]{l72118.pdf}}%
    \put(-0.00662291,0.46599008){\makebox(0,0)[lt]{\lineheight{1.45000005}\smash{\begin{tabular}[t]{l}$\epsilon$\end{tabular}}}}%
    \put(0.13784475,0.29732621){\makebox(0,0)[lt]{\lineheight{1.45000005}\smash{\begin{tabular}[t]{l}$\eta$\end{tabular}}}}%
    \put(0,0){\includegraphics[width=\unitlength,page=3]{l72118.pdf}}%
    \put(0.54588177,0.46250042){\makebox(0,0)[lt]{\lineheight{1.45000005}\smash{\begin{tabular}[t]{l}$\epsilon$\end{tabular}}}}%
    \put(0,0){\includegraphics[width=\unitlength,page=4]{l72118.pdf}}%
  \end{picture}%
\endgroup%
}=0.$
    Thus, when $p=0$, we have $D=0$.
    When $p\geq 1$, by Lemma \ref{l723} (\ref{l7234}), we have
    \begin{gather*}
     \begin{split}
       D&=\scalebox{0.95}{$\centre{
\begingroup%
  \makeatletter%
  \providecommand\color[2][]{%
    \errmessage{(Inkscape) Color is used for the text in Inkscape, but the package 'color.sty' is not loaded}%
    \renewcommand\color[2][]{}%
  }%
  \providecommand\transparent[1]{%
    \errmessage{(Inkscape) Transparency is used (non-zero) for the text in Inkscape, but the package 'transparent.sty' is not loaded}%
    \renewcommand\transparent[1]{}%
  }%
  \providecommand\rotatebox[2]{#2}%
  \newcommand*\fsize{\dimexpr\f@size pt\relax}%
  \newcommand*\lineheight[1]{\fontsize{\fsize}{#1\fsize}\selectfont}%
  \ifx\svgwidth\undefined%
    \setlength{\unitlength}{98.35891973bp}%
    \ifx\svgscale\undefined%
      \relax%
    \else%
      \setlength{\unitlength}{\unitlength * \real{\svgscale}}%
    \fi%
  \else%
    \setlength{\unitlength}{\svgwidth}%
  \fi%
  \global\let\svgwidth\undefined%
  \global\let\svgscale\undefined%
  \makeatother%
  \begin{picture}(1,1.14533211)%
    \lineheight{1}%
    \setlength\tabcolsep{0pt}%
    \put(0,0){\includegraphics[width=\unitlength,page=1]{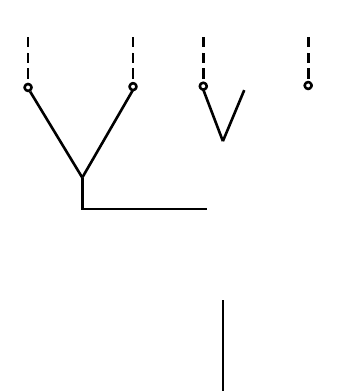}}%
    \put(-0.00307785,0.87785214){\makebox(0,0)[lt]{\lineheight{1.45000005}\smash{\begin{tabular}[t]{l}$i$\end{tabular}}}}%
    \put(0.93987268,0.8806785){\makebox(0,0)[lt]{\lineheight{1.45000005}\smash{\begin{tabular}[t]{l}$i$\end{tabular}}}}%
    \put(0,0){\includegraphics[width=\unitlength,page=2]{l721110.pdf}}%
    \put(0.62655596,1.11189044){\makebox(0,0)[lt]{\lineheight{1.45000005}\smash{\begin{tabular}[t]{l}$p$\end{tabular}}}}%
    \put(0.84671807,1.11403731){\makebox(0,0)[lt]{\lineheight{1.45000005}\smash{\begin{tabular}[t]{l}$q$\end{tabular}}}}%
    \put(0,0){\includegraphics[width=\unitlength,page=3]{l721110.pdf}}%
    \put(0.18762471,1.1075757){\makebox(0,0)[lt]{\lineheight{1.45000005}\smash{\begin{tabular}[t]{l}$m$\end{tabular}}}}%
    \put(0,0){\includegraphics[width=\unitlength,page=4]{l721110.pdf}}%
    \put(0.4201634,0.44447938){\makebox(0,0)[lt]{\lineheight{1.45000005}\smash{\begin{tabular}[t]{l}$ad_H$\end{tabular}}}}%
    \put(0,0){\includegraphics[width=\unitlength,page=5]{l721110.pdf}}%
  \end{picture}%
\endgroup%
}$}
       =\scalebox{0.95}{$\centre{
\begingroup%
  \makeatletter%
  \providecommand\color[2][]{%
    \errmessage{(Inkscape) Color is used for the text in Inkscape, but the package 'color.sty' is not loaded}%
    \renewcommand\color[2][]{}%
  }%
  \providecommand\transparent[1]{%
    \errmessage{(Inkscape) Transparency is used (non-zero) for the text in Inkscape, but the package 'transparent.sty' is not loaded}%
    \renewcommand\transparent[1]{}%
  }%
  \providecommand\rotatebox[2]{#2}%
  \newcommand*\fsize{\dimexpr\f@size pt\relax}%
  \newcommand*\lineheight[1]{\fontsize{\fsize}{#1\fsize}\selectfont}%
  \ifx\svgwidth\undefined%
    \setlength{\unitlength}{98.35891973bp}%
    \ifx\svgscale\undefined%
      \relax%
    \else%
      \setlength{\unitlength}{\unitlength * \real{\svgscale}}%
    \fi%
  \else%
    \setlength{\unitlength}{\svgwidth}%
  \fi%
  \global\let\svgwidth\undefined%
  \global\let\svgscale\undefined%
  \makeatother%
  \begin{picture}(1,1.14679808)%
    \lineheight{1}%
    \setlength\tabcolsep{0pt}%
    \put(0,0){\includegraphics[width=\unitlength,page=1]{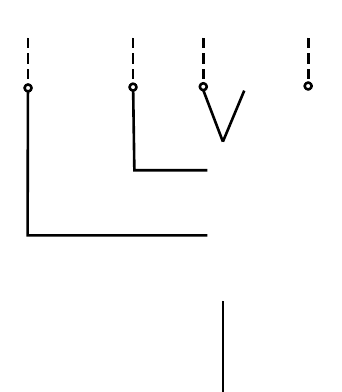}}%
    \put(-0.00307785,0.87785214){\makebox(0,0)[lt]{\lineheight{1.45000005}\smash{\begin{tabular}[t]{l}$i$\end{tabular}}}}%
    \put(0.93987268,0.8806785){\makebox(0,0)[lt]{\lineheight{1.45000005}\smash{\begin{tabular}[t]{l}$i$\end{tabular}}}}%
    \put(0,0){\includegraphics[width=\unitlength,page=2]{l721111.pdf}}%
    \put(0.62242947,1.11550328){\makebox(0,0)[lt]{\lineheight{1.45000005}\smash{\begin{tabular}[t]{l}$p$\end{tabular}}}}%
    \put(0.83954884,1.11522562){\makebox(0,0)[lt]{\lineheight{1.45000005}\smash{\begin{tabular}[t]{l}$q$\end{tabular}}}}%
    \put(0,0){\includegraphics[width=\unitlength,page=3]{l721111.pdf}}%
    \put(0.18762471,1.1075757){\makebox(0,0)[lt]{\lineheight{1.45000005}\smash{\begin{tabular}[t]{l}$m$\end{tabular}}}}%
    \put(0,0){\includegraphics[width=\unitlength,page=4]{l721111.pdf}}%
    \put(0.40414899,0.38421708){\makebox(0,0)[lt]{\lineheight{1.45000005}\smash{\begin{tabular}[t]{l}$ad_H$\end{tabular}}}}%
    \put(0,0){\includegraphics[width=\unitlength,page=5]{l721111.pdf}}%
    \put(0.40497962,0.57440181){\makebox(0,0)[lt]{\lineheight{1.45000005}\smash{\begin{tabular}[t]{l}$ad_H$\end{tabular}}}}%
  \end{picture}%
\endgroup%
}$}
       =\scalebox{0.95}{$\centre{
\begingroup%
  \makeatletter%
  \providecommand\color[2][]{%
    \errmessage{(Inkscape) Color is used for the text in Inkscape, but the package 'color.sty' is not loaded}%
    \renewcommand\color[2][]{}%
  }%
  \providecommand\transparent[1]{%
    \errmessage{(Inkscape) Transparency is used (non-zero) for the text in Inkscape, but the package 'transparent.sty' is not loaded}%
    \renewcommand\transparent[1]{}%
  }%
  \providecommand\rotatebox[2]{#2}%
  \newcommand*\fsize{\dimexpr\f@size pt\relax}%
  \newcommand*\lineheight[1]{\fontsize{\fsize}{#1\fsize}\selectfont}%
  \ifx\svgwidth\undefined%
    \setlength{\unitlength}{98.35891973bp}%
    \ifx\svgscale\undefined%
      \relax%
    \else%
      \setlength{\unitlength}{\unitlength * \real{\svgscale}}%
    \fi%
  \else%
    \setlength{\unitlength}{\svgwidth}%
  \fi%
  \global\let\svgwidth\undefined%
  \global\let\svgscale\undefined%
  \makeatother%
  \begin{picture}(1,1.1388705)%
    \lineheight{1}%
    \setlength\tabcolsep{0pt}%
    \put(0,0){\includegraphics[width=\unitlength,page=1]{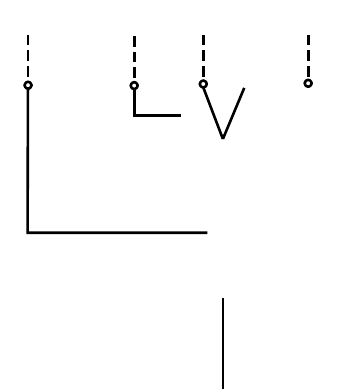}}%
    \put(-0.00307785,0.87785214){\makebox(0,0)[lt]{\lineheight{1.45000005}\smash{\begin{tabular}[t]{l}$i$\end{tabular}}}}%
    \put(0.93987268,0.8806785){\makebox(0,0)[lt]{\lineheight{1.45000005}\smash{\begin{tabular}[t]{l}$i$\end{tabular}}}}%
    \put(0,0){\includegraphics[width=\unitlength,page=2]{l721112.pdf}}%
    \put(0.62923734,1.0984835){\makebox(0,0)[lt]{\lineheight{1.45000005}\smash{\begin{tabular}[t]{l}$p$\end{tabular}}}}%
    \put(0.84805876,1.10331177){\makebox(0,0)[lt]{\lineheight{1.45000005}\smash{\begin{tabular}[t]{l}$q$\end{tabular}}}}%
    \put(0,0){\includegraphics[width=\unitlength,page=3]{l721112.pdf}}%
    \put(0.18762471,1.1075757){\makebox(0,0)[lt]{\lineheight{1.45000005}\smash{\begin{tabular}[t]{l}$m$\end{tabular}}}}%
    \put(0,0){\includegraphics[width=\unitlength,page=4]{l721112.pdf}}%
    \put(0.40510677,0.38892739){\makebox(0,0)[lt]{\lineheight{1.45000005}\smash{\begin{tabular}[t]{l}$ad_H$\end{tabular}}}}%
    \put(0,0){\includegraphics[width=\unitlength,page=5]{l721112.pdf}}%
    \put(0.4035716,0.71807355){\makebox(0,0)[lt]{\lineheight{1.45000005}\smash{\begin{tabular}[t]{l}$ad_H$\end{tabular}}}}%
    \put(0,0){\includegraphics[width=\unitlength,page=6]{l721112.pdf}}%
    \put(0.40516333,0.53725503){\makebox(0,0)[lt]{\lineheight{1.45000005}\smash{\begin{tabular}[t]{l}$ad_H$\end{tabular}}}}%
  \end{picture}%
\endgroup%
}$}
       \\
       &=\scalebox{0.95}{$\centre{
\begingroup%
  \makeatletter%
  \providecommand\color[2][]{%
    \errmessage{(Inkscape) Color is used for the text in Inkscape, but the package 'color.sty' is not loaded}%
    \renewcommand\color[2][]{}%
  }%
  \providecommand\transparent[1]{%
    \errmessage{(Inkscape) Transparency is used (non-zero) for the text in Inkscape, but the package 'transparent.sty' is not loaded}%
    \renewcommand\transparent[1]{}%
  }%
  \providecommand\rotatebox[2]{#2}%
  \newcommand*\fsize{\dimexpr\f@size pt\relax}%
  \newcommand*\lineheight[1]{\fontsize{\fsize}{#1\fsize}\selectfont}%
  \ifx\svgwidth\undefined%
    \setlength{\unitlength}{98.35891973bp}%
    \ifx\svgscale\undefined%
      \relax%
    \else%
      \setlength{\unitlength}{\unitlength * \real{\svgscale}}%
    \fi%
  \else%
    \setlength{\unitlength}{\svgwidth}%
  \fi%
  \global\let\svgwidth\undefined%
  \global\let\svgscale\undefined%
  \makeatother%
  \begin{picture}(1,1.14484468)%
    \lineheight{1}%
    \setlength\tabcolsep{0pt}%
    \put(0,0){\includegraphics[width=\unitlength,page=1]{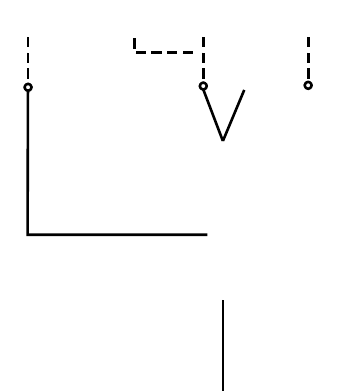}}%
    \put(-0.00307785,0.87785214){\makebox(0,0)[lt]{\lineheight{1.45000005}\smash{\begin{tabular}[t]{l}$i$\end{tabular}}}}%
    \put(0.93987268,0.8806785){\makebox(0,0)[lt]{\lineheight{1.45000005}\smash{\begin{tabular}[t]{l}$i$\end{tabular}}}}%
    \put(0,0){\includegraphics[width=\unitlength,page=2]{l721113.pdf}}%
    \put(0.61558653,1.11042797){\makebox(0,0)[lt]{\lineheight{1.45000005}\smash{\begin{tabular}[t]{l}$p$\end{tabular}}}}%
    \put(0.83440795,1.11354988){\makebox(0,0)[lt]{\lineheight{1.45000005}\smash{\begin{tabular}[t]{l}$q$\end{tabular}}}}%
    \put(0,0){\includegraphics[width=\unitlength,page=3]{l721113.pdf}}%
    \put(0.18762471,1.1075757){\makebox(0,0)[lt]{\lineheight{1.45000005}\smash{\begin{tabular}[t]{l}$m$\end{tabular}}}}%
    \put(0,0){\includegraphics[width=\unitlength,page=4]{l721113.pdf}}%
    \put(0.39933115,0.37687453){\makebox(0,0)[lt]{\lineheight{1.45000005}\smash{\begin{tabular}[t]{l}$ad_H$\end{tabular}}}}%
    \put(0,0){\includegraphics[width=\unitlength,page=5]{l721113.pdf}}%
    \put(0.40018663,0.53576544){\makebox(0,0)[lt]{\lineheight{1.45000005}\smash{\begin{tabular}[t]{l}$ad_H$\end{tabular}}}}%
    \put(0,0){\includegraphics[width=\unitlength,page=6]{l721113.pdf}}%
  \end{picture}%
\endgroup%
}$}
       =\cdots
       =\scalebox{0.95}{$\centre{
\begingroup%
  \makeatletter%
  \providecommand\color[2][]{%
    \errmessage{(Inkscape) Color is used for the text in Inkscape, but the package 'color.sty' is not loaded}%
    \renewcommand\color[2][]{}%
  }%
  \providecommand\transparent[1]{%
    \errmessage{(Inkscape) Transparency is used (non-zero) for the text in Inkscape, but the package 'transparent.sty' is not loaded}%
    \renewcommand\transparent[1]{}%
  }%
  \providecommand\rotatebox[2]{#2}%
  \newcommand*\fsize{\dimexpr\f@size pt\relax}%
  \newcommand*\lineheight[1]{\fontsize{\fsize}{#1\fsize}\selectfont}%
  \ifx\svgwidth\undefined%
    \setlength{\unitlength}{77.11321935bp}%
    \ifx\svgscale\undefined%
      \relax%
    \else%
      \setlength{\unitlength}{\unitlength * \real{\svgscale}}%
    \fi%
  \else%
    \setlength{\unitlength}{\svgwidth}%
  \fi%
  \global\let\svgwidth\undefined%
  \global\let\svgscale\undefined%
  \makeatother%
  \begin{picture}(1,1.32460301)%
    \lineheight{1}%
    \setlength\tabcolsep{0pt}%
    \put(0,0){\includegraphics[width=\unitlength,page=1]{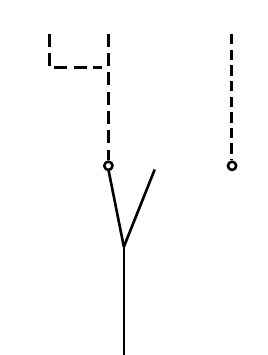}}%
    \put(0.90917594,0.68325336){\makebox(0,0)[lt]{\lineheight{1.45000005}\smash{\begin{tabular}[t]{l}$i$\end{tabular}}}}%
    \put(0,0){\includegraphics[width=\unitlength,page=2]{l721114.pdf}}%
    \put(0.45694022,1.27189789){\makebox(0,0)[lt]{\lineheight{1.45000005}\smash{\begin{tabular}[t]{l}$p$\end{tabular}}}}%
    \put(0.76923146,1.27497101){\makebox(0,0)[lt]{\lineheight{1.45000005}\smash{\begin{tabular}[t]{l}$q$\end{tabular}}}}%
    \put(0,0){\includegraphics[width=\unitlength,page=3]{l721114.pdf}}%
    \put(0.28416956,0.68268432){\makebox(0,0)[lt]{\lineheight{1.45000005}\smash{\begin{tabular}[t]{l}$i$\end{tabular}}}}%
    \put(0,0){\includegraphics[width=\unitlength,page=4]{l721114.pdf}}%
    \put(0.042509,1.2773313){\makebox(0,0)[lt]{\lineheight{1.45000005}\smash{\begin{tabular}[t]{l}$m$\end{tabular}}}}%
  \end{picture}%
\endgroup%
}$}.
     \end{split}
    \end{gather*}
    Note that the last term is a $\Z$-linear sum of unions of tree diagrams with $m$ trivalent vertices.
    Therefore, the first equality of (\ref{l7211}) follows.
    If $m=n=1$, then the equality follows from the case where $m=p=1, q=0$.
    The second equality of (\ref{l7211}) follows similarly.

    The first equality of (\ref{l7212}) follows from $\centre{
\begingroup%
  \makeatletter%
  \providecommand\color[2][]{%
    \errmessage{(Inkscape) Color is used for the text in Inkscape, but the package 'color.sty' is not loaded}%
    \renewcommand\color[2][]{}%
  }%
  \providecommand\transparent[1]{%
    \errmessage{(Inkscape) Transparency is used (non-zero) for the text in Inkscape, but the package 'transparent.sty' is not loaded}%
    \renewcommand\transparent[1]{}%
  }%
  \providecommand\rotatebox[2]{#2}%
  \newcommand*\fsize{\dimexpr\f@size pt\relax}%
  \newcommand*\lineheight[1]{\fontsize{\fsize}{#1\fsize}\selectfont}%
  \ifx\svgwidth\undefined%
    \setlength{\unitlength}{57.55368174bp}%
    \ifx\svgscale\undefined%
      \relax%
    \else%
      \setlength{\unitlength}{\unitlength * \real{\svgscale}}%
    \fi%
  \else%
    \setlength{\unitlength}{\svgwidth}%
  \fi%
  \global\let\svgwidth\undefined%
  \global\let\svgscale\undefined%
  \makeatother%
  \begin{picture}(1,0.84182288)%
    \lineheight{1}%
    \setlength\tabcolsep{0pt}%
    \put(0,0){\includegraphics[width=\unitlength,page=1]{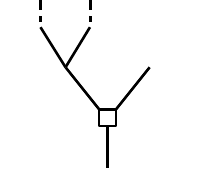}}%
    \put(0.06983426,0.70760373){\makebox(0,0)[lt]{\lineheight{1.45000005}\smash{\begin{tabular}[t]{l}$i$\end{tabular}}}}%
    \put(0.47000421,0.70311628){\makebox(0,0)[lt]{\lineheight{1.45000005}\smash{\begin{tabular}[t]{l}$i$\end{tabular}}}}%
    \put(-0.00526003,0.37990078){\makebox(0,0)[lt]{\lineheight{1.45000005}\smash{\begin{tabular}[t]{l}$\mu^{[m]}$\end{tabular}}}}%
    \put(0,0){\includegraphics[width=\unitlength,page=2]{l72123.pdf}}%
    \put(0.77412409,0.45315119){\makebox(0,0)[lt]{\lineheight{1.45000005}\smash{\begin{tabular}[t]{l}$\eta$\end{tabular}}}}%
    \put(0,0){\includegraphics[width=\unitlength,page=3]{l72123.pdf}}%
  \end{picture}%
\endgroup%
}=\centre{}=0.$
    The second equality follows similarly.

    \vspace{0.1in}
    We have (\ref{l7213}) because  $\centre{}=\centre{
\begingroup%
  \makeatletter%
  \providecommand\color[2][]{%
    \errmessage{(Inkscape) Color is used for the text in Inkscape, but the package 'color.sty' is not loaded}%
    \renewcommand\color[2][]{}%
  }%
  \providecommand\transparent[1]{%
    \errmessage{(Inkscape) Transparency is used (non-zero) for the text in Inkscape, but the package 'transparent.sty' is not loaded}%
    \renewcommand\transparent[1]{}%
  }%
  \providecommand\rotatebox[2]{#2}%
  \newcommand*\fsize{\dimexpr\f@size pt\relax}%
  \newcommand*\lineheight[1]{\fontsize{\fsize}{#1\fsize}\selectfont}%
  \ifx\svgwidth\undefined%
    \setlength{\unitlength}{37.79130196bp}%
    \ifx\svgscale\undefined%
      \relax%
    \else%
      \setlength{\unitlength}{\unitlength * \real{\svgscale}}%
    \fi%
  \else%
    \setlength{\unitlength}{\svgwidth}%
  \fi%
  \global\let\svgwidth\undefined%
  \global\let\svgscale\undefined%
  \makeatother%
  \begin{picture}(1,0.83289445)%
    \lineheight{1}%
    \setlength\tabcolsep{0pt}%
    \put(0,0){\includegraphics[width=\unitlength,page=1]{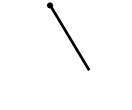}}%
    \put(-0.00801068,0.75144389){\makebox(0,0)[lt]{\lineheight{1.45000005}\smash{\begin{tabular}[t]{l}$\eta$\end{tabular}}}}%
    \put(0,0){\includegraphics[width=\unitlength,page=2]{l72132.pdf}}%
  \end{picture}%
\endgroup%
}=\quad\centre{}.$
  \end{proof}

 \subsection{Proof of Theorem \ref{th731}}\label{ss43}
  In this subsection, we prove Theorem \ref{th731}.

  For any $y_1,\cdots,y_r\in F_n$, we call $[y_1,\cdots,[y_{r-1},y_r]]\in\Gamma_r$ an \emph{$r$-fold commutator}.

  For $i\in [n]$,
  define $d_i\in \End(F_n)=\F^{\op}(n,n)$ by
  $$d_i(x_i)=[y_1,\cdots,[y_{r},y_{r+1}]]^{\epsilon},\quad d_i(x_j)=1 \quad (j\neq i)$$
  for $y_1,\cdots,y_{r+1}\in F_n, \epsilon\in\{\pm 1\}$, which we call an \emph{$(r+1)$-fold commutator at $i$}.
  Via the isomorphism $\K\F^{\op}(n,n)\cong \A_0(n,n)$, we identify $d_i\in\F^{\op}(n,n)$ with a morphism of the following form
  $$\centre{
\begingroup%
  \makeatletter%
  \providecommand\color[2][]{%
    \errmessage{(Inkscape) Color is used for the text in Inkscape, but the package 'color.sty' is not loaded}%
    \renewcommand\color[2][]{}%
  }%
  \providecommand\transparent[1]{%
    \errmessage{(Inkscape) Transparency is used (non-zero) for the text in Inkscape, but the package 'transparent.sty' is not loaded}%
    \renewcommand\transparent[1]{}%
  }%
  \providecommand\rotatebox[2]{#2}%
  \newcommand*\fsize{\dimexpr\f@size pt\relax}%
  \newcommand*\lineheight[1]{\fontsize{\fsize}{#1\fsize}\selectfont}%
  \ifx\svgwidth\undefined%
    \setlength{\unitlength}{98.34342073bp}%
    \ifx\svgscale\undefined%
      \relax%
    \else%
      \setlength{\unitlength}{\unitlength * \real{\svgscale}}%
    \fi%
  \else%
    \setlength{\unitlength}{\svgwidth}%
  \fi%
  \global\let\svgwidth\undefined%
  \global\let\svgscale\undefined%
  \makeatother%
  \begin{picture}(1,1.14548943)%
    \lineheight{1}%
    \setlength\tabcolsep{0pt}%
    \put(0,0){\includegraphics[width=\unitlength,page=1]{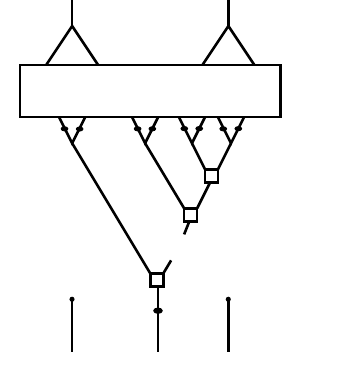}}%
    \put(0.14561517,0.30031685){\makebox(0,0)[lt]{\lineheight{1.45000005}\smash{\begin{tabular}[t]{l}$\eta$\end{tabular}}}}%
    \put(0.67480749,0.29810764){\makebox(0,0)[lt]{\lineheight{1.45000005}\smash{\begin{tabular}[t]{l}$\eta$\end{tabular}}}}%
    \put(0,0){\includegraphics[width=\unitlength,page=2]{d_i.pdf}}%
    \put(0.44813837,0.00716955){\makebox(0,0)[lt]{\lineheight{1.45000005}\smash{\begin{tabular}[t]{l}$i$\end{tabular}}}}%
    \put(0,0){\includegraphics[width=\unitlength,page=3]{d_i.pdf}}%
    \put(0.67917577,0.43407139){\makebox(0,0)[lt]{\lineheight{1.45000005}\smash{\begin{tabular}[t]{l}$r$\end{tabular}}}}%
    \put(0.20698859,0.84847741){\makebox(0,0)[lt]{\lineheight{1.45000005}\smash{\begin{tabular}[t]{l}$\text{permutation}$\end{tabular}}}}%
    \put(0.18260712,0.01004992){\makebox(0,0)[lt]{\lineheight{1.45000005}\smash{\begin{tabular}[t]{l}$1$\end{tabular}}}}%
    \put(0.64461785,0.01111461){\makebox(0,0)[lt]{\lineheight{1.45000005}\smash{\begin{tabular}[t]{l}$n$\end{tabular}}}}%
    \put(-0.00069511,1.04965253){\makebox(0,0)[lt]{\lineheight{1.45000005}\smash{\begin{tabular}[t]{l}\small{$\Delta^{[q_1]}$}\end{tabular}}}}%
    \put(0.69533619,1.05321135){\makebox(0,0)[lt]{\lineheight{1.45000005}\smash{\begin{tabular}[t]{l}\small{$\Delta^{[q_n]}$}\end{tabular}}}}%
    \put(-0.00069511,0.68536584){\makebox(0,0)[lt]{\lineheight{1.45000005}\smash{\begin{tabular}[t]{l}\small{$\mu^{[p_1]}$}\end{tabular}}}}%
    \put(0.69432076,0.68790926){\makebox(0,0)[lt]{\lineheight{1.45000005}\smash{\begin{tabular}[t]{l}\small{$\mu^{[p_{r+1}]}$}\end{tabular}}}}%
  \end{picture}%
\endgroup%
}\in\A_0(n,n),$$
  which we also call an \emph{$(r+1)$-fold commutator at $i$},
  where each $\centre{
\begingroup%
  \makeatletter%
  \providecommand\color[2][]{%
    \errmessage{(Inkscape) Color is used for the text in Inkscape, but the package 'color.sty' is not loaded}%
    \renewcommand\color[2][]{}%
  }%
  \providecommand\transparent[1]{%
    \errmessage{(Inkscape) Transparency is used (non-zero) for the text in Inkscape, but the package 'transparent.sty' is not loaded}%
    \renewcommand\transparent[1]{}%
  }%
  \providecommand\rotatebox[2]{#2}%
  \newcommand*\fsize{\dimexpr\f@size pt\relax}%
  \newcommand*\lineheight[1]{\fontsize{\fsize}{#1\fsize}\selectfont}%
  \ifx\svgwidth\undefined%
    \setlength{\unitlength}{3.37417307bp}%
    \ifx\svgscale\undefined%
      \relax%
    \else%
      \setlength{\unitlength}{\unitlength * \real{\svgscale}}%
    \fi%
  \else%
    \setlength{\unitlength}{\svgwidth}%
  \fi%
  \global\let\svgwidth\undefined%
  \global\let\svgscale\undefined%
  \makeatother%
  \begin{picture}(1,6.66830034)%
    \lineheight{1}%
    \setlength\tabcolsep{0pt}%
    \put(0,0){\includegraphics[width=\unitlength,page=1]{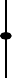}}%
  \end{picture}%
\endgroup%
}$ depicts $S$ or $\id_H$, and $q_k,p_l\geq 0$ satisfy $\sum_{k=1}^{n}q_k=\sum_{l=1}^{r+1}p_l$.

  \begin{claim}\label{c731}
   An element $g\in\jfE_r(n)$ can be written as a convolution product
   $$g= d_{1,1}\ast\cdots\ast d_{1,l_1}\ast\cdots\ast d_{n,1}\ast\cdots\ast d_{n,l_n}\ast\id_{H^{\otimes n}},$$
   where $d_{i,j}$ is an $(r+1)$-fold commutator at $i$ for $i\in [n]$ ($l_i\geq 0, 1\leq j\leq l_i$).
  \end{claim}
  \begin{proof}
   Let $g\in\jfE_r(n)$.
   Since $\Gamma_{r+1}$ is generated by $(r+1)$-fold commutators, $g(x_i)x_i^{-1}$ is a product of $(r+1)$-fold commutators or their inverses for any $i\in [n]$.
   Thus, we can decompose $g$ into a convolution product of $(r+1)$-fold commutators and $\id_{H^{\otimes n}}$.
  \end{proof}

  \begin{proof}[Proof of Theorem \ref{th731}]
   We show that $[A_{d,k}(n),\jfE_{r}(n)]\subset A_{d,k+r}(n)$.
   We can write an element of $A_{d,k}(n)\subset \A^{L}(I,H^{\otimes n})$ as a linear sum of the following diagrams:
   $$u=\centre{
\begingroup%
  \makeatletter%
  \providecommand\color[2][]{%
    \errmessage{(Inkscape) Color is used for the text in Inkscape, but the package 'color.sty' is not loaded}%
    \renewcommand\color[2][]{}%
  }%
  \providecommand\transparent[1]{%
    \errmessage{(Inkscape) Transparency is used (non-zero) for the text in Inkscape, but the package 'transparent.sty' is not loaded}%
    \renewcommand\transparent[1]{}%
  }%
  \providecommand\rotatebox[2]{#2}%
  \newcommand*\fsize{\dimexpr\f@size pt\relax}%
  \newcommand*\lineheight[1]{\fontsize{\fsize}{#1\fsize}\selectfont}%
  \ifx\svgwidth\undefined%
    \setlength{\unitlength}{64.49999809bp}%
    \ifx\svgscale\undefined%
      \relax%
    \else%
      \setlength{\unitlength}{\unitlength * \real{\svgscale}}%
    \fi%
  \else%
    \setlength{\unitlength}{\svgwidth}%
  \fi%
  \global\let\svgwidth\undefined%
  \global\let\svgscale\undefined%
  \makeatother%
  \begin{picture}(1,0.70582251)%
    \lineheight{1}%
    \setlength\tabcolsep{0pt}%
    \put(0,0){\includegraphics[width=\unitlength,page=1]{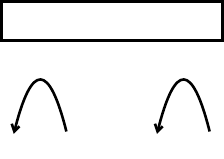}}%
    \put(0.43291744,0.56380374){\makebox(0,0)[lt]{\lineheight{1.45000005}\smash{\begin{tabular}[t]{l}$D$\end{tabular}}}}%
    \put(0.7785644,0.0184876){\makebox(0,0)[lt]{\lineheight{1.45000005}\smash{\begin{tabular}[t]{l}$n$\end{tabular}}}}%
    \put(0.1377383,0.01093144){\makebox(0,0)[lt]{\lineheight{1.45000005}\smash{\begin{tabular}[t]{l}$1$\end{tabular}}}}%
    \put(0,0){\includegraphics[width=\unitlength,page=2]{u2.pdf}}%
  \end{picture}%
\endgroup%
}=\centre{
\begingroup%
  \makeatletter%
  \providecommand\color[2][]{%
    \errmessage{(Inkscape) Color is used for the text in Inkscape, but the package 'color.sty' is not loaded}%
    \renewcommand\color[2][]{}%
  }%
  \providecommand\transparent[1]{%
    \errmessage{(Inkscape) Transparency is used (non-zero) for the text in Inkscape, but the package 'transparent.sty' is not loaded}%
    \renewcommand\transparent[1]{}%
  }%
  \providecommand\rotatebox[2]{#2}%
  \newcommand*\fsize{\dimexpr\f@size pt\relax}%
  \newcommand*\lineheight[1]{\fontsize{\fsize}{#1\fsize}\selectfont}%
  \ifx\svgwidth\undefined%
    \setlength{\unitlength}{65.38865842bp}%
    \ifx\svgscale\undefined%
      \relax%
    \else%
      \setlength{\unitlength}{\unitlength * \real{\svgscale}}%
    \fi%
  \else%
    \setlength{\unitlength}{\svgwidth}%
  \fi%
  \global\let\svgwidth\undefined%
  \global\let\svgscale\undefined%
  \makeatother%
  \begin{picture}(1,0.78097065)%
    \lineheight{1}%
    \setlength\tabcolsep{0pt}%
    \put(0,0){\includegraphics[width=\unitlength,page=1]{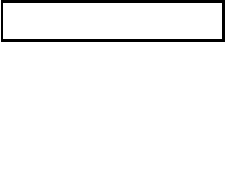}}%
    \put(0.42974986,0.64309589){\makebox(0,0)[lt]{\lineheight{1.45000005}\smash{\begin{tabular}[t]{l}$D$\end{tabular}}}}%
    \put(0.77244646,0.01078287){\makebox(0,0)[lt]{\lineheight{1.45000005}\smash{\begin{tabular}[t]{l}$n$\end{tabular}}}}%
    \put(0.13136111,0.01288464){\makebox(0,0)[lt]{\lineheight{1.45000005}\smash{\begin{tabular}[t]{l}$1$\end{tabular}}}}%
    \put(0,0){\includegraphics[width=\unitlength,page=2]{u2A.pdf}}%
    \put(-0.00462977,0.34143971){\makebox(0,0)[lt]{\lineheight{1.45000005}\smash{\begin{tabular}[t]{l}$i$\end{tabular}}}}%
    \put(0.90955529,0.34143971){\makebox(0,0)[lt]{\lineheight{1.45000005}\smash{\begin{tabular}[t]{l}$i$\end{tabular}}}}%
    \put(0,0){\includegraphics[width=\unitlength,page=3]{u2A.pdf}}%
  \end{picture}%
\endgroup%
},$$
   where $D$ is a Jacobi diagram with at least $k$ trivalent vertices.
   Let $g\in \jfE_{r}(n)$.
   By Claim \ref{c731}, we can write $g$ as a convolution product
   $$g=d_{1,1}\ast\cdots\ast d_{1,l_1}\ast\cdots\ast d_{n,1}\ast\cdots\ast d_{n,l_n}\ast \id_{H^{\otimes n}},$$
   where $d_{i,j}\in \A_0(n,n)$ is an $(r+1)$-fold commutator at $i$.
   Let $l=1+\sum_{i=1}^{n}{l_i}$.

   Now, we have
   $$ u\cdot g=\centre{
\begingroup%
  \makeatletter%
  \providecommand\color[2][]{%
    \errmessage{(Inkscape) Color is used for the text in Inkscape, but the package 'color.sty' is not loaded}%
    \renewcommand\color[2][]{}%
  }%
  \providecommand\transparent[1]{%
    \errmessage{(Inkscape) Transparency is used (non-zero) for the text in Inkscape, but the package 'transparent.sty' is not loaded}%
    \renewcommand\transparent[1]{}%
  }%
  \providecommand\rotatebox[2]{#2}%
  \newcommand*\fsize{\dimexpr\f@size pt\relax}%
  \newcommand*\lineheight[1]{\fontsize{\fsize}{#1\fsize}\selectfont}%
  \ifx\svgwidth\undefined%
    \setlength{\unitlength}{131.31157619bp}%
    \ifx\svgscale\undefined%
      \relax%
    \else%
      \setlength{\unitlength}{\unitlength * \real{\svgscale}}%
    \fi%
  \else%
    \setlength{\unitlength}{\svgwidth}%
  \fi%
  \global\let\svgwidth\undefined%
  \global\let\svgscale\undefined%
  \makeatother%
  \begin{picture}(1,0.89904978)%
    \lineheight{1}%
    \setlength\tabcolsep{0pt}%
    \put(0,0){\includegraphics[width=\unitlength,page=1]{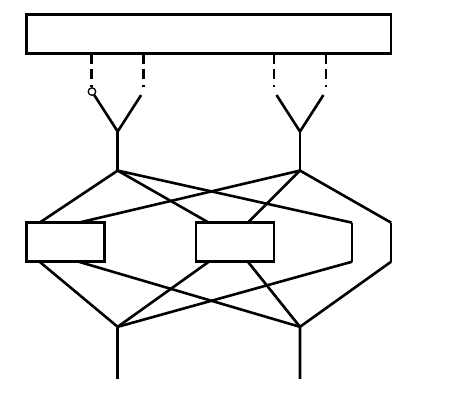}}%
    \put(0.4171057,0.80120121){\makebox(0,0)[lt]{\lineheight{1.45000005}\smash{\begin{tabular}[t]{l}$D$\end{tabular}}}}%
    \put(0.17250832,0.57584319){\makebox(0,0)[lt]{\lineheight{1.45000005}\smash{\begin{tabular}[t]{l}$1$\end{tabular}}}}%
    \put(0.6944751,0.58050254){\makebox(0,0)[lt]{\lineheight{1.45000005}\smash{\begin{tabular}[t]{l}$n$\end{tabular}}}}%
    \put(-0.00104118,0.41892649){\makebox(0,0)[lt]{\lineheight{1.45000005}\smash{\begin{tabular}[t]{l}\tiny$1$\end{tabular}}}}%
    \put(0.36374982,0.40759125){\makebox(0,0)[lt]{\lineheight{1.45000005}\smash{\begin{tabular}[t]{l}\tiny$1$\end{tabular}}}}%
    \put(0.7014147,0.38123919){\makebox(0,0)[lt]{\lineheight{1.45000005}\smash{\begin{tabular}[t]{l}\tiny$1$\end{tabular}}}}%
    \put(0.24942861,0.40977847){\makebox(0,0)[lt]{\lineheight{1.45000005}\smash{\begin{tabular}[t]{l}\tiny$n$\end{tabular}}}}%
    \put(0.5788497,0.41457758){\makebox(0,0)[lt]{\lineheight{1.45000005}\smash{\begin{tabular}[t]{l}\tiny$n$\end{tabular}}}}%
    \put(0.86397011,0.38302986){\makebox(0,0)[lt]{\lineheight{1.45000005}\smash{\begin{tabular}[t]{l}\tiny$n$\end{tabular}}}}%
    \put(0.62465508,0.00549591){\makebox(0,0)[lt]{\lineheight{1.45000005}\smash{\begin{tabular}[t]{l}$n$\end{tabular}}}}%
    \put(0.23032575,0.0053695){\makebox(0,0)[lt]{\lineheight{1.45000005}\smash{\begin{tabular}[t]{l}$1$\end{tabular}}}}%
    \put(0.09351131,0.35239164){\makebox(0,0)[lt]{\lineheight{1.45000005}\smash{\begin{tabular}[t]{l}\tiny$d_{1,1}$\end{tabular}}}}%
    \put(0.44425296,0.35479119){\makebox(0,0)[lt]{\lineheight{1.45000005}\smash{\begin{tabular}[t]{l}\tiny$d_{n,l_n}$\end{tabular}}}}%
    \put(0,0){\includegraphics[width=\unitlength,page=2]{ug1.pdf}}%
    \put(0.97203691,0.72484579){\makebox(0,0)[lt]{\lineheight{1.45000005}\smash{\begin{tabular}[t]{l}$u$\end{tabular}}}}%
    \put(0.97203691,0.29647534){\makebox(0,0)[lt]{\lineheight{1.45000005}\smash{\begin{tabular}[t]{l}$g$\end{tabular}}}}%
    \put(0,0){\includegraphics[width=\unitlength,page=3]{ug1.pdf}}%
  \end{picture}%
\endgroup%
}=\centre{
\begingroup%
  \makeatletter%
  \providecommand\color[2][]{%
    \errmessage{(Inkscape) Color is used for the text in Inkscape, but the package 'color.sty' is not loaded}%
    \renewcommand\color[2][]{}%
  }%
  \providecommand\transparent[1]{%
    \errmessage{(Inkscape) Transparency is used (non-zero) for the text in Inkscape, but the package 'transparent.sty' is not loaded}%
    \renewcommand\transparent[1]{}%
  }%
  \providecommand\rotatebox[2]{#2}%
  \newcommand*\fsize{\dimexpr\f@size pt\relax}%
  \newcommand*\lineheight[1]{\fontsize{\fsize}{#1\fsize}\selectfont}%
  \ifx\svgwidth\undefined%
    \setlength{\unitlength}{130.19145204bp}%
    \ifx\svgscale\undefined%
      \relax%
    \else%
      \setlength{\unitlength}{\unitlength * \real{\svgscale}}%
    \fi%
  \else%
    \setlength{\unitlength}{\svgwidth}%
  \fi%
  \global\let\svgwidth\undefined%
  \global\let\svgscale\undefined%
  \makeatother%
  \begin{picture}(1,0.87798117)%
    \lineheight{1}%
    \setlength\tabcolsep{0pt}%
    \put(0,0){\includegraphics[width=\unitlength,page=1]{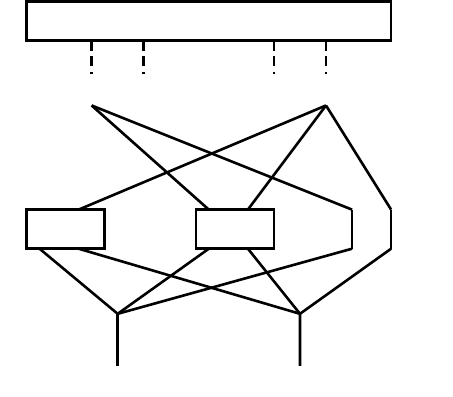}}%
    \put(0.42069434,0.80809448){\makebox(0,0)[lt]{\lineheight{1.45000005}\smash{\begin{tabular}[t]{l}$D$\end{tabular}}}}%
    \put(0.11621892,0.66312493){\makebox(0,0)[lt]{\lineheight{1.45000005}\smash{\begin{tabular}[t]{l}$1$\end{tabular}}}}%
    \put(0.52131282,0.67156028){\makebox(0,0)[lt]{\lineheight{1.45000005}\smash{\begin{tabular}[t]{l}$n$\end{tabular}}}}%
    \put(-0.00105014,0.42253079){\makebox(0,0)[lt]{\lineheight{1.45000005}\smash{\begin{tabular}[t]{l}\tiny$1$\end{tabular}}}}%
    \put(0.36687941,0.41109803){\makebox(0,0)[lt]{\lineheight{1.45000005}\smash{\begin{tabular}[t]{l}\tiny$1$\end{tabular}}}}%
    \put(0.70744943,0.38451925){\makebox(0,0)[lt]{\lineheight{1.45000005}\smash{\begin{tabular}[t]{l}\tiny$1$\end{tabular}}}}%
    \put(0.25157462,0.41330407){\makebox(0,0)[lt]{\lineheight{1.45000005}\smash{\begin{tabular}[t]{l}\tiny$n$\end{tabular}}}}%
    \put(0.58382993,0.41814447){\makebox(0,0)[lt]{\lineheight{1.45000005}\smash{\begin{tabular}[t]{l}\tiny$n$\end{tabular}}}}%
    \put(0.87140343,0.38632532){\makebox(0,0)[lt]{\lineheight{1.45000005}\smash{\begin{tabular}[t]{l}\tiny$n$\end{tabular}}}}%
    \put(0.6300294,0.00554319){\makebox(0,0)[lt]{\lineheight{1.45000005}\smash{\begin{tabular}[t]{l}$n$\end{tabular}}}}%
    \put(0.2323074,0.0054157){\makebox(0,0)[lt]{\lineheight{1.45000005}\smash{\begin{tabular}[t]{l}$1$\end{tabular}}}}%
    \put(0.09431585,0.3554235){\makebox(0,0)[lt]{\lineheight{1.45000005}\smash{\begin{tabular}[t]{l}\tiny$d_{1,1}$\end{tabular}}}}%
    \put(0.44807517,0.35784369){\makebox(0,0)[lt]{\lineheight{1.45000005}\smash{\begin{tabular}[t]{l}\tiny$d_{n,l_n}$\end{tabular}}}}%
    \put(0,0){\includegraphics[width=\unitlength,page=2]{ug2.pdf}}%
    \put(0.74096992,0.62855456){\makebox(0,0)[lt]{\lineheight{1.45000005}\smash{\begin{tabular}[t]{l}$\Delta^{[l]}$\end{tabular}}}}%
    \put(0,0){\includegraphics[width=\unitlength,page=3]{ug2.pdf}}%
  \end{picture}%
\endgroup%
}.$$
   Here, each $\centre{}$ is once connected to all of the diagrams $d_{1,1}, \cdots, d_{n,l_n}$ and $\id_{H^{\otimes n}}$.
   Since we have $\centre{
\begingroup%
  \makeatletter%
  \providecommand\color[2][]{%
    \errmessage{(Inkscape) Color is used for the text in Inkscape, but the package 'color.sty' is not loaded}%
    \renewcommand\color[2][]{}%
  }%
  \providecommand\transparent[1]{%
    \errmessage{(Inkscape) Transparency is used (non-zero) for the text in Inkscape, but the package 'transparent.sty' is not loaded}%
    \renewcommand\transparent[1]{}%
  }%
  \providecommand\rotatebox[2]{#2}%
  \newcommand*\fsize{\dimexpr\f@size pt\relax}%
  \newcommand*\lineheight[1]{\fontsize{\fsize}{#1\fsize}\selectfont}%
  \ifx\svgwidth\undefined%
    \setlength{\unitlength}{42.51968504bp}%
    \ifx\svgscale\undefined%
      \relax%
    \else%
      \setlength{\unitlength}{\unitlength * \real{\svgscale}}%
    \fi%
  \else%
    \setlength{\unitlength}{\svgwidth}%
  \fi%
  \global\let\svgwidth\undefined%
  \global\let\svgscale\undefined%
  \makeatother%
  \begin{picture}(1,1)%
    \lineheight{1}%
    \setlength\tabcolsep{0pt}%
    \put(0,0){\includegraphics[width=\unitlength,page=1]{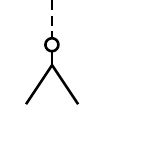}}%
    \put(0.42903395,0.66629947){\makebox(0,0)[lt]{\lineheight{1.45000005}\smash{\begin{tabular}[t]{l}$i$\end{tabular}}}}%
    \put(0.37668098,0.50300989){\makebox(0,0)[lt]{\lineheight{1.45000005}\smash{\begin{tabular}[t]{l}\tiny$\Delta^{[l]}$\end{tabular}}}}%
    \put(0.08805019,0.08807596){\makebox(0,0)[lt]{\lineheight{1.45000005}\smash{\begin{tabular}[t]{l}$1$\end{tabular}}}}%
    \put(0.47670196,0.08810482){\makebox(0,0)[lt]{\lineheight{1.45000005}\smash{\begin{tabular}[t]{l}$l$\end{tabular}}}}%
    \put(0,0){\includegraphics[width=\unitlength,page=2]{multidelta.pdf}}%
  \end{picture}%
\endgroup%
}=\sum_{j=1}^l\centre{
\begingroup%
  \makeatletter%
  \providecommand\color[2][]{%
    \errmessage{(Inkscape) Color is used for the text in Inkscape, but the package 'color.sty' is not loaded}%
    \renewcommand\color[2][]{}%
  }%
  \providecommand\transparent[1]{%
    \errmessage{(Inkscape) Transparency is used (non-zero) for the text in Inkscape, but the package 'transparent.sty' is not loaded}%
    \renewcommand\transparent[1]{}%
  }%
  \providecommand\rotatebox[2]{#2}%
  \newcommand*\fsize{\dimexpr\f@size pt\relax}%
  \newcommand*\lineheight[1]{\fontsize{\fsize}{#1\fsize}\selectfont}%
  \ifx\svgwidth\undefined%
    \setlength{\unitlength}{52.42418096bp}%
    \ifx\svgscale\undefined%
      \relax%
    \else%
      \setlength{\unitlength}{\unitlength * \real{\svgscale}}%
    \fi%
  \else%
    \setlength{\unitlength}{\svgwidth}%
  \fi%
  \global\let\svgwidth\undefined%
  \global\let\svgscale\undefined%
  \makeatother%
  \begin{picture}(1,0.79222913)%
    \lineheight{1}%
    \setlength\tabcolsep{0pt}%
    \put(0,0){\includegraphics[width=\unitlength,page=1]{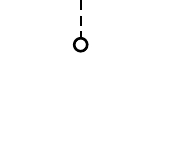}}%
    \put(0.50612279,0.52157462){\makebox(0,0)[lt]{\lineheight{1.45000005}\smash{\begin{tabular}[t]{l}$i$\end{tabular}}}}%
    \put(0.08319923,0.01442542){\makebox(0,0)[lt]{\lineheight{1.45000005}\smash{\begin{tabular}[t]{l}$1$\end{tabular}}}}%
    \put(0.66195685,0.01344947){\makebox(0,0)[lt]{\lineheight{1.45000005}\smash{\begin{tabular}[t]{l}$l$\end{tabular}}}}%
    \put(0,0){\includegraphics[width=\unitlength,page=2]{imultidelta.pdf}}%
    \put(0.38092178,0.01829466){\makebox(0,0)[lt]{\lineheight{1.45000005}\smash{\begin{tabular}[t]{l}$j$\end{tabular}}}}%
    \put(0,0){\includegraphics[width=\unitlength,page=3]{imultidelta.pdf}}%
    \put(-0.0057747,0.44585022){\makebox(0,0)[lt]{\lineheight{1.45000005}\smash{\begin{tabular}[t]{l}$\eta$\end{tabular}}}}%
    \put(0.75202301,0.44426388){\makebox(0,0)[lt]{\lineheight{1.45000005}\smash{\begin{tabular}[t]{l}$\eta$\end{tabular}}}}%
    \put(0,0){\includegraphics[width=\unitlength,page=4]{imultidelta.pdf}}%
  \end{picture}%
\endgroup%
}$ by Lemma \ref{l723} (\ref{l7232}),
   the element $u\cdot g$ is a linear sum of diagrams of shape $\centre{
\begingroup%
  \makeatletter%
  \providecommand\color[2][]{%
    \errmessage{(Inkscape) Color is used for the text in Inkscape, but the package 'color.sty' is not loaded}%
    \renewcommand\color[2][]{}%
  }%
  \providecommand\transparent[1]{%
    \errmessage{(Inkscape) Transparency is used (non-zero) for the text in Inkscape, but the package 'transparent.sty' is not loaded}%
    \renewcommand\transparent[1]{}%
  }%
  \providecommand\rotatebox[2]{#2}%
  \newcommand*\fsize{\dimexpr\f@size pt\relax}%
  \newcommand*\lineheight[1]{\fontsize{\fsize}{#1\fsize}\selectfont}%
  \ifx\svgwidth\undefined%
    \setlength{\unitlength}{130.37348203bp}%
    \ifx\svgscale\undefined%
      \relax%
    \else%
      \setlength{\unitlength}{\unitlength * \real{\svgscale}}%
    \fi%
  \else%
    \setlength{\unitlength}{\svgwidth}%
  \fi%
  \global\let\svgwidth\undefined%
  \global\let\svgscale\undefined%
  \makeatother%
  \begin{picture}(1,0.81922828)%
    \lineheight{1}%
    \setlength\tabcolsep{0pt}%
    \put(0,0){\includegraphics[width=\unitlength,page=1]{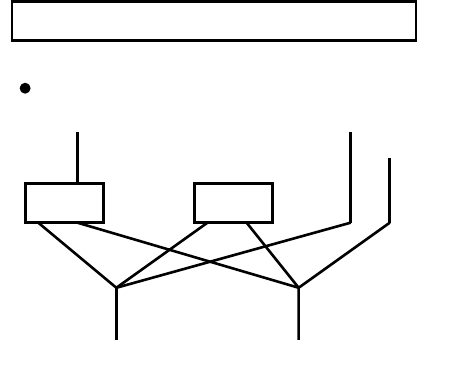}}%
    \put(-0.00104867,0.47848517){\makebox(0,0)[lt]{\lineheight{1.45000005}\smash{\begin{tabular}[t]{l}\tiny$1$\end{tabular}}}}%
    \put(0.18125858,0.47800166){\makebox(0,0)[lt]{\lineheight{1.45000005}\smash{\begin{tabular}[t]{l}\tiny$n$\end{tabular}}}}%
    \put(0.62625003,0.00553545){\makebox(0,0)[lt]{\lineheight{1.45000005}\smash{\begin{tabular}[t]{l}$n$\end{tabular}}}}%
    \put(0.22908333,0.00540814){\makebox(0,0)[lt]{\lineheight{1.45000005}\smash{\begin{tabular}[t]{l}$1$\end{tabular}}}}%
    \put(0.09128445,0.35492725){\makebox(0,0)[lt]{\lineheight{1.45000005}\smash{\begin{tabular}[t]{l}\tiny$d_{1,1}$\end{tabular}}}}%
    \put(0.44454984,0.35734407){\makebox(0,0)[lt]{\lineheight{1.45000005}\smash{\begin{tabular}[t]{l}\tiny$d_{n,l_n}$\end{tabular}}}}%
    \put(0,0){\includegraphics[width=\unitlength,page=2]{ug3.pdf}}%
    \put(0.37287707,0.47848517){\makebox(0,0)[lt]{\lineheight{1.45000005}\smash{\begin{tabular}[t]{l}\tiny$1$\end{tabular}}}}%
    \put(0.55518431,0.47800166){\makebox(0,0)[lt]{\lineheight{1.45000005}\smash{\begin{tabular}[t]{l}\tiny$n$\end{tabular}}}}%
    \put(0,0){\includegraphics[width=\unitlength,page=3]{ug3.pdf}}%
    \put(0.6892757,0.47848517){\makebox(0,0)[lt]{\lineheight{1.45000005}\smash{\begin{tabular}[t]{l}\tiny$1$\end{tabular}}}}%
    \put(0.87158298,0.47800166){\makebox(0,0)[lt]{\lineheight{1.45000005}\smash{\begin{tabular}[t]{l}\tiny$n$\end{tabular}}}}%
    \put(0,0){\includegraphics[width=\unitlength,page=4]{ug3.pdf}}%
  \end{picture}%
\endgroup%
}$, where $\centre{
\begingroup%
  \makeatletter%
  \providecommand\color[2][]{%
    \errmessage{(Inkscape) Color is used for the text in Inkscape, but the package 'color.sty' is not loaded}%
    \renewcommand\color[2][]{}%
  }%
  \providecommand\transparent[1]{%
    \errmessage{(Inkscape) Transparency is used (non-zero) for the text in Inkscape, but the package 'transparent.sty' is not loaded}%
    \renewcommand\transparent[1]{}%
  }%
  \providecommand\rotatebox[2]{#2}%
  \newcommand*\fsize{\dimexpr\f@size pt\relax}%
  \newcommand*\lineheight[1]{\fontsize{\fsize}{#1\fsize}\selectfont}%
  \ifx\svgwidth\undefined%
    \setlength{\unitlength}{4.50118095bp}%
    \ifx\svgscale\undefined%
      \relax%
    \else%
      \setlength{\unitlength}{\unitlength * \real{\svgscale}}%
    \fi%
  \else%
    \setlength{\unitlength}{\svgwidth}%
  \fi%
  \global\let\svgwidth\undefined%
  \global\let\svgscale\undefined%
  \makeatother%
  \begin{picture}(1,6.66491736)%
    \lineheight{1}%
    \setlength\tabcolsep{0pt}%
    \put(0,0){\includegraphics[width=\unitlength,page=1]{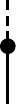}}%
  \end{picture}%
\endgroup%
}$ denotes $\centre{}$ or $\centre{}$.
   If all $\centre{}$ that are connected to $\id_{H^{\otimes n}}$ are $\centre{}$, then it is easily checked that the corresponding summand is just $u$ by using Lemma \ref{l721} (\ref{l7213}).
   Otherwise, at least one of $\centre{}$ that are connected to diagrams $d_{1,1},\cdots,d_{n,l_n}$ are $\centre{}$.
   By using Lemma \ref{l721}, it follows that each summand is a linear sum of diagrams with at least $k+r$ trivalent vertices.
   Therefore, we have $[u,g]=u\cdot g-u \in A_{d,k+r}(n)$.
  \end{proof}


\section{Contraction map}\label{s5}
 Recall that $H=\Lie_1(n)= \bigoplus_{i=1}^n \Z \bar{x}_i$ and $H^{\ast}= \bigoplus_{i=1}^n \Z v_i$.
 In what follows, we identify $H^{\ast}\otimes\Lie_{r+1}(n)$ with $T_{r}(n)$ as we remarked in Section \ref{ss315}.

 \subsection{Contraction map}\label{ss51}
  We define a contraction map
  $$c:B_{d,k}(n)\otimes T_r(n) \rightarrow B_{d,k+r}(n),$$
  which is an analogue of the contraction map defined in Appendix B of \cite{FH}.

  Let $p\geq q$.
  For $I=(i_1,\cdots,i_{q})$ such that $i_1,\cdots,i_{q}$ are distinct elements of $[p]$, define a contraction map
  $$
    c^I:V_n^{\otimes p}\otimes (V_n^{\ast})^{\otimes q}\rightarrow V_n^{\otimes (p-q)}
  $$
  by
  $$
    c^I((w_1\otimes\cdots\otimes w_p)\otimes(y_1\otimes\cdots\otimes y_q))=
    \left(\prod_{j=1}^{q}\langle w_{i_{j}},y_j\rangle\right)
    w_1\otimes\cdots\hat{w}_{i_1}\cdots\hat{w}_{i_q}\cdots\otimes w_{p},
  $$
  where $\hat{w}_{i_1}\cdots\hat{w}_{i_q}$ denotes the omission of $w_{i_1},\cdots,w_{i_q}$ and where $\langle -,-\rangle:V_n\otimes V_n^{\ast}\rightarrow \K$ denotes the dual pairing.
  (See \cite{FH} for details.)

  We next consider a diagrammatic version of the above contraction map $c^I$. Let $2d-k\geq r+1$. For $I=(i_1,\cdots,i_{r+1})\in [2d-k]^{r+1}$ such that $i_1,\cdots,i_{r+1}$ are distinct, we define a linear map
  $$c^I:B_{d,k}(n)\otimes T_r(n) \rightarrow B_{d,k+r}(n)$$
  by contracting colorings of a Jacobi diagram and leaves of a rooted trivalent tree; that is,
  \begin{gather*}
   \begin{split}
     c^I&\left(\centre{
\begingroup%
  \makeatletter%
  \providecommand\color[2][]{%
    \errmessage{(Inkscape) Color is used for the text in Inkscape, but the package 'color.sty' is not loaded}%
    \renewcommand\color[2][]{}%
  }%
  \providecommand\transparent[1]{%
    \errmessage{(Inkscape) Transparency is used (non-zero) for the text in Inkscape, but the package 'transparent.sty' is not loaded}%
    \renewcommand\transparent[1]{}%
  }%
  \providecommand\rotatebox[2]{#2}%
  \newcommand*\fsize{\dimexpr\f@size pt\relax}%
  \newcommand*\lineheight[1]{\fontsize{\fsize}{#1\fsize}\selectfont}%
  \ifx\svgwidth\undefined%
    \setlength{\unitlength}{63.61197298bp}%
    \ifx\svgscale\undefined%
      \relax%
    \else%
      \setlength{\unitlength}{\unitlength * \real{\svgscale}}%
    \fi%
  \else%
    \setlength{\unitlength}{\svgwidth}%
  \fi%
  \global\let\svgwidth\undefined%
  \global\let\svgscale\undefined%
  \makeatother%
  \begin{picture}(1,0.47976077)%
    \lineheight{1}%
    \setlength\tabcolsep{0pt}%
    \put(0,0){\includegraphics[width=\unitlength,page=1]{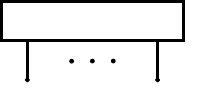}}%
    \put(0.36214569,0.3658384){\makebox(0,0)[lt]{\lineheight{1.45000005}\smash{\begin{tabular}[t]{l}$u$\end{tabular}}}}%
    \put(0.01966035,0.01108404){\makebox(0,0)[lt]{\lineheight{1.45000005}\smash{\begin{tabular}[t]{l}$w_1$\end{tabular}}}}%
    \put(0.60797503,0.01407652){\makebox(0,0)[lt]{\lineheight{1.45000005}\smash{\begin{tabular}[t]{l}$w_{2d-k}$\end{tabular}}}}%
  \end{picture}%
\endgroup%
} \otimes \centre{
\begingroup%
  \makeatletter%
  \providecommand\color[2][]{%
    \errmessage{(Inkscape) Color is used for the text in Inkscape, but the package 'color.sty' is not loaded}%
    \renewcommand\color[2][]{}%
  }%
  \providecommand\transparent[1]{%
    \errmessage{(Inkscape) Transparency is used (non-zero) for the text in Inkscape, but the package 'transparent.sty' is not loaded}%
    \renewcommand\transparent[1]{}%
  }%
  \providecommand\rotatebox[2]{#2}%
  \newcommand*\fsize{\dimexpr\f@size pt\relax}%
  \newcommand*\lineheight[1]{\fontsize{\fsize}{#1\fsize}\selectfont}%
  \ifx\svgwidth\undefined%
    \setlength{\unitlength}{66.49334872bp}%
    \ifx\svgscale\undefined%
      \relax%
    \else%
      \setlength{\unitlength}{\unitlength * \real{\svgscale}}%
    \fi%
  \else%
    \setlength{\unitlength}{\svgwidth}%
  \fi%
  \global\let\svgwidth\undefined%
  \global\let\svgscale\undefined%
  \makeatother%
  \begin{picture}(1,0.73595273)%
    \lineheight{1}%
    \setlength\tabcolsep{0pt}%
    \put(0,0){\includegraphics[width=\unitlength,page=1]{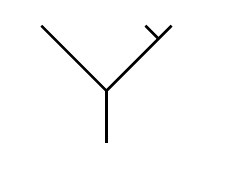}}%
    \put(-0.00455285,0.68966054){\makebox(0,0)[lt]{\lineheight{1.45000005}\smash{\begin{tabular}[t]{l}$y_1$\end{tabular}}}}%
    \put(0.38268163,0.01060373){\makebox(0,0)[lt]{\lineheight{1.45000005}\smash{\begin{tabular}[t]{l}$w$\end{tabular}}}}%
    \put(0.48493984,0.68736319){\makebox(0,0)[lt]{\lineheight{1.45000005}\smash{\begin{tabular}[t]{l}$y_r$\end{tabular}}}}%
    \put(0,0){\includegraphics[width=\unitlength,page=2]{rootedtree2.pdf}}%
    \put(0.66408789,0.68896926){\makebox(0,0)[lt]{\lineheight{1.45000005}\smash{\begin{tabular}[t]{l}$y_{r+1}$\end{tabular}}}}%
  \end{picture}%
\endgroup%
}\right)\\
     &=\left(\prod_{j=1}^{r+1}\langle w_{i_{j}},y_j\rangle\right)
     \centre{
\begingroup%
  \makeatletter%
  \providecommand\color[2][]{%
    \errmessage{(Inkscape) Color is used for the text in Inkscape, but the package 'color.sty' is not loaded}%
    \renewcommand\color[2][]{}%
  }%
  \providecommand\transparent[1]{%
    \errmessage{(Inkscape) Transparency is used (non-zero) for the text in Inkscape, but the package 'transparent.sty' is not loaded}%
    \renewcommand\transparent[1]{}%
  }%
  \providecommand\rotatebox[2]{#2}%
  \newcommand*\fsize{\dimexpr\f@size pt\relax}%
  \newcommand*\lineheight[1]{\fontsize{\fsize}{#1\fsize}\selectfont}%
  \ifx\svgwidth\undefined%
    \setlength{\unitlength}{142.37769723bp}%
    \ifx\svgscale\undefined%
      \relax%
    \else%
      \setlength{\unitlength}{\unitlength * \real{\svgscale}}%
    \fi%
  \else%
    \setlength{\unitlength}{\svgwidth}%
  \fi%
  \global\let\svgwidth\undefined%
  \global\let\svgscale\undefined%
  \makeatother%
  \begin{picture}(1,0.39946787)%
    \lineheight{1}%
    \setlength\tabcolsep{0pt}%
    \put(0,0){\includegraphics[width=\unitlength,page=1]{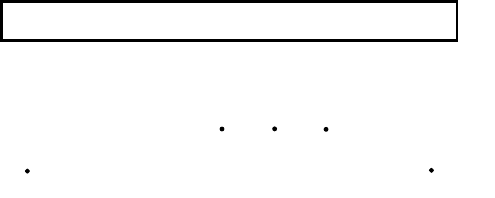}}%
    \put(0.42770167,0.34581546){\makebox(0,0)[lt]{\lineheight{1.45000005}\smash{\begin{tabular}[t]{l}$u$\end{tabular}}}}%
    \put(0.20610356,0.00513507){\makebox(0,0)[lt]{\lineheight{1.45000005}\smash{\begin{tabular}[t]{l}$w_1$\end{tabular}}}}%
    \put(0.8248498,0.00607728){\makebox(0,0)[lt]{\lineheight{1.45000005}\smash{\begin{tabular}[t]{l}$w_{2d-k}$\end{tabular}}}}%
    \put(0,0){\includegraphics[width=\unitlength,page=2]{contra.pdf}}%
    \put(0.0261359,0.00745737){\makebox(0,0)[lt]{\lineheight{1.45000005}\smash{\begin{tabular}[t]{l}$w$\end{tabular}}}}%
    \put(0,0){\includegraphics[width=\unitlength,page=3]{contra.pdf}}%
    \put(0.37794933,0.00537742){\makebox(0,0)[lt]{\lineheight{1.45000005}\smash{\begin{tabular}[t]{l}$\hat{w}_{i_1}\cdots\hat{w}_{i_{r+1}}$\end{tabular}}}}%
    \put(0,0){\includegraphics[width=\unitlength,page=4]{contra.pdf}}%
    \put(0.42197447,0.22722021){\makebox(0,0)[lt]{\lineheight{1.45000005}\smash{\begin{tabular}[t]{l}$\sigma$\end{tabular}}}}%
    \put(0,0){\includegraphics[width=\unitlength,page=5]{contra.pdf}}%
  \end{picture}%
\endgroup%
},
   \end{split}
  \end{gather*}
  where
  $\sigma^{-1}= \begin{pmatrix}
     1&\cdots  &r+1    &r+2&\hdotsfor{5}          &2d-k\\
     i_1&\cdots&i_{r+1}&1&\cdots&\hat{i}_1&\cdots&\hat{i}_{r+1}&\cdots&2d-k
   \end{pmatrix}.$
  We define a \emph{contraction map}
  $$
   c:B_{d,k}(n)\otimes T_r(n) \rightarrow B_{d,k+r}(n)
  $$
  by $c=\sum_{I=(i_1,\cdots,i_{r+1})\in [2d-k]^{r+1}:\text{distinct}}c^I$.
  By using the contraction map $c$, we define a map
  $$\gamma_{d,k}^r:T_r(n)\rightarrow \Hom(B_{d,k}(n),B_{d,k+r}(n))$$
  by $\gamma_{d,k}^r(g)(u):=c(u\otimes g)$ for $g\in T_r(n), u\in B_{d,k}(n)$.

 \subsection{Preliminaries to computation}\label{ss515}
  Let $N\geq 1$.
  We briefly review the construction of the irreducible representations of the symmetric group $\gpS_N$. See Fulton--Harris \cite{FH} and Sagan \cite{Sagan} for basic facts of representation theory of $\gpS_N$.
  Let $\lambda=(\lambda_1,\cdots,\lambda_l)$ be a partition of $N$ and write $\lambda\vdash N$.
  A \emph{Young diagram} of $\lambda$ consists of $\lambda_i$ boxes in the $i$-th row for $i\in[l]$, such that the rows of boxes are lined up on the left.
  A \emph{$\lambda$-tableau} is a numbering of the boxes by the integers in $[N]$. We call a $\lambda$-tableau \emph{standard} if the numbering increases in each row and in each column. The \emph{canonical} $\lambda$-tableau is a standard tableau whose numbering starts from the first row from left to right and then the second row from left to right and so on.

  Let $t_0$ be the canonical $\lambda$-tableau.
  Define $R_{t_0}$ (resp. $C_{t_0}$) to be the subgroup of $\gpS_N$ that preserves each row (resp. column) of $t_0$.
  We define
  $$a_{\lambda}:=\sum_{\sigma\in R_{t_0}}\sigma,\quad b_{\lambda}:=\sum_{\sigma\in C_{t_0}}\sgn(\sigma)\sigma\in\K\gpS_N.$$

  For each $\lambda\vdash N$, the \emph{Young symmetrizer} $c_{\lambda}$ is defined by
  \begin{equation}\label{youngsym}
    c_{\lambda}=b_{\lambda}a_{\lambda}\in\K\gpS_{N}.
  \end{equation}
  The Specht module $S^{\lambda}$, which is an irreducible representation of $\gpS_N$ corresponding to $\lambda$, can be constructed as
  $$
    S^{\lambda}=\K\gpS_N\cdot c_{\lambda}.
  $$

  \begin{lemma}\label{l822}
   We have the following decomposition of $\K\gpS_N$-bimodules
   \begin{gather*}
    \begin{split}
     \K\gpS_{N}=\bigoplus_{\lambda\vdash N}\K\gpS_{N}\cdot c_{\lambda}\cdot\K\gpS_{N}.
    \end{split}
   \end{gather*}
  \end{lemma}
  \begin{proof}
   This follows from basic facts of representation theory. The reader is referred to \cite{FH} and \cite{Sagan}.
  \end{proof}

  For $N',N''\geq 0$, let $N=N'+N''$.
  For $\mu\vdash N',\nu\vdash N''$, let $S^{\mu}\diamond S^{\nu}$ denote the representation of $\gpS_N$ induced from the tensor product representation $S^{\mu}\boxtimes S^{\nu}$ of $\gpS_{N'}\times \gpS_{N''}$ by the inclusion of $\gpS_{N'}\times \gpS_{N''}$ in $\gpS_N$.
  By the Littlewood--Richardson rule, we have
  $$
    S^{\mu}\diamond S^{\nu}=\bigoplus_{\lambda\vdash N} (S^{\lambda})^{LR_{\mu,\nu}^{\lambda}},
  $$
  where $LR_{\mu,\nu}^{\lambda}$ denotes the Littlewood--Richardson coefficient.
  We have the following lemma by using basic facts of representation theory of $\gpS_N$.
  \begin{lemma}\label{l824}
   Let $N=N'+N''$ for $N',N''\geq 0$.
   Let $\lambda\vdash N,\mu\vdash N',\nu\vdash N''$, respectively.
   We have $$\dim_{\K}((c_{\mu}\diamond c_{\nu})\cdot \K\gpS_N\cdot c_{\lambda})= LR_{\mu,\nu}^{\lambda}.$$
   In particular, if the Littlewood--Richardson coefficient $LR_{\mu,\nu}^{\lambda}=0$, then we have
   $$
     (c_{\mu}\diamond c_{\nu})\cdot \K\gpS_N\cdot c_{\lambda}=0.
   $$
  \end{lemma}

 \subsection{Vanishing conditions for the contraction map}\label{ss52}
  Here, we observe that the contraction map vanishes under certain specific conditions.

  We have an isomorphism of $\GL(V_n)$-modules
  \begin{equation}\label{eqBD}
    B_{d,k}(n)\cong V_n^{\otimes 2d-k}\otimes_{\K\gpS_{2d-k}} D_{d,k},
  \end{equation}
  where $D_{d,k}$ is the $\K$-vector space spanned by $[2d-k]$-colored open Jacobi diagrams of degree $d$ such that the map $\{\text{univalent vertices of $D$}\}\rightarrow T$ that gives the coloring of $D$ is a bijection.
  Thus, any element of $B_{d,k}(n)$ can be written in the form
  $$
    u(w_1,\cdots,w_{2d-k}):=(w_1\otimes\cdots\otimes w_{2d-k})\otimes u
  $$
  for $u\in D_{d,k}$ and $w_1,\cdots, w_{2d-k}\in V_n$.

  For $\lambda\vdash 2d-k$, let $B_{d,k}(n)_\lambda$ be the \emph{isotypic component} of $B_{d,k}(n)$ corresponding to $\lambda$; that is,
  $$B_{d,k}(n)_\lambda\cong V_n^{\otimes 2d-k}\otimes_{\K\gpS_{2d-k}}\K\gpS_{2d-k}c_{\lambda}D_{d,k}.$$
  We have $B_{d,k}(n)=\bigoplus_{\lambda\vdash 2d-k}B_{d,k}(n)_{\lambda}$.

  For $r\geq 0$, a trivalent tree is called a \emph{based trivalent tree of degree $r$} if it has one distinguished univalent vertex with no coloring (called a \emph{base}) and $r+1$ univalent vertices (called \emph{leaves}) that are colored by distinct elements of $[r+1]$. (Note that a based trivalent tree is different from a rooted trivalent tree.)
  Let $L_{r}$ denote the $\Z$-module spanned by based trivalent trees of degree $r$ modulo the AS and IHX relations.
  The symmetric group $\gpS_{r+1}$ acts on the $\Z$-module $L_{r}$ by the action on colorings of based trivalent trees.
  Then we have
  $$\Lie_{r+1}(n)\cong H^{\otimes (r+1)}\otimes_{\Z\gpS_{r+1}} L_{r}.$$

  On the other hand, $\Lie_{r+1}(n)$ has a $\GL(n;\Z)$-module structure by the standard action on each factor. (See \cite{F} for representation theory of $\GL(n;\Z)$.)
  For $\mu\vdash r+1$, let $\Lie_{r+1}(n)_{\mu}$ denote the isotypic component of $\Lie_{r+1}(n)$ corresponding to $\mu$; that is,
  $$\Lie_{r+1}(n)_{\mu}\cong H^{\otimes (r+1)}\otimes_{\Z\gpS_{r+1}} \Z\gpS_{r+1}c_{\mu} L_{r}.$$
  We have $\Lie_{r+1}(n)=\bigoplus_{\mu\vdash r+1}\Lie_{r+1}(n)_{\mu}$.

  For partitions $\lambda$ and $\mu$, we write $\lambda\nsupseteq\mu$ if the Young diagram of $\lambda$ does not contain that of $\mu$.

  \begin{proposition}\label{p822}
   For $2d-k\geq r+1$, let $\lambda\vdash 2d-k$ and $\mu\vdash r+1$.
   We have
   $$
      c(B_{d,k}(n)_{\lambda}\otimes (H^{\ast}\otimes \Lie_{r+1}(n)_{\mu})) \subset \bigoplus_{\rho:LR^{\lambda}_{\mu,\nu}LR^{\rho}_{\nu,(1)}\neq0 \text{ for some }\nu} B_{d,k+r}(n)_{\rho}.
   $$
   In particular, if $\lambda\nsupseteq\mu$, then we have $$c(B_{d,k}(n)_{\lambda}\otimes (H^{\ast}\otimes \Lie_{r+1}(n)_{\mu}))=0.$$
  \end{proposition}
  \begin{proof}
   Any element of $B_{d,k}(n)_{\lambda}$ is a linear sum of $(c_{\lambda}\cdot u)(w_1,\cdots,w_{2d-k})$, where $u(w_1,\cdots,w_{2d-k})\in B_{d,k}(n)$.
   Any element of $L_{r}$ is a linear sum of
   $$L=\pi^{-1}\cdot \centre{
\begingroup%
  \makeatletter%
  \providecommand\color[2][]{%
    \errmessage{(Inkscape) Color is used for the text in Inkscape, but the package 'color.sty' is not loaded}%
    \renewcommand\color[2][]{}%
  }%
  \providecommand\transparent[1]{%
    \errmessage{(Inkscape) Transparency is used (non-zero) for the text in Inkscape, but the package 'transparent.sty' is not loaded}%
    \renewcommand\transparent[1]{}%
  }%
  \providecommand\rotatebox[2]{#2}%
  \newcommand*\fsize{\dimexpr\f@size pt\relax}%
  \newcommand*\lineheight[1]{\fontsize{\fsize}{#1\fsize}\selectfont}%
  \ifx\svgwidth\undefined%
    \setlength{\unitlength}{40.85642127bp}%
    \ifx\svgscale\undefined%
      \relax%
    \else%
      \setlength{\unitlength}{\unitlength * \real{\svgscale}}%
    \fi%
  \else%
    \setlength{\unitlength}{\svgwidth}%
  \fi%
  \global\let\svgwidth\undefined%
  \global\let\svgscale\undefined%
  \makeatother%
  \begin{picture}(1,0.7343828)%
    \lineheight{1}%
    \setlength\tabcolsep{0pt}%
    \put(0,0){\includegraphics[width=\unitlength,page=1]{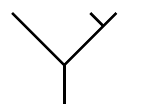}}%
    \put(-0.00167316,0.69671275){\makebox(0,0)[lt]{\lineheight{1.45000005}\smash{\begin{tabular}[t]{l}\tiny$1$\end{tabular}}}}%
    \put(0.53877545,0.69297388){\makebox(0,0)[lt]{\lineheight{1.45000005}\smash{\begin{tabular}[t]{l}\tiny$r$\end{tabular}}}}%
    \put(0,0){\includegraphics[width=\unitlength,page=2]{basedtree.pdf}}%
    \put(0.73391942,0.69283297){\makebox(0,0)[lt]{\lineheight{1.45000005}\smash{\begin{tabular}[t]{l}\tiny$r+1$\end{tabular}}}}%
  \end{picture}%
\endgroup%
}$$
   for $\pi\in \gpS_{r+1}$.
   Thus, any element of $H^{\ast}\otimes \Lie_{r+1}(n)_{\mu}$ is a linear sum of
   $w\otimes ((y_1\otimes \cdots\otimes y_{r+1})\otimes c_{\mu}\cdot L)$ for
   $w\in H^{\ast}, y_1,\cdots,y_{r+1}\in H$.

   For any $I=(i_1,\cdots,i_{r+1})\in[2d-k]^{r+1}$ such that $i_1,\cdots,i_{r+1}$ are distinct, we have
   $$
     c^I((c_{\lambda}\cdot u)(w_1,\cdots,w_{2d-k})\otimes (w\otimes ((y_1 \otimes \cdots\otimes y_{r+1})\otimes c_{\mu}\cdot L)))
     =\prod_{j=1}^{r+1}\langle w_{i_j}, y_j \rangle D,
   $$
   where
   $$D=\centre{
\begingroup%
  \makeatletter%
  \providecommand\color[2][]{%
    \errmessage{(Inkscape) Color is used for the text in Inkscape, but the package 'color.sty' is not loaded}%
    \renewcommand\color[2][]{}%
  }%
  \providecommand\transparent[1]{%
    \errmessage{(Inkscape) Transparency is used (non-zero) for the text in Inkscape, but the package 'transparent.sty' is not loaded}%
    \renewcommand\transparent[1]{}%
  }%
  \providecommand\rotatebox[2]{#2}%
  \newcommand*\fsize{\dimexpr\f@size pt\relax}%
  \newcommand*\lineheight[1]{\fontsize{\fsize}{#1\fsize}\selectfont}%
  \ifx\svgwidth\undefined%
    \setlength{\unitlength}{94.50058921bp}%
    \ifx\svgscale\undefined%
      \relax%
    \else%
      \setlength{\unitlength}{\unitlength * \real{\svgscale}}%
    \fi%
  \else%
    \setlength{\unitlength}{\svgwidth}%
  \fi%
  \global\let\svgwidth\undefined%
  \global\let\svgscale\undefined%
  \makeatother%
  \begin{picture}(1,1.04261735)%
    \lineheight{1}%
    \setlength\tabcolsep{0pt}%
    \put(0,0){\includegraphics[width=\unitlength,page=1]{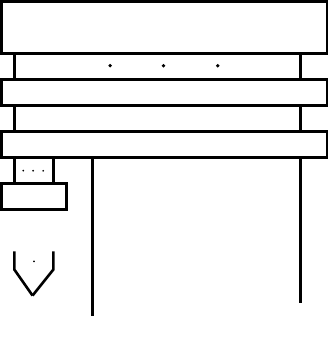}}%
    \put(0.45415485,0.93148726){\makebox(0,0)[lt]{\lineheight{1.45000005}\smash{\begin{tabular}[t]{l}$u$\end{tabular}}}}%
    \put(0.43620022,0.73807784){\makebox(0,0)[lt]{\lineheight{1.45000005}\smash{\begin{tabular}[t]{l}\tiny$c_{\lambda}$\end{tabular}}}}%
    \put(0.45409008,0.57777662){\makebox(0,0)[lt]{\lineheight{1.45000005}\smash{\begin{tabular}[t]{l}$\sigma$\end{tabular}}}}%
    \put(0.01950071,0.42910499){\makebox(0,0)[lt]{\lineheight{1.45000005}\smash{\begin{tabular}[t]{l}\tiny$a_{\mu}b_{\mu}$\end{tabular}}}}%
    \put(0,0){\includegraphics[width=\unitlength,page=2]{D11.pdf}}%
    \put(0.06138215,0.00425758){\makebox(0,0)[lt]{\lineheight{1.45000005}\smash{\begin{tabular}[t]{l}\tiny$w$\end{tabular}}}}%
    \put(0.27216852,0.00616071){\makebox(0,0)[lt]{\lineheight{1.45000005}\smash{\begin{tabular}[t]{l}\tiny$w_1\; \hat{w}_{i_1}\cdots\hat{w}_{i_{r+1}}\; w_{2d-k}$\end{tabular}}}}%
    \put(0,0){\includegraphics[width=\unitlength,page=3]{D11.pdf}}%
    \put(0.07677795,0.29985293){\makebox(0,0)[lt]{\lineheight{1.45000005}\smash{\begin{tabular}[t]{l}\tiny$\pi$\end{tabular}}}}%
  \end{picture}%
\endgroup%
},$$
   $$\sigma^{-1}=
     \begin{pmatrix}
       1&\cdots  &r+1    &r+2&\hdotsfor{5}          &2d-k\\
       i_1&\cdots&i_{r+1}&1&\cdots&\hat{i}_1&\cdots&\hat{i}_{r+1}&\cdots&2d-k
     \end{pmatrix}.
   $$
   Let $l=2d-k-r-1$.
   By Lemma \ref{l822}, we have
   $$\id_l=\sum_{\nu\vdash l,1\leq i\leq \dim{S^{\nu}}}\tau_{i,1} c_{\nu}\tau_{i,2},$$
   where $\tau_{i,1},\tau_{i,2}\in \K\gpS_l$.
   Thus, we have
   $$D=\sum_{\nu\vdash l,1\leq i\leq \dim{S^{\nu}}} \centre{
\begingroup%
  \makeatletter%
  \providecommand\color[2][]{%
    \errmessage{(Inkscape) Color is used for the text in Inkscape, but the package 'color.sty' is not loaded}%
    \renewcommand\color[2][]{}%
  }%
  \providecommand\transparent[1]{%
    \errmessage{(Inkscape) Transparency is used (non-zero) for the text in Inkscape, but the package 'transparent.sty' is not loaded}%
    \renewcommand\transparent[1]{}%
  }%
  \providecommand\rotatebox[2]{#2}%
  \newcommand*\fsize{\dimexpr\f@size pt\relax}%
  \newcommand*\lineheight[1]{\fontsize{\fsize}{#1\fsize}\selectfont}%
  \ifx\svgwidth\undefined%
    \setlength{\unitlength}{93.21019643bp}%
    \ifx\svgscale\undefined%
      \relax%
    \else%
      \setlength{\unitlength}{\unitlength * \real{\svgscale}}%
    \fi%
  \else%
    \setlength{\unitlength}{\svgwidth}%
  \fi%
  \global\let\svgwidth\undefined%
  \global\let\svgscale\undefined%
  \makeatother%
  \begin{picture}(1,1.01107055)%
    \lineheight{1}%
    \setlength\tabcolsep{0pt}%
    \put(0,0){\includegraphics[width=\unitlength,page=1]{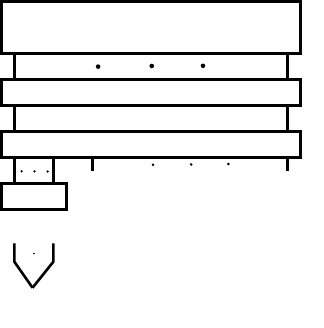}}%
    \put(0.4277086,0.89754057){\makebox(0,0)[lt]{\lineheight{1.45000005}\smash{\begin{tabular}[t]{l}$u$\end{tabular}}}}%
    \put(0.42907172,0.70806867){\makebox(0,0)[lt]{\lineheight{1.45000005}\smash{\begin{tabular}[t]{l}$c_{\lambda}$\end{tabular}}}}%
    \put(0.43077474,0.54312567){\makebox(0,0)[lt]{\lineheight{1.45000005}\smash{\begin{tabular}[t]{l}$\sigma$\end{tabular}}}}%
    \put(0.01977068,0.38906474){\makebox(0,0)[lt]{\lineheight{1.45000005}\smash{\begin{tabular}[t]{l}\tiny$a_{\mu}b_{\mu}$\end{tabular}}}}%
    \put(0,0){\includegraphics[width=\unitlength,page=2]{D12.pdf}}%
    \put(0.05971707,0.00448253){\makebox(0,0)[lt]{\lineheight{1.45000005}\smash{\begin{tabular}[t]{l}\tiny$w$\end{tabular}}}}%
    \put(0.26095254,0.00431652){\makebox(0,0)[lt]{\lineheight{1.45000005}\smash{\begin{tabular}[t]{l}\tiny$w_1\; \hat{w}_{i_1}\cdots\hat{w}_{i_{r+1}} \; w_{2d-k}$\end{tabular}}}}%
    \put(0,0){\includegraphics[width=\unitlength,page=3]{D12.pdf}}%
    \put(0.54404592,0.42203751){\makebox(0,0)[lt]{\lineheight{1.45000005}\smash{\begin{tabular}[t]{l}$\tau_{i,2}$\end{tabular}}}}%
    \put(0,0){\includegraphics[width=\unitlength,page=4]{D12.pdf}}%
    \put(0.54237526,0.29967184){\makebox(0,0)[lt]{\lineheight{1.45000005}\smash{\begin{tabular}[t]{l}$c_{\nu}$\end{tabular}}}}%
    \put(0,0){\includegraphics[width=\unitlength,page=5]{D12.pdf}}%
    \put(0.54254365,0.17914534){\makebox(0,0)[lt]{\lineheight{1.45000005}\smash{\begin{tabular}[t]{l}$\tau_{i,1}$\end{tabular}}}}%
    \put(0,0){\includegraphics[width=\unitlength,page=6]{D12.pdf}}%
    \put(0.08713574,0.2780427){\makebox(0,0)[lt]{\lineheight{1.45000005}\smash{\begin{tabular}[t]{l}\tiny$\pi$\end{tabular}}}}%
  \end{picture}%
\endgroup%
}.$$

   If $LR^{\lambda}_{\mu,\nu}=0$ for any $\nu\vdash l$, then we have $D=0$ by Lemma \ref{l824}.
   Otherwise, since we have
   $$
     \id_1\otimes c_{\nu}\in \bigoplus_{\rho\vdash l+1} (S^{\rho})^{LR^{\rho}_{\nu,(1)}},
   $$
   by the Littlewood--Richardson rule, it follows that
   $$
     D\in\bigoplus_{\rho:LR^{\lambda}_{\mu,\nu}LR^{\rho}_{\nu,(1)}\neq0 \text{ for some }\nu}\left( B_{d,k+r}(n)\right)_{\rho}.
   $$

   If $\lambda\nsupseteq\mu$, then $LR_{\mu,\nu}^{\lambda}=0$ for any $\nu\vdash l$. Thus, we have
   $$
     c(B_{d,k}(n)_{\lambda}\otimes (H^{\ast}\otimes \Lie_{r+1}(n)_{\mu}))=0.
   $$
  \end{proof}

  \begin{remark}\label{rem52}
    Note that we have $\Lie_2(n)=\Lie_2(n)_{(1^2)}$. Thus, the restriction
    $$c:B_{d,k}(n)_{\lambda}\otimes (H^{\ast}\otimes \Lie_2(n)_{(1^2)})\rightarrow B_{d,k+1}(n)_{\rho}$$
    of the contraction map vanishes unless $\rho$ can be obtained from $\lambda$ by taking away one box from each of two different rows of $\lambda$ and then by adding one box.
  \end{remark}


\section{Correspondence between the map $\ti{\beta}_{d,k}^r$ and the
map $\gamma_{d,k}^r$}\label{s6}
 In this section, we prove that the map $\ti{\beta}_{d,k}^r$ defined in Section \ref{s4} can be identified with the map $\gamma_{d,k}^r$ defined in Section \ref{s5} via the Johnson homomorphism of $\End(F_n)$ defined in Section \ref{s3}.
 \begin{theorem}\label{th91}
   We have $\ti{\beta}_{d,k}^r=(-1)^r\cdot \gamma_{d,k}^r\circ \ti{\tau}_r$. That is, we have the following commutative diagram (up to sign):
   \begin{gather*}
    \xymatrix{
         \opegr^r(\jfE_{\ast}(n))\ar[d]_{\cong}^{\ti{\tau}_r}\ar[rr]^-{\ti{\beta}_{d,k}^r}
         &&\Hom(B_{d,k}(n),B_{d,k+r}(n))
         \\
         H^{\ast}\otimes\Lie_{r+1}(n)
         \ar[rru]_-{\gamma_{d,k}^r}.
       }
   \end{gather*}
 \end{theorem}

 \begin{proof}
  The $\Z$-module $H^{\ast}\otimes \Lie_{r+1}(n)$ is spanned by
  $v_i\otimes [\bar{x}_{i_1},\cdots,[\bar{x}_{i_r},\bar{x}_{i_{r+1}}]\cdots]$ for $i, i_1,\cdots,i_{r+1}\in [n]$.
  Define $\phi\in\End(F_n)$ by
  $$
    \phi(x_i)=[x_{i_1},\cdots,[x_{i_r},x_{i_{r+1}}]\cdots]\cdot x_i,\quad \phi(x_j)= x_j \:(j\neq i).
  $$
  It is easily checked that $\phi\in\jfE_{r}(n)$ and that $\ti{\tau}_r([\phi]_r)=v_i\otimes [\bar{x}_{i_1},\cdots,[\bar{x}_{i_r},\bar{x}_{i_{r+1}}]\cdots]$,
  where $[\phi]_r\in \opegr^r(\jfE_{\ast}(n))$ denotes the image of $\phi$ under the projection.

  Any element of $B_{d,k}(n)$ can be written as a linear sum of
  $u=\centre{
\begingroup%
  \makeatletter%
  \providecommand\color[2][]{%
    \errmessage{(Inkscape) Color is used for the text in Inkscape, but the package 'color.sty' is not loaded}%
    \renewcommand\color[2][]{}%
  }%
  \providecommand\transparent[1]{%
    \errmessage{(Inkscape) Transparency is used (non-zero) for the text in Inkscape, but the package 'transparent.sty' is not loaded}%
    \renewcommand\transparent[1]{}%
  }%
  \providecommand\rotatebox[2]{#2}%
  \newcommand*\fsize{\dimexpr\f@size pt\relax}%
  \newcommand*\lineheight[1]{\fontsize{\fsize}{#1\fsize}\selectfont}%
  \ifx\svgwidth\undefined%
    \setlength{\unitlength}{68.91362872bp}%
    \ifx\svgscale\undefined%
      \relax%
    \else%
      \setlength{\unitlength}{\unitlength * \real{\svgscale}}%
    \fi%
  \else%
    \setlength{\unitlength}{\svgwidth}%
  \fi%
  \global\let\svgwidth\undefined%
  \global\let\svgscale\undefined%
  \makeatother%
  \begin{picture}(1,0.52113134)%
    \lineheight{1}%
    \setlength\tabcolsep{0pt}%
    \put(0,0){\includegraphics[width=\unitlength,page=1]{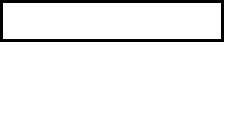}}%
    \put(0.40519088,0.39030895){\makebox(0,0)[lt]{\lineheight{1.45000005}\smash{\begin{tabular}[t]{l}$D$\end{tabular}}}}%
    \put(0,0){\includegraphics[width=\unitlength,page=2]{D2.pdf}}%
    \put(0.06859371,0.00599426){\makebox(0,0)[lt]{\lineheight{1.45000005}\smash{\begin{tabular}[t]{l}$v_{j_1}$\end{tabular}}}}%
    \put(0.76034298,0.00828944){\makebox(0,0)[lt]{\lineheight{1.45000005}\smash{\begin{tabular}[t]{l}$v_{j_{2d-k}}$\end{tabular}}}}%
  \end{picture}%
\endgroup%
}$,
  where $1\leq j_1\leq \cdots\leq j_{2d-k}\leq n$,
  by arranging the univalent vertices according to the order of indices of the colorings from left to right.
  We have
  \begin{gather*}
    \begin{split}
      &\gamma_{d,k}^r\circ  \ti{\tau}_r([\phi]_r)(u)\\
      &=c(u\otimes (v_i\otimes [\bar{x}_{i_1},\cdots,[\bar{x}_{i_r},\bar{x}_{i_{r+1}}]\cdots]))\\
      &=\sum_{(\alpha_l)\in[2d-k]^{r+1}:\text{ distinct}}
      \left(\prod_{l=1}^{r+1}\langle v_{j_{\alpha_l}},\bar{x}_{i_l}\rangle\right) \centre{
\begingroup%
  \makeatletter%
  \providecommand\color[2][]{%
    \errmessage{(Inkscape) Color is used for the text in Inkscape, but the package 'color.sty' is not loaded}%
    \renewcommand\color[2][]{}%
  }%
  \providecommand\transparent[1]{%
    \errmessage{(Inkscape) Transparency is used (non-zero) for the text in Inkscape, but the package 'transparent.sty' is not loaded}%
    \renewcommand\transparent[1]{}%
  }%
  \providecommand\rotatebox[2]{#2}%
  \newcommand*\fsize{\dimexpr\f@size pt\relax}%
  \newcommand*\lineheight[1]{\fontsize{\fsize}{#1\fsize}\selectfont}%
  \ifx\svgwidth\undefined%
    \setlength{\unitlength}{152.50079675bp}%
    \ifx\svgscale\undefined%
      \relax%
    \else%
      \setlength{\unitlength}{\unitlength * \real{\svgscale}}%
    \fi%
  \else%
    \setlength{\unitlength}{\svgwidth}%
  \fi%
  \global\let\svgwidth\undefined%
  \global\let\svgscale\undefined%
  \makeatother%
  \begin{picture}(1,0.56880366)%
    \lineheight{1}%
    \setlength\tabcolsep{0pt}%
    \put(0,0){\includegraphics[width=\unitlength,page=1]{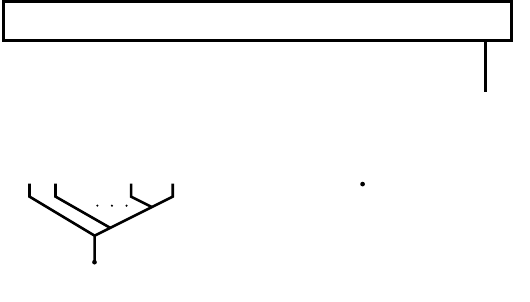}}%
    \put(0.41101439,0.00270875){\makebox(0,0)[lt]{\lineheight{1.45000005}\smash{\begin{tabular}[t]{l}$v_{j_1}$\end{tabular}}}}%
    \put(0.15408072,0.00519449){\makebox(0,0)[lt]{\lineheight{1.45000005}\smash{\begin{tabular}[t]{l}$v_i$\end{tabular}}}}%
    \put(0.89170132,0.00302942){\makebox(0,0)[lt]{\lineheight{1.45000005}\smash{\begin{tabular}[t]{l}$v_{j_{2d-k}}$\end{tabular}}}}%
    \put(0,0){\includegraphics[width=\unitlength,page=2]{gammaim.pdf}}%
    \put(0.47041854,0.35802049){\makebox(0,0)[lt]{\lineheight{1.45000005}\smash{\begin{tabular}[t]{l}$\tau$\end{tabular}}}}%
    \put(0,0){\includegraphics[width=\unitlength,page=3]{gammaim.pdf}}%
    \put(0.17253767,0.23243008){\makebox(0,0)[lt]{\lineheight{1.45000005}\smash{\begin{tabular}[t]{l}$\sigma$\end{tabular}}}}%
    \put(0,0){\includegraphics[width=\unitlength,page=4]{gammaim.pdf}}%
    \put(0.50916882,0.00325838){\makebox(0,0)[lt]{\lineheight{1.45000005}\smash{\begin{tabular}[t]{l}$\hat{v}_{j_{\alpha_1}}\cdots \hat{v}_{j_{\alpha_{r+1}}}$\end{tabular}}}}%
    \put(0,0){\includegraphics[width=\unitlength,page=5]{gammaim.pdf}}%
    \put(0.46824738,0.51770738){\makebox(0,0)[lt]{\lineheight{1.45000005}\smash{\begin{tabular}[t]{l}$D$\end{tabular}}}}%
  \end{picture}%
\endgroup%
},
    \end{split}
  \end{gather*}
  where $\tau^{-1}\in \gpS_{2d-k}$ is the $(r+1,2d-k-r-1)$-shuffle that maps $[r+1]\subset [2d-k]$ to $\{\alpha_l\}$, and $\sigma\in \gpS_{r+1}$ satisfies $\sigma^{-1}(l)=\tau(\alpha_l)$ for any $l\in [r+1]$.

  Let $\ti{u}=\centre{} \in A_{d,k}(n)$, which can be obtained from $u$ by replacing univalent vertices with $\centre{}$ and combining solid lines whose corresponding colorings of $u$ are the same.
  Then $\ti{u}$ is a lift of $u$; that is, we have $\theta_{d,n,k}(\ti{u})=u$.
  By the definition of $\ti{\beta}_{d,k}^r$, we have
  \begin{equation*}
    \ti{\beta}_{d,k}^r([\phi]_r)(u)=[u,[\phi]_r]=\theta_{d,n,k+r}([\ti{u},\phi]).
  \end{equation*}
  We have
  \begin{gather*}
   \begin{split}
     [\ti{u},\phi]&=\ti{u}\cdot \phi- \ti{u}\\
     &=\ti{u}\cdot \left(\centre{
\begingroup%
  \makeatletter%
  \providecommand\color[2][]{%
    \errmessage{(Inkscape) Color is used for the text in Inkscape, but the package 'color.sty' is not loaded}%
    \renewcommand\color[2][]{}%
  }%
  \providecommand\transparent[1]{%
    \errmessage{(Inkscape) Transparency is used (non-zero) for the text in Inkscape, but the package 'transparent.sty' is not loaded}%
    \renewcommand\transparent[1]{}%
  }%
  \providecommand\rotatebox[2]{#2}%
  \newcommand*\fsize{\dimexpr\f@size pt\relax}%
  \newcommand*\lineheight[1]{\fontsize{\fsize}{#1\fsize}\selectfont}%
  \ifx\svgwidth\undefined%
    \setlength{\unitlength}{78.22326224bp}%
    \ifx\svgscale\undefined%
      \relax%
    \else%
      \setlength{\unitlength}{\unitlength * \real{\svgscale}}%
    \fi%
  \else%
    \setlength{\unitlength}{\svgwidth}%
  \fi%
  \global\let\svgwidth\undefined%
  \global\let\svgscale\undefined%
  \makeatother%
  \begin{picture}(1,1.49503295)%
    \lineheight{1}%
    \setlength\tabcolsep{0pt}%
    \put(0,0){\includegraphics[width=\unitlength,page=1]{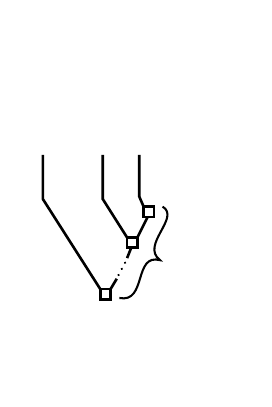}}%
    \put(0.60721419,0.51974855){\makebox(0,0)[lt]{\lineheight{1.45000005}\smash{\begin{tabular}[t]{l}$r$\end{tabular}}}}%
    \put(0,0){\includegraphics[width=\unitlength,page=2]{phi.pdf}}%
    \put(0.39014131,1.19508456){\makebox(0,0)[lt]{\lineheight{1.45000005}\smash{\begin{tabular}[t]{l}$\rho$\end{tabular}}}}%
    \put(0,0){\includegraphics[width=\unitlength,page=3]{phi.pdf}}%
    \put(0.30292676,0.95181191){\makebox(0,0)[lt]{\lineheight{1.45000005}\smash{\begin{tabular}[t]{l}$\pi$\end{tabular}}}}%
    \put(0,0){\includegraphics[width=\unitlength,page=4]{phi.pdf}}%
    \put(0.83205133,0.81136218){\makebox(0,0)[lt]{\lineheight{1.45000005}\smash{\begin{tabular}[t]{l}$\epsilon$\end{tabular}}}}%
    \put(0,0){\includegraphics[width=\unitlength,page=5]{phi.pdf}}%
    \put(0.09826244,0.00901366){\makebox(0,0)[lt]{\lineheight{1.45000005}\smash{\begin{tabular}[t]{l}$1$\end{tabular}}}}%
    \put(0.34883796,0.01022834){\makebox(0,0)[lt]{\lineheight{1.45000005}\smash{\begin{tabular}[t]{l}$i$\end{tabular}}}}%
    \put(0.77188498,0.01140598){\makebox(0,0)[lt]{\lineheight{1.45000005}\smash{\begin{tabular}[t]{l}$n$\end{tabular}}}}%
    \put(0,0){\includegraphics[width=\unitlength,page=6]{phi.pdf}}%
    \put(0.83380915,0.38225211){\makebox(0,0)[lt]{\lineheight{1.45000005}\smash{\begin{tabular}[t]{l}$\eta$\end{tabular}}}}%
    \put(0.04571354,0.37910446){\makebox(0,0)[lt]{\lineheight{1.45000005}\smash{\begin{tabular}[t]{l}$\eta$\end{tabular}}}}%
  \end{picture}%
\endgroup%
}\ast \id\right) - \ti{u}\\
     &=\centre{
\begingroup%
  \makeatletter%
  \providecommand\color[2][]{%
    \errmessage{(Inkscape) Color is used for the text in Inkscape, but the package 'color.sty' is not loaded}%
    \renewcommand\color[2][]{}%
  }%
  \providecommand\transparent[1]{%
    \errmessage{(Inkscape) Transparency is used (non-zero) for the text in Inkscape, but the package 'transparent.sty' is not loaded}%
    \renewcommand\transparent[1]{}%
  }%
  \providecommand\rotatebox[2]{#2}%
  \newcommand*\fsize{\dimexpr\f@size pt\relax}%
  \newcommand*\lineheight[1]{\fontsize{\fsize}{#1\fsize}\selectfont}%
  \ifx\svgwidth\undefined%
    \setlength{\unitlength}{154.49999809bp}%
    \ifx\svgscale\undefined%
      \relax%
    \else%
      \setlength{\unitlength}{\unitlength * \real{\svgscale}}%
    \fi%
  \else%
    \setlength{\unitlength}{\svgwidth}%
  \fi%
  \global\let\svgwidth\undefined%
  \global\let\svgscale\undefined%
  \makeatother%
  \begin{picture}(1,1.03284592)%
    \lineheight{1}%
    \setlength\tabcolsep{0pt}%
    \put(0,0){\includegraphics[width=\unitlength,page=1]{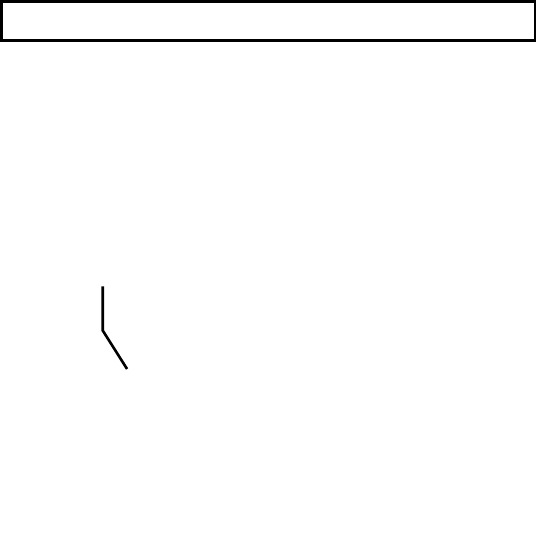}}%
    \put(0.49295342,0.97904106){\makebox(0,0)[lt]{\lineheight{1.45000005}\smash{\begin{tabular}[t]{l}$D$\end{tabular}}}}%
    \put(0.91004304,0.00899021){\makebox(0,0)[lt]{\lineheight{1.45000005}\smash{\begin{tabular}[t]{l}$n$\end{tabular}}}}%
    \put(0.66253689,0.00456361){\makebox(0,0)[lt]{\lineheight{1.45000005}\smash{\begin{tabular}[t]{l}$1$\end{tabular}}}}%
    \put(0,0){\includegraphics[width=\unitlength,page=2]{uphi.pdf}}%
    \put(0.3074322,0.29391435){\makebox(0,0)[lt]{\lineheight{1.45000005}\smash{\begin{tabular}[t]{l}$r$\end{tabular}}}}%
    \put(0.77917477,0.00717504){\makebox(0,0)[lt]{\lineheight{1.45000005}\smash{\begin{tabular}[t]{l}$i$\end{tabular}}}}%
    \put(0,0){\includegraphics[width=\unitlength,page=3]{uphi.pdf}}%
    \put(0.06289012,0.81768336){\makebox(0,0)[lt]{\lineheight{1.45000005}\smash{\begin{tabular}[t]{l}\tiny$1$\end{tabular}}}}%
    \put(0.88813743,0.81734884){\makebox(0,0)[lt]{\lineheight{1.45000005}\smash{\begin{tabular}[t]{l}\tiny$n$\end{tabular}}}}%
    \put(0,0){\includegraphics[width=\unitlength,page=4]{uphi.pdf}}%
    \put(0.1618264,0.63583659){\makebox(0,0)[lt]{\lineheight{1.45000005}\smash{\begin{tabular}[t]{l}$\rho$\end{tabular}}}}%
    \put(0,0){\includegraphics[width=\unitlength,page=5]{uphi.pdf}}%
    \put(0.15337165,0.51266778){\makebox(0,0)[lt]{\lineheight{1.45000005}\smash{\begin{tabular}[t]{l}$\pi$\end{tabular}}}}%
    \put(0,0){\includegraphics[width=\unitlength,page=6]{uphi.pdf}}%
    \put(0.64541441,0.81768336){\makebox(0,0)[lt]{\lineheight{1.45000005}\smash{\begin{tabular}[t]{l}\tiny$i$\end{tabular}}}}%
    \put(0,0){\includegraphics[width=\unitlength,page=7]{uphi.pdf}}%
    \put(0.39699526,0.44155816){\makebox(0,0)[lt]{\lineheight{1.45000005}\smash{\begin{tabular}[t]{l}$\epsilon$\end{tabular}}}}%
    \put(0.00444504,0.77509652){\makebox(0,0)[lt]{\lineheight{1.45000005}\smash{\begin{tabular}[t]{l}$\Delta$\end{tabular}}}}%
    \put(0.82248416,0.1067922){\makebox(0,0)[lt]{\lineheight{1.45000005}\smash{\begin{tabular}[t]{l}$\mu$\end{tabular}}}}%
    \put(0,0){\includegraphics[width=\unitlength,page=8]{uphi.pdf}}%
  \end{picture}%
\endgroup%
}-\centre{}\in A_{d,k+r}(n),
   \end{split}
  \end{gather*}
  where $\rho^{-1}\in\gpS_{n}$ is the $(r+1,n-r-1)$-shuffle that maps $[r+1]\subset [n]$ to $\{i_1,\cdots,i_{r+1}\}$ and $\pi\in \gpS_{r+1}$ satisfies $\pi^{-1}(j)=\rho(i_j)$ for any $j\in[r+1]$.
  By using Lemma \ref{l721}, we have for $\beta_1,\cdots,\beta_{r+1}\geq 0$,
  \begin{gather}\label{commcomp}
    \centre{
\begingroup%
  \makeatletter%
  \providecommand\color[2][]{%
    \errmessage{(Inkscape) Color is used for the text in Inkscape, but the package 'color.sty' is not loaded}%
    \renewcommand\color[2][]{}%
  }%
  \providecommand\transparent[1]{%
    \errmessage{(Inkscape) Transparency is used (non-zero) for the text in Inkscape, but the package 'transparent.sty' is not loaded}%
    \renewcommand\transparent[1]{}%
  }%
  \providecommand\rotatebox[2]{#2}%
  \newcommand*\fsize{\dimexpr\f@size pt\relax}%
  \newcommand*\lineheight[1]{\fontsize{\fsize}{#1\fsize}\selectfont}%
  \ifx\svgwidth\undefined%
    \setlength{\unitlength}{64.359327bp}%
    \ifx\svgscale\undefined%
      \relax%
    \else%
      \setlength{\unitlength}{\unitlength * \real{\svgscale}}%
    \fi%
  \else%
    \setlength{\unitlength}{\svgwidth}%
  \fi%
  \global\let\svgwidth\undefined%
  \global\let\svgscale\undefined%
  \makeatother%
  \begin{picture}(1,1.15703306)%
    \lineheight{1}%
    \setlength\tabcolsep{0pt}%
    \put(0,0){\includegraphics[width=\unitlength,page=1]{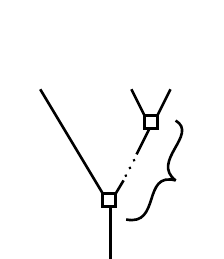}}%
    \put(0.83361105,0.28397921){\makebox(0,0)[lt]{\lineheight{1.45000005}\smash{\begin{tabular}[t]{l}$r$\end{tabular}}}}%
    \put(0,0){\includegraphics[width=\unitlength,page=2]{d_izero.pdf}}%
    \put(-0.00470381,0.84790352){\makebox(0,0)[lt]{\lineheight{1.45000005}\smash{\begin{tabular}[t]{l}$i$\end{tabular}}}}%
    \put(0.85928639,0.85119216){\makebox(0,0)[lt]{\lineheight{1.45000005}\smash{\begin{tabular}[t]{l}$i$\end{tabular}}}}%
    \put(0.09102075,1.12790978){\makebox(0,0)[lt]{\lineheight{1.45000005}\smash{\begin{tabular}[t]{l}\tiny$\beta_1$\end{tabular}}}}%
    \put(0.68888039,1.13311945){\makebox(0,0)[lt]{\lineheight{1.45000005}\smash{\begin{tabular}[t]{l}\tiny$\beta_{r+1}$\end{tabular}}}}%
    \put(0.49795332,1.1321876){\makebox(0,0)[lt]{\lineheight{1.45000005}\smash{\begin{tabular}[t]{l}\tiny$\beta_r$\end{tabular}}}}%
    \put(0,0){\includegraphics[width=\unitlength,page=3]{d_izero.pdf}}%
  \end{picture}%
\endgroup%
}=
    \begin{cases}
      \centre{} & (\beta_j=0 \text{ for all }j\in [r+1])\\
      (-1)^r \centre{
\begingroup%
  \makeatletter%
  \providecommand\color[2][]{%
    \errmessage{(Inkscape) Color is used for the text in Inkscape, but the package 'color.sty' is not loaded}%
    \renewcommand\color[2][]{}%
  }%
  \providecommand\transparent[1]{%
    \errmessage{(Inkscape) Transparency is used (non-zero) for the text in Inkscape, but the package 'transparent.sty' is not loaded}%
    \renewcommand\transparent[1]{}%
  }%
  \providecommand\rotatebox[2]{#2}%
  \newcommand*\fsize{\dimexpr\f@size pt\relax}%
  \newcommand*\lineheight[1]{\fontsize{\fsize}{#1\fsize}\selectfont}%
  \ifx\svgwidth\undefined%
    \setlength{\unitlength}{39.45260866bp}%
    \ifx\svgscale\undefined%
      \relax%
    \else%
      \setlength{\unitlength}{\unitlength * \real{\svgscale}}%
    \fi%
  \else%
    \setlength{\unitlength}{\svgwidth}%
  \fi%
  \global\let\svgwidth\undefined%
  \global\let\svgscale\undefined%
  \makeatother%
  \begin{picture}(1,1.33496126)%
    \lineheight{1}%
    \setlength\tabcolsep{0pt}%
    \put(0,0){\includegraphics[width=\unitlength,page=1]{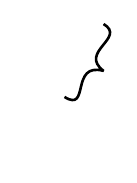}}%
    \put(0.83465142,0.74331404){\makebox(0,0)[lt]{\lineheight{1.45000005}\smash{\begin{tabular}[t]{l}$r$\end{tabular}}}}%
    \put(0,0){\includegraphics[width=\unitlength,page=2]{treer.pdf}}%
    \put(0.43584326,0.33894294){\makebox(0,0)[lt]{\lineheight{1.45000005}\smash{\begin{tabular}[t]{l}$i$\end{tabular}}}}%
  \end{picture}%
\endgroup%
} & (\beta_j=1 \text{ for all }j\in [r+1])\\
      \text{a linear sum of diagrams}\\
      \text{ with at least $r+1$ trivalent vertices}& (\text{otherwise}).
    \end{cases}
  \end{gather}
  In the last case, the corresponding term of $[\ti{u},\phi]$ is included in $A_{d,k+r+1}(n)$.

  Thus, by \eqref{commcomp} and Lemma \ref{l723} \eqref{l7232}, we have
  \begin{gather*}
    \begin{split}
     &\ti{\beta}_{d,k}^r([\phi]_r)(u)\\
     &=\sum_{(\alpha_l)\in[2d-k]^{r+1}}
     (-1)^r
     \left(\prod_{j=1}^{r+1} \langle v_{j_{\alpha_l}},\bar{x}_{i_l}\rangle\right) \theta_{d,n,k+r}\left(\centre{
\begingroup%
  \makeatletter%
  \providecommand\color[2][]{%
    \errmessage{(Inkscape) Color is used for the text in Inkscape, but the package 'color.sty' is not loaded}%
    \renewcommand\color[2][]{}%
  }%
  \providecommand\transparent[1]{%
    \errmessage{(Inkscape) Transparency is used (non-zero) for the text in Inkscape, but the package 'transparent.sty' is not loaded}%
    \renewcommand\transparent[1]{}%
  }%
  \providecommand\rotatebox[2]{#2}%
  \newcommand*\fsize{\dimexpr\f@size pt\relax}%
  \newcommand*\lineheight[1]{\fontsize{\fsize}{#1\fsize}\selectfont}%
  \ifx\svgwidth\undefined%
    \setlength{\unitlength}{135.74999903bp}%
    \ifx\svgscale\undefined%
      \relax%
    \else%
      \setlength{\unitlength}{\unitlength * \real{\svgscale}}%
    \fi%
  \else%
    \setlength{\unitlength}{\svgwidth}%
  \fi%
  \global\let\svgwidth\undefined%
  \global\let\svgscale\undefined%
  \makeatother%
  \begin{picture}(1,0.89299293)%
    \lineheight{1}%
    \setlength\tabcolsep{0pt}%
    \put(0,0){\includegraphics[width=\unitlength,page=1]{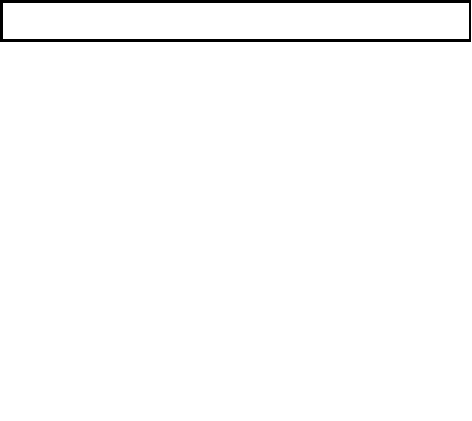}}%
    \put(0.46189061,0.83175646){\makebox(0,0)[lt]{\lineheight{1.45000005}\smash{\begin{tabular}[t]{l}$D$\end{tabular}}}}%
    \put(0.87172458,0.00769992){\makebox(0,0)[lt]{\lineheight{1.45000005}\smash{\begin{tabular}[t]{l}$n$\end{tabular}}}}%
    \put(0.39989459,0.00791324){\makebox(0,0)[lt]{\lineheight{1.45000005}\smash{\begin{tabular}[t]{l}$1$\end{tabular}}}}%
    \put(0,0){\includegraphics[width=\unitlength,page=2]{betaimi.pdf}}%
    \put(0.6157188,0.00519394){\makebox(0,0)[lt]{\lineheight{1.45000005}\smash{\begin{tabular}[t]{l}$i$\end{tabular}}}}%
    \put(0,0){\includegraphics[width=\unitlength,page=3]{betaimi.pdf}}%
    \put(0.46243522,0.6592956){\makebox(0,0)[lt]{\lineheight{1.45000005}\smash{\begin{tabular}[t]{l}$\tau$\end{tabular}}}}%
    \put(0,0){\includegraphics[width=\unitlength,page=4]{betaimi.pdf}}%
    \put(0.1209442,0.49051855){\makebox(0,0)[lt]{\lineheight{1.45000005}\smash{\begin{tabular}[t]{l}$\sigma$\end{tabular}}}}%
    \put(0,0){\includegraphics[width=\unitlength,page=5]{betaimi.pdf}}%
    \put(0.24501844,0.32773414){\makebox(0,0)[lt]{\lineheight{1.45000005}\smash{\begin{tabular}[t]{l}$r$\end{tabular}}}}%
  \end{picture}%
\endgroup%
}\right)\\
     &=(-1)^r \sum_{(\alpha_l)\in[2d-k]^{r+1}}
     \left(\prod_{l=1}^{r+1}\langle v_{j_{\alpha_l}},\bar{x}_{i_l}\rangle\right) \centre{}\\
     &=(-1)^r \gamma_{d,k}^r\cdot \ti{\tau}_r([\phi]_r)(u).
   \end{split}
  \end{gather*}

 \end{proof}


\section{The $\GL(V_n)$-module structure of $B_d(n)$}\label{s7}
 In this section, we consider the $\GL(V_n)$-module structure of $B_d(n)$
 and give a decomposition of $B_d(n)$ with respect to connected parts. Moreover, we compute the irreducible decomposition of $B_d(n)$ for $d=3,4,5$ and that of $B_{d,0}(n),B_{d,1}(n)$ for any $d$.
 Lastly, we show that the surjectivity of the bracket map which we defined in Section \ref{s4}.

 Let $B_{d,k}^c(n)\subset B_{d,k}(n)$ denote the \emph{connected part} of $B_{d,k}(n)$, which is spanned by connected $V_n$-colored open Jacobi diagrams. Let $D_{d,k}^c\subset D_{d,k}$ denote the \emph{connected part} of $D_{d,k}$, which is spanned by connected $[2d-k]$-colored open Jacobi diagrams.
 We have an isomorphism of $\GL(V_n)$-modules
 \begin{equation*}
   B_{d,k}^c(n)\cong V_n^{\otimes 2d-k}\otimes_{\K\gpS_{2d-k}} D_{d,k}^c,
 \end{equation*}
 which is the connected version of \eqref{eqBD}.

 The direct sum $\bigoplus_{d\geq 0}B_d(n)$ has the following coalgebra structure. This is an analogue of the coalgebra structure of the space of open Jacobi diagrams colored by one element \cite{BN_v}.
 Let $C=\bigcup_{i\in I}C_i$ be a presentation of a diagram $C\in \bigoplus_{d\geq 0}B_d(n)$ as the disjoint union of its connected components.
 The comultiplication $\Delta$ is defined by
 $$\Delta(C)=\sum_{J\subset I}\left(\bigcup_{i\in J}C_i\right)\otimes \left(\bigcup_{i\in I\setminus J}C_i\right).$$
 Note that the connected part $\bigoplus_{d,k\geq 0}B_{d,k}^c(n)$ coincides with the \emph{primitive part} of the coalgebra $\bigoplus_{d\geq 0}B_d(n)$.

 \subsection{Decomposition of $B_d(n)$ with respect to connected parts}\label{ss71}
  Note that $D_{d,k}^c\neq 0$ if and only if $d-1\leq k\leq 2d-2$ because each element of $D_{d,k}^c$ has at least two univalent vertex and is connected.
  For $d\geq 1,k\geq 0$, the pair $(d,k)$ is called a \emph{good pair} if $d-1\leq k\leq 2d-2$.
  We consider the following decomposition of a pair $(d,k)$ to consider the decomposition of an element of $D_{d,k}$ into the connected parts.
  \begin{definition}
   Let $d,\ k\geq 0$.
   A \emph{decomposition of $(d,k)$ into good pairs} is a sequence of triples of integers
   $$
     \pi = ((a_1,d_1, k_1),\cdots, (a_l,d_l, k_l))
   $$
   such that $(d_i,k_i)$ are good pairs, $a_i\geq 1$,
   $$
     \sum_{i=1}^l a_i d_i = d, \quad \sum_{i=1}^l a_i k_i = k,
   $$
   and
   $$
     \quad(d_1,k_1)>(d_2,k_2)>\cdots>(d_l,k_l)
   $$
   in the lexicographical order.

   Let $\Pi(d,k)$ be the set of all decompositions of $(d,k)$ into good pairs.
  \end{definition}
  For example, we have
  \begin{equation}\label{42}
    \Pi(4,2)=\{((1,3,2),(1,1,0)),((1,2,2),(2,1,0)),((2,2,1))\}.
  \end{equation}

  For any diagram $K\in D_{d,k}$, we can assign a decomposition of $(d,k)$ into good pairs such that $d_i$ and $k_i$ correspond to the degree and the number of trivalent vertices of each connected component of $K$, respectively, and $a_i$ corresponds to the multiplicity of $(d_i,k_i)$.
  We call a coloring of $K=\bigsqcup_{1\leq i\leq l, 1\leq j\leq a_i} K_{i}^{(j)}\in D_{d,k}$ \emph{standard}
  if the set of colorings of $K_{i}^{(j)}\in D_{d_i,k_i}^c$ is
  $$\{\sum_{p=1}^{i-1}(2d_p-k_p)a_p +(j-1)(2d_i-k_i)+1 ,\cdots, \sum_{p=1}^{i-1}(2d_p-k_p)a_p +j(2d_i-k_i)\}$$
  for each $i\in[l],j\in[a_i]$.

  \begin{theorem}\label{th521}
   For $d,k,n\geq 0$, we have an isomorphism of $\GL(V_n)$-modules
   \begin{equation}\label{eqB}
     B_{d,k}(n) \cong
        \bigoplus_{\pi =((a_1,d_1, k_1),\cdots,(a_l,d_l, k_l))\in\Pi(d,k)}
        \left(\bigotimes_{i=1}^l \Sym^{a_i}(B_{d_i,k_i}^c(n))\right).
   \end{equation}
  \end{theorem}

  To prove this, we need the following proposition.

  \begin{proposition}\label{p521}
   Let $d,k\geq 0$. We have an isomorphism of $\gpS_{2d-k}$-modules
   \begin{equation}\label{eqD}
     D_{d,k}\cong
     \bigoplus_{\pi =((a_1,d_1,k_1),\cdots,(a_l,d_l,k_l))\in\Pi(d,k)}
      \Ind_{\prod_{i=1}^l(\gpS_{2d_i-k_i}\wr\gpS_{a_i})}^{\gpS_{2d-k}}
         (\bigotimes_{i=1}^l(D_{d_i,k_i}^c)^{\otimes a_i}),
   \end{equation}
   where $\gpS_{2d_i-k_i}\wr\gpS_{a_i}={\gpS_{2d_i-k_i}}^{a_i}\rtimes\gpS_{a_i}\subset \gpS_{(2d_i-k_i)a_i}$ is the wreath product.
  \end{proposition}

  For example, we have an isomorphism of $\gpS_6$-modules for $(d,k)=(4,2)$, which corresponds to \eqref{42},
  $$D_{4,2}\cong
   \Ind_{\gpS_{4}\times \gpS_{2}}^{\gpS_{6}} (D_{3,2}^c \otimes D_{1,0}^c)
   \oplus
   \Ind_{\gpS_{2}\times (\gpS_{2}\wr\gpS_{2})}^{\gpS_{6}} (D_{2,2}^c \otimes (D_{1,0}^c)^{\otimes 2})
   \oplus
   \Ind_{\gpS_{3}\wr\gpS_{2}}^{\gpS_{6}} (D_{2,1}^c)^{\otimes 2}.
  $$
  For example, $$\centre{
\begingroup%
  \makeatletter%
  \providecommand\color[2][]{%
    \errmessage{(Inkscape) Color is used for the text in Inkscape, but the package 'color.sty' is not loaded}%
    \renewcommand\color[2][]{}%
  }%
  \providecommand\transparent[1]{%
    \errmessage{(Inkscape) Transparency is used (non-zero) for the text in Inkscape, but the package 'transparent.sty' is not loaded}%
    \renewcommand\transparent[1]{}%
  }%
  \providecommand\rotatebox[2]{#2}%
  \newcommand*\fsize{\dimexpr\f@size pt\relax}%
  \newcommand*\lineheight[1]{\fontsize{\fsize}{#1\fsize}\selectfont}%
  \ifx\svgwidth\undefined%
    \setlength{\unitlength}{35.26304341bp}%
    \ifx\svgscale\undefined%
      \relax%
    \else%
      \setlength{\unitlength}{\unitlength * \real{\svgscale}}%
    \fi%
  \else%
    \setlength{\unitlength}{\svgwidth}%
  \fi%
  \global\let\svgwidth\undefined%
  \global\let\svgscale\undefined%
  \makeatother%
  \begin{picture}(1,0.63318251)%
    \lineheight{1}%
    \setlength\tabcolsep{0pt}%
    \put(0,0){\includegraphics[width=\unitlength,page=1]{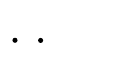}}%
    \put(-0.00858503,0.02210553){\makebox(0,0)[lt]{\lineheight{1.45000005}\smash{\begin{tabular}[t]{l}$1$\end{tabular}}}}%
    \put(0.25845581,0.02555271){\makebox(0,0)[lt]{\lineheight{1.45000005}\smash{\begin{tabular}[t]{l}$3$\end{tabular}}}}%
    \put(0,0){\includegraphics[width=\unitlength,page=2]{ex11.pdf}}%
    \put(0.5242601,0.0199948){\makebox(0,0)[lt]{\lineheight{1.45000005}\smash{\begin{tabular}[t]{l}$2$\end{tabular}}}}%
    \put(0.79130078,0.02344198){\makebox(0,0)[lt]{\lineheight{1.45000005}\smash{\begin{tabular}[t]{l}$4$\end{tabular}}}}%
  \end{picture}%
\endgroup%
}\otimes \centre{
\begingroup%
  \makeatletter%
  \providecommand\color[2][]{%
    \errmessage{(Inkscape) Color is used for the text in Inkscape, but the package 'color.sty' is not loaded}%
    \renewcommand\color[2][]{}%
  }%
  \providecommand\transparent[1]{%
    \errmessage{(Inkscape) Transparency is used (non-zero) for the text in Inkscape, but the package 'transparent.sty' is not loaded}%
    \renewcommand\transparent[1]{}%
  }%
  \providecommand\rotatebox[2]{#2}%
  \newcommand*\fsize{\dimexpr\f@size pt\relax}%
  \newcommand*\lineheight[1]{\fontsize{\fsize}{#1\fsize}\selectfont}%
  \ifx\svgwidth\undefined%
    \setlength{\unitlength}{33.36463348bp}%
    \ifx\svgscale\undefined%
      \relax%
    \else%
      \setlength{\unitlength}{\unitlength * \real{\svgscale}}%
    \fi%
  \else%
    \setlength{\unitlength}{\svgwidth}%
  \fi%
  \global\let\svgwidth\undefined%
  \global\let\svgscale\undefined%
  \makeatother%
  \begin{picture}(1,0.1168453)%
    \lineheight{1}%
    \setlength\tabcolsep{0pt}%
    \put(0,0){\includegraphics[width=\unitlength,page=1]{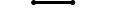}}%
    \put(-0.0090735,0.0230201){\makebox(0,0)[lt]{\lineheight{1.45000005}\smash{\begin{tabular}[t]{l}$1$\end{tabular}}}}%
    \put(0.77942603,0.02113248){\makebox(0,0)[lt]{\lineheight{1.45000005}\smash{\begin{tabular}[t]{l}$2$\end{tabular}}}}%
  \end{picture}%
\endgroup%
} \in \Ind_{\gpS_{4}\times \gpS_{2}}^{\gpS_{6}} (D_{3,2}^c \otimes D_{1,0}^c),$$
  $$\centre{
\begingroup%
  \makeatletter%
  \providecommand\color[2][]{%
    \errmessage{(Inkscape) Color is used for the text in Inkscape, but the package 'color.sty' is not loaded}%
    \renewcommand\color[2][]{}%
  }%
  \providecommand\transparent[1]{%
    \errmessage{(Inkscape) Transparency is used (non-zero) for the text in Inkscape, but the package 'transparent.sty' is not loaded}%
    \renewcommand\transparent[1]{}%
  }%
  \providecommand\rotatebox[2]{#2}%
  \newcommand*\fsize{\dimexpr\f@size pt\relax}%
  \newcommand*\lineheight[1]{\fontsize{\fsize}{#1\fsize}\selectfont}%
  \ifx\svgwidth\undefined%
    \setlength{\unitlength}{44.71747486bp}%
    \ifx\svgscale\undefined%
      \relax%
    \else%
      \setlength{\unitlength}{\unitlength * \real{\svgscale}}%
    \fi%
  \else%
    \setlength{\unitlength}{\svgwidth}%
  \fi%
  \global\let\svgwidth\undefined%
  \global\let\svgscale\undefined%
  \makeatother%
  \begin{picture}(1,0.18451798)%
    \lineheight{1}%
    \setlength\tabcolsep{0pt}%
    \put(0,0){\includegraphics[width=\unitlength,page=1]{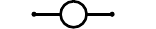}}%
    \put(-0.00676993,0.03460368){\makebox(0,0)[lt]{\lineheight{1.45000005}\smash{\begin{tabular}[t]{l}$1$\end{tabular}}}}%
    \put(0.8354252,0.03477331){\makebox(0,0)[lt]{\lineheight{1.45000005}\smash{\begin{tabular}[t]{l}$2$\end{tabular}}}}%
  \end{picture}%
\endgroup%
} \otimes \centre{} \otimes \centre{}\in \Ind_{\gpS_{2}\times (\gpS_{2}\wr\gpS_{2})}^{\gpS_{6}} (D_{2,2}^c \otimes (D_{1,0}^c)^{\otimes 2})$$
  and $$\centre{
\begingroup%
  \makeatletter%
  \providecommand\color[2][]{%
    \errmessage{(Inkscape) Color is used for the text in Inkscape, but the package 'color.sty' is not loaded}%
    \renewcommand\color[2][]{}%
  }%
  \providecommand\transparent[1]{%
    \errmessage{(Inkscape) Transparency is used (non-zero) for the text in Inkscape, but the package 'transparent.sty' is not loaded}%
    \renewcommand\transparent[1]{}%
  }%
  \providecommand\rotatebox[2]{#2}%
  \newcommand*\fsize{\dimexpr\f@size pt\relax}%
  \newcommand*\lineheight[1]{\fontsize{\fsize}{#1\fsize}\selectfont}%
  \ifx\svgwidth\undefined%
    \setlength{\unitlength}{24.22929078bp}%
    \ifx\svgscale\undefined%
      \relax%
    \else%
      \setlength{\unitlength}{\unitlength * \real{\svgscale}}%
    \fi%
  \else%
    \setlength{\unitlength}{\svgwidth}%
  \fi%
  \global\let\svgwidth\undefined%
  \global\let\svgscale\undefined%
  \makeatother%
  \begin{picture}(1,0.93831901)%
    \lineheight{1}%
    \setlength\tabcolsep{0pt}%
    \put(0,0){\includegraphics[width=\unitlength,page=1]{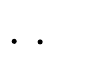}}%
    \put(-0.01249455,0.03968702){\makebox(0,0)[lt]{\lineheight{1.45000005}\smash{\begin{tabular}[t]{l}$1$\end{tabular}}}}%
    \put(0.69626146,0.02910022){\makebox(0,0)[lt]{\lineheight{1.45000005}\smash{\begin{tabular}[t]{l}$3$\end{tabular}}}}%
    \put(0,0){\includegraphics[width=\unitlength,page=2]{ex12.pdf}}%
    \put(0.35162214,0.03380223){\makebox(0,0)[lt]{\lineheight{1.45000005}\smash{\begin{tabular}[t]{l}$2$\end{tabular}}}}%
  \end{picture}%
\endgroup%
}\otimes \centre{} \in \Ind_{\gpS_{3}\wr\gpS_{2}}^{\gpS_{6}} (D_{2,1}^c)^{\otimes 2}.$$
  Via the above isomorphism, the element
  $$(2,3)(4,5)\cdot(\centre{} \otimes \centre{} \otimes \centre{})\in \Ind_{\gpS_{2}\times (\gpS_{2}\wr\gpS_{2})}^{\gpS_{6}} (D_{2,2}^c \otimes (D_{1,0}^c)^{\otimes 2})$$
  corresponds to the element
  $$\centre{
\begingroup%
  \makeatletter%
  \providecommand\color[2][]{%
    \errmessage{(Inkscape) Color is used for the text in Inkscape, but the package 'color.sty' is not loaded}%
    \renewcommand\color[2][]{}%
  }%
  \providecommand\transparent[1]{%
    \errmessage{(Inkscape) Transparency is used (non-zero) for the text in Inkscape, but the package 'transparent.sty' is not loaded}%
    \renewcommand\transparent[1]{}%
  }%
  \providecommand\rotatebox[2]{#2}%
  \newcommand*\fsize{\dimexpr\f@size pt\relax}%
  \newcommand*\lineheight[1]{\fontsize{\fsize}{#1\fsize}\selectfont}%
  \ifx\svgwidth\undefined%
    \setlength{\unitlength}{127.75237937bp}%
    \ifx\svgscale\undefined%
      \relax%
    \else%
      \setlength{\unitlength}{\unitlength * \real{\svgscale}}%
    \fi%
  \else%
    \setlength{\unitlength}{\svgwidth}%
  \fi%
  \global\let\svgwidth\undefined%
  \global\let\svgscale\undefined%
  \makeatother%
  \begin{picture}(1,0.06458727)%
    \lineheight{1}%
    \setlength\tabcolsep{0pt}%
    \put(0,0){\includegraphics[width=\unitlength,page=1]{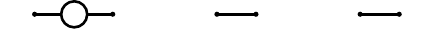}}%
    \put(-0.00236969,0.01256311){\makebox(0,0)[lt]{\lineheight{1.45000005}\smash{\begin{tabular}[t]{l}$1$\end{tabular}}}}%
    \put(0.29457536,0.01184174){\makebox(0,0)[lt]{\lineheight{1.45000005}\smash{\begin{tabular}[t]{l}$3$\end{tabular}}}}%
    \put(0.41091223,0.0122128){\makebox(0,0)[lt]{\lineheight{1.45000005}\smash{\begin{tabular}[t]{l}$2$\end{tabular}}}}%
    \put(0.61611152,0.01282443){\makebox(0,0)[lt]{\lineheight{1.45000005}\smash{\begin{tabular}[t]{l}$5$\end{tabular}}}}%
    \put(0.73591405,0.01349109){\makebox(0,0)[lt]{\lineheight{1.45000005}\smash{\begin{tabular}[t]{l}$4$\end{tabular}}}}%
    \put(0.94239348,0.01371453){\makebox(0,0)[lt]{\lineheight{1.45000005}\smash{\begin{tabular}[t]{l}$6$\end{tabular}}}}%
  \end{picture}%
\endgroup%
}=(2,3)(4,5)\cdot(\centre{
\begingroup%
  \makeatletter%
  \providecommand\color[2][]{%
    \errmessage{(Inkscape) Color is used for the text in Inkscape, but the package 'color.sty' is not loaded}%
    \renewcommand\color[2][]{}%
  }%
  \providecommand\transparent[1]{%
    \errmessage{(Inkscape) Transparency is used (non-zero) for the text in Inkscape, but the package 'transparent.sty' is not loaded}%
    \renewcommand\transparent[1]{}%
  }%
  \providecommand\rotatebox[2]{#2}%
  \newcommand*\fsize{\dimexpr\f@size pt\relax}%
  \newcommand*\lineheight[1]{\fontsize{\fsize}{#1\fsize}\selectfont}%
  \ifx\svgwidth\undefined%
    \setlength{\unitlength}{127.75237937bp}%
    \ifx\svgscale\undefined%
      \relax%
    \else%
      \setlength{\unitlength}{\unitlength * \real{\svgscale}}%
    \fi%
  \else%
    \setlength{\unitlength}{\svgwidth}%
  \fi%
  \global\let\svgwidth\undefined%
  \global\let\svgscale\undefined%
  \makeatother%
  \begin{picture}(1,0.06458727)%
    \lineheight{1}%
    \setlength\tabcolsep{0pt}%
    \put(0,0){\includegraphics[width=\unitlength,page=1]{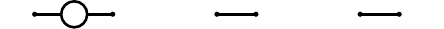}}%
    \put(-0.00236969,0.01256311){\makebox(0,0)[lt]{\lineheight{1.45000005}\smash{\begin{tabular}[t]{l}$1$\end{tabular}}}}%
    \put(0.29457536,0.01184174){\makebox(0,0)[lt]{\lineheight{1.45000005}\smash{\begin{tabular}[t]{l}$2$\end{tabular}}}}%
    \put(0.41091223,0.0122128){\makebox(0,0)[lt]{\lineheight{1.45000005}\smash{\begin{tabular}[t]{l}$3$\end{tabular}}}}%
    \put(0.61611152,0.01282443){\makebox(0,0)[lt]{\lineheight{1.45000005}\smash{\begin{tabular}[t]{l}$4$\end{tabular}}}}%
    \put(0.73591405,0.01349109){\makebox(0,0)[lt]{\lineheight{1.45000005}\smash{\begin{tabular}[t]{l}$5$\end{tabular}}}}%
    \put(0.94239348,0.01371453){\makebox(0,0)[lt]{\lineheight{1.45000005}\smash{\begin{tabular}[t]{l}$6$\end{tabular}}}}%
  \end{picture}%
\endgroup%
})\in D_{4,2}.$$

  \begin{proof}[Proof of Proposition \ref{p521}]
   Let $D'_{d,k}$ denote the right-hand side of \eqref{eqD}.

   For any coset $\sigma\in \gpS_{2d-k}/\prod_{i=1}^l (\gpS_{2d_i-k_i}\wr\gpS_{a_i})$, we fix a representative $\ti{\sigma}\in \gpS_{2d-k}$ of $\sigma$.

   Any element of $D'_{d,k}$ can be written uniquely as a linear sum of
   $$
     K=\ti{\sigma}\cdot\bigotimes_{1\leq i\leq l, 1\leq j\leq a_i}
     K_{i}^{(j)},
   $$
   where $K_{i}^{(j)}\in D_{d_i,k_i}^c$.
   We assign $\bigsqcup_{1\leq i\leq l, 1\leq j\leq a_i}K_{i}^{(j)}$ a standard coloring in $[2d-k]$ according to the order of the colorings in $\bigsqcup_{i=1}^{l} [2d_i-k_i]^{a_i}$ of $\bigotimes_{1\leq i\leq l, 1\leq j\leq a_i} K_{i}^{(j)}$.
   For example, if
   $$\bigotimes_{1\leq i\leq l, 1\leq j\leq a_i} K_{i}^{(j)}= \centre{
\begingroup%
  \makeatletter%
  \providecommand\color[2][]{%
    \errmessage{(Inkscape) Color is used for the text in Inkscape, but the package 'color.sty' is not loaded}%
    \renewcommand\color[2][]{}%
  }%
  \providecommand\transparent[1]{%
    \errmessage{(Inkscape) Transparency is used (non-zero) for the text in Inkscape, but the package 'transparent.sty' is not loaded}%
    \renewcommand\transparent[1]{}%
  }%
  \providecommand\rotatebox[2]{#2}%
  \newcommand*\fsize{\dimexpr\f@size pt\relax}%
  \newcommand*\lineheight[1]{\fontsize{\fsize}{#1\fsize}\selectfont}%
  \ifx\svgwidth\undefined%
    \setlength{\unitlength}{45.29298517bp}%
    \ifx\svgscale\undefined%
      \relax%
    \else%
      \setlength{\unitlength}{\unitlength * \real{\svgscale}}%
    \fi%
  \else%
    \setlength{\unitlength}{\svgwidth}%
  \fi%
  \global\let\svgwidth\undefined%
  \global\let\svgscale\undefined%
  \makeatother%
  \begin{picture}(1,0.80466442)%
    \lineheight{1}%
    \setlength\tabcolsep{0pt}%
    \put(0,0){\includegraphics[width=\unitlength,page=1]{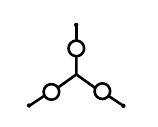}}%
    \put(0.40313431,0.73670417){\makebox(0,0)[lt]{\lineheight{1.45000005}\smash{\begin{tabular}[t]{l}$1$\end{tabular}}}}%
    \put(-0.00668391,0.01556704){\makebox(0,0)[lt]{\lineheight{1.45000005}\smash{\begin{tabular}[t]{l}$2$\end{tabular}}}}%
    \put(0.83751635,0.0185679){\makebox(0,0)[lt]{\lineheight{1.45000005}\smash{\begin{tabular}[t]{l}$3$\end{tabular}}}}%
  \end{picture}%
\endgroup%
}\otimes \centre{
\begingroup%
  \makeatletter%
  \providecommand\color[2][]{%
    \errmessage{(Inkscape) Color is used for the text in Inkscape, but the package 'color.sty' is not loaded}%
    \renewcommand\color[2][]{}%
  }%
  \providecommand\transparent[1]{%
    \errmessage{(Inkscape) Transparency is used (non-zero) for the text in Inkscape, but the package 'transparent.sty' is not loaded}%
    \renewcommand\transparent[1]{}%
  }%
  \providecommand\rotatebox[2]{#2}%
  \newcommand*\fsize{\dimexpr\f@size pt\relax}%
  \newcommand*\lineheight[1]{\fontsize{\fsize}{#1\fsize}\selectfont}%
  \ifx\svgwidth\undefined%
    \setlength{\unitlength}{59.88938044bp}%
    \ifx\svgscale\undefined%
      \relax%
    \else%
      \setlength{\unitlength}{\unitlength * \real{\svgscale}}%
    \fi%
  \else%
    \setlength{\unitlength}{\svgwidth}%
  \fi%
  \global\let\svgwidth\undefined%
  \global\let\svgscale\undefined%
  \makeatother%
  \begin{picture}(1,0.82657007)%
    \lineheight{1}%
    \setlength\tabcolsep{0pt}%
    \put(0,0){\includegraphics[width=\unitlength,page=1]{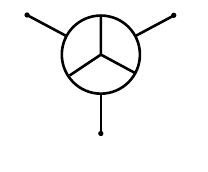}}%
    \put(0.44557471,0.011773){\makebox(0,0)[lt]{\lineheight{1.45000005}\smash{\begin{tabular}[t]{l}$1$\end{tabular}}}}%
    \put(-0.00505489,0.77517327){\makebox(0,0)[lt]{\lineheight{1.45000005}\smash{\begin{tabular}[t]{l}$2$\end{tabular}}}}%
    \put(0.87711729,0.77285091){\makebox(0,0)[lt]{\lineheight{1.45000005}\smash{\begin{tabular}[t]{l}$3$\end{tabular}}}}%
  \end{picture}%
\endgroup%
}\otimes \centre{
\begingroup%
  \makeatletter%
  \providecommand\color[2][]{%
    \errmessage{(Inkscape) Color is used for the text in Inkscape, but the package 'color.sty' is not loaded}%
    \renewcommand\color[2][]{}%
  }%
  \providecommand\transparent[1]{%
    \errmessage{(Inkscape) Transparency is used (non-zero) for the text in Inkscape, but the package 'transparent.sty' is not loaded}%
    \renewcommand\transparent[1]{}%
  }%
  \providecommand\rotatebox[2]{#2}%
  \newcommand*\fsize{\dimexpr\f@size pt\relax}%
  \newcommand*\lineheight[1]{\fontsize{\fsize}{#1\fsize}\selectfont}%
  \ifx\svgwidth\undefined%
    \setlength{\unitlength}{51.95463606bp}%
    \ifx\svgscale\undefined%
      \relax%
    \else%
      \setlength{\unitlength}{\unitlength * \real{\svgscale}}%
    \fi%
  \else%
    \setlength{\unitlength}{\svgwidth}%
  \fi%
  \global\let\svgwidth\undefined%
  \global\let\svgscale\undefined%
  \makeatother%
  \begin{picture}(1,0.54787743)%
    \lineheight{1}%
    \setlength\tabcolsep{0pt}%
    \put(0,0){\includegraphics[width=\unitlength,page=1]{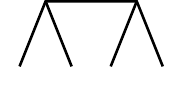}}%
    \put(-0.00582689,0.01552318){\makebox(0,0)[lt]{\lineheight{1.45000005}\smash{\begin{tabular}[t]{l}$1$\end{tabular}}}}%
    \put(0.28253797,0.01357102){\makebox(0,0)[lt]{\lineheight{1.45000005}\smash{\begin{tabular}[t]{l}$3$\end{tabular}}}}%
    \put(0.54718944,0.01472011){\makebox(0,0)[lt]{\lineheight{1.45000005}\smash{\begin{tabular}[t]{l}$4$\end{tabular}}}}%
    \put(0.85835009,0.01560685){\makebox(0,0)[lt]{\lineheight{1.45000005}\smash{\begin{tabular}[t]{l}$2$\end{tabular}}}}%
    \put(0,0){\includegraphics[width=\unitlength,page=2]{D32.pdf}}%
  \end{picture}%
\endgroup%
}\otimes \centre{
\begingroup%
  \makeatletter%
  \providecommand\color[2][]{%
    \errmessage{(Inkscape) Color is used for the text in Inkscape, but the package 'color.sty' is not loaded}%
    \renewcommand\color[2][]{}%
  }%
  \providecommand\transparent[1]{%
    \errmessage{(Inkscape) Transparency is used (non-zero) for the text in Inkscape, but the package 'transparent.sty' is not loaded}%
    \renewcommand\transparent[1]{}%
  }%
  \providecommand\rotatebox[2]{#2}%
  \newcommand*\fsize{\dimexpr\f@size pt\relax}%
  \newcommand*\lineheight[1]{\fontsize{\fsize}{#1\fsize}\selectfont}%
  \ifx\svgwidth\undefined%
    \setlength{\unitlength}{44.00624799bp}%
    \ifx\svgscale\undefined%
      \relax%
    \else%
      \setlength{\unitlength}{\unitlength * \real{\svgscale}}%
    \fi%
  \else%
    \setlength{\unitlength}{\svgwidth}%
  \fi%
  \global\let\svgwidth\undefined%
  \global\let\svgscale\undefined%
  \makeatother%
  \begin{picture}(1,0.18750015)%
    \lineheight{1}%
    \setlength\tabcolsep{0pt}%
    \put(0,0){\includegraphics[width=\unitlength,page=1]{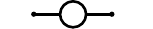}}%
    \put(-0.00687934,0.03329945){\makebox(0,0)[lt]{\lineheight{1.45000005}\smash{\begin{tabular}[t]{l}$1$\end{tabular}}}}%
    \put(0.83276535,0.03329945){\makebox(0,0)[lt]{\lineheight{1.45000005}\smash{\begin{tabular}[t]{l}$2$\end{tabular}}}}%
  \end{picture}%
\endgroup%
},$$
   then the corresponding coloring of $\bigsqcup_{1\leq i\leq l, 1\leq j\leq a_i}K_{i}^{(j)}$ is
   $$\centre{
\begingroup%
  \makeatletter%
  \providecommand\color[2][]{%
    \errmessage{(Inkscape) Color is used for the text in Inkscape, but the package 'color.sty' is not loaded}%
    \renewcommand\color[2][]{}%
  }%
  \providecommand\transparent[1]{%
    \errmessage{(Inkscape) Transparency is used (non-zero) for the text in Inkscape, but the package 'transparent.sty' is not loaded}%
    \renewcommand\transparent[1]{}%
  }%
  \providecommand\rotatebox[2]{#2}%
  \newcommand*\fsize{\dimexpr\f@size pt\relax}%
  \newcommand*\lineheight[1]{\fontsize{\fsize}{#1\fsize}\selectfont}%
  \ifx\svgwidth\undefined%
    \setlength{\unitlength}{227.48867172bp}%
    \ifx\svgscale\undefined%
      \relax%
    \else%
      \setlength{\unitlength}{\unitlength * \real{\svgscale}}%
    \fi%
  \else%
    \setlength{\unitlength}{\svgwidth}%
  \fi%
  \global\let\svgwidth\undefined%
  \global\let\svgscale\undefined%
  \makeatother%
  \begin{picture}(1,0.21760544)%
    \lineheight{1}%
    \setlength\tabcolsep{0pt}%
    \put(0,0){\includegraphics[width=\unitlength,page=1]{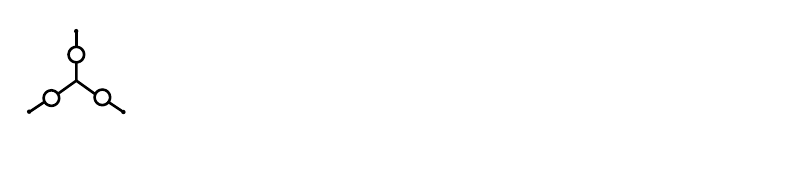}}%
    \put(0.08026403,0.19625209){\makebox(0,0)[lt]{\lineheight{1.45000005}\smash{\begin{tabular}[t]{l}$1$\end{tabular}}}}%
    \put(-0.00133077,0.05267371){\makebox(0,0)[lt]{\lineheight{1.45000005}\smash{\begin{tabular}[t]{l}$2$\end{tabular}}}}%
    \put(0.16674948,0.05327119){\makebox(0,0)[lt]{\lineheight{1.45000005}\smash{\begin{tabular}[t]{l}$3$\end{tabular}}}}%
    \put(0,0){\includegraphics[width=\unitlength,page=2]{D1518.pdf}}%
    \put(0.33350517,0.0030994){\makebox(0,0)[lt]{\lineheight{1.45000005}\smash{\begin{tabular}[t]{l}$4$\end{tabular}}}}%
    \put(0.21487101,0.20407456){\makebox(0,0)[lt]{\lineheight{1.45000005}\smash{\begin{tabular}[t]{l}$5$\end{tabular}}}}%
    \put(0.44711442,0.20346317){\makebox(0,0)[lt]{\lineheight{1.45000005}\smash{\begin{tabular}[t]{l}$6$\end{tabular}}}}%
    \put(0,0){\includegraphics[width=\unitlength,page=3]{D1518.pdf}}%
    \put(0.49911622,0.02854753){\makebox(0,0)[lt]{\lineheight{1.45000005}\smash{\begin{tabular}[t]{l}$7$\end{tabular}}}}%
    \put(0.57317001,0.02810169){\makebox(0,0)[lt]{\lineheight{1.45000005}\smash{\begin{tabular}[t]{l}$9$\end{tabular}}}}%
    \put(0.62541599,0.02836412){\makebox(0,0)[lt]{\lineheight{1.45000005}\smash{\begin{tabular}[t]{l}$10$\end{tabular}}}}%
    \put(0.69647993,0.02856666){\makebox(0,0)[lt]{\lineheight{1.45000005}\smash{\begin{tabular}[t]{l}$8$\end{tabular}}}}%
    \put(0,0){\includegraphics[width=\unitlength,page=4]{D1518.pdf}}%
    \put(0.7831018,0.10381753){\makebox(0,0)[lt]{\lineheight{1.45000005}\smash{\begin{tabular}[t]{l}$11$\end{tabular}}}}%
    \put(0.95645391,0.10381753){\makebox(0,0)[lt]{\lineheight{1.45000005}\smash{\begin{tabular}[t]{l}$12$\end{tabular}}}}%
  \end{picture}%
\endgroup%
}.$$

   Define a map $\Psi:D'_{d,k}\rightarrow D_{d,k}$ by
   $$
     \Psi(K)=\ti{\sigma}\cdot \bigsqcup_{1\leq i\leq l, 1\leq j\leq a_i} K_{i}^{(j)},
   $$
   where $\ti{\sigma}\in \gpS_{2d-k}$ acts on the colorings in $[2d-k]$.
   We can check that the map $\Psi$ is an $\gpS_{2d-k}$-module map.

   We need to check that $\Psi$ is bijective.
   If we have $\Psi(K)=\Psi(L)$
   for $K=\ti{\sigma}\cdot\bigotimes_{1\leq i\leq l, 1\leq j\leq a_i}
   K_{i}^{(j)}$, $L=\ti{\tau}\cdot\bigotimes_{1\leq i\leq l, 1\leq j\leq a_i}
   L_{i}^{(j)}$,
   then we have $\sigma=\tau$ by looking at the set of colorings of each connected component.
   Since we fix the representatives of cosets of $\gpS_{2d-k}/\prod_{i=1}^l (\gpS_{2d_i-k_i}\wr\gpS_{a_i})$, we have $\ti{\sigma}=\ti{\tau}$. Thus, we have $K=L$ and $\Psi$ is injective.
   For any element $K\in D_{d,k}$, we can take $\sigma\in \gpS_{2d-k}/\prod_{i=1}^l (\gpS_{2d_i-k_i}\wr\gpS_{a_i})$ such that $K=\ti{\sigma}\cdot \bigsqcup_{1\leq i\leq l, 1\leq j\leq a_i} K_{i}^{(j)}$, where $K_{i}^{(j)}\in D_{(d_i,k_i)}^c$
   and $\bigsqcup_{1\leq i\leq l, 1\leq j\leq a_i} K_{i}^{(j)}$ has a standard coloring.
   Therefore, $\Psi$ is surjective.
  \end{proof}

  \begin{proof}[Proof of Theorem \ref{th521}]
   By Proposition \ref{p521}, we have
   \begin{gather*}
    \begin{split}
     B_{d,k}(n)
      &\cong V_n^{\otimes 2d-k}\otimes_{\K\gpS_{2d-k}}D_{d,k} \\
      &\cong\bigoplus_{\pi\in \Pi(d,k)}
       \left(V_n^{\otimes 2d-k}\otimes_{\K\gpS_{2d-k}}
          \Ind_{\prod_{i=1}^l (\gpS_{2d_i-k_i}\wr\gpS_{a_i})}^{\gpS_{2d-k}} \left(\bigotimes_{i=1}^l(D_{d_i,k_i}^c)^{\otimes a_i} \right)
         \right).
    \end{split}
   \end{gather*}
   Moreover, we can check \eqref{eqB} as follows.
   \begin{gather*}
    \begin{split}
      &V_n^{\otimes 2d-k}\otimes_{\K\gpS_{2d-k}}
         \Ind_{\prod_{i=1}^l (\gpS_{2d_i-k_i}\wr\gpS_{a_i})}^{\gpS_{2d-k}} \left(\bigotimes_{i=1}^l(D_{d_i,k_i}^c)^{\otimes a_i} \right)\\
      &\cong
         V_n^{\otimes 2d-k}\otimes_{\K\gpS_{2d-k}}
         \Ind_{\prod_{i=1}^l \gpS_{a_i(2d_i-k_i)}}^{\gpS_{2d-k}}
          \left(
         \Ind_{\prod_{i=1}^l (\gpS_{2d_i-k_i}\wr\gpS_{a_i})}^{\prod_{i=1}^l \gpS_{a_i(2d_i-k_i)}}
          \left(\bigotimes_{i=1}^l(D_{d_i,k_i}^c)^{\otimes a_i} \right)
          \right)\\
      &\cong
         V_n^{\otimes 2d-k}\otimes_{\K\gpS_{2d-k}}
         \Ind_{\prod_{i=1}^l \gpS_{a_i(2d_i-k_i)}}^{\gpS_{2d-k}}
          \left(
          \bigotimes_{i=1}^l
          \Ind_{\gpS_{2d_i-k_i}\wr\gpS_{a_i}}^{\gpS_{a_i(2d_i-k_i)}}
           \left((D_{d_i,k_i}^c)^{\otimes a_i} \right)
           \right)\\
      &\cong
         V_n^{\otimes 2d-k}\otimes_{\K(\prod_{i=1}^l \gpS_{a_i(2d_i-k_i)})}
         \left(
         \bigotimes_{i=1}^l
         \Ind_{\gpS_{2d_i-k_i}\wr\gpS_{a_i}}^{\gpS_{a_i(2d_i-k_i)}}
          \left((D_{d_i,k_i}^c)^{\otimes a_i} \right)
          \right)\\
      &\cong
         \bigotimes_{i=1}^l
         \left(
          V_n^{\otimes a_i(2d_i-k_i)}\otimes_{\K\gpS_{a_i(2d_i-k_i)}}
           \left(
          \Ind_{\gpS_{2d_i-k_i}\wr\gpS_{a_i}}^{\gpS_{a_i(2d_i-k_i)}}
           \left((D_{d_i,k_i}^c)^{\otimes a_i} \right)
           \right)
         \right)\\
      &\cong
          \bigotimes_{i=1}^l
          \left(
           V_n^{\otimes a_i(2d_i-k_i)}\otimes_{\K(\gpS_{2d_i-k_i}\wr\gpS_{a_i})}
           \left((D_{d_i,k_i}^c)^{\otimes a_i} \right)
          \right)\\
      &\cong
          \bigotimes_{i=1}^l
          \Sym^{a_i}\left(
            V_n^{\otimes(2d_i-k_i)}\otimes_{\K\gpS_{2d_i-k_i}}
            D_{d_i,k_i}^c
          \right)\\
      &\cong
          \bigotimes_{i=1}^l
          \Sym^{a_i}\left(B_{d_i,k_i}^c(n)\right).
    \end{split}
  \end{gather*}
  \end{proof}

 \subsection{Irreducible decomposition of $B_d(n)$ as $\GL(V_n)$-modules}\label{ss72}
  In this subsection, for simplicity, we write $V=V_n$, $B_{d,k}=B_{d,k}(n)$ and $B_{d,k}^c=B_{d,k}^c(n)$.

  Let $N$ be a nonnegative integer and $\lambda\vdash N$.
  Recall from Section \ref{ss515} that $S^\lambda$ denotes the Specht module, which is an irreducible representation of $\gpS_N$ corresponding to $\lambda$.
  Let $V_{\lambda}=\repS_{\lambda}V$ denote the image of $V$ under the Schur functor $\repS_{\lambda}$.
  Note that $V_{\lambda}$ is a simple $\GL(V)$-module if $n\geq r(\lambda)$ and that $V_{\lambda}=0$ if $n<r(\lambda)$, where $r(\lambda)$ is the number of rows of $\lambda$.

  We use the Littlewood-Richardson rule, plethysms and results by Bar-Natan \cite{BN} to compute the irreducible decompositions of the $\GL(V)$-modules $B_d$.

  \begin{proposition}[Bar-Natan \cite{BN}]\label{fbn}
   As $\gpS_{2d-k}$-modules, we have isomorphisms
   $$D^c_{1,0}\cong S^{(2)},$$
   $$D^c_{2,1}\cong S^{(1^3)},\quad D^c_{2,2}\cong S^{(2)},$$
   $$D^c_{3,2}\cong S^{(2^2)},\quad D^c_{3,3}\cong S^{(1^3)},\quad D^c_{3,4}\cong S^{(2)},$$
   $$D^c_{4,3}\cong S^{(3,1^2)},\quad D^c_{4,4}\cong S^{(4)}\oplus S^{(2^2)} ,\quad D^c_{4,5}\cong S^{(1^3)} ,\quad D^c_{4,6}\cong S^{(2)},$$
   $$D^c_{5,4}\cong S^{(4,2)}\oplus S^{(2^3)}\oplus S^{(3,1^3)},\quad D^c_{5,5}\cong ({S^{(3,1^2)}})^{\oplus 2},$$
   $$D^c_{5,6}\cong S^{(4)}\oplus ({S^{(2^2)}})^{\oplus 2},\quad D^c_{5,7}\cong ({S^{(1^3)}})^{\oplus 2},\quad D^c_{5,8}\cong ({S^{(2)}})^{\oplus 2}.$$
  \end{proposition}

  \begin{lemma}\label{lbn}
   We have the following isomorphisms of the $\GL(V)$-modules:
   $$B^c_{1,0}\cong V_{(2)},$$
   $$B^c_{2,1}\cong V_{(1^3)},\quad B^c_{2,2}\cong V_{(2)},$$
   $$B^c_{3,2}\cong V_{(2^2)},\quad B^c_{3,3}\cong V_{(1^3)},\quad B^c_{3,4}\cong V_{(2)},$$
   $$B^c_{4,3}\cong V_{(3,1^2)},\quad B^c_{4,4}\cong V_{(4)}\oplus V_{(2^2)} ,\quad B^c_{4,5}\cong V_{(1^3)},\quad B^c_{4,6}\cong V_{(2)},$$
   $$B^c_{5,4}\cong V_{(4,2)}\oplus V_{(2^3)}\oplus V_{(3,1^3)},\quad B^c_{5,5}\cong (V_{(3,1^2)})^{\oplus 2},$$
   $$B^c_{5,6}\cong V_{(4)}\oplus (V_{(2^2)})^{\oplus 2},\quad B^c_{5,7}\cong (V_{(1^3)})^{\oplus 2},\quad B^c_{5,8}\cong (V_{(2)})^{\oplus 2}.$$
  \end{lemma}
  \begin{proof}
   These follow from Proposition \ref{fbn}.
  \end{proof}

  \begin{proposition}\label{p531}
    For $d=3,4,5$, we have the following irreducible decompositions of the $\GL(V)$-modules $B_d$.
    \begin{enumerate}
     \item\label{eqB_3} We have $B_3=B_{3,0}\oplus\cdots\oplus B_{3,4},$
     where
     \begin{gather*}
       \begin{split}
        B_{3,0}&\cong V_{(6)}\oplus V_{(4,2)}\oplus V_{(2^3)},\\
        B_{3,1}&\cong V_{(3,1^2)}\oplus V_{(2,1^3)},\\
        B_{3,2}&\cong V_{(4)}\oplus V_{(3,1)}\oplus (V_{(2^2)})^{\oplus 2},\\
        B_{3,3}&=B^c_{3,3}\cong V_{(1^3)},\\
        B_{3,4}&=B^c_{3,4}\cong V_{(2)}.
       \end{split}
     \end{gather*}
     \item\label{eqB_4} We have $B_4=B_{4,0}\oplus\cdots\oplus B_{4,6},$
     where
     \begin{gather*}
       \begin{split}
         B_{4,0}&\cong V_{(8)}\oplus V_{(6,2)}\oplus V_{(4^2)}\oplus V_{(4,2^2)}\oplus V_{(2^4)},\\
         B_{4,1}&\cong V_{(5,1^2)}\oplus V_{(4,1^3)}\oplus V_{(3^2,1)} \oplus V_{(3,2,1^2)}\oplus V_{(2^2,1^3)},\\
         B_{4,2}&\cong V_{(6)}\oplus V_{(5,1)}\oplus(V_{(4,2)})^{\oplus 3}\oplus(V_{(3,2,1)})^{\oplus 2}\oplus(V_{(2^3)})^{\oplus 3}\oplus V_{(2,1^4)},\\
         B_{4,3}&\cong(V_{(3,1^2)})^{\oplus 3}\oplus(V_{(2,1^3)})^{\oplus 2},\\
         B_{4,4}&\cong(V_{(4)})^{\oplus 3}\oplus V_{(3,1)}\oplus(V_{(2^2)})^{\oplus 3},\\
         B_{4,5}&\cong V_{(1^3)},\\
         B_{4,6}&\cong V_{(2)}.
       \end{split}
     \end{gather*}
     \item \label{eqB_5} We have $B_5=B_{5,0}\oplus\cdots\oplus B_{5,8},$
     where
      \begin{gather*}
       \begin{split}
        B_{5,0}&\cong V_{(10)}\oplus V_{(8,2)}\oplus V_{(6,4)}\oplus V_{(6,2^2)}\oplus V_{(4^2,2)}\oplus V_{(4,2^3)} \oplus V_{(2^5)},\\
        B_{5,1}&\cong V_{(7,1^2)}\oplus V_{(6,1^3)}\oplus V_{(5,3,1)}\oplus V_{(5,2,1^2)}\oplus V_{(4,3,1^2)}\oplus V_{(4,2,1^3)}\\
        &\oplus V_{(3^3)}\oplus V_{(3^2,2,1)}\oplus V_{(3,2^2,1^2)}\oplus V_{(2^3,1^3)},\\
        B_{5,2}&\cong V_{(8)}\oplus V_{(7,1)}\oplus(V_{(6,2)})^{\oplus 3}\oplus V_{(5,3)}\oplus(V_{(5,2,1)})^{\oplus 2}\oplus(V_{(4^2)})^{\oplus 2}\\
        &\oplus(V_{(4,3,1)})^{\oplus 2}\oplus(V_{(4,2^2)})^{\oplus 5}\oplus V_{(4,1^4)}\oplus V_{(3^2,1^2)}\oplus(V_{(3,2^2,1)})^{\oplus 3}\\
        &\oplus V_{(3,2,1^3)}\oplus V_{(3,1^5)}\oplus(V_{(2^4)})^{\oplus 3}\oplus V_{(2^2,1^4)},\\
        B_{5,3}&\cong (V_{(5,1^2)})^{\oplus 3}\oplus(V_{(4,2,1)})^{\oplus 2}\oplus(V_{(4,1^3)})^{\oplus 4}\oplus(V_{(3^2,1)})^{\oplus 4}\oplus(V_{(3,2,1^2)})^{\oplus 5}\\
        &\oplus V_{(3,1^4)}\oplus(V_{(2^2,1^3)})^{\oplus 3},\\
        B_{5,4}&\cong (V_{(6)})^{\oplus 3}\oplus(V_{(5,1)})^{\oplus 3}\oplus(V_{(4,2)})^{\oplus 8}\oplus(V_{(3,2,1)})^{\oplus 4}\oplus V_{(3,1^3)}\oplus(V_{(2^3)})^{\oplus 6}\\
        &\oplus V_{(2^2,1^2)}\oplus V_{(2,1^4)}\oplus V_{(1^6)},\\
        B_{5,5}&\cong (V_{(3,1^2)})^{\oplus 5}\oplus(V_{(2,1^3)})^{\oplus 3},\\
        B_{5,6}&\cong (V_{(4)})^{\oplus 3}\oplus(V_{(3,1)})^{\oplus 2}\oplus(V_{(2^2)})^{\oplus 4},\\
        B_{5,7}&\cong (V_{(1^3)})^{\oplus 2},\\
        B_{5,8}&\cong (V_{(2)})^{\oplus 2}.
       \end{split}
      \end{gather*}
    \end{enumerate}
  \end{proposition}
  \begin{proof}
   By using Theorem \ref{th521}, Lemma \ref{lbn} and plethysm, we have
   $$
     B_{3,0}\cong\Sym^3(B^c_{1,0})\cong\repS_{(3)}(\repS_{(2)} V) \cong V_{(6)}\oplus V_{(4,2)}\oplus V_{(2^3)}.
   $$
   By using Theorem \ref{th521}, Lemma \ref{lbn}, and the Littlewood--Richardson rule, we have
   $$
     B_{3,1}\cong B^c_{2,1}\otimes B^c_{1,0}\cong V_{(1^3)}\otimes V_{(2)}\cong V_{(3,1^2)}\oplus V_{(2,1^3)},
   $$
   and
   $$
     B_{3,2}\cong B^c_{3,2}\oplus(B^c_{2,2}\otimes B^c_{1,0})\cong V_{(2^2)}\oplus (V_{(4)}\oplus V_{(3,1)} \oplus V_{(2^2)}).
   $$
   The other isomorphisms of (\ref{eqB_3}) follow from Lemma \ref{lbn}.

   The irreducible decompositions (\ref{eqB_4}) and (\ref{eqB_5}) follow in a similar way.
  \end{proof}

  We need the irreducible decompositions of $B_{d,0}$ and $B_{d,1}$ to study the $\Aut(F_n)$-module structure of $A_d(n)$.
  For $\lambda=(\lambda_1,\cdots,\lambda_r)\vdash N$, let $2\lambda$ denote the partition $(2\lambda_1,\cdots,2\lambda_r)$ of $2N$.
  \begin{proposition}\label{pdecompgen}
   For any $d\geq 0$, we have
   $$
     B_{d,0}\cong \bigoplus_{\lambda\vdash d} V_{2\lambda}.
   $$
   For any $d\geq 2$, we have
   $$
     B_{d,1}\cong \bigoplus_{\lambda\vdash 2d-1 \text{ with exactly $3$ odd parts}}V_{\lambda}.
   $$
  \end{proposition}
  \begin{proof}
   By Theorem 5.4.23 in \cite{JamesKerber}, we have
   $$\repS_{(d)}(\repS_{(2)} V)\cong\bigoplus_{\lambda\vdash d} V_{2\lambda}.$$
   Therefore, by Theorem \ref{th521} and Lemma \ref{lbn}, we have
   $$B_{d,0}\cong\Sym^{d}(B^c_{1,0})\cong\repS_{(d)}(\repS_{(2)} V)\cong \bigoplus_{\lambda\vdash d} V_{2\lambda}.$$
   By Theorem \ref{th521}, Lemma \ref{lbn}, plethysm and the Littlewood--Richardson rule, we have
   $$B_{d,1}\cong B^c_{2,1}\otimes \Sym^{d-2}(B^c_{1,0})\cong V_{(1^3)}\otimes \bigoplus_{\mu\vdash d-2} V_{2\mu}
   \cong \bigoplus_{\lambda\vdash 2d-1\text{ with exactly $3$ odd parts}}V_{\lambda}.$$
  \end{proof}

 \subsection{Surjectivity of the bracket map $[\cdot,\cdot]:B_{d,k}(n)\otimes \opegr^1(\IA(n))\rightarrow B_{d,k+1}(n)$}\label{ss73}
  Here, we show that the bracket map $[\cdot,\cdot]: B_{d,k}(n)\otimes \opegr^1(\IA(n))\rightarrow B_{d,k+1}(n)$ is surjective for $n\geq 2d$.
  Since we have abelian group isomorphisms \eqref{bracketeq}, the bracket map of $\opegr^1(\IA(n))$ coincides with that of $\opegr^1(\jfE_{\ast}(n))$.
  Thus, we can compute the bracket map by using the contraction map $c$ defined in Section \ref{s5}.

  Define $K_{i,j}, K_{i,j,k}\in \IA(n)$ by
  $$K_{i,j}(x_i)=x_jx_ix_j^{-1},\quad K_{i,j}(x_l)=x_l \quad(l\neq i),$$
  \begin{equation}\label{Kijk}
    K_{i,j,k}(x_i)=x_i[x_j,x_k],\quad K_{i,j,k}(x_l)=x_l \quad(l\neq i).
  \end{equation}

  \begin{proposition}\label{pbracketsurj}
    For $n\geq 2d-k$, the bracket map
    $$[\cdot,\cdot]:B_{d,k}(n)\otimes \opegr^1(\IA(n))\rightarrow B_{d,k+1}(n)$$
    is surjective.
  \end{proposition}
  \begin{proof}
    Any element of $B_{d,k+1}(n)$ is a linear sum of $u=\centre{
\begingroup%
  \makeatletter%
  \providecommand\color[2][]{%
    \errmessage{(Inkscape) Color is used for the text in Inkscape, but the package 'color.sty' is not loaded}%
    \renewcommand\color[2][]{}%
  }%
  \providecommand\transparent[1]{%
    \errmessage{(Inkscape) Transparency is used (non-zero) for the text in Inkscape, but the package 'transparent.sty' is not loaded}%
    \renewcommand\transparent[1]{}%
  }%
  \providecommand\rotatebox[2]{#2}%
  \newcommand*\fsize{\dimexpr\f@size pt\relax}%
  \newcommand*\lineheight[1]{\fontsize{\fsize}{#1\fsize}\selectfont}%
  \ifx\svgwidth\undefined%
    \setlength{\unitlength}{53.28929028bp}%
    \ifx\svgscale\undefined%
      \relax%
    \else%
      \setlength{\unitlength}{\unitlength * \real{\svgscale}}%
    \fi%
  \else%
    \setlength{\unitlength}{\svgwidth}%
  \fi%
  \global\let\svgwidth\undefined%
  \global\let\svgscale\undefined%
  \makeatother%
  \begin{picture}(1,0.66319365)%
    \lineheight{1}%
    \setlength\tabcolsep{0pt}%
    \put(0,0){\includegraphics[width=\unitlength,page=1]{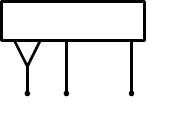}}%
    \put(0.32205571,0.50298645){\makebox(0,0)[lt]{\lineheight{1.45000005}\smash{\begin{tabular}[t]{l}$D$\end{tabular}}}}%
    \put(0.04450978,0.01370462){\makebox(0,0)[lt]{\lineheight{1.45000005}\smash{\begin{tabular}[t]{l}$v_{i_1}$\end{tabular}}}}%
    \put(0.28873857,0.01194441){\makebox(0,0)[lt]{\lineheight{1.45000005}\smash{\begin{tabular}[t]{l}$v_{i_2}$\end{tabular}}}}%
    \put(0.65166527,0.00718367){\makebox(0,0)[lt]{\lineheight{1.45000005}\smash{\begin{tabular}[t]{l}$v_{i_{2d-k-1}}$\end{tabular}}}}%
    \put(0,0){\includegraphics[width=\unitlength,page=2]{utriv.pdf}}%
  \end{picture}%
\endgroup%
}$, where $i_1,\cdots,i_{2d-k-1}\in[n]$.
    Since $n\geq 2d-k$, we can take $\ti{u}=\centre{
\begingroup%
  \makeatletter%
  \providecommand\color[2][]{%
    \errmessage{(Inkscape) Color is used for the text in Inkscape, but the package 'color.sty' is not loaded}%
    \renewcommand\color[2][]{}%
  }%
  \providecommand\transparent[1]{%
    \errmessage{(Inkscape) Transparency is used (non-zero) for the text in Inkscape, but the package 'transparent.sty' is not loaded}%
    \renewcommand\transparent[1]{}%
  }%
  \providecommand\rotatebox[2]{#2}%
  \newcommand*\fsize{\dimexpr\f@size pt\relax}%
  \newcommand*\lineheight[1]{\fontsize{\fsize}{#1\fsize}\selectfont}%
  \ifx\svgwidth\undefined%
    \setlength{\unitlength}{53.82668171bp}%
    \ifx\svgscale\undefined%
      \relax%
    \else%
      \setlength{\unitlength}{\unitlength * \real{\svgscale}}%
    \fi%
  \else%
    \setlength{\unitlength}{\svgwidth}%
  \fi%
  \global\let\svgwidth\undefined%
  \global\let\svgscale\undefined%
  \makeatother%
  \begin{picture}(1,0.65747109)%
    \lineheight{1}%
    \setlength\tabcolsep{0pt}%
    \put(0,0){\includegraphics[width=\unitlength,page=1]{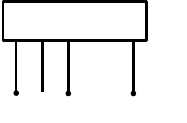}}%
    \put(0.32882412,0.49886335){\makebox(0,0)[lt]{\lineheight{1.45000005}\smash{\begin{tabular}[t]{l}$D$\end{tabular}}}}%
    \put(-0.00281212,0.00814845){\makebox(0,0)[lt]{\lineheight{1.45000005}\smash{\begin{tabular}[t]{l}$v_{i}$\end{tabular}}}}%
    \put(0.29162768,0.00711195){\makebox(0,0)[lt]{\lineheight{1.45000005}\smash{\begin{tabular}[t]{l}$v_{i_2}$\end{tabular}}}}%
    \put(0.65514295,0.00801054){\makebox(0,0)[lt]{\lineheight{1.45000005}\smash{\begin{tabular}[t]{l}$v_{i_{2d-k-1}}$\end{tabular}}}}%
    \put(0,0){\includegraphics[width=\unitlength,page=2]{ucontra.pdf}}%
    \put(0.13701707,0.00839715){\makebox(0,0)[lt]{\lineheight{1.45000005}\smash{\begin{tabular}[t]{l}$v_{j}$\end{tabular}}}}%
  \end{picture}%
\endgroup%
} \in B_{d,k}(n)$, where $i,j\in [n]\setminus\{i_2,\cdots,i_{2d-k-1}\}$ are distinct.
    We have $[\ti{u},K_{i_1,j,i}]=u$ and therefore, the bracket map is surjective.
  \end{proof}

  As in Section \ref{ss52}, for $\lambda\vdash 2d-k$, let $B_{d,k}(n)_{\lambda}$ denote the isotypic component of $\GL(n;\Z)$-module $B_{d,k}(n)$ corresponding to $\lambda$.

  In Proposition \ref{pdecompgen}, we computed a decomposition of $B_{d,0}(n)$.
  Since the Young diagram of $(2d)$ does not contain that of $(1^2)$, by Remark \ref{rem52}, we have the following corollary.
  \begin{corollary}\label{c82}
    The restriction of the bracket map
    $$[\cdot,\cdot]:\bigoplus_{\lambda\vdash d,\lambda\neq (d)} B_{d,0}(n)_{2\lambda}\otimes \opegr^1(\IA(n))\rightarrow B_{d,1}(n)$$
    is surjective for $n\geq 2d$.
  \end{corollary}

  Lastly, we consider the condition for $\lambda\vdash 2d-k$ that the isotypic component $B_{d,k}(n)_{\lambda}$ of $B_{d,k}(n)$ does not vanish.
  Let $o(\lambda)$ be the number of odd parts of $\lambda$.
  We have $$o(\lambda)\equiv 2d-k\equiv k \quad(\bmod\: 2).$$
  In Proposition \ref{pdecompgen}, we observed that $o(\lambda)=0\; (k=0)$ and $o(\lambda)=3 \;(k=1)$.
  Moreover, by Proposition \ref{pbracketsurj} and Remark \ref{rem52}, we have $o(\lambda)\leq 3k$ if $B_{d,k}(n)_{\lambda}\neq 0$.


\section{The $\Aut(F_n)$-module structure of $A_d(n)$}\label{s8}
 In this section, we study the $\Aut(F_n)$-module structure of $A_d(n)$.
 We have $A_0(n)=\K$ for any $n\geq 0$ and we studied the cases where $d=1, 2$ in \cite{Mai1}.
 Note that we have $A_d(0)=0$ for $d\geq 1$. Thus, we have only to consider $n\geq 1$.
 Here, we construct a direct decomposition of $A_d(n)$ as $\Aut(F_n)$-modules for any $d\geq 3,n\geq 1$, which is indecomposable for $n\geq 2d$. Moreover, we study the degree $3$ case in detail.
 \subsection{A direct decomposition of $A_d(n)$}\label{ss81}
  Here, we give a direct decomposition of the $\Aut(F_n)$-module $A_d(n)$.

  Let $c=\centre{
\begingroup%
  \makeatletter%
  \providecommand\color[2][]{%
    \errmessage{(Inkscape) Color is used for the text in Inkscape, but the package 'color.sty' is not loaded}%
    \renewcommand\color[2][]{}%
  }%
  \providecommand\transparent[1]{%
    \errmessage{(Inkscape) Transparency is used (non-zero) for the text in Inkscape, but the package 'transparent.sty' is not loaded}%
    \renewcommand\transparent[1]{}%
  }%
  \providecommand\rotatebox[2]{#2}%
  \newcommand*\fsize{\dimexpr\f@size pt\relax}%
  \newcommand*\lineheight[1]{\fontsize{\fsize}{#1\fsize}\selectfont}%
  \ifx\svgwidth\undefined%
    \setlength{\unitlength}{23.975989bp}%
    \ifx\svgscale\undefined%
      \relax%
    \else%
      \setlength{\unitlength}{\unitlength * \real{\svgscale}}%
    \fi%
  \else%
    \setlength{\unitlength}{\svgwidth}%
  \fi%
  \global\let\svgwidth\undefined%
  \global\let\svgscale\undefined%
  \makeatother%
  \begin{picture}(1,1.08854031)%
    \lineheight{1}%
    \setlength\tabcolsep{0pt}%
    \put(0,0){\includegraphics[width=\unitlength,page=1]{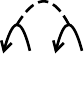}}%
    \put(0.14029882,0.01470385){\makebox(0,0)[lt]{\lineheight{1.45000005}\smash{\begin{tabular}[t]{l}$1$\end{tabular}}}}%
    \put(0.76480238,0.01778309){\makebox(0,0)[lt]{\lineheight{1.45000005}\smash{\begin{tabular}[t]{l}$2$\end{tabular}}}}%
  \end{picture}%
\endgroup%
}\in A_1(2)=\A_1(0,2)$ and depict it as
  $\centre{
\begingroup%
  \makeatletter%
  \providecommand\color[2][]{%
    \errmessage{(Inkscape) Color is used for the text in Inkscape, but the package 'color.sty' is not loaded}%
    \renewcommand\color[2][]{}%
  }%
  \providecommand\transparent[1]{%
    \errmessage{(Inkscape) Transparency is used (non-zero) for the text in Inkscape, but the package 'transparent.sty' is not loaded}%
    \renewcommand\transparent[1]{}%
  }%
  \providecommand\rotatebox[2]{#2}%
  \newcommand*\fsize{\dimexpr\f@size pt\relax}%
  \newcommand*\lineheight[1]{\fontsize{\fsize}{#1\fsize}\selectfont}%
  \ifx\svgwidth\undefined%
    \setlength{\unitlength}{15.72738029bp}%
    \ifx\svgscale\undefined%
      \relax%
    \else%
      \setlength{\unitlength}{\unitlength * \real{\svgscale}}%
    \fi%
  \else%
    \setlength{\unitlength}{\svgwidth}%
  \fi%
  \global\let\svgwidth\undefined%
  \global\let\svgscale\undefined%
  \makeatother%
  \begin{picture}(1,1.05203066)%
    \lineheight{1}%
    \setlength\tabcolsep{0pt}%
    \put(0,0){\includegraphics[width=\unitlength,page=1]{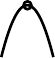}}%
  \end{picture}%
\endgroup%
}\:.$
  Here, we use the same graphical notation of morphisms $\mu,\eta,\Delta,\epsilon,S$ in the category $\A$ as in the category $\A^{L}$.
  As in Section \ref{ss425}, we can define the iterated multiplications $\mu^{[q]}\in \A(q,1)$ for $q\geq 0$.
  For $m\geq 0$, there is a group homomorphism
  $$
    \gpS_{m}\rightarrow \A(m,m),\quad \sigma \mapsto P_{\sigma},
  $$
  where $P_{\sigma}$ is the symmetry in $\A$ corresponding to $\sigma$.
  Set
  $$\centre{
\begingroup%
  \makeatletter%
  \providecommand\color[2][]{%
    \errmessage{(Inkscape) Color is used for the text in Inkscape, but the package 'color.sty' is not loaded}%
    \renewcommand\color[2][]{}%
  }%
  \providecommand\transparent[1]{%
    \errmessage{(Inkscape) Transparency is used (non-zero) for the text in Inkscape, but the package 'transparent.sty' is not loaded}%
    \renewcommand\transparent[1]{}%
  }%
  \providecommand\rotatebox[2]{#2}%
  \newcommand*\fsize{\dimexpr\f@size pt\relax}%
  \newcommand*\lineheight[1]{\fontsize{\fsize}{#1\fsize}\selectfont}%
  \ifx\svgwidth\undefined%
    \setlength{\unitlength}{30.75117737bp}%
    \ifx\svgscale\undefined%
      \relax%
    \else%
      \setlength{\unitlength}{\unitlength * \real{\svgscale}}%
    \fi%
  \else%
    \setlength{\unitlength}{\svgwidth}%
  \fi%
  \global\let\svgwidth\undefined%
  \global\let\svgscale\undefined%
  \makeatother%
  \begin{picture}(1,1.09751899)%
    \lineheight{1}%
    \setlength\tabcolsep{0pt}%
    \put(0,0){\includegraphics[width=\unitlength,page=1]{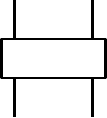}}%
    \put(0.07929683,0.50062537){\makebox(0,0)[lt]{\lineheight{1.45000005}\smash{\begin{tabular}[t]{l}$sym_m$\end{tabular}}}}%
    \put(0,0){\includegraphics[width=\unitlength,page=2]{symn.pdf}}%
  \end{picture}%
\endgroup%
}:= \sum_{\sigma\in\gpS_m}P_{\sigma},\quad
  \centre{
\begingroup%
  \makeatletter%
  \providecommand\color[2][]{%
    \errmessage{(Inkscape) Color is used for the text in Inkscape, but the package 'color.sty' is not loaded}%
    \renewcommand\color[2][]{}%
  }%
  \providecommand\transparent[1]{%
    \errmessage{(Inkscape) Transparency is used (non-zero) for the text in Inkscape, but the package 'transparent.sty' is not loaded}%
    \renewcommand\transparent[1]{}%
  }%
  \providecommand\rotatebox[2]{#2}%
  \newcommand*\fsize{\dimexpr\f@size pt\relax}%
  \newcommand*\lineheight[1]{\fontsize{\fsize}{#1\fsize}\selectfont}%
  \ifx\svgwidth\undefined%
    \setlength{\unitlength}{30.75117737bp}%
    \ifx\svgscale\undefined%
      \relax%
    \else%
      \setlength{\unitlength}{\unitlength * \real{\svgscale}}%
    \fi%
  \else%
    \setlength{\unitlength}{\svgwidth}%
  \fi%
  \global\let\svgwidth\undefined%
  \global\let\svgscale\undefined%
  \makeatother%
  \begin{picture}(1,1.09751899)%
    \lineheight{1}%
    \setlength\tabcolsep{0pt}%
    \put(0,0){\includegraphics[width=\unitlength,page=1]{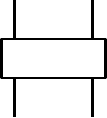}}%
    \put(0.22078711,0.48330002){\makebox(0,0)[lt]{\lineheight{1.45000005}\smash{\begin{tabular}[t]{l}$alt_m$\end{tabular}}}}%
    \put(0,0){\includegraphics[width=\unitlength,page=2]{altn.pdf}}%
  \end{picture}%
\endgroup%
}:= \sum_{\sigma\in\gpS_m}\sgn(\sigma)P_{\sigma}\in\A(m,m).$$

  By Habiro--Massuyeau \cite[Lemma 5.16]{HM_k}, every element of $A_d(n)$ is a linear combination of morphisms of the form
  $$
    (\mu^{[q_1]}\otimes\cdots\otimes\mu^{[q_n]})\circ P_{\sigma}\circ c^{\otimes d}=\centre{
\begingroup%
  \makeatletter%
  \providecommand\color[2][]{%
    \errmessage{(Inkscape) Color is used for the text in Inkscape, but the package 'color.sty' is not loaded}%
    \renewcommand\color[2][]{}%
  }%
  \providecommand\transparent[1]{%
    \errmessage{(Inkscape) Transparency is used (non-zero) for the text in Inkscape, but the package 'transparent.sty' is not loaded}%
    \renewcommand\transparent[1]{}%
  }%
  \providecommand\rotatebox[2]{#2}%
  \newcommand*\fsize{\dimexpr\f@size pt\relax}%
  \newcommand*\lineheight[1]{\fontsize{\fsize}{#1\fsize}\selectfont}%
  \ifx\svgwidth\undefined%
    \setlength{\unitlength}{98.18362578bp}%
    \ifx\svgscale\undefined%
      \relax%
    \else%
      \setlength{\unitlength}{\unitlength * \real{\svgscale}}%
    \fi%
  \else%
    \setlength{\unitlength}{\svgwidth}%
  \fi%
  \global\let\svgwidth\undefined%
  \global\let\svgscale\undefined%
  \makeatother%
  \begin{picture}(1,0.34928896)%
    \lineheight{1}%
    \setlength\tabcolsep{0pt}%
    \put(0,0){\includegraphics[width=\unitlength,page=1]{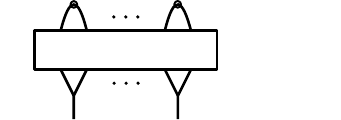}}%
    \put(0.27584011,0.15897959){\makebox(0,0)[lt]{\lineheight{1.45000005}\smash{\begin{tabular}[t]{l}\small$P_{\sigma}$\end{tabular}}}}%
    \put(-0.00139248,0.03357592){\makebox(0,0)[lt]{\lineheight{1.45000005}\smash{\begin{tabular}[t]{l}\small$\mu^{[q_1]}$\end{tabular}}}}%
    \put(0.54143674,0.03628154){\makebox(0,0)[lt]{\lineheight{1.45000005}\smash{\begin{tabular}[t]{l}\small$\mu^{[q_n]}$\end{tabular}}}}%
  \end{picture}%
\endgroup%
}
  $$
  for $\sigma\in\gpS_{2d}$ and $q_1,\cdots, q_n\geq 0$ such that $q_1+\cdots+q_n=2d$. The following lemma easily follows.
  \begin{lemma}\label{l1021}
   For $n\geq 0$, we have
   $$
     A_d(n)=\Span_\K\{A_d(f)(c^{\otimes d})\mid f\in\F^\op(2d,n)\}.
   $$
  \end{lemma}

  For $X\in A_d(m)$, let
  $$A_d X:\F^{\op}\rightarrow \fVect$$
  denote the subfunctor of $A_d$ generated by $X$.
  That is, for any $n\in\N$, $A_d X(n)$ is the $\Aut(F_n)$-submodule of $A_d(n)$ defined by
  $$A_d X (n):=\Span_\K\{A_d(f)(X)\mid f\in\F^\op(m,n)\}.$$

  Set
  $$P=\centre{},\quad Q=\centre{}\in A_d(2d).$$
  Note that we have $A_1 Q=0$.
  \begin{theorem}\label{decompositionofAd}
   We have
   \begin{equation}\label{decomp}
     A_d(n)=A_d P(n)\oplus A_d Q(n).
   \end{equation}
  \end{theorem}

  \begin{proof}
   By Lemma \ref{l1021}, any element of $A_d(n)$ is a linear sum of $A_d(f)(c^{\otimes d})$ for $f\in \F^{\op}(2d,n)$.
   Define an $\Aut(F_n)$-module map
   $$e_n:A_d(n)\rightarrow A_d(n)$$
   by $e_n(A_d(f)(c^{\otimes d}))=\frac{1}{(2d)!} A_d(f)(P)$ for $f\in \F^{\op}(2d,n)$.
   This is well defined because the $4$T relation is sent to $0$.
   Since $A_d P$ is generated by $P$, we have $\im(e_n)=A_d P(n)$.

   Since we have $e_n(A_d(f)(P))=A_d(f)(P)$ for any $f\in \F^{\op}(2d,n)$, the $\Aut(F_n)$-endomorphism $e_n$ is an idempotent in $\End(A_d(n))$, where we consider $A_d(n)$ as a right $\Aut(F_n)$-module.
   Therefore, we have
   $$A_d(n)= \im(e_n)\oplus \ker(e_n),\quad \ker(e_n)=\im(1-e_n).$$

   Since $\centre{
\begingroup%
  \makeatletter%
  \providecommand\color[2][]{%
    \errmessage{(Inkscape) Color is used for the text in Inkscape, but the package 'color.sty' is not loaded}%
    \renewcommand\color[2][]{}%
  }%
  \providecommand\transparent[1]{%
    \errmessage{(Inkscape) Transparency is used (non-zero) for the text in Inkscape, but the package 'transparent.sty' is not loaded}%
    \renewcommand\transparent[1]{}%
  }%
  \providecommand\rotatebox[2]{#2}%
  \newcommand*\fsize{\dimexpr\f@size pt\relax}%
  \newcommand*\lineheight[1]{\fontsize{\fsize}{#1\fsize}\selectfont}%
  \ifx\svgwidth\undefined%
    \setlength{\unitlength}{53.25118286bp}%
    \ifx\svgscale\undefined%
      \relax%
    \else%
      \setlength{\unitlength}{\unitlength * \real{\svgscale}}%
    \fi%
  \else%
    \setlength{\unitlength}{\svgwidth}%
  \fi%
  \global\let\svgwidth\undefined%
  \global\let\svgscale\undefined%
  \makeatother%
  \begin{picture}(1,0.74646229)%
    \lineheight{1}%
    \setlength\tabcolsep{0pt}%
    \put(0,0){\includegraphics[width=\unitlength,page=1]{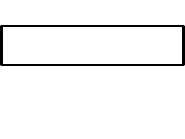}}%
    \put(0.21831612,0.46477843){\makebox(0,0)[lt]{\lineheight{1.45000005}\smash{\begin{tabular}[t]{l}$\sym_{2d}$\end{tabular}}}}%
    \put(0,0){\includegraphics[width=\unitlength,page=2]{symalt2.pdf}}%
    \put(0.08244809,0.1611395){\makebox(0,0)[lt]{\lineheight{1.45000005}\smash{\begin{tabular}[t]{l}$\alt_{2}$\end{tabular}}}}%
    \put(0,0){\includegraphics[width=\unitlength,page=3]{symalt2.pdf}}%
  \end{picture}%
\endgroup%
}=0$, we have $A_d Q(n)\subset \ker(e_n)$.
   Finally, we need to check that $\im(1-e_n)\subset A_d Q(n)$.
   Since we have for $f\in \F^{\op}(2d,n)$,
   \begin{gather*}
     \begin{split}
       (1-e_n)(A_d(f)(c^{\otimes d}))&=
        A_d(f)(c^{\otimes d})-\frac{1}{(2d)!} A_d(f)(P)\\
        &=\frac{1}{(2d)!}\sum_{\sigma\in\gpS_{2d}}A_d(f)(c^{\otimes d}- \sigma c^{\otimes d}),
     \end{split}
   \end{gather*}
   we need to show that for any $\sigma\in \gpS_{2d}$, there exists $\tau\in\K\gpS_{2d}$ such that
   \begin{gather}\label{c^d}
     c^{\otimes d}- \sigma c^{\otimes d}=\tau Q\in A_d Q(2d).
   \end{gather}

   It suffices to show the existence of $\tau$ satisfying \eqref{c^d} when $\sigma$ is an adjacent transposition, because any permutation is generated by adjacent transpositions and we have such $\tau$ by inductively using
   $$c^{\otimes d}- \sigma\rho c^{\otimes d}=c^{\otimes d}-\sigma c^{\otimes d}+ \sigma(c^{\otimes d}-\rho c^{\otimes d}).$$

   If $\sigma$ is an adjacent transposition $(2i,2i+1)$ for $i\in[n-1]$, then we set
   $$\tau=
   \begin{pmatrix}
     1&2&3&4&5&\cdots&2d\\
     2i-1&2i+2&2i+1&2i&1&\cdots \widehat{2i-1}\cdots \widehat{2i+2}\cdots&2d\\
   \end{pmatrix}.$$
   If $\sigma$ is an adjacent transposition $(2i-1,2i)$ for $i\in[n]$,
   then we set $\tau=0$.
   The proof is complete.
  \end{proof}

  \begin{lemma}\label{l81}
    The $\Aut(F_n)$-module $A_d P(n)$ is irreducible and thus indecomposable.
  \end{lemma}
  \begin{proof}
    Since $\centre{
\begingroup%
  \makeatletter%
  \providecommand\color[2][]{%
    \errmessage{(Inkscape) Color is used for the text in Inkscape, but the package 'color.sty' is not loaded}%
    \renewcommand\color[2][]{}%
  }%
  \providecommand\transparent[1]{%
    \errmessage{(Inkscape) Transparency is used (non-zero) for the text in Inkscape, but the package 'transparent.sty' is not loaded}%
    \renewcommand\transparent[1]{}%
  }%
  \providecommand\rotatebox[2]{#2}%
  \newcommand*\fsize{\dimexpr\f@size pt\relax}%
  \newcommand*\lineheight[1]{\fontsize{\fsize}{#1\fsize}\selectfont}%
  \ifx\svgwidth\undefined%
    \setlength{\unitlength}{53.25118286bp}%
    \ifx\svgscale\undefined%
      \relax%
    \else%
      \setlength{\unitlength}{\unitlength * \real{\svgscale}}%
    \fi%
  \else%
    \setlength{\unitlength}{\svgwidth}%
  \fi%
  \global\let\svgwidth\undefined%
  \global\let\svgscale\undefined%
  \makeatother%
  \begin{picture}(1,0.72842959)%
    \lineheight{1}%
    \setlength\tabcolsep{0pt}%
    \put(0,0){\includegraphics[width=\unitlength,page=1]{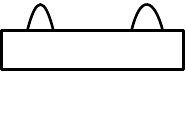}}%
    \put(0.21831612,0.42252585){\makebox(0,0)[lt]{\lineheight{1.45000005}\smash{\begin{tabular}[t]{l}$\sym_{2d}$\end{tabular}}}}%
    \put(0,0){\includegraphics[width=\unitlength,page=2]{symtriv.pdf}}%
  \end{picture}%
\endgroup%
}=0$, we have $\theta_{d,n}(A_d P(n))=B_{d,0}(n)_{(2d)}$ by the PBW map. Therefore, $A_d P(n)$ is an  irreducible $\Aut(F_n)$-module.
  \end{proof}

  For $\lambda\vdash d$, set $Q_{\lambda}=\centre{
\begingroup%
  \makeatletter%
  \providecommand\color[2][]{%
    \errmessage{(Inkscape) Color is used for the text in Inkscape, but the package 'color.sty' is not loaded}%
    \renewcommand\color[2][]{}%
  }%
  \providecommand\transparent[1]{%
    \errmessage{(Inkscape) Transparency is used (non-zero) for the text in Inkscape, but the package 'transparent.sty' is not loaded}%
    \renewcommand\transparent[1]{}%
  }%
  \providecommand\rotatebox[2]{#2}%
  \newcommand*\fsize{\dimexpr\f@size pt\relax}%
  \newcommand*\lineheight[1]{\fontsize{\fsize}{#1\fsize}\selectfont}%
  \ifx\svgwidth\undefined%
    \setlength{\unitlength}{45.98981849bp}%
    \ifx\svgscale\undefined%
      \relax%
    \else%
      \setlength{\unitlength}{\unitlength * \real{\svgscale}}%
    \fi%
  \else%
    \setlength{\unitlength}{\svgwidth}%
  \fi%
  \global\let\svgwidth\undefined%
  \global\let\svgscale\undefined%
  \makeatother%
  \begin{picture}(1,0.59892514)%
    \lineheight{1}%
    \setlength\tabcolsep{0pt}%
    \put(0,0){\includegraphics[width=\unitlength,page=1]{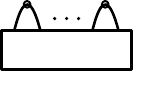}}%
    \put(0.26767131,0.233){\makebox(0,0)[lt]{\lineheight{1.45000005}\smash{\begin{tabular}[t]{l}$c_{2\lambda}$\\\end{tabular}}}}%
    \put(0,0){\includegraphics[width=\unitlength,page=2]{Q_lambda.pdf}}%
  \end{picture}%
\endgroup%
}$,
  where $c_{2\lambda}\in\K\gpS_{2d}$ is the Young symmetrizer.
  Note that we have $Q_{(d)}=P$.

  \begin{lemma}\label{pQlambda}
   For $\lambda\vdash d, \lambda\neq (d)$, we have $Q_{\lambda}\in A_d Q(2d)$.
  \end{lemma}
  \begin{proof}
   For $\lambda=(\lambda_1,\cdots,\lambda_r)\neq (d)$, we have $r\geq 2$.
   By expanding $a_{\lambda}$ and $b_{\lambda}$ except for the first column, we can write $Q_{\lambda}$ as a linear sum of
   $$\scalebox{0.9}{$\centre{
\begingroup%
  \makeatletter%
  \providecommand\color[2][]{%
    \errmessage{(Inkscape) Color is used for the text in Inkscape, but the package 'color.sty' is not loaded}%
    \renewcommand\color[2][]{}%
  }%
  \providecommand\transparent[1]{%
    \errmessage{(Inkscape) Transparency is used (non-zero) for the text in Inkscape, but the package 'transparent.sty' is not loaded}%
    \renewcommand\transparent[1]{}%
  }%
  \providecommand\rotatebox[2]{#2}%
  \newcommand*\fsize{\dimexpr\f@size pt\relax}%
  \newcommand*\lineheight[1]{\fontsize{\fsize}{#1\fsize}\selectfont}%
  \ifx\svgwidth\undefined%
    \setlength{\unitlength}{158.25118165bp}%
    \ifx\svgscale\undefined%
      \relax%
    \else%
      \setlength{\unitlength}{\unitlength * \real{\svgscale}}%
    \fi%
  \else%
    \setlength{\unitlength}{\svgwidth}%
  \fi%
  \global\let\svgwidth\undefined%
  \global\let\svgscale\undefined%
  \makeatother%
  \begin{picture}(1,0.52950287)%
    \lineheight{1}%
    \setlength\tabcolsep{0pt}%
    \put(0,0){\includegraphics[width=\unitlength,page=1]{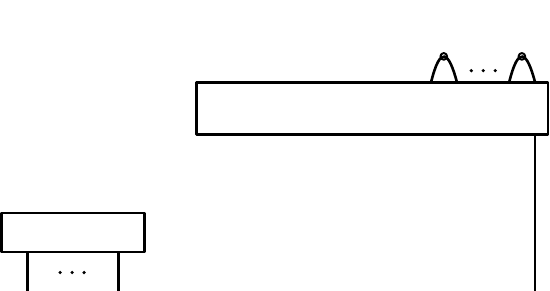}}%
    \put(0.63302828,0.31461035){\makebox(0,0)[lt]{\lineheight{1.45000005}\smash{\begin{tabular}[t]{l}$\sigma$\end{tabular}}}}%
    \put(0.06866326,0.08823317){\makebox(0,0)[lt]{\lineheight{1.45000005}\smash{\begin{tabular}[t]{l}$\alt_r$\end{tabular}}}}%
    \put(0,0){\includegraphics[width=\unitlength,page=2]{Q_lambda1.pdf}}%
  \end{picture}%
\endgroup%
}$}=\frac{1}{2}\centre{
\begingroup%
  \makeatletter%
  \providecommand\color[2][]{%
    \errmessage{(Inkscape) Color is used for the text in Inkscape, but the package 'color.sty' is not loaded}%
    \renewcommand\color[2][]{}%
  }%
  \providecommand\transparent[1]{%
    \errmessage{(Inkscape) Transparency is used (non-zero) for the text in Inkscape, but the package 'transparent.sty' is not loaded}%
    \renewcommand\transparent[1]{}%
  }%
  \providecommand\rotatebox[2]{#2}%
  \newcommand*\fsize{\dimexpr\f@size pt\relax}%
  \newcommand*\lineheight[1]{\fontsize{\fsize}{#1\fsize}\selectfont}%
  \ifx\svgwidth\undefined%
    \setlength{\unitlength}{170.1499878bp}%
    \ifx\svgscale\undefined%
      \relax%
    \else%
      \setlength{\unitlength}{\unitlength * \real{\svgscale}}%
    \fi%
  \else%
    \setlength{\unitlength}{\svgwidth}%
  \fi%
  \global\let\svgwidth\undefined%
  \global\let\svgscale\undefined%
  \makeatother%
  \begin{picture}(1,0.55318841)%
    \lineheight{1}%
    \setlength\tabcolsep{0pt}%
    \put(0,0){\includegraphics[width=\unitlength,page=1]{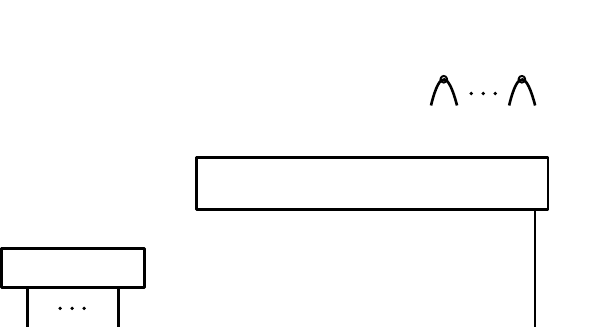}}%
    \put(0.5887598,0.22649113){\makebox(0,0)[lt]{\lineheight{1.45000005}\smash{\begin{tabular}[t]{l}$\sigma$\end{tabular}}}}%
    \put(0.06386155,0.08206292){\makebox(0,0)[lt]{\lineheight{1.45000005}\smash{\begin{tabular}[t]{l}$\alt_r$\end{tabular}}}}%
    \put(0,0){\includegraphics[width=\unitlength,page=2]{Q_lambda2.pdf}}%
    \put(0.02580344,0.43391639){\makebox(0,0)[lt]{\lineheight{1.45000005}\smash{\begin{tabular}[t]{l}$\alt_{2}$\end{tabular}}}}%
    \put(0,0){\includegraphics[width=\unitlength,page=3]{Q_lambda2.pdf}}%
    \put(0.95321227,0.44078758){\makebox(0,0)[lt]{\lineheight{1.45000005}\smash{\begin{tabular}[t]{l}$Q$\end{tabular}}}}%
  \end{picture}%
\endgroup%
},$$
   where $\sigma\in \gpS_{2d-r}$.
   The latter diagram is obtained from $Q$ by composing a morphism of $\K\F^{\op}(2d,2d)$, so is included in $A_d Q(2d)$.
  \end{proof}

  By Lemma \ref{pQlambda}, we have $A_d Q(n)\supset \sum_{\lambda\vdash d, \lambda\neq (d)}A_d Q_{\lambda}(n)$. Moreover, we have the following corollary.
  \begin{corollary}
   The $\Aut(F_n)$-module $A_d Q(n)$ is generated by $\{Q_{\lambda}\mid \lambda\vdash d, \lambda\neq (d)\}$ for $n\geq 2d$.
   That is, we have $A_d Q(n)=\sum_{\lambda\vdash d, \lambda\neq (d)}A_d Q_{\lambda}(n)$.
  \end{corollary}
  \begin{proof}
   For simplicity, let $A$ denote $\sum_{\lambda\vdash d, \lambda\neq (d)}A_d Q_{\lambda}(n)$.
   By Lemma \ref{l81}, we have $\theta_{d,n}(A_d P(n))= B_{d,0}(n)_{(2d)}$.
   Thus, by Theorem \ref{decompositionofAd}, we have
   $$\theta_{d,n}(A_d Q(n))= \left(\bigoplus_{\lambda\vdash d, \lambda\neq (d)}B_{d,0}(n)_{2\lambda}\right)\oplus \left(\bigoplus_{k\geq 1}B_{d,k}(n)\right).$$
   On the other hand, by the PBW theorem, we have
   $$\theta_{d,n}(A)\supset \left(\bigoplus_{\lambda\vdash d, \lambda\neq (d)}B_{d,0}(n)_{2\lambda}\right).$$
   By Corollary \ref{c82} and Proposition \ref{pbracketsurj}, we have
   $$\theta_{d,n}(A)\supset \left(\bigoplus_{\lambda\vdash d, \lambda\neq (d)}B_{d,0}(n)_{2\lambda}\right)\oplus \left(\bigoplus_{k\geq 1}B_{d,k}(n)\right).$$
   Therefore, we have $A_d Q(n)\subset A.$
   Hence, we have $A_d Q(n)= A$.
  \end{proof}

 \subsection{Radical filtration of $A_d(n)$}
  For an $\Aut(F_n)$-module $M$, let $\Rad(M)$ denote the \emph{radical} of $M$; that is,
  \begin{gather*}
    \begin{split}
      \Rad(M)
      =\bigcap\{K\subset M\mid K \text{ is maximal in }M\}.
    \end{split}
  \end{gather*}
  We have a radical filtration of $A_d(n)$
  $$A_d(n)\supset \Rad(A_d(n))\supset \Rad^2(A_d(n))=\Rad(\Rad(A_d(n)))\supset\cdots.$$

  \begin{theorem}\label{rad}
   Let $n\geq 2d$.
   Then, the filtration of $A_d(n)$ by the number of trivalent vertices coincides with the radical filtration.
   That is, we have $\Rad(A_{d,k}(n))=A_{d,k+1}(n)$ for any $k\geq 0$.
  \end{theorem}
  \begin{proof}
    For $\lambda\vdash 2d-k$, we have $B_{d,k}(n)_{\lambda}\cong \bigoplus_{i=1}^{r_{\lambda}}(V_{\lambda})_{i}$ as $\GL(n;\Z)$-modules.
    Let $B_{d,k}(n)_{\lambda, i}\subset B_{d,k}(n)_{\lambda}$ be a $\GL(n;\Z)$-submodule corresponding to $(V_{\lambda})_{i}$.
    Let $A_{d,k}(n)_{\lambda, i}\subset A_{d,k}(n)$ be the $\Aut(F_n)$-submodule generated by $\theta_{d,n}^{-1}(B_{d,k}(n)_{\lambda, i})$.
    For each $\lambda\vdash 2d-k, i\in [r_{\lambda}]$, we have a maximal submodule
    $$R_{\lambda, i}= \left(\sum_{(\mu,j)\neq (\lambda,i)} A_{d,k}(n)_{\mu, j}\right)+ A_{d,k+1}(n).$$
    Since we have $\bigcap_{(\lambda,i)} R_{\lambda, i}=A_{d,k+1}(n)$, it follows that $\Rad(A_{d,k}(n))\subset A_{d,k+1}(n)$.

    For any maximal submodule $K$ of $A_{d,k}(n)$, the quotient $A_{d,k}(n)/K$ is an irreducible $\Aut(F_n)$-module, which factors through an irreducible polynomial $\GL(n;\Z)$-module.
    It follows that $\theta_{d,n}(A_{d,k}(n))/\theta_{d,n}(K)$ is isomorphic to one of the irreducible components of the $\GL(n;\Z)$-module $\bigoplus_{i\geq k}B_{d,i}(n)$.
    If $B_{d,k}(n)\subset \theta_{d,n}(K)$, then by Proposition \ref{pbracketsurj}, we have $K= A_{d,k}(n)$, which contradicts to the maximality of $K$.
    Therefore, $\theta_{d,n}(A_{d,k}(n))/\theta_{d,n}(K)$ is isomorphic to one of the irreducible components of $B_{d,k}(n)$, and we have $K\supset A_{d,k+1}(n)$.
    This implies that $\Rad(A_{d,k}(n))\supset A_{d,k+1}(n)$ and the proof is complete.
  \end{proof}
  It is possible that Theorem \ref{rad} holds for some $n<2d$.
  However, it does not hold for all $n$. (See Remark \ref{rad3}.)

 \subsection{Indecomposability of the decomposition of $A_d(n)$}\label{ss825}
  Here, we consider the indecomposability of the decomposition \eqref{decomp} of $A_d(n)$.

  In Proposition \ref{pdecompgen}, we observed that $$B_{d,0}(n)\cong \bigoplus_{\lambda\vdash d}B_{d,0}(n)_{2\lambda},\quad B_{d,1}(n)\cong \bigoplus_{\mu\vdash 2d-1 \text{ with exactly $3$ odd parts}}B_{d,1}(n)_{\mu}.$$
  In order to study the indecomposability of \eqref{decomp}, we observe certain connectivity at the level of partitions.

  Let $X_d=\{ 2\lambda \mid \lambda\vdash d, \lambda \neq (d)\}$ and $Y_d=\{\mu\vdash 2d-1\mid  \mu \text{ has exactly $3$ odd parts}\}$.
  We consider the bipartite graph $G_d$ with vertex sets $X_d$ and $Y_d$ and with an edge between each pair of vertices $2\lambda$ and $\mu$ if $\mu$ is obtained from $2\lambda$ by taking away one box from each of two different rows of $2\lambda$ and then by adding one box to another row.
  For example,
  $G_2$ is
  $$\xymatrix@C=30pt@R=0pt{
    &(2^2)\ar@{-}[r]&(1^3), \\
    }
  $$
  $G_3$ is
  $$\xymatrix@C=30pt@R=0pt{
    &(4,2)\ar@{-}[r]&(3,1^2) \\
    &(2^3)\ar@{-}[ru]\ar@{-}[r]&(2,1^3)\\
    }
  $$
  and $G_4$ is
  $$\xymatrix@C=30pt@R=0pt{
    &(6,2)\ar@{-}[r]&(5,1^2) \\
    &(4^2)\ar@{-}[rd]&(4,1^3)\\
    &(4,2^2)\ar@{-}[ruu]\ar@{-}[ru]\ar@{-}[r]\ar@{-}[rd]&(3^2,1)\\
    &(2^4)\ar@{-}[r]\ar@{-}[rd]&(3,2,1^2)\\
    &&(2^2,1^3).\\
    }
  $$

  \begin{proposition}\label{p825}
   The graph $G_d$ is path-connected.
  \end{proposition}
  \begin{proof}
   For $\lambda\vdash d, \lambda\neq (d)$, let $r(\lambda)$ be the number of rows of $\lambda$. We write $\lambda=(\lambda_1^{a_1},\lambda_2^{a_2},\cdots,\lambda_l^{a_l})$, where $\lambda_1>\lambda_2>\cdots>\lambda_l$, $\sum_{i=1}^l a_i=r(\lambda)$, $a_i\geq 1$.

   We show that for $\lambda\vdash d$ such that $r(\lambda)<d$, there is a path between $2\lambda$ and some $2\lambda'\in X_d$ such that $r(\lambda')=r(\lambda)+1$.
   Then, since $(2^d)$ is the only partition that has $d$ rows, it follows  by induction on $k=r(\lambda)$ that all vertices in $X_d$ are path-connected.

   If $a_1=k$, then we have $2\lambda=((2\lambda_1)^k)$ and $2\lambda_1\geq 4$ because we assume that $k<d$.
   Thus, we have
   $$
      2\lambda \edge \mu',
   $$
   where $\mu'$ is obtained from $2\lambda$ by taking away a box from each of the $(k-1)$-st and $k$-th row and adding one box to the $(k+1)$-st row,
   and
   $$
     \mu' \edge 2\lambda',
   $$
   where $2\lambda'$ is obtained from $\mu'$ by taking away a box from the $k$-th row and adding a box to each of the $(k-1)$-st and $(k+1)$-st row.
   Therefore, we have a path between $2\lambda$ and $2\lambda'$ such that $r(\lambda')=k+1$.

   If $a_1< k$, then we have
   $$
      2\lambda \edge \mu'',
   $$
   where $\mu''$ is obtained from $2\lambda$ by taking away a box from each of the $a_1$-th and $(a_1+a_2)$-th row, and adding a box to the $(k+1)$-st row, and
   $$
     \mu''\edge 2\lambda'',
   $$
   where $2\lambda''$ is obtained from $\mu''$ by taking away a box from the $a_1$-th row
   and adding a box to each of the $(a_1+a_2)$-th and $(k+1)$-st row.
   Therefore, we have a path between $2\lambda$ and $2\lambda''$ such that $r(\lambda'')=k+1$.

   Lastly, we will show that each vertex of $Y_d$ is connected to a vertex of $X_d$. Any element $\mu\in Y_d$ is a partition of $2d-1$ and has three odd parts.
   Therefore, by taking away a box from the last odd row and then adding one box to each of the other two odd rows,
   we obtain a partition of $2d$ with only even parts, which is a vertex of $X_d$.
   The proof is complete.
  \end{proof}

  If $n\geq d$, then for any $2\lambda\in X_d$,
  $B_{d,0}(n)_{2\lambda}$ is a nonzero $\GL(n;\Z)$-submodule of $B_d(n)$.
  If $n\geq d$, then for any $\mu\in Y_d$ (except $\mu=(2^{d-2},1^3)$ if $n=d$), $B_{d,1}(n)_{\mu}$ is a nonzero $\GL(n;\Z)$-submodule of $B_d(n)$.

  Let $\pi_{\mu}:B_{d,1}(n)\twoheadrightarrow B_{d,1}(n)_{\mu}$ be the projection.
  \begin{proposition}\label{conj}
   Let $n\geq 2d$.
   Let $2\lambda\in X_d$, $\mu\in Y_d$ be two endpoints of an edge of the bipartite graph $G_d$.
   Then the composition of the bracket map and the projection $\pi_{\mu}$
   \begin{equation}\label{conjbra}
     B_{d,0}(n)_{2\lambda}\otimes \opegr^1(\IA(n))\xrightarrow{[\cdot,\cdot]} B_{d,1}(n)\xrightarrow{\pi_{\mu}} B_{d,1}(n)_{\mu}
   \end{equation}
   does not vanish.
  \end{proposition}
  Note that this proposition holds for $d=1,2$ because we have $X_1=Y_1=\emptyset, X_2=\{(2^2)\}, Y_2=\{(1^3)\}$ and by Lemma 6.7 in \cite{Mai1}.

  We introduce an intermediate vector space $B_d'(n)$ between $B_{d,0}(n)$ and $B_{d,1}(n)$.
  Define $B_d'(n)$ by
  $$\frac{\Span_{\K}\{\centre{
\begingroup%
  \makeatletter%
  \providecommand\color[2][]{%
    \errmessage{(Inkscape) Color is used for the text in Inkscape, but the package 'color.sty' is not loaded}%
    \renewcommand\color[2][]{}%
  }%
  \providecommand\transparent[1]{%
    \errmessage{(Inkscape) Transparency is used (non-zero) for the text in Inkscape, but the package 'transparent.sty' is not loaded}%
    \renewcommand\transparent[1]{}%
  }%
  \providecommand\rotatebox[2]{#2}%
  \newcommand*\fsize{\dimexpr\f@size pt\relax}%
  \newcommand*\lineheight[1]{\fontsize{\fsize}{#1\fsize}\selectfont}%
  \ifx\svgwidth\undefined%
    \setlength{\unitlength}{127.18011075bp}%
    \ifx\svgscale\undefined%
      \relax%
    \else%
      \setlength{\unitlength}{\unitlength * \real{\svgscale}}%
    \fi%
  \else%
    \setlength{\unitlength}{\svgwidth}%
  \fi%
  \global\let\svgwidth\undefined%
  \global\let\svgscale\undefined%
  \makeatother%
  \begin{picture}(1,0.27900956)%
    \lineheight{1}%
    \setlength\tabcolsep{0pt}%
    \put(0,0){\includegraphics[width=\unitlength,page=1]{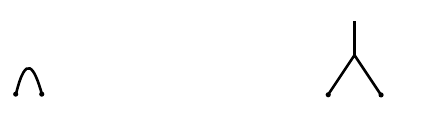}}%
    \put(0.07746383,0.0048321){\makebox(0,0)[lt]{\lineheight{1.45000005}\smash{\begin{tabular}[t]{l}\tiny$w_2$\end{tabular}}}}%
    \put(0.86430407,0.00316358){\makebox(0,0)[lt]{\lineheight{1.45000005}\smash{\begin{tabular}[t]{l}\tiny$w_{2d-2}$\end{tabular}}}}%
    \put(0.63880979,0.00455803){\makebox(0,0)[lt]{\lineheight{1.45000005}\smash{\begin{tabular}[t]{l}\tiny$w_{2d-3}$\end{tabular}}}}%
    \put(-0.0005375,0.00330051){\makebox(0,0)[lt]{\lineheight{1.45000005}\smash{\begin{tabular}[t]{l}\tiny$w_1$\end{tabular}}}}%
    \put(0.77844789,0.26690811){\makebox(0,0)[lt]{\lineheight{1.45000005}\smash{\begin{tabular}[t]{l}$\ast_1$\end{tabular}}}}%
    \put(0,0){\includegraphics[width=\unitlength,page=2]{Bdbase.pdf}}%
    \put(0.43518044,0.00747734){\makebox(0,0)[lt]{\lineheight{1.45000005}\smash{\begin{tabular}[t]{l}\tiny$w_{2d-4}$\end{tabular}}}}%
    \put(0.25475196,0.00522391){\makebox(0,0)[lt]{\lineheight{1.45000005}\smash{\begin{tabular}[t]{l}\tiny$w_{2d-5}$\end{tabular}}}}%
    \put(0,0){\includegraphics[width=\unitlength,page=3]{Bdbase.pdf}}%
  \end{picture}%
\endgroup%
} \mid w_1,\cdots, w_{2d-2}\in V_n\}}{\text{AS relation and the multilinearity}}\:,$$
  where $\centre{
\begingroup%
  \makeatletter%
  \providecommand\color[2][]{%
    \errmessage{(Inkscape) Color is used for the text in Inkscape, but the package 'color.sty' is not loaded}%
    \renewcommand\color[2][]{}%
  }%
  \providecommand\transparent[1]{%
    \errmessage{(Inkscape) Transparency is used (non-zero) for the text in Inkscape, but the package 'transparent.sty' is not loaded}%
    \renewcommand\transparent[1]{}%
  }%
  \providecommand\rotatebox[2]{#2}%
  \newcommand*\fsize{\dimexpr\f@size pt\relax}%
  \newcommand*\lineheight[1]{\fontsize{\fsize}{#1\fsize}\selectfont}%
  \ifx\svgwidth\undefined%
    \setlength{\unitlength}{45.86785417bp}%
    \ifx\svgscale\undefined%
      \relax%
    \else%
      \setlength{\unitlength}{\unitlength * \real{\svgscale}}%
    \fi%
  \else%
    \setlength{\unitlength}{\svgwidth}%
  \fi%
  \global\let\svgwidth\undefined%
  \global\let\svgscale\undefined%
  \makeatother%
  \begin{picture}(1,0.77362389)%
    \lineheight{1}%
    \setlength\tabcolsep{0pt}%
    \put(0,0){\includegraphics[width=\unitlength,page=1]{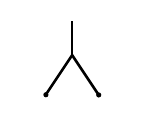}}%
    \put(0.62374906,0.00877181){\makebox(0,0)[lt]{\lineheight{1.45000005}\smash{\begin{tabular}[t]{l}\tiny$w_{2d-2}$\end{tabular}}}}%
    \put(-0.00149036,0.01263827){\makebox(0,0)[lt]{\lineheight{1.45000005}\smash{\begin{tabular}[t]{l}\tiny$w_{2d-3}$\end{tabular}}}}%
    \put(0.38569123,0.74006959){\makebox(0,0)[lt]{\lineheight{1.45000005}\smash{\begin{tabular}[t]{l}$\ast_1$\end{tabular}}}}%
    \put(0,0){\includegraphics[width=\unitlength,page=2]{base.pdf}}%
  \end{picture}%
\endgroup%
}$ is a based trivalent tree of degree $1$.
  Then, $B_d'(n)$ is a $\GL(n;\Z)$-module and we have an irreducible decomposition
  $$B_d'(n)\cong \repS_{(d-2)}(\repS_{(2)} V_n)\otimes \repS_{(1^2)} V_n \cong \bigoplus_{\nu\vdash 2d-2 \text{ with exactly $2$ odd parts}} V_{\nu}$$
  in a way similar to Proposition \ref{pdecompgen}.
  Let $B_d'(n)_{\nu}$ be the isotypic component of $B_d'(n)$ corresponding to $\nu$.

  \begin{proof}[Proof of Proposition \ref{conj}]
    Let $2\lambda=(2\lambda_1,\cdots,2\lambda_r)\vdash 2d \in X_d$.
    Any vertex $\mu\in Y_d$ that is connected to $2\lambda$ by an edge of $G_d$ is obtained from $2\lambda$ by taking away a box from each of the $i$-th and $j$-th row of $2\lambda$ and adding a box to the $k$-th row of $2\lambda$ for some $i,j\in [r], i<j, k\in[r+1], k\neq i,j$.
    We write $\mu=(\mu_1,\cdots,\mu_s)$. Then we have $\mu_i=2\lambda_i-1, \mu_j=2\lambda_j-1, \mu_k=2\lambda_k+1$ and $\mu_l=2\lambda_l$ for $l\in [s], l\neq i,j,k$.

    Since we have $\opegr^1(\IA(n))\cong H^{\ast}\otimes \Lie_2(n)$, we can write \eqref{conjbra} by
    $$h_{\lambda,\mu}:B_{d,0}(n)_{2\lambda}\otimes H^{\ast}\otimes \Lie_2(n)\rightarrow B_{d,1}(n) \xrightarrow{\pi_{\mu}} B_{d,1}(n)_{\mu}.$$
    We will show that $h_{\lambda,\mu}$ does not vanish.

    Let $\nu\vdash 2d-2$ be the partition that is obtained from $2\lambda$ by taking away a box from each of the $i$-th and $j$-th row of $2\lambda$.
    We decompose $h_{\lambda,\mu}$ into the composition
    $$h_{\lambda,\mu}=h_{\nu,\mu} h_{\lambda,\nu},$$
    where $h_{\nu,\mu}$ and $h_{\lambda,\nu}$ are $\GL(n;\Z)$-module maps defined as follows.

    Let $$h'_{\lambda}:B_{d,0}(n)_{2\lambda}\otimes\Lie_2(n)\rightarrow B_d'(n)$$
    be a $\GL(n;\Z)$-module map defined in a way similar to the contraction map in Section \ref{ss51}.
    Define
    $$h_{\lambda}:B_{d,0}(n)_{2\lambda}\otimes H^{\ast}\otimes \Lie_2(n)\rightarrow B_d'(n)\otimes H^{\ast}$$
    by $h_{\lambda}(x\otimes y\otimes z)=h'_{\lambda}(x\otimes z)\otimes y$ for $x\in B_{d,0}(n)_{2\lambda}, y\in H^{\ast}, z\in \Lie_2(n)$.
    We also define a $\GL(n;\Z)$-module map
    $$h:B_d'(n)\otimes H^{\ast}\rightarrow B_{d,1}(n)$$
    by connecting two bases $\ast_1,\ast_2$, that is, for $w_1,\cdots, w_{2d-2}\in V_n, v\in H^{\ast}$,
    $$h(\centre{} \otimes \centre{
\begingroup%
  \makeatletter%
  \providecommand\color[2][]{%
    \errmessage{(Inkscape) Color is used for the text in Inkscape, but the package 'color.sty' is not loaded}%
    \renewcommand\color[2][]{}%
  }%
  \providecommand\transparent[1]{%
    \errmessage{(Inkscape) Transparency is used (non-zero) for the text in Inkscape, but the package 'transparent.sty' is not loaded}%
    \renewcommand\transparent[1]{}%
  }%
  \providecommand\rotatebox[2]{#2}%
  \newcommand*\fsize{\dimexpr\f@size pt\relax}%
  \newcommand*\lineheight[1]{\fontsize{\fsize}{#1\fsize}\selectfont}%
  \ifx\svgwidth\undefined%
    \setlength{\unitlength}{8.74707543bp}%
    \ifx\svgscale\undefined%
      \relax%
    \else%
      \setlength{\unitlength}{\unitlength * \real{\svgscale}}%
    \fi%
  \else%
    \setlength{\unitlength}{\svgwidth}%
  \fi%
  \global\let\svgwidth\undefined%
  \global\let\svgscale\undefined%
  \makeatother%
  \begin{picture}(1,4.08114966)%
    \lineheight{1}%
    \setlength\tabcolsep{0pt}%
    \put(0,0){\includegraphics[width=\unitlength,page=1]{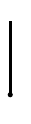}}%
    \put(0.01924798,0.04030367){\makebox(0,0)[lt]{\lineheight{1.45000005}\smash{\begin{tabular}[t]{l}$v$\end{tabular}}}}%
    \put(-0.0173049,3.90519788){\makebox(0,0)[lt]{\lineheight{1.45000005}\smash{\begin{tabular}[t]{l}$\ast_2$\end{tabular}}}}%
    \put(0,0){\includegraphics[width=\unitlength,page=2]{Hdual.pdf}}%
  \end{picture}%
\endgroup%
})= \centre{
\begingroup%
  \makeatletter%
  \providecommand\color[2][]{%
    \errmessage{(Inkscape) Color is used for the text in Inkscape, but the package 'color.sty' is not loaded}%
    \renewcommand\color[2][]{}%
  }%
  \providecommand\transparent[1]{%
    \errmessage{(Inkscape) Transparency is used (non-zero) for the text in Inkscape, but the package 'transparent.sty' is not loaded}%
    \renewcommand\transparent[1]{}%
  }%
  \providecommand\rotatebox[2]{#2}%
  \newcommand*\fsize{\dimexpr\f@size pt\relax}%
  \newcommand*\lineheight[1]{\fontsize{\fsize}{#1\fsize}\selectfont}%
  \ifx\svgwidth\undefined%
    \setlength{\unitlength}{127.18011075bp}%
    \ifx\svgscale\undefined%
      \relax%
    \else%
      \setlength{\unitlength}{\unitlength * \real{\svgscale}}%
    \fi%
  \else%
    \setlength{\unitlength}{\svgwidth}%
  \fi%
  \global\let\svgwidth\undefined%
  \global\let\svgscale\undefined%
  \makeatother%
  \begin{picture}(1,0.28267802)%
    \lineheight{1}%
    \setlength\tabcolsep{0pt}%
    \put(0,0){\includegraphics[width=\unitlength,page=1]{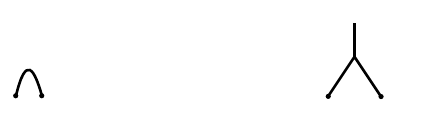}}%
    \put(0.07746383,0.0048321){\makebox(0,0)[lt]{\lineheight{1.45000005}\smash{\begin{tabular}[t]{l}\tiny$w_2$\end{tabular}}}}%
    \put(0.86430407,0.00316358){\makebox(0,0)[lt]{\lineheight{1.45000005}\smash{\begin{tabular}[t]{l}\tiny$w_{2d-2}$\end{tabular}}}}%
    \put(0.63880979,0.00455803){\makebox(0,0)[lt]{\lineheight{1.45000005}\smash{\begin{tabular}[t]{l}\tiny$w_{2d-3}$\end{tabular}}}}%
    \put(-0.0005375,0.00330051){\makebox(0,0)[lt]{\lineheight{1.45000005}\smash{\begin{tabular}[t]{l}\tiny$w_1$\end{tabular}}}}%
    \put(0.7867019,0.27057657){\makebox(0,0)[lt]{\lineheight{1.45000005}\smash{\begin{tabular}[t]{l}$v$\end{tabular}}}}%
    \put(0,0){\includegraphics[width=\unitlength,page=2]{connect.pdf}}%
    \put(0.43518044,0.00747734){\makebox(0,0)[lt]{\lineheight{1.45000005}\smash{\begin{tabular}[t]{l}\tiny$w_{2d-4}$\end{tabular}}}}%
    \put(0.25475196,0.00522391){\makebox(0,0)[lt]{\lineheight{1.45000005}\smash{\begin{tabular}[t]{l}\tiny$w_{2d-5}$\end{tabular}}}}%
    \put(0,0){\includegraphics[width=\unitlength,page=3]{connect.pdf}}%
  \end{picture}%
\endgroup%
}.$$
    Let $\pi_{\nu}: B_d'(n)\otimes H^{\ast}\rightarrow B_d'(n)_{\nu}\otimes H^{\ast}$ be the tensor product of the projection and $\id_{H^{\ast}}$.
    Then we have two $\GL(n;\Z)$-module maps
    $$h_{\lambda, \nu}:B_{d,0}(n)_{2\lambda}\otimes H^{\ast}\otimes \Lie_2(n)\xrightarrow{h_{\lambda}} B_d'(n)\otimes H^{\ast}\xrightarrow{\pi_{\nu}} B_d'(n)_{\nu}\otimes H^{\ast}$$
    and
    $$h_{\nu,\mu}:B_d'(n)_{\nu}\otimes H^{\ast}\xrightarrow{h} B_{d,1}(n)\xrightarrow{\pi_\mu} B_{d,1}(n)_{\mu}.$$

    Since $h_{\lambda,\nu}$ and $h_{\nu,\mu}$ are $\GL(n;\Z)$-module maps and since $B_{d,0}(n)_{2\lambda}$ and $B_d'(n)_{\nu}$ are irreducible, it suffices to prove that $h_{\lambda,\nu}\neq 0$ and $h_{\nu,\mu}\neq 0$.

    We will prove that $h_{\lambda,\nu}$ does not vanish.
    Let
    $$u= \centre{
\begingroup%
  \makeatletter%
  \providecommand\color[2][]{%
    \errmessage{(Inkscape) Color is used for the text in Inkscape, but the package 'color.sty' is not loaded}%
    \renewcommand\color[2][]{}%
  }%
  \providecommand\transparent[1]{%
    \errmessage{(Inkscape) Transparency is used (non-zero) for the text in Inkscape, but the package 'transparent.sty' is not loaded}%
    \renewcommand\transparent[1]{}%
  }%
  \providecommand\rotatebox[2]{#2}%
  \newcommand*\fsize{\dimexpr\f@size pt\relax}%
  \newcommand*\lineheight[1]{\fontsize{\fsize}{#1\fsize}\selectfont}%
  \ifx\svgwidth\undefined%
    \setlength{\unitlength}{222.80727726bp}%
    \ifx\svgscale\undefined%
      \relax%
    \else%
      \setlength{\unitlength}{\unitlength * \real{\svgscale}}%
    \fi%
  \else%
    \setlength{\unitlength}{\svgwidth}%
  \fi%
  \global\let\svgwidth\undefined%
  \global\let\svgscale\undefined%
  \makeatother%
  \begin{picture}(1,0.27898008)%
    \lineheight{1}%
    \setlength\tabcolsep{0pt}%
    \put(0,0){\includegraphics[width=\unitlength,page=1]{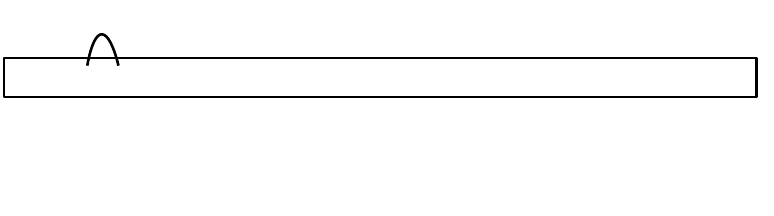}}%
    \put(0.84219177,0.17257127){\makebox(0,0)[lt]{\lineheight{1.45000005}\smash{\begin{tabular}[t]{l}\tiny$a_{2\lambda}$\end{tabular}}}}%
    \put(0,0){\includegraphics[width=\unitlength,page=2]{Bd02lambda.pdf}}%
    \put(0.01302547,0.0026817){\makebox(0,0)[lt]{\lineheight{1.45000005}\smash{\begin{tabular}[t]{l}\tiny$v_1$\end{tabular}}}}%
    \put(0.05447536,0.00325722){\makebox(0,0)[lt]{\lineheight{1.45000005}\smash{\begin{tabular}[t]{l}\tiny$v_2$\end{tabular}}}}%
    \put(0.73371616,0.00205974){\makebox(0,0)[lt]{\lineheight{1.45000005}\smash{\begin{tabular}[t]{l}\tiny$v_{2d-2}$\end{tabular}}}}%
    \put(0,0){\includegraphics[width=\unitlength,page=3]{Bd02lambda.pdf}}%
    \put(0.84470311,0.00191089){\makebox(0,0)[lt]{\lineheight{1.45000005}\smash{\begin{tabular}[t]{l}\tiny$v_{2d-1}$\end{tabular}}}}%
    \put(0.9343778,0.00303497){\makebox(0,0)[lt]{\lineheight{1.45000005}\smash{\begin{tabular}[t]{l}\tiny$v_{2d}$\end{tabular}}}}%
    \put(0,0){\includegraphics[width=\unitlength,page=4]{Bd02lambda.pdf}}%
    \put(0.36838959,0.08050605){\makebox(0,0)[lt]{\lineheight{1.45000005}\smash{\begin{tabular}[t]{l}\tiny$b_{\nu}$\end{tabular}}}}%
    \put(0,0){\includegraphics[width=\unitlength,page=5]{Bd02lambda.pdf}}%
    \put(0.8596151,0.07978357){\makebox(0,0)[lt]{\lineheight{1.45000005}\smash{\begin{tabular}[t]{l}\tiny$b_{(1^2)}$\end{tabular}}}}%
    \put(0,0){\includegraphics[width=\unitlength,page=6]{Bd02lambda.pdf}}%
    \put(0.06229753,0.27207248){\makebox(0,0)[lt]{\lineheight{1.45000005}\smash{\begin{tabular}[t]{l}\tiny$\lambda_1$\end{tabular}}}}%
    \put(0.26581926,0.27029277){\makebox(0,0)[lt]{\lineheight{1.45000005}\smash{\begin{tabular}[t]{l}\tiny$\lambda_i$\end{tabular}}}}%
    \put(0.46787532,0.27118263){\makebox(0,0)[lt]{\lineheight{1.45000005}\smash{\begin{tabular}[t]{l}\tiny$\lambda_j$\end{tabular}}}}%
    \put(0.66883116,0.2713105){\makebox(0,0)[lt]{\lineheight{1.45000005}\smash{\begin{tabular}[t]{l}\tiny$\lambda_r$\end{tabular}}}}%
    \put(0.29797765,0.00180579){\makebox(0,0)[lt]{\lineheight{1.45000005}\smash{\begin{tabular}[t]{l}\tiny$v_{\bar{i}}$\end{tabular}}}}%
    \put(0.50126903,0.00260279){\makebox(0,0)[lt]{\lineheight{1.45000005}\smash{\begin{tabular}[t]{l}\tiny$v_{\bar{j}}$\end{tabular}}}}%
    \put(0,0){\includegraphics[width=\unitlength,page=7]{Bd02lambda.pdf}}%
    \put(0.02633311,0.17772303){\makebox(0,0)[lt]{\lineheight{1.45000005}\smash{\begin{tabular}[t]{l}\tiny$\sym_{2\lambda_1}$\end{tabular}}}}%
    \put(0,0){\includegraphics[width=\unitlength,page=8]{Bd02lambda.pdf}}%
    \put(0.22904689,0.17697751){\makebox(0,0)[lt]{\lineheight{1.45000005}\smash{\begin{tabular}[t]{l}\tiny$\sym_{2\lambda_i}$\end{tabular}}}}%
    \put(0,0){\includegraphics[width=\unitlength,page=9]{Bd02lambda.pdf}}%
    \put(0.63447448,0.17846855){\makebox(0,0)[lt]{\lineheight{1.45000005}\smash{\begin{tabular}[t]{l}\tiny$\sym_{2\lambda_r}$\end{tabular}}}}%
    \put(0,0){\includegraphics[width=\unitlength,page=10]{Bd02lambda.pdf}}%
    \put(0.43176071,0.17846855){\makebox(0,0)[lt]{\lineheight{1.45000005}\smash{\begin{tabular}[t]{l}\tiny$\sym_{2\lambda_j}$\end{tabular}}}}%
    \put(0,0){\includegraphics[width=\unitlength,page=11]{Bd02lambda.pdf}}%
  \end{picture}%
\endgroup%
} \in B_{d,0}(n),$$
    where $\bar{i}=\sum_{l=1}^{i}2\lambda_l-1, \bar{j}=\sum_{l=1}^{j}2\lambda_l-2$.
    Since we have
    $$c_{\nu}\diamond c_{(1^2)}\in S^{\nu}\diamond S^{(1^2)}=\bigoplus_{\rho\vdash 2d}(S^{\rho})^{LR^{\rho}_{\nu,(1^2)}}$$
    and
    $$\{\rho\vdash 2d \mid LR^{\rho}_{\nu,(1^2)}\neq 0 \}\cap X_d=\{2\lambda\},$$
    we have $u\in B_{d,0}(n)_{2\lambda}$.
    Moreover, we have
    \begin{equation}\label{hlambda}
      h_{\lambda}\left(u\otimes \centre{
\begingroup%
  \makeatletter%
  \providecommand\color[2][]{%
    \errmessage{(Inkscape) Color is used for the text in Inkscape, but the package 'color.sty' is not loaded}%
    \renewcommand\color[2][]{}%
  }%
  \providecommand\transparent[1]{%
    \errmessage{(Inkscape) Transparency is used (non-zero) for the text in Inkscape, but the package 'transparent.sty' is not loaded}%
    \renewcommand\transparent[1]{}%
  }%
  \providecommand\rotatebox[2]{#2}%
  \newcommand*\fsize{\dimexpr\f@size pt\relax}%
  \newcommand*\lineheight[1]{\fontsize{\fsize}{#1\fsize}\selectfont}%
  \ifx\svgwidth\undefined%
    \setlength{\unitlength}{8.74707543bp}%
    \ifx\svgscale\undefined%
      \relax%
    \else%
      \setlength{\unitlength}{\unitlength * \real{\svgscale}}%
    \fi%
  \else%
    \setlength{\unitlength}{\svgwidth}%
  \fi%
  \global\let\svgwidth\undefined%
  \global\let\svgscale\undefined%
  \makeatother%
  \begin{picture}(1,4.08114966)%
    \lineheight{1}%
    \setlength\tabcolsep{0pt}%
    \put(0,0){\includegraphics[width=\unitlength,page=1]{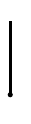}}%
    \put(0.01924798,0.04030367){\makebox(0,0)[lt]{\lineheight{1.45000005}\smash{\begin{tabular}[t]{l}$v_1$\end{tabular}}}}%
    \put(-0.0173049,3.90519788){\makebox(0,0)[lt]{\lineheight{1.45000005}\smash{\begin{tabular}[t]{l}$\ast_2$\end{tabular}}}}%
    \put(0,0){\includegraphics[width=\unitlength,page=2]{Hdual1.pdf}}%
  \end{picture}%
\endgroup%
}\otimes  \centre{
\begingroup%
  \makeatletter%
  \providecommand\color[2][]{%
    \errmessage{(Inkscape) Color is used for the text in Inkscape, but the package 'color.sty' is not loaded}%
    \renewcommand\color[2][]{}%
  }%
  \providecommand\transparent[1]{%
    \errmessage{(Inkscape) Transparency is used (non-zero) for the text in Inkscape, but the package 'transparent.sty' is not loaded}%
    \renewcommand\transparent[1]{}%
  }%
  \providecommand\rotatebox[2]{#2}%
  \newcommand*\fsize{\dimexpr\f@size pt\relax}%
  \newcommand*\lineheight[1]{\fontsize{\fsize}{#1\fsize}\selectfont}%
  \ifx\svgwidth\undefined%
    \setlength{\unitlength}{31.13256704bp}%
    \ifx\svgscale\undefined%
      \relax%
    \else%
      \setlength{\unitlength}{\unitlength * \real{\svgscale}}%
    \fi%
  \else%
    \setlength{\unitlength}{\svgwidth}%
  \fi%
  \global\let\svgwidth\undefined%
  \global\let\svgscale\undefined%
  \makeatother%
  \begin{picture}(1,1.12344744)%
    \lineheight{1}%
    \setlength\tabcolsep{0pt}%
    \put(0,0){\includegraphics[width=\unitlength,page=1]{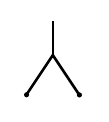}}%
    \put(0.65670975,0.00912599){\makebox(0,0)[lt]{\lineheight{1.45000005}\smash{\begin{tabular}[t]{l}\tiny$\bar{x}_{2d}$\end{tabular}}}}%
    \put(-0.00109788,0.00646179){\makebox(0,0)[lt]{\lineheight{1.45000005}\smash{\begin{tabular}[t]{l}\tiny$\bar{x}_{2d-1}$\end{tabular}}}}%
    \put(0.38958445,1.07401164){\makebox(0,0)[lt]{\lineheight{1.45000005}\smash{\begin{tabular}[t]{l}$\ast_1$\end{tabular}}}}%
    \put(0,0){\includegraphics[width=\unitlength,page=2]{based.pdf}}%
  \end{picture}%
\endgroup%
}\right)=\centre{
\begingroup%
  \makeatletter%
  \providecommand\color[2][]{%
    \errmessage{(Inkscape) Color is used for the text in Inkscape, but the package 'color.sty' is not loaded}%
    \renewcommand\color[2][]{}%
  }%
  \providecommand\transparent[1]{%
    \errmessage{(Inkscape) Transparency is used (non-zero) for the text in Inkscape, but the package 'transparent.sty' is not loaded}%
    \renewcommand\transparent[1]{}%
  }%
  \providecommand\rotatebox[2]{#2}%
  \newcommand*\fsize{\dimexpr\f@size pt\relax}%
  \newcommand*\lineheight[1]{\fontsize{\fsize}{#1\fsize}\selectfont}%
  \ifx\svgwidth\undefined%
    \setlength{\unitlength}{218.25117512bp}%
    \ifx\svgscale\undefined%
      \relax%
    \else%
      \setlength{\unitlength}{\unitlength * \real{\svgscale}}%
    \fi%
  \else%
    \setlength{\unitlength}{\svgwidth}%
  \fi%
  \global\let\svgwidth\undefined%
  \global\let\svgscale\undefined%
  \makeatother%
  \begin{picture}(1,0.3165036)%
    \lineheight{1}%
    \setlength\tabcolsep{0pt}%
    \put(0,0){\includegraphics[width=\unitlength,page=1]{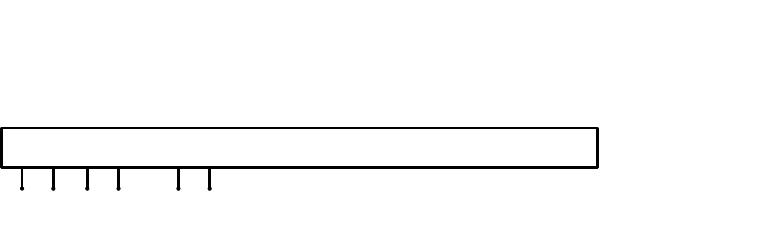}}%
    \put(0.37607991,0.11388633){\makebox(0,0)[lt]{\lineheight{1.45000005}\smash{\begin{tabular}[t]{l}\tiny$b_{\nu}$\end{tabular}}}}%
    \put(0,0){\includegraphics[width=\unitlength,page=2]{Bdpramu.pdf}}%
    \put(0.87755999,0.11314876){\makebox(0,0)[lt]{\lineheight{1.45000005}\smash{\begin{tabular}[t]{l}\tiny$b_{(1^2)}$\end{tabular}}}}%
    \put(0,0){\includegraphics[width=\unitlength,page=3]{Bdpramu.pdf}}%
    \put(0.90236444,0.00161529){\makebox(0,0)[lt]{\lineheight{1.45000005}\smash{\begin{tabular}[t]{l}$\ast_1$\end{tabular}}}}%
    \put(0.00938993,0.03403564){\makebox(0,0)[lt]{\lineheight{1.45000005}\smash{\begin{tabular}[t]{l}\tiny$v_1$\end{tabular}}}}%
    \put(0.05170511,0.03462317){\makebox(0,0)[lt]{\lineheight{1.45000005}\smash{\begin{tabular}[t]{l}\tiny$v_2$\end{tabular}}}}%
    \put(0.74512542,0.03340069){\makebox(0,0)[lt]{\lineheight{1.45000005}\smash{\begin{tabular}[t]{l}\tiny$v_{2d-2}$\end{tabular}}}}%
    \put(0.30029062,0.03314145){\makebox(0,0)[lt]{\lineheight{1.45000005}\smash{\begin{tabular}[t]{l}\tiny$v_{\bar{i}}$\end{tabular}}}}%
    \put(0.50782583,0.03395508){\makebox(0,0)[lt]{\lineheight{1.45000005}\smash{\begin{tabular}[t]{l}\tiny$v_{\bar{j}}$\end{tabular}}}}%
    \put(0,0){\includegraphics[width=\unitlength,page=4]{Bdpramu.pdf}}%
    \put(0.85977304,0.20787345){\makebox(0,0)[lt]{\lineheight{1.45000005}\smash{\begin{tabular}[t]{l}\tiny$a_{2\lambda}$\end{tabular}}}}%
    \put(0,0){\includegraphics[width=\unitlength,page=5]{Bdpramu.pdf}}%
    \put(0.06359805,0.3094518){\makebox(0,0)[lt]{\lineheight{1.45000005}\smash{\begin{tabular}[t]{l}\tiny$\lambda_1$\end{tabular}}}}%
    \put(0.27136839,0.30763493){\makebox(0,0)[lt]{\lineheight{1.45000005}\smash{\begin{tabular}[t]{l}\tiny$\lambda_i$\end{tabular}}}}%
    \put(0.47764245,0.30854337){\makebox(0,0)[lt]{\lineheight{1.45000005}\smash{\begin{tabular}[t]{l}\tiny$\lambda_j$\end{tabular}}}}%
    \put(0.68279335,0.30867392){\makebox(0,0)[lt]{\lineheight{1.45000005}\smash{\begin{tabular}[t]{l}\tiny$\lambda_r$\end{tabular}}}}%
    \put(0,0){\includegraphics[width=\unitlength,page=6]{Bdpramu.pdf}}%
    \put(0.02688285,0.21313275){\makebox(0,0)[lt]{\lineheight{1.45000005}\smash{\begin{tabular}[t]{l}\tiny$\sym_{2\lambda_1}$\end{tabular}}}}%
    \put(0,0){\includegraphics[width=\unitlength,page=7]{Bdpramu.pdf}}%
    \put(0.23382838,0.21237168){\makebox(0,0)[lt]{\lineheight{1.45000005}\smash{\begin{tabular}[t]{l}\tiny$\sym_{2\lambda_i}$\end{tabular}}}}%
    \put(0,0){\includegraphics[width=\unitlength,page=8]{Bdpramu.pdf}}%
    \put(0.64771945,0.21389384){\makebox(0,0)[lt]{\lineheight{1.45000005}\smash{\begin{tabular}[t]{l}\tiny$\sym_{2\lambda_r}$\end{tabular}}}}%
    \put(0,0){\includegraphics[width=\unitlength,page=9]{Bdpramu.pdf}}%
    \put(0.44077393,0.21389384){\makebox(0,0)[lt]{\lineheight{1.45000005}\smash{\begin{tabular}[t]{l}\tiny$\sym_{2\lambda_j}$\end{tabular}}}}%
    \put(0,0){\includegraphics[width=\unitlength,page=10]{Bdpramu.pdf}}%
  \end{picture}%
\endgroup%
}\otimes \centre{}.
    \end{equation}
    By pulling $\ast_1$ toward the upper side, since we have $\centre{
\begingroup%
  \makeatletter%
  \providecommand\color[2][]{%
    \errmessage{(Inkscape) Color is used for the text in Inkscape, but the package 'color.sty' is not loaded}%
    \renewcommand\color[2][]{}%
  }%
  \providecommand\transparent[1]{%
    \errmessage{(Inkscape) Transparency is used (non-zero) for the text in Inkscape, but the package 'transparent.sty' is not loaded}%
    \renewcommand\transparent[1]{}%
  }%
  \providecommand\rotatebox[2]{#2}%
  \newcommand*\fsize{\dimexpr\f@size pt\relax}%
  \newcommand*\lineheight[1]{\fontsize{\fsize}{#1\fsize}\selectfont}%
  \ifx\svgwidth\undefined%
    \setlength{\unitlength}{60.75117863bp}%
    \ifx\svgscale\undefined%
      \relax%
    \else%
      \setlength{\unitlength}{\unitlength * \real{\svgscale}}%
    \fi%
  \else%
    \setlength{\unitlength}{\svgwidth}%
  \fi%
  \global\let\svgwidth\undefined%
  \global\let\svgscale\undefined%
  \makeatother%
  \begin{picture}(1,0.73634129)%
    \lineheight{1}%
    \setlength\tabcolsep{0pt}%
    \put(0,0){\includegraphics[width=\unitlength,page=1]{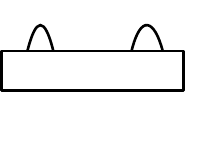}}%
    \put(0.19136404,0.37036319){\makebox(0,0)[lt]{\lineheight{1.45000005}\smash{\begin{tabular}[t]{l}$\sym_{2l}$\end{tabular}}}}%
    \put(0,0){\includegraphics[width=\unitlength,page=2]{symm2l.pdf}}%
    \put(0.41035241,0.71100741){\makebox(0,0)[lt]{\lineheight{1.45000005}\smash{\begin{tabular}[t]{l}\tiny$l$\end{tabular}}}}%
  \end{picture}%
\endgroup%
}=2l\centre{
\begingroup%
  \makeatletter%
  \providecommand\color[2][]{%
    \errmessage{(Inkscape) Color is used for the text in Inkscape, but the package 'color.sty' is not loaded}%
    \renewcommand\color[2][]{}%
  }%
  \providecommand\transparent[1]{%
    \errmessage{(Inkscape) Transparency is used (non-zero) for the text in Inkscape, but the package 'transparent.sty' is not loaded}%
    \renewcommand\transparent[1]{}%
  }%
  \providecommand\rotatebox[2]{#2}%
  \newcommand*\fsize{\dimexpr\f@size pt\relax}%
  \newcommand*\lineheight[1]{\fontsize{\fsize}{#1\fsize}\selectfont}%
  \ifx\svgwidth\undefined%
    \setlength{\unitlength}{53.25118286bp}%
    \ifx\svgscale\undefined%
      \relax%
    \else%
      \setlength{\unitlength}{\unitlength * \real{\svgscale}}%
    \fi%
  \else%
    \setlength{\unitlength}{\svgwidth}%
  \fi%
  \global\let\svgwidth\undefined%
  \global\let\svgscale\undefined%
  \makeatother%
  \begin{picture}(1,0.84505164)%
    \lineheight{1}%
    \setlength\tabcolsep{0pt}%
    \put(0,0){\includegraphics[width=\unitlength,page=1]{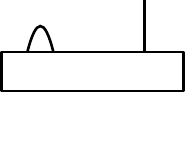}}%
    \put(0.21831611,0.42252583){\makebox(0,0)[lt]{\lineheight{1.45000005}\smash{\begin{tabular}[t]{l}$\sym_{2l-1}$\end{tabular}}}}%
    \put(0,0){\includegraphics[width=\unitlength,page=2]{symm2l-1.pdf}}%
    \put(0.29600788,0.80885649){\makebox(0,0)[lt]{\lineheight{1.45000005}\smash{\begin{tabular}[t]{l}\tiny$l-1$\end{tabular}}}}%
  \end{picture}%
\endgroup%
}$, the right-hand side of \eqref{hlambda} is
    $$u'=(-2)(2\lambda_i)(2\lambda_j)\centre{
\begingroup%
  \makeatletter%
  \providecommand\color[2][]{%
    \errmessage{(Inkscape) Color is used for the text in Inkscape, but the package 'color.sty' is not loaded}%
    \renewcommand\color[2][]{}%
  }%
  \providecommand\transparent[1]{%
    \errmessage{(Inkscape) Transparency is used (non-zero) for the text in Inkscape, but the package 'transparent.sty' is not loaded}%
    \renewcommand\transparent[1]{}%
  }%
  \providecommand\rotatebox[2]{#2}%
  \newcommand*\fsize{\dimexpr\f@size pt\relax}%
  \newcommand*\lineheight[1]{\fontsize{\fsize}{#1\fsize}\selectfont}%
  \ifx\svgwidth\undefined%
    \setlength{\unitlength}{192.0011803bp}%
    \ifx\svgscale\undefined%
      \relax%
    \else%
      \setlength{\unitlength}{\unitlength * \real{\svgscale}}%
    \fi%
  \else%
    \setlength{\unitlength}{\svgwidth}%
  \fi%
  \global\let\svgwidth\undefined%
  \global\let\svgscale\undefined%
  \makeatother%
  \begin{picture}(1,0.41570471)%
    \lineheight{1}%
    \setlength\tabcolsep{0pt}%
    \put(0,0){\includegraphics[width=\unitlength,page=1]{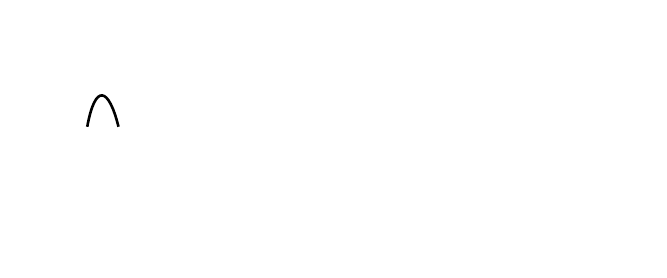}}%
    \put(0.91647105,0.20225132){\makebox(0,0)[lt]{\lineheight{1.45000005}\smash{\begin{tabular}[t]{l}\tiny$a_{\nu}$\end{tabular}}}}%
    \put(0,0){\includegraphics[width=\unitlength,page=2]{Bdpramu2.pdf}}%
    \put(0.42749676,0.09374051){\makebox(0,0)[lt]{\lineheight{1.45000005}\smash{\begin{tabular}[t]{l}\tiny$b_{\nu}$\end{tabular}}}}%
    \put(0,0){\includegraphics[width=\unitlength,page=3]{Bdpramu2.pdf}}%
    \put(0.07229302,0.31604329){\makebox(0,0)[lt]{\lineheight{1.45000005}\smash{\begin{tabular}[t]{l}\tiny$\lambda_1$\end{tabular}}}}%
    \put(0.77614341,0.31515905){\makebox(0,0)[lt]{\lineheight{1.45000005}\smash{\begin{tabular}[t]{l}\tiny$\lambda_r$\end{tabular}}}}%
    \put(0.46923643,0.4076888){\makebox(0,0)[lt]{\lineheight{1.45000005}\smash{\begin{tabular}[t]{l}$\ast_1$\end{tabular}}}}%
    \put(0.01685026,0.00311197){\makebox(0,0)[lt]{\lineheight{1.45000005}\smash{\begin{tabular}[t]{l}\tiny$v_1$\end{tabular}}}}%
    \put(0.06495067,0.00377983){\makebox(0,0)[lt]{\lineheight{1.45000005}\smash{\begin{tabular}[t]{l}\tiny$v_2$\end{tabular}}}}%
    \put(0.85317395,0.00239022){\makebox(0,0)[lt]{\lineheight{1.45000005}\smash{\begin{tabular}[t]{l}\tiny$v_{2d-2}$\end{tabular}}}}%
    \put(0.34752228,0.00209553){\makebox(0,0)[lt]{\lineheight{1.45000005}\smash{\begin{tabular}[t]{l}\tiny$v_{\bar{i}}$\end{tabular}}}}%
    \put(0.58343127,0.0030204){\makebox(0,0)[lt]{\lineheight{1.45000005}\smash{\begin{tabular}[t]{l}\tiny$v_{\bar{j}}$\end{tabular}}}}%
    \put(0,0){\includegraphics[width=\unitlength,page=4]{Bdpramu2.pdf}}%
    \put(0.03055822,0.20655569){\makebox(0,0)[lt]{\lineheight{1.45000005}\smash{\begin{tabular}[t]{l}\tiny$\sym_{2\lambda_1}$\end{tabular}}}}%
    \put(0,0){\includegraphics[width=\unitlength,page=5]{Bdpramu2.pdf}}%
    \put(0.2423595,0.20569056){\makebox(0,0)[lt]{\lineheight{1.45000005}\smash{\begin{tabular}[t]{l}\tiny$\sym_{2\lambda_i-1}$\end{tabular}}}}%
    \put(0,0){\includegraphics[width=\unitlength,page=6]{Bdpramu2.pdf}}%
    \put(0.73627434,0.20742084){\makebox(0,0)[lt]{\lineheight{1.45000005}\smash{\begin{tabular}[t]{l}\tiny$\sym_{2\lambda_r}$\end{tabular}}}}%
    \put(0,0){\includegraphics[width=\unitlength,page=7]{Bdpramu2.pdf}}%
    \put(0.47759829,0.20742084){\makebox(0,0)[lt]{\lineheight{1.45000005}\smash{\begin{tabular}[t]{l}\tiny$\sym_{2\lambda_j-1}$\end{tabular}}}}%
    \put(0,0){\includegraphics[width=\unitlength,page=8]{Bdpramu2.pdf}}%
  \end{picture}%
\endgroup%
}\otimes \centre{}\in B_d'(n)_{\nu}\otimes H^{\ast}.$$
    We will look at the coefficient in $u'$ of $u_0=\centre{
\begingroup%
  \makeatletter%
  \providecommand\color[2][]{%
    \errmessage{(Inkscape) Color is used for the text in Inkscape, but the package 'color.sty' is not loaded}%
    \renewcommand\color[2][]{}%
  }%
  \providecommand\transparent[1]{%
    \errmessage{(Inkscape) Transparency is used (non-zero) for the text in Inkscape, but the package 'transparent.sty' is not loaded}%
    \renewcommand\transparent[1]{}%
  }%
  \providecommand\rotatebox[2]{#2}%
  \newcommand*\fsize{\dimexpr\f@size pt\relax}%
  \newcommand*\lineheight[1]{\fontsize{\fsize}{#1\fsize}\selectfont}%
  \ifx\svgwidth\undefined%
    \setlength{\unitlength}{119.52560238bp}%
    \ifx\svgscale\undefined%
      \relax%
    \else%
      \setlength{\unitlength}{\unitlength * \real{\svgscale}}%
    \fi%
  \else%
    \setlength{\unitlength}{\svgwidth}%
  \fi%
  \global\let\svgwidth\undefined%
  \global\let\svgscale\undefined%
  \makeatother%
  \begin{picture}(1,0.29963527)%
    \lineheight{1}%
    \setlength\tabcolsep{0pt}%
    \put(0,0){\includegraphics[width=\unitlength,page=1]{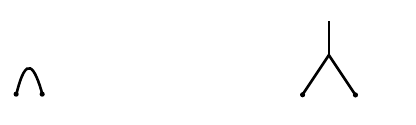}}%
    \put(0.09079134,0.00487478){\makebox(0,0)[lt]{\lineheight{1.45000005}\smash{\begin{tabular}[t]{l}\tiny$v_2$\end{tabular}}}}%
    \put(0.83649544,0.00803888){\makebox(0,0)[lt]{\lineheight{1.45000005}\smash{\begin{tabular}[t]{l}\tiny$v_{\bar{j}}$\end{tabular}}}}%
    \put(0.69876861,0.00624746){\makebox(0,0)[lt]{\lineheight{1.45000005}\smash{\begin{tabular}[t]{l}\tiny$v_{\bar{i}}$\end{tabular}}}}%
    \put(-0.00057192,0.00626961){\makebox(0,0)[lt]{\lineheight{1.45000005}\smash{\begin{tabular}[t]{l}\tiny$v_1$\end{tabular}}}}%
    \put(0.77114499,0.28675883){\makebox(0,0)[lt]{\lineheight{1.45000005}\smash{\begin{tabular}[t]{l}$\ast_1$\end{tabular}}}}%
    \put(0,0){\includegraphics[width=\unitlength,page=2]{basedij.pdf}}%
    \put(0.51458194,0.0035683){\makebox(0,0)[lt]{\lineheight{1.45000005}\smash{\begin{tabular}[t]{l}\tiny$v_{2d-2}$\end{tabular}}}}%
    \put(0.33532755,0.00336617){\makebox(0,0)[lt]{\lineheight{1.45000005}\smash{\begin{tabular}[t]{l}\tiny$v_{2d-3}$\end{tabular}}}}%
  \end{picture}%
\endgroup%
}$ to show that $u'$ does not vanish.
    Note that $b_{\nu}a_{\nu}=\sum_{\tau\in C_{t_0},\rho\in R_{t_0}} \sgn(\tau)\tau\rho$, where $t_0$ is the canonical $\nu$-tableau.
    If $\lambda_i\neq \lambda_j$, then there is no $\tau\in C_{t_0}$ such that $\tau(\bar{i})=\bar{j}, \tau(\bar{j})=\bar{i}$.
    Thus, the diagram $u_0$ appears only when $\tau$ is an even permutation which fixes $\bar{i}$ and $\bar{j}$.
    Then, the coefficient of $u_0$ in $u'$ is negative.
    If $\lambda_i= \lambda_j$, then the diagram $u_0$ appears
    when $\tau$ preserves the subset $\{\bar{i},\bar{j}\}$ and the parity of $\tau$ coincides with that of the restriction of $\tau$ to $\{\bar{i},\bar{j}\}$.
    Hence, by the AS relation, the coefficient of $u_0$ in $u'$ is negative.
    Therefore, $h_{\lambda, \nu}$ does not vanish.

    We will prove that $h_{\nu,\mu}$ does not vanish.
    Let $N\in \N$. Set $c'_{\rho}=a_{\rho} b_{\rho}\in\K\gpS_{N}$ for $\rho \vdash N$.
    From basic facts of representation theory, we have of $\K\gpS_{N}$-modules
    $$\K\gpS_{N}c_{\rho} \cong \K\gpS_{N}c'_{\rho}.$$
    In what follows, we use $c'_{\rho}$ instead of $c_{\rho}$ as the Young symmetrizer.
    Let
    \begin{gather*}
     \begin{split}
       Z_{\mu}&=\left(c'_{\mu}\sigma\cdot \centre{
\begingroup%
  \makeatletter%
  \providecommand\color[2][]{%
    \errmessage{(Inkscape) Color is used for the text in Inkscape, but the package 'color.sty' is not loaded}%
    \renewcommand\color[2][]{}%
  }%
  \providecommand\transparent[1]{%
    \errmessage{(Inkscape) Transparency is used (non-zero) for the text in Inkscape, but the package 'transparent.sty' is not loaded}%
    \renewcommand\transparent[1]{}%
  }%
  \providecommand\rotatebox[2]{#2}%
  \newcommand*\fsize{\dimexpr\f@size pt\relax}%
  \newcommand*\lineheight[1]{\fontsize{\fsize}{#1\fsize}\selectfont}%
  \ifx\svgwidth\undefined%
    \setlength{\unitlength}{119.8200084bp}%
    \ifx\svgscale\undefined%
      \relax%
    \else%
      \setlength{\unitlength}{\unitlength * \real{\svgscale}}%
    \fi%
  \else%
    \setlength{\unitlength}{\svgwidth}%
  \fi%
  \global\let\svgwidth\undefined%
  \global\let\svgscale\undefined%
  \makeatother%
  \begin{picture}(1,0.29739414)%
    \lineheight{1}%
    \setlength\tabcolsep{0pt}%
    \put(0,0){\includegraphics[width=\unitlength,page=1]{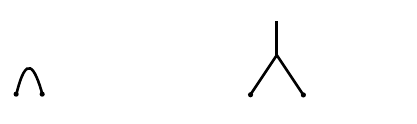}}%
    \put(0.09056825,0.0033579){\makebox(0,0)[lt]{\lineheight{1.45000005}\smash{\begin{tabular}[t]{l}\tiny$2$\end{tabular}}}}%
    \put(0.6790817,0.00651422){\makebox(0,0)[lt]{\lineheight{1.45000005}\smash{\begin{tabular}[t]{l}\tiny$2d-2$\end{tabular}}}}%
    \put(0.48436904,0.00623575){\makebox(0,0)[lt]{\lineheight{1.45000005}\smash{\begin{tabular}[t]{l}\tiny$2d-3$\end{tabular}}}}%
    \put(-0.00057052,0.0047493){\makebox(0,0)[lt]{\lineheight{1.45000005}\smash{\begin{tabular}[t]{l}\tiny$1$\end{tabular}}}}%
    \put(0.64406247,0.28454935){\makebox(0,0)[lt]{\lineheight{1.45000005}\smash{\begin{tabular}[t]{l}$\ast_1$\end{tabular}}}}%
    \put(0,0){\includegraphics[width=\unitlength,page=2]{based2.pdf}}%
    \put(0.92573492,0.28454935){\makebox(0,0)[lt]{\lineheight{1.45000005}\smash{\begin{tabular}[t]{l}$\ast_2$\end{tabular}}}}%
    \put(0,0){\includegraphics[width=\unitlength,page=3]{based2.pdf}}%
    \put(0.86333031,0.00607459){\makebox(0,0)[lt]{\lineheight{1.45000005}\smash{\begin{tabular}[t]{l}\tiny$2d-1$\end{tabular}}}}%
    \put(0,0){\includegraphics[width=\unitlength,page=4]{based2.pdf}}%
  \end{picture}%
\endgroup%
}\right)
       (v_1^{\otimes \mu_1}\otimes \cdots \otimes v_s^{\otimes \mu_s})\\
       &=\centre{
\begingroup%
  \makeatletter%
  \providecommand\color[2][]{%
    \errmessage{(Inkscape) Color is used for the text in Inkscape, but the package 'color.sty' is not loaded}%
    \renewcommand\color[2][]{}%
  }%
  \providecommand\transparent[1]{%
    \errmessage{(Inkscape) Transparency is used (non-zero) for the text in Inkscape, but the package 'transparent.sty' is not loaded}%
    \renewcommand\transparent[1]{}%
  }%
  \providecommand\rotatebox[2]{#2}%
  \newcommand*\fsize{\dimexpr\f@size pt\relax}%
  \newcommand*\lineheight[1]{\fontsize{\fsize}{#1\fsize}\selectfont}%
  \ifx\svgwidth\undefined%
    \setlength{\unitlength}{94.77853823bp}%
    \ifx\svgscale\undefined%
      \relax%
    \else%
      \setlength{\unitlength}{\unitlength * \real{\svgscale}}%
    \fi%
  \else%
    \setlength{\unitlength}{\svgwidth}%
  \fi%
  \global\let\svgwidth\undefined%
  \global\let\svgscale\undefined%
  \makeatother%
  \begin{picture}(1,0.76491741)%
    \lineheight{1}%
    \setlength\tabcolsep{0pt}%
    \put(0,0){\includegraphics[width=\unitlength,page=1]{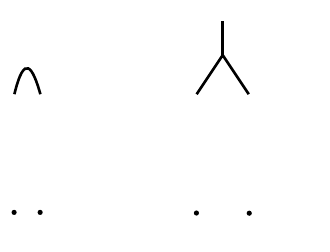}}%
    \put(-0.00072125,0.0046338){\makebox(0,0)[lt]{\lineheight{1.45000005}\smash{\begin{tabular}[t]{l}\tiny$v_1$\end{tabular}}}}%
    \put(0.64987967,0.74867889){\makebox(0,0)[lt]{\lineheight{1.45000005}\smash{\begin{tabular}[t]{l}$\ast_1$\end{tabular}}}}%
    \put(0,0){\includegraphics[width=\unitlength,page=2]{zmu.pdf}}%
    \put(0.88727513,0.74867889){\makebox(0,0)[lt]{\lineheight{1.45000005}\smash{\begin{tabular}[t]{l}$\ast_2$\end{tabular}}}}%
    \put(0,0){\includegraphics[width=\unitlength,page=3]{zmu.pdf}}%
    \put(0.88810257,0.0042451){\makebox(0,0)[lt]{\lineheight{1.45000005}\smash{\begin{tabular}[t]{l}\tiny$v_s$\end{tabular}}}}%
    \put(0,0){\includegraphics[width=\unitlength,page=4]{zmu.pdf}}%
    \put(0.44174804,0.42016383){\makebox(0,0)[lt]{\lineheight{1.45000005}\smash{\begin{tabular}[t]{l}\tiny$\sigma$\end{tabular}}}}%
    \put(0,0){\includegraphics[width=\unitlength,page=5]{zmu.pdf}}%
    \put(0.42074338,0.27085036){\makebox(0,0)[lt]{\lineheight{1.45000005}\smash{\begin{tabular}[t]{l}\tiny$c'_{\mu}$\end{tabular}}}}%
  \end{picture}%
\endgroup%
}\in B_d'(n)\otimes H^{\ast},
     \end{split}
    \end{gather*}
    where $\sigma\in \gpS_{2d-1}$ is defined by
    $$\sigma=
    \begin{pmatrix}
      1&\cdots  &2d-3 &2d-2 & 2d-1\\
      1&\cdots  &i'   &j'   & k'
    \end{pmatrix}
    \text{ for } i'=\sum_{l=1}^{i}\mu_{l},\: j'=\sum_{l=1}^{j} \mu_{l},\: k'=\sum_{l=1}^{k} \mu_{l}.
    $$
    We will show that $h(\pi_{\nu}(Z_{\mu}))\in B_{d,1}(n)_{\mu}$ and that $h(\pi_{\nu}(Z_{\mu}))\neq 0$.

    If the diagram that is obtained from $\mu$ by taking away a box from the $i$-th (resp. $j$-th) row of $\mu$ is a partition of $2d-2$, then write it $\nu_i$ (resp. $\nu_j$).
    Since any partition $\rho\vdash 2d-2$ with exactly two odd parts other than $\nu,\nu_i,\nu_j$ is not included in $\mu$, it follows that
    $$Z_\mu \in (B_d'(n)_{\nu}\otimes H^{\ast}) \oplus  (B_d'(n)_{\nu_i}\otimes H^{\ast}) \oplus(B_d'(n)_{\nu_j}\otimes H^{\ast}).$$
    By using an argument similar to Proposition \ref{p822}, we have
    \begin{gather*}
      h(B_d'(n)_{\nu}\otimes H^{\ast})\subset \bigoplus_{\alpha=\nu\sqcup \square} B_{d,1}(n)_{\alpha},\quad
      h(B_d'(n)_{\nu_i}\otimes H^{\ast})\subset \bigoplus_{\alpha=\nu_i\sqcup \square} B_{d,1}(n)_{\alpha},\\
      h(B_d'(n)_{\nu_j}\otimes H^{\ast})\subset \bigoplus_{\alpha=\nu_j\sqcup \square} B_{d,1}(n)_{\alpha}.
    \end{gather*}
    Since $\{\nu\sqcup \square\}\cap \{\nu_i\sqcup \square\}\cap \{\nu_j\sqcup \square\}=\{\mu\}$ and since $h(Z_{\mu})\in B_{d,1}(n)_{\mu}$, we have $h(\pi_{\nu}(Z_{\mu}))\in B_{d,1}(n)_{\mu}$.

    In order to prove that $h(\pi_{\nu}(Z_{\mu}))\neq 0$, we will look at the coefficient in $h(\pi_{\nu}(Z_{\mu}))$
    \vspace{0.05in}\\
    of
    \begin{gather*}
     \begin{split}
       z&=h\left(\left(\sigma\cdot \centre{}\right) (v_1^{\otimes \mu_1}\otimes \cdots \otimes v_s^{\otimes \mu_s})\right)\\
       &=\centre{
\begingroup%
  \makeatletter%
  \providecommand\color[2][]{%
    \errmessage{(Inkscape) Color is used for the text in Inkscape, but the package 'color.sty' is not loaded}%
    \renewcommand\color[2][]{}%
  }%
  \providecommand\transparent[1]{%
    \errmessage{(Inkscape) Transparency is used (non-zero) for the text in Inkscape, but the package 'transparent.sty' is not loaded}%
    \renewcommand\transparent[1]{}%
  }%
  \providecommand\rotatebox[2]{#2}%
  \newcommand*\fsize{\dimexpr\f@size pt\relax}%
  \newcommand*\lineheight[1]{\fontsize{\fsize}{#1\fsize}\selectfont}%
  \ifx\svgwidth\undefined%
    \setlength{\unitlength}{94.77853823bp}%
    \ifx\svgscale\undefined%
      \relax%
    \else%
      \setlength{\unitlength}{\unitlength * \real{\svgscale}}%
    \fi%
  \else%
    \setlength{\unitlength}{\svgwidth}%
  \fi%
  \global\let\svgwidth\undefined%
  \global\let\svgscale\undefined%
  \makeatother%
  \begin{picture}(1,0.54733959)%
    \lineheight{1}%
    \setlength\tabcolsep{0pt}%
    \put(0,0){\includegraphics[width=\unitlength,page=1]{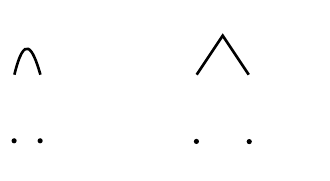}}%
    \put(-0.00072125,0.0046338){\makebox(0,0)[lt]{\lineheight{1.45000005}\smash{\begin{tabular}[t]{l}\tiny$v_1$\end{tabular}}}}%
    \put(0,0){\includegraphics[width=\unitlength,page=2]{z.pdf}}%
    \put(0.88810257,0.0042451){\makebox(0,0)[lt]{\lineheight{1.45000005}\smash{\begin{tabular}[t]{l}\tiny$v_s$\end{tabular}}}}%
    \put(0,0){\includegraphics[width=\unitlength,page=3]{z.pdf}}%
    \put(0.44174804,0.26190017){\makebox(0,0)[lt]{\lineheight{1.45000005}\smash{\begin{tabular}[t]{l}\tiny$\sigma$\end{tabular}}}}%
    \put(0,0){\includegraphics[width=\unitlength,page=4]{z.pdf}}%
  \end{picture}%
\endgroup%
}.
     \end{split}
    \end{gather*}
    Note that $c'_{\mu}=\sum_{\rho\in R_{s_0},\tau\in C_{s_0}} \sgn(\tau)\rho\tau$, where $s_0$ is the canonical $\mu$-tableau.

    Firstly, we consider the case where $\mu_i,\mu_j,\mu_k$ are distinct.
    Then $z$ appears only when $\tau$ is an even permutation which fixes $i'$, $j'$ and $k'$.
    Therefore, the coefficient of $z$ in $h(Z_{\mu})$ is positive.
    Moreover, the linear sum of terms in $Z_{\mu}$ such that $\ast_2$ is connected to $v_k$ lies in $\pi_{\nu}(Z_{\mu})$,
    so the coefficient of $z$ in $h(\pi_{\nu}(Z_{\mu}))$ is equal to that of $z$ in $h(Z_{\mu})$, which is nonzero.

    The other cases, where at least two of $\mu_i,\mu_j$ and $\mu_k$ are equal, follow in a similar argument.
    The only thing that differs from the above case is that $z$ appears when $\tau$ preserves the subset $\{i',j',k'\}\subset [2d-1]$ and the parity of $\tau$ coincides with that of the restriction of $\tau$ to $\{i',j',k'\}$.
    Since we have the AS relation, the sign due to the permutation of $\{i',j',k'\}$ is cancelled.
    Therefore, the coefficient of $z$ in $h(Z_{\mu})$ is positive in any case.
    The proof is complete.
  \end{proof}

  \begin{theorem}\label{propindecomp}
   Let $d\geq 2$. The direct decomposition $$A_d(n)=A_d P(n)\oplus A_d Q(n)$$ of $\Aut(F_n)$-modules is indecomposable for $n\geq 2d$.
  \end{theorem}
  \begin{proof}
   By Lemma \ref{l81}, it suffices to show that $A_d Q(n)$ is indecomposable.
   Since the radical preserves the direct sum, we have only to show that $A_d Q(n)/\Rad^2(A_d Q(n))$ is indecomposable.
   Suppose that we have a nontrivial decomposition of $\Aut(F_n)$-modules
   \begin{gather*}
    \begin{split}
      A_d Q(n)/\Rad^2(A_d Q(n))&=A_d Q(n)/A_{d,2}(n)\\
      &=(M_1+A_{d,2}(n))/ A_{d,2}(n)\oplus (M_2+A_{d,2}(n))/A_{d,2}(n),
    \end{split}
   \end{gather*}
   where $M_i$ is an $\Aut(F_n)$-submodule of $A_d Q(n)$ for $i=1,2$.
   Let $$N_i=\theta_{d,n}(M_i+A_{d,2}(n))/\theta_{d,n}(A_{d,2}(n))$$ for $i=1,2$.
   We have
   $$N_1\oplus N_2=\theta_{d,n}(A_d Q(n))/\theta_{d,n}(A_{d,2}(n))=\left(\bigoplus_{\lambda\vdash d, \lambda\neq (d)} B_{d,0}(n)_{2\lambda}\right) \oplus B_{d,1}(n).$$
   For any $2\lambda\in X_d$, there uniquely exists $i\in \{1,2\}$ such that
   $N_i$ includes a $\GL(n;\Z)$-submodule $(N_i)_{2\lambda}\cong V_{2\lambda}$.
   Let $x\in (N_i)_{2\lambda}$ be a generator of the irreducible $\GL(n;\Z)$-module $(N_i)_{2\lambda}$.
   Then, the image $x'$ of $x$ under the composition of $\GL(n;\Z)$-module maps
   $$(N_i)_{2\lambda}\hookrightarrow N_i \hookrightarrow B_{d,0}(n)\oplus B_{d,1}(n)\twoheadrightarrow B_{d,0}(n)$$
   is an element of $B_{d,0}(n)_{2\lambda}$.
   For any $\mu\in Y_d$ that is connected to $2\lambda$ by an edge of $G_d$, by Proposition \ref{conj}, there exists $g\in \opegr^1(\IA(n))$ such that $[x',g]\neq 0\in B_{d,1}(n)_{\mu}$.
   Therefore, we have
   $$[x,g]=[x',g]+[x-x',g]=[x',g]\neq 0\in B_{d,1}(n)_{\mu}.$$
   It follows that $N_i$ includes a $\GL(n;\Z)$-submodule $(N_i)_{\mu}$ that is isomorphic to $V_{\mu}$ for any $\mu \in Y_d$ that is connected to $2\lambda$ by an edge of $G_d$.
   Hence, by Proposition \ref{p825}, we have $N_1\cap N_2\neq \{0\}$, a contradiction.
   Therefore, $A_d Q(n)$ is indecomposable.
  \end{proof}

  Note that the assumption $n\geq 2d$ is needed for the surjectivity of the bracket map and the nontriviality of the bracket map for each pair of nonzero irreducible $\GL(n;\Z)$-submodules. Thus, if we have the surjectivity and the nontriviality of the bracket map for some $n<2d$, we can loose the assumption.

 \subsection{The $\Aut(F_n)$-module structure of $A_3(n)$}\label{ss83}
  Here, we consider the $\Aut(F_n)$-module structure of $A_3(n)$ in detail.

  In degree $3$, the restrictions of the bracket map to each isotypic component induce $\GL(n;\Z)$-module homomorphisms
  \begin{gather*}
   \begin{split}
     \rho_1:B_{3,0}(n)_{(4,2)}&\rightarrow \Hom(\opegr^1(\IA(n)), B_{3,1}(n)_{(3,1^2)}),\\
     \rho_2:B_{3,0}(n)_{(2^3)}&\rightarrow \Hom(\opegr^1(\IA(n)), B_{3,1}(n)_{(3,1^2)}),\\
     \rho_3:B_{3,0}(n)_{(2^3)}&\rightarrow \Hom(\opegr^1(\IA(n)), B_{3,1}(n)_{(2,1^3)}).
   \end{split}
  \end{gather*}

  \begin{proposition}\label{propd3}
   The $\GL(n;\Z)$-module homomorphisms $\rho_1$ and $\rho_2$ are injective for $n\geq 3$ and $\rho_3$ for $n\geq 4$.
  \end{proposition}
  \begin{proof}
    Recall that $c_{\lambda}$ denotes the Young symmetrizer defined in \eqref{youngsym} and that $K_{i,j,k}\in\IA(n)$ is defined by \eqref{Kijk}.
    For $n\geq 3$, we have
    $$\rho_1(u)(K_{3,2,1})=[u,K_{3,2,1}]=-10 w\neq 0 \in B_{3,1}(n)_{(3,1^2)},$$
    where
    $$u= \frac{1}{64}\centre{
\begingroup%
  \makeatletter%
  \providecommand\color[2][]{%
    \errmessage{(Inkscape) Color is used for the text in Inkscape, but the package 'color.sty' is not loaded}%
    \renewcommand\color[2][]{}%
  }%
  \providecommand\transparent[1]{%
    \errmessage{(Inkscape) Transparency is used (non-zero) for the text in Inkscape, but the package 'transparent.sty' is not loaded}%
    \renewcommand\transparent[1]{}%
  }%
  \providecommand\rotatebox[2]{#2}%
  \newcommand*\fsize{\dimexpr\f@size pt\relax}%
  \newcommand*\lineheight[1]{\fontsize{\fsize}{#1\fsize}\selectfont}%
  \ifx\svgwidth\undefined%
    \setlength{\unitlength}{56.25118034bp}%
    \ifx\svgscale\undefined%
      \relax%
    \else%
      \setlength{\unitlength}{\unitlength * \real{\svgscale}}%
    \fi%
  \else%
    \setlength{\unitlength}{\svgwidth}%
  \fi%
  \global\let\svgwidth\undefined%
  \global\let\svgscale\undefined%
  \makeatother%
  \begin{picture}(1,0.61186316)%
    \lineheight{1}%
    \setlength\tabcolsep{0pt}%
    \put(0,0){\includegraphics[width=\unitlength,page=1]{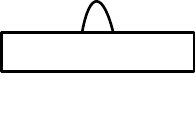}}%
    \put(0.35176061,0.31609032){\makebox(0,0)[lt]{\lineheight{1.45000005}\smash{\begin{tabular}[t]{l}$c_{(4,2)}$\end{tabular}}}}%
    \put(0,0){\includegraphics[width=\unitlength,page=2]{rho1u3.pdf}}%
    \put(0.03826003,0.00626723){\makebox(0,0)[lt]{\lineheight{1.45000005}\smash{\begin{tabular}[t]{l}$v_1$\end{tabular}}}}%
    \put(0.20244032,0.00854682){\makebox(0,0)[lt]{\lineheight{1.45000005}\smash{\begin{tabular}[t]{l}$v_1$\end{tabular}}}}%
    \put(0.36843795,0.00878529){\makebox(0,0)[lt]{\lineheight{1.45000005}\smash{\begin{tabular}[t]{l}$v_1$\end{tabular}}}}%
    \put(0.52619889,0.00646523){\makebox(0,0)[lt]{\lineheight{1.45000005}\smash{\begin{tabular}[t]{l}$v_1$\end{tabular}}}}%
    \put(0,0){\includegraphics[width=\unitlength,page=3]{rho1u3.pdf}}%
    \put(0.68209209,0.0077282){\makebox(0,0)[lt]{\lineheight{1.45000005}\smash{\begin{tabular}[t]{l}$v_2$\end{tabular}}}}%
    \put(0.84375766,0.01000779){\makebox(0,0)[lt]{\lineheight{1.45000005}\smash{\begin{tabular}[t]{l}$v_2$\end{tabular}}}}%
    \put(0,0){\includegraphics[width=\unitlength,page=4]{rho1u3.pdf}}%
  \end{picture}%
\endgroup%
}= \centre{
\begingroup%
  \makeatletter%
  \providecommand\color[2][]{%
    \errmessage{(Inkscape) Color is used for the text in Inkscape, but the package 'color.sty' is not loaded}%
    \renewcommand\color[2][]{}%
  }%
  \providecommand\transparent[1]{%
    \errmessage{(Inkscape) Transparency is used (non-zero) for the text in Inkscape, but the package 'transparent.sty' is not loaded}%
    \renewcommand\transparent[1]{}%
  }%
  \providecommand\rotatebox[2]{#2}%
  \newcommand*\fsize{\dimexpr\f@size pt\relax}%
  \newcommand*\lineheight[1]{\fontsize{\fsize}{#1\fsize}\selectfont}%
  \ifx\svgwidth\undefined%
    \setlength{\unitlength}{52.38792046bp}%
    \ifx\svgscale\undefined%
      \relax%
    \else%
      \setlength{\unitlength}{\unitlength * \real{\svgscale}}%
    \fi%
  \else%
    \setlength{\unitlength}{\svgwidth}%
  \fi%
  \global\let\svgwidth\undefined%
  \global\let\svgscale\undefined%
  \makeatother%
  \begin{picture}(1,0.30363651)%
    \lineheight{1}%
    \setlength\tabcolsep{0pt}%
    \put(0,0){\includegraphics[width=\unitlength,page=1]{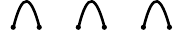}}%
    \put(-0.00288935,0.0067294){\makebox(0,0)[lt]{\lineheight{1.45000005}\smash{\begin{tabular}[t]{l}$v_1$\end{tabular}}}}%
    \put(0.16116657,0.00684365){\makebox(0,0)[lt]{\lineheight{1.45000005}\smash{\begin{tabular}[t]{l}$v_1$\end{tabular}}}}%
    \put(0.88285997,0.0093188){\makebox(0,0)[lt]{\lineheight{1.45000005}\smash{\begin{tabular}[t]{l}$v_2$\end{tabular}}}}%
    \put(0.35869189,0.00956298){\makebox(0,0)[lt]{\lineheight{1.45000005}\smash{\begin{tabular}[t]{l}$v_1$\end{tabular}}}}%
    \put(0.71665635,0.01027122){\makebox(0,0)[lt]{\lineheight{1.45000005}\smash{\begin{tabular}[t]{l}$v_2$\end{tabular}}}}%
    \put(0.52943186,0.00935714){\makebox(0,0)[lt]{\lineheight{1.45000005}\smash{\begin{tabular}[t]{l}$v_1$\end{tabular}}}}%
  \end{picture}%
\endgroup%
}- \centre{
\begingroup%
  \makeatletter%
  \providecommand\color[2][]{%
    \errmessage{(Inkscape) Color is used for the text in Inkscape, but the package 'color.sty' is not loaded}%
    \renewcommand\color[2][]{}%
  }%
  \providecommand\transparent[1]{%
    \errmessage{(Inkscape) Transparency is used (non-zero) for the text in Inkscape, but the package 'transparent.sty' is not loaded}%
    \renewcommand\transparent[1]{}%
  }%
  \providecommand\rotatebox[2]{#2}%
  \newcommand*\fsize{\dimexpr\f@size pt\relax}%
  \newcommand*\lineheight[1]{\fontsize{\fsize}{#1\fsize}\selectfont}%
  \ifx\svgwidth\undefined%
    \setlength{\unitlength}{52.38792046bp}%
    \ifx\svgscale\undefined%
      \relax%
    \else%
      \setlength{\unitlength}{\unitlength * \real{\svgscale}}%
    \fi%
  \else%
    \setlength{\unitlength}{\svgwidth}%
  \fi%
  \global\let\svgwidth\undefined%
  \global\let\svgscale\undefined%
  \makeatother%
  \begin{picture}(1,0.30363651)%
    \lineheight{1}%
    \setlength\tabcolsep{0pt}%
    \put(0,0){\includegraphics[width=\unitlength,page=1]{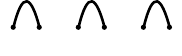}}%
    \put(-0.00288935,0.0067294){\makebox(0,0)[lt]{\lineheight{1.45000005}\smash{\begin{tabular}[t]{l}$v_1$\end{tabular}}}}%
    \put(0.16116657,0.00684365){\makebox(0,0)[lt]{\lineheight{1.45000005}\smash{\begin{tabular}[t]{l}$v_1$\end{tabular}}}}%
    \put(0.88285997,0.0093188){\makebox(0,0)[lt]{\lineheight{1.45000005}\smash{\begin{tabular}[t]{l}$v_2$\end{tabular}}}}%
    \put(0.35869189,0.00956298){\makebox(0,0)[lt]{\lineheight{1.45000005}\smash{\begin{tabular}[t]{l}$v_1$\end{tabular}}}}%
    \put(0.71665635,0.01027122){\makebox(0,0)[lt]{\lineheight{1.45000005}\smash{\begin{tabular}[t]{l}$v_1$\end{tabular}}}}%
    \put(0.52943186,0.00935714){\makebox(0,0)[lt]{\lineheight{1.45000005}\smash{\begin{tabular}[t]{l}$v_2$\end{tabular}}}}%
  \end{picture}%
\endgroup%
} \in B_{3,0}(n)_{(4,2)}$$
    and
    $$w= \frac{1}{20}\centre{
\begingroup%
  \makeatletter%
  \providecommand\color[2][]{%
    \errmessage{(Inkscape) Color is used for the text in Inkscape, but the package 'color.sty' is not loaded}%
    \renewcommand\color[2][]{}%
  }%
  \providecommand\transparent[1]{%
    \errmessage{(Inkscape) Transparency is used (non-zero) for the text in Inkscape, but the package 'transparent.sty' is not loaded}%
    \renewcommand\transparent[1]{}%
  }%
  \providecommand\rotatebox[2]{#2}%
  \newcommand*\fsize{\dimexpr\f@size pt\relax}%
  \newcommand*\lineheight[1]{\fontsize{\fsize}{#1\fsize}\selectfont}%
  \ifx\svgwidth\undefined%
    \setlength{\unitlength}{56.25118034bp}%
    \ifx\svgscale\undefined%
      \relax%
    \else%
      \setlength{\unitlength}{\unitlength * \real{\svgscale}}%
    \fi%
  \else%
    \setlength{\unitlength}{\svgwidth}%
  \fi%
  \global\let\svgwidth\undefined%
  \global\let\svgscale\undefined%
  \makeatother%
  \begin{picture}(1,0.6539475)%
    \lineheight{1}%
    \setlength\tabcolsep{0pt}%
    \put(0,0){\includegraphics[width=\unitlength,page=1]{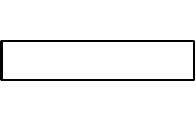}}%
    \put(0.29944908,0.31428646){\makebox(0,0)[lt]{\lineheight{1.45000005}\smash{\begin{tabular}[t]{l}$c_{(3,1^2)}$\end{tabular}}}}%
    \put(0,0){\includegraphics[width=\unitlength,page=2]{rho1v2.pdf}}%
    \put(0.03826003,0.00626723){\makebox(0,0)[lt]{\lineheight{1.45000005}\smash{\begin{tabular}[t]{l}$v_1$\end{tabular}}}}%
    \put(0.20244032,0.00854682){\makebox(0,0)[lt]{\lineheight{1.45000005}\smash{\begin{tabular}[t]{l}$v_1$\end{tabular}}}}%
    \put(0.36843795,0.00878529){\makebox(0,0)[lt]{\lineheight{1.45000005}\smash{\begin{tabular}[t]{l}$v_1$\end{tabular}}}}%
    \put(0.56407963,0.00646523){\makebox(0,0)[lt]{\lineheight{1.45000005}\smash{\begin{tabular}[t]{l}$v_2$\end{tabular}}}}%
    \put(0,0){\includegraphics[width=\unitlength,page=3]{rho1v2.pdf}}%
    \put(0.76146131,0.0077282){\makebox(0,0)[lt]{\lineheight{1.45000005}\smash{\begin{tabular}[t]{l}$v_3$\end{tabular}}}}%
    \put(0,0){\includegraphics[width=\unitlength,page=4]{rho1v2.pdf}}%
  \end{picture}%
\endgroup%
}= \centre{
\begingroup%
  \makeatletter%
  \providecommand\color[2][]{%
    \errmessage{(Inkscape) Color is used for the text in Inkscape, but the package 'color.sty' is not loaded}%
    \renewcommand\color[2][]{}%
  }%
  \providecommand\transparent[1]{%
    \errmessage{(Inkscape) Transparency is used (non-zero) for the text in Inkscape, but the package 'transparent.sty' is not loaded}%
    \renewcommand\transparent[1]{}%
  }%
  \providecommand\rotatebox[2]{#2}%
  \newcommand*\fsize{\dimexpr\f@size pt\relax}%
  \newcommand*\lineheight[1]{\fontsize{\fsize}{#1\fsize}\selectfont}%
  \ifx\svgwidth\undefined%
    \setlength{\unitlength}{47.88484651bp}%
    \ifx\svgscale\undefined%
      \relax%
    \else%
      \setlength{\unitlength}{\unitlength * \real{\svgscale}}%
    \fi%
  \else%
    \setlength{\unitlength}{\svgwidth}%
  \fi%
  \global\let\svgwidth\undefined%
  \global\let\svgscale\undefined%
  \makeatother%
  \begin{picture}(1,0.42019778)%
    \lineheight{1}%
    \setlength\tabcolsep{0pt}%
    \put(0,0){\includegraphics[width=\unitlength,page=1]{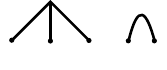}}%
    \put(0.87184417,0.00816378){\makebox(0,0)[lt]{\lineheight{1.45000005}\smash{\begin{tabular}[t]{l}$v_1$\end{tabular}}}}%
    \put(0.4628079,0.00862174){\makebox(0,0)[lt]{\lineheight{1.45000005}\smash{\begin{tabular}[t]{l}$v_3$\end{tabular}}}}%
    \put(-0.00316107,0.00939551){\makebox(0,0)[lt]{\lineheight{1.45000005}\smash{\begin{tabular}[t]{l}$v_1$\end{tabular}}}}%
    \put(0.69283711,0.00736223){\makebox(0,0)[lt]{\lineheight{1.45000005}\smash{\begin{tabular}[t]{l}$v_1$\end{tabular}}}}%
    \put(0.23052756,0.00764354){\makebox(0,0)[lt]{\lineheight{1.45000005}\smash{\begin{tabular}[t]{l}$v_2$\end{tabular}}}}%
  \end{picture}%
\endgroup%
}\neq 0 \in B_{3,1}(n)_{(3,1^2)}.$$
    Thus, we have $\rho_1\neq 0$ for $n\geq 3$. Since $B_{3,0}(n)_{(4,2)}$ is irreducible, $\rho_1$ is injective.

    Let $$x= \frac{1}{48}\centre{
\begingroup%
  \makeatletter%
  \providecommand\color[2][]{%
    \errmessage{(Inkscape) Color is used for the text in Inkscape, but the package 'color.sty' is not loaded}%
    \renewcommand\color[2][]{}%
  }%
  \providecommand\transparent[1]{%
    \errmessage{(Inkscape) Transparency is used (non-zero) for the text in Inkscape, but the package 'transparent.sty' is not loaded}%
    \renewcommand\transparent[1]{}%
  }%
  \providecommand\rotatebox[2]{#2}%
  \newcommand*\fsize{\dimexpr\f@size pt\relax}%
  \newcommand*\lineheight[1]{\fontsize{\fsize}{#1\fsize}\selectfont}%
  \ifx\svgwidth\undefined%
    \setlength{\unitlength}{56.25118034bp}%
    \ifx\svgscale\undefined%
      \relax%
    \else%
      \setlength{\unitlength}{\unitlength * \real{\svgscale}}%
    \fi%
  \else%
    \setlength{\unitlength}{\svgwidth}%
  \fi%
  \global\let\svgwidth\undefined%
  \global\let\svgscale\undefined%
  \makeatother%
  \begin{picture}(1,0.61186316)%
    \lineheight{1}%
    \setlength\tabcolsep{0pt}%
    \put(0,0){\includegraphics[width=\unitlength,page=1]{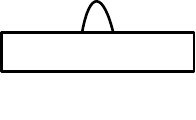}}%
    \put(0.35176061,0.31609032){\makebox(0,0)[lt]{\lineheight{1.45000005}\smash{\begin{tabular}[t]{l}$c_{(2^3)}$\end{tabular}}}}%
    \put(0,0){\includegraphics[width=\unitlength,page=2]{rhox.pdf}}%
    \put(0.03826003,0.00626723){\makebox(0,0)[lt]{\lineheight{1.45000005}\smash{\begin{tabular}[t]{l}$v_1$\end{tabular}}}}%
    \put(0.20244032,0.00854682){\makebox(0,0)[lt]{\lineheight{1.45000005}\smash{\begin{tabular}[t]{l}$v_1$\end{tabular}}}}%
    \put(0.36843795,0.00878529){\makebox(0,0)[lt]{\lineheight{1.45000005}\smash{\begin{tabular}[t]{l}$v_2$\end{tabular}}}}%
    \put(0.52619889,0.00646523){\makebox(0,0)[lt]{\lineheight{1.45000005}\smash{\begin{tabular}[t]{l}$v_2$\end{tabular}}}}%
    \put(0,0){\includegraphics[width=\unitlength,page=3]{rhox.pdf}}%
    \put(0.68209209,0.0077282){\makebox(0,0)[lt]{\lineheight{1.45000005}\smash{\begin{tabular}[t]{l}$v_3$\end{tabular}}}}%
    \put(0.84375766,0.01000779){\makebox(0,0)[lt]{\lineheight{1.45000005}\smash{\begin{tabular}[t]{l}$v_3$\end{tabular}}}}%
    \put(0,0){\includegraphics[width=\unitlength,page=4]{rhox.pdf}}%
  \end{picture}%
\endgroup%
}=\centre{
\begingroup%
  \makeatletter%
  \providecommand\color[2][]{%
    \errmessage{(Inkscape) Color is used for the text in Inkscape, but the package 'color.sty' is not loaded}%
    \renewcommand\color[2][]{}%
  }%
  \providecommand\transparent[1]{%
    \errmessage{(Inkscape) Transparency is used (non-zero) for the text in Inkscape, but the package 'transparent.sty' is not loaded}%
    \renewcommand\transparent[1]{}%
  }%
  \providecommand\rotatebox[2]{#2}%
  \newcommand*\fsize{\dimexpr\f@size pt\relax}%
  \newcommand*\lineheight[1]{\fontsize{\fsize}{#1\fsize}\selectfont}%
  \ifx\svgwidth\undefined%
    \setlength{\unitlength}{63.65976125bp}%
    \ifx\svgscale\undefined%
      \relax%
    \else%
      \setlength{\unitlength}{\unitlength * \real{\svgscale}}%
    \fi%
  \else%
    \setlength{\unitlength}{\svgwidth}%
  \fi%
  \global\let\svgwidth\undefined%
  \global\let\svgscale\undefined%
  \makeatother%
  \begin{picture}(1,0.65880512)%
    \lineheight{1}%
    \setlength\tabcolsep{0pt}%
    \put(0,0){\includegraphics[width=\unitlength,page=1]{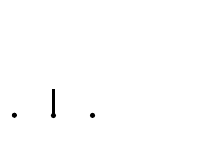}}%
    \put(0.36147511,0.0125764){\makebox(0,0)[lt]{\lineheight{1.45000005}\smash{\begin{tabular}[t]{l}$v_3$\end{tabular}}}}%
    \put(0.01097346,0.01315843){\makebox(0,0)[lt]{\lineheight{1.45000005}\smash{\begin{tabular}[t]{l}$v_1$\end{tabular}}}}%
    \put(0.18675393,0.01184056){\makebox(0,0)[lt]{\lineheight{1.45000005}\smash{\begin{tabular}[t]{l}$v_2$\end{tabular}}}}%
    \put(0.90360123,0.0062737){\makebox(0,0)[lt]{\lineheight{1.45000005}\smash{\begin{tabular}[t]{l}$v_3$\end{tabular}}}}%
    \put(0.55309954,0.00685573){\makebox(0,0)[lt]{\lineheight{1.45000005}\smash{\begin{tabular}[t]{l}$v_1$\end{tabular}}}}%
    \put(0.72888002,0.00553787){\makebox(0,0)[lt]{\lineheight{1.45000005}\smash{\begin{tabular}[t]{l}$v_2$\end{tabular}}}}%
    \put(0,0){\includegraphics[width=\unitlength,page=2]{rho23x.pdf}}%
    \put(0.11588898,0.29902858){\makebox(0,0)[lt]{\lineheight{1.45000005}\smash{\begin{tabular}[t]{l}$\alt_{3}$\end{tabular}}}}%
    \put(0,0){\includegraphics[width=\unitlength,page=3]{rho23x.pdf}}%
  \end{picture}%
\endgroup%
}\in B_{3,0}(n)_{(2^3)}.$$
    We have
    $$\rho_2(x)(K_{1,3,2})=[x,K_{1,3,2}]= -6 w\neq 0 \in  B_{3,1}(n)_{(3,1^2)}.$$
    Thus, we have $\rho_2\neq 0$ for $n\geq 3$. Since $B_{3,0}(n)_{(2^3)}$ is irreducible, $\rho_2$ is injective.

    For $n\geq 4$, we have
    $$[x,K_{4,3,2}]=-\frac{6}{5}y-\frac{24}{5}z$$
    and thus
    $$\rho_3(x)(K_{4,3,2})=-\frac{24}{5}z\neq 0 \in B_{3,1}(n)_{(2,1^3)},$$
    where
    \begin{gather*}
     \begin{split}
       y= \frac{1}{4}\centre{
\begingroup%
  \makeatletter%
  \providecommand\color[2][]{%
    \errmessage{(Inkscape) Color is used for the text in Inkscape, but the package 'color.sty' is not loaded}%
    \renewcommand\color[2][]{}%
  }%
  \providecommand\transparent[1]{%
    \errmessage{(Inkscape) Transparency is used (non-zero) for the text in Inkscape, but the package 'transparent.sty' is not loaded}%
    \renewcommand\transparent[1]{}%
  }%
  \providecommand\rotatebox[2]{#2}%
  \newcommand*\fsize{\dimexpr\f@size pt\relax}%
  \newcommand*\lineheight[1]{\fontsize{\fsize}{#1\fsize}\selectfont}%
  \ifx\svgwidth\undefined%
    \setlength{\unitlength}{56.25118034bp}%
    \ifx\svgscale\undefined%
      \relax%
    \else%
      \setlength{\unitlength}{\unitlength * \real{\svgscale}}%
    \fi%
  \else%
    \setlength{\unitlength}{\svgwidth}%
  \fi%
  \global\let\svgwidth\undefined%
  \global\let\svgscale\undefined%
  \makeatother%
  \begin{picture}(1,0.6539475)%
    \lineheight{1}%
    \setlength\tabcolsep{0pt}%
    \put(0,0){\includegraphics[width=\unitlength,page=1]{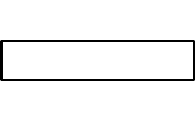}}%
    \put(0.29944908,0.31428646){\makebox(0,0)[lt]{\lineheight{1.45000005}\smash{\begin{tabular}[t]{l}$c_{(3,1^2)}$\end{tabular}}}}%
    \put(0,0){\includegraphics[width=\unitlength,page=2]{rhoy.pdf}}%
    \put(0.03826003,0.00626723){\makebox(0,0)[lt]{\lineheight{1.45000005}\smash{\begin{tabular}[t]{l}$v_1$\end{tabular}}}}%
    \put(0.20244032,0.00854682){\makebox(0,0)[lt]{\lineheight{1.45000005}\smash{\begin{tabular}[t]{l}$v_1$\end{tabular}}}}%
    \put(0.36843795,0.00878529){\makebox(0,0)[lt]{\lineheight{1.45000005}\smash{\begin{tabular}[t]{l}$v_4$\end{tabular}}}}%
    \put(0.56407963,0.00646523){\makebox(0,0)[lt]{\lineheight{1.45000005}\smash{\begin{tabular}[t]{l}$v_2$\end{tabular}}}}%
    \put(0,0){\includegraphics[width=\unitlength,page=3]{rhoy.pdf}}%
    \put(0.76146131,0.0077282){\makebox(0,0)[lt]{\lineheight{1.45000005}\smash{\begin{tabular}[t]{l}$v_3$\end{tabular}}}}%
    \put(0,0){\includegraphics[width=\unitlength,page=4]{rhoy.pdf}}%
  \end{picture}%
\endgroup%
}=\centre{
\begingroup%
  \makeatletter%
  \providecommand\color[2][]{%
    \errmessage{(Inkscape) Color is used for the text in Inkscape, but the package 'color.sty' is not loaded}%
    \renewcommand\color[2][]{}%
  }%
  \providecommand\transparent[1]{%
    \errmessage{(Inkscape) Transparency is used (non-zero) for the text in Inkscape, but the package 'transparent.sty' is not loaded}%
    \renewcommand\transparent[1]{}%
  }%
  \providecommand\rotatebox[2]{#2}%
  \newcommand*\fsize{\dimexpr\f@size pt\relax}%
  \newcommand*\lineheight[1]{\fontsize{\fsize}{#1\fsize}\selectfont}%
  \ifx\svgwidth\undefined%
    \setlength{\unitlength}{47.75353501bp}%
    \ifx\svgscale\undefined%
      \relax%
    \else%
      \setlength{\unitlength}{\unitlength * \real{\svgscale}}%
    \fi%
  \else%
    \setlength{\unitlength}{\svgwidth}%
  \fi%
  \global\let\svgwidth\undefined%
  \global\let\svgscale\undefined%
  \makeatother%
  \begin{picture}(1,0.42107114)%
    \lineheight{1}%
    \setlength\tabcolsep{0pt}%
    \put(0,0){\includegraphics[width=\unitlength,page=1]{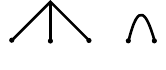}}%
    \put(0.87149177,0.00832418){\makebox(0,0)[lt]{\lineheight{1.45000005}\smash{\begin{tabular}[t]{l}$v_1$\end{tabular}}}}%
    \put(0.46408052,0.00836336){\makebox(0,0)[lt]{\lineheight{1.45000005}\smash{\begin{tabular}[t]{l}$v_4$\end{tabular}}}}%
    \put(-0.00316976,0.00913925){\makebox(0,0)[lt]{\lineheight{1.45000005}\smash{\begin{tabular}[t]{l}$v_2$\end{tabular}}}}%
    \put(0.69583652,0.00825546){\makebox(0,0)[lt]{\lineheight{1.45000005}\smash{\begin{tabular}[t]{l}$v_1$\end{tabular}}}}%
    \put(0.23116146,0.00738247){\makebox(0,0)[lt]{\lineheight{1.45000005}\smash{\begin{tabular}[t]{l}$v_3$\end{tabular}}}}%
  \end{picture}%
\endgroup%
}&-\centre{
\begingroup%
  \makeatletter%
  \providecommand\color[2][]{%
    \errmessage{(Inkscape) Color is used for the text in Inkscape, but the package 'color.sty' is not loaded}%
    \renewcommand\color[2][]{}%
  }%
  \providecommand\transparent[1]{%
    \errmessage{(Inkscape) Transparency is used (non-zero) for the text in Inkscape, but the package 'transparent.sty' is not loaded}%
    \renewcommand\transparent[1]{}%
  }%
  \providecommand\rotatebox[2]{#2}%
  \newcommand*\fsize{\dimexpr\f@size pt\relax}%
  \newcommand*\lineheight[1]{\fontsize{\fsize}{#1\fsize}\selectfont}%
  \ifx\svgwidth\undefined%
    \setlength{\unitlength}{47.75353501bp}%
    \ifx\svgscale\undefined%
      \relax%
    \else%
      \setlength{\unitlength}{\unitlength * \real{\svgscale}}%
    \fi%
  \else%
    \setlength{\unitlength}{\svgwidth}%
  \fi%
  \global\let\svgwidth\undefined%
  \global\let\svgscale\undefined%
  \makeatother%
  \begin{picture}(1,0.42107114)%
    \lineheight{1}%
    \setlength\tabcolsep{0pt}%
    \put(0,0){\includegraphics[width=\unitlength,page=1]{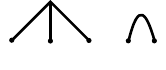}}%
    \put(0.87149177,0.00832418){\makebox(0,0)[lt]{\lineheight{1.45000005}\smash{\begin{tabular}[t]{l}$v_2$\end{tabular}}}}%
    \put(0.46408052,0.00836336){\makebox(0,0)[lt]{\lineheight{1.45000005}\smash{\begin{tabular}[t]{l}$v_4$\end{tabular}}}}%
    \put(-0.00316976,0.00913925){\makebox(0,0)[lt]{\lineheight{1.45000005}\smash{\begin{tabular}[t]{l}$v_1$\end{tabular}}}}%
    \put(0.69583652,0.00825546){\makebox(0,0)[lt]{\lineheight{1.45000005}\smash{\begin{tabular}[t]{l}$v_1$\end{tabular}}}}%
    \put(0.23116146,0.00738247){\makebox(0,0)[lt]{\lineheight{1.45000005}\smash{\begin{tabular}[t]{l}$v_3$\end{tabular}}}}%
  \end{picture}%
\endgroup%
} +\centre{
\begingroup%
  \makeatletter%
  \providecommand\color[2][]{%
    \errmessage{(Inkscape) Color is used for the text in Inkscape, but the package 'color.sty' is not loaded}%
    \renewcommand\color[2][]{}%
  }%
  \providecommand\transparent[1]{%
    \errmessage{(Inkscape) Transparency is used (non-zero) for the text in Inkscape, but the package 'transparent.sty' is not loaded}%
    \renewcommand\transparent[1]{}%
  }%
  \providecommand\rotatebox[2]{#2}%
  \newcommand*\fsize{\dimexpr\f@size pt\relax}%
  \newcommand*\lineheight[1]{\fontsize{\fsize}{#1\fsize}\selectfont}%
  \ifx\svgwidth\undefined%
    \setlength{\unitlength}{47.75353501bp}%
    \ifx\svgscale\undefined%
      \relax%
    \else%
      \setlength{\unitlength}{\unitlength * \real{\svgscale}}%
    \fi%
  \else%
    \setlength{\unitlength}{\svgwidth}%
  \fi%
  \global\let\svgwidth\undefined%
  \global\let\svgscale\undefined%
  \makeatother%
  \begin{picture}(1,0.42107114)%
    \lineheight{1}%
    \setlength\tabcolsep{0pt}%
    \put(0,0){\includegraphics[width=\unitlength,page=1]{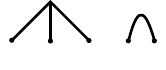}}%
    \put(0.87149177,0.00832418){\makebox(0,0)[lt]{\lineheight{1.45000005}\smash{\begin{tabular}[t]{l}$v_3$\end{tabular}}}}%
    \put(0.46408052,0.00836336){\makebox(0,0)[lt]{\lineheight{1.45000005}\smash{\begin{tabular}[t]{l}$v_4$\end{tabular}}}}%
    \put(-0.00316976,0.00913925){\makebox(0,0)[lt]{\lineheight{1.45000005}\smash{\begin{tabular}[t]{l}$v_1$\end{tabular}}}}%
    \put(0.69583652,0.00825546){\makebox(0,0)[lt]{\lineheight{1.45000005}\smash{\begin{tabular}[t]{l}$v_1$\end{tabular}}}}%
    \put(0.23116146,0.00738247){\makebox(0,0)[lt]{\lineheight{1.45000005}\smash{\begin{tabular}[t]{l}$v_2$\end{tabular}}}}%
  \end{picture}%
\endgroup%
}\\
       &+4\centre{
\begingroup%
  \makeatletter%
  \providecommand\color[2][]{%
    \errmessage{(Inkscape) Color is used for the text in Inkscape, but the package 'color.sty' is not loaded}%
    \renewcommand\color[2][]{}%
  }%
  \providecommand\transparent[1]{%
    \errmessage{(Inkscape) Transparency is used (non-zero) for the text in Inkscape, but the package 'transparent.sty' is not loaded}%
    \renewcommand\transparent[1]{}%
  }%
  \providecommand\rotatebox[2]{#2}%
  \newcommand*\fsize{\dimexpr\f@size pt\relax}%
  \newcommand*\lineheight[1]{\fontsize{\fsize}{#1\fsize}\selectfont}%
  \ifx\svgwidth\undefined%
    \setlength{\unitlength}{47.75353501bp}%
    \ifx\svgscale\undefined%
      \relax%
    \else%
      \setlength{\unitlength}{\unitlength * \real{\svgscale}}%
    \fi%
  \else%
    \setlength{\unitlength}{\svgwidth}%
  \fi%
  \global\let\svgwidth\undefined%
  \global\let\svgscale\undefined%
  \makeatother%
  \begin{picture}(1,0.42107114)%
    \lineheight{1}%
    \setlength\tabcolsep{0pt}%
    \put(0,0){\includegraphics[width=\unitlength,page=1]{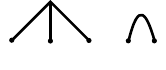}}%
    \put(0.87149177,0.00832418){\makebox(0,0)[lt]{\lineheight{1.45000005}\smash{\begin{tabular}[t]{l}$v_4$\end{tabular}}}}%
    \put(0.46408052,0.00836336){\makebox(0,0)[lt]{\lineheight{1.45000005}\smash{\begin{tabular}[t]{l}$v_3$\end{tabular}}}}%
    \put(-0.00316976,0.00913925){\makebox(0,0)[lt]{\lineheight{1.45000005}\smash{\begin{tabular}[t]{l}$v_1$\end{tabular}}}}%
    \put(0.69583652,0.00825546){\makebox(0,0)[lt]{\lineheight{1.45000005}\smash{\begin{tabular}[t]{l}$v_1$\end{tabular}}}}%
    \put(0.23116146,0.00738247){\makebox(0,0)[lt]{\lineheight{1.45000005}\smash{\begin{tabular}[t]{l}$v_2$\end{tabular}}}}%
  \end{picture}%
\endgroup%
}
       \in B_{3,1}(n)_{(3,1^2)}
     \end{split}
    \end{gather*}
    and
    \begin{gather*}
     \begin{split}
       z= \frac{1}{12}\centre{
\begingroup%
  \makeatletter%
  \providecommand\color[2][]{%
    \errmessage{(Inkscape) Color is used for the text in Inkscape, but the package 'color.sty' is not loaded}%
    \renewcommand\color[2][]{}%
  }%
  \providecommand\transparent[1]{%
    \errmessage{(Inkscape) Transparency is used (non-zero) for the text in Inkscape, but the package 'transparent.sty' is not loaded}%
    \renewcommand\transparent[1]{}%
  }%
  \providecommand\rotatebox[2]{#2}%
  \newcommand*\fsize{\dimexpr\f@size pt\relax}%
  \newcommand*\lineheight[1]{\fontsize{\fsize}{#1\fsize}\selectfont}%
  \ifx\svgwidth\undefined%
    \setlength{\unitlength}{56.25118034bp}%
    \ifx\svgscale\undefined%
      \relax%
    \else%
      \setlength{\unitlength}{\unitlength * \real{\svgscale}}%
    \fi%
  \else%
    \setlength{\unitlength}{\svgwidth}%
  \fi%
  \global\let\svgwidth\undefined%
  \global\let\svgscale\undefined%
  \makeatother%
  \begin{picture}(1,0.6539475)%
    \lineheight{1}%
    \setlength\tabcolsep{0pt}%
    \put(0,0){\includegraphics[width=\unitlength,page=1]{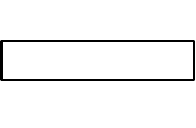}}%
    \put(0.29944908,0.31428646){\makebox(0,0)[lt]{\lineheight{1.45000005}\smash{\begin{tabular}[t]{l}$c_{(2,1^3)}$\end{tabular}}}}%
    \put(0,0){\includegraphics[width=\unitlength,page=2]{rhoz.pdf}}%
    \put(0.03826003,0.00626723){\makebox(0,0)[lt]{\lineheight{1.45000005}\smash{\begin{tabular}[t]{l}$v_1$\end{tabular}}}}%
    \put(0.20244032,0.00854682){\makebox(0,0)[lt]{\lineheight{1.45000005}\smash{\begin{tabular}[t]{l}$v_1$\end{tabular}}}}%
    \put(0.36843795,0.00878529){\makebox(0,0)[lt]{\lineheight{1.45000005}\smash{\begin{tabular}[t]{l}$v_2$\end{tabular}}}}%
    \put(0.56407963,0.00646523){\makebox(0,0)[lt]{\lineheight{1.45000005}\smash{\begin{tabular}[t]{l}$v_3$\end{tabular}}}}%
    \put(0,0){\includegraphics[width=\unitlength,page=3]{rhoz.pdf}}%
    \put(0.76146131,0.0077282){\makebox(0,0)[lt]{\lineheight{1.45000005}\smash{\begin{tabular}[t]{l}$v_4$\end{tabular}}}}%
    \put(0,0){\includegraphics[width=\unitlength,page=4]{rhoz.pdf}}%
  \end{picture}%
\endgroup%
}=
       \centre{}&-\centre{} +\centre{}\\
       &-\centre{}\neq 0
       \in B_{3,1}(n)_{(2,1^3)}.
     \end{split}
    \end{gather*}
    Therefore, we have $\rho_3\neq 0$ for $n\geq 4$. Since $B_{3,0}(n)_{(2^3)}$ is irreducible, $\rho_3$ is injective.
  \end{proof}

  \begin{remark}\label{compd3}
   We consider a restriction of the bracket map
   \begin{equation}\label{bra}
     [\cdot,\cdot]:V_{\lambda}\otimes \opegr^1(\IA(n))\rightarrow V_{\mu}
   \end{equation}
   for each irreducible $\GL(n;\Z)$-submodule $V_{\lambda}$ (resp. $V_{\mu}$) of $B_{d,k}(n)$ (resp. $B_{d,k+1}(n)$).
   We write a \emph{wavy arrow} $$V_{\lambda}\rightsquigarrow V_{\mu}$$ if the restriction map \eqref{bra} does not vanish.
   Then, we have the following diagram for $n\geq 4$:
   $$\xymatrix@C=9pt@R=0pt{
     B_3(n)&=B_{3,0}(n)\quad\oplus & B_{3,1}(n)\quad\oplus & B_{3,2}(n)\quad\oplus & B_{3,3}(n)\quad\oplus & B_{3,4}(n),\\
     &\rotatebox{90}{$\cong$}&\rotatebox{90}{$\cong$}&\rotatebox{90}{$\cong$} &\rotatebox{90}{$\cong$}&\rotatebox{90}{$\cong$}\\
     &B_{3,0}(n)_{(6)} \\
     &\oplus \\
     &B_{3,0}(n)_{(4,2)}\ar@{~>}[r]
     &B_{3,1}(n)_{(3,1^2)}\ar@{~>}[r]\ar@{~>}[rdd]\ar@{~>}[rdddd]\ar@{~>}[rdddddd] &B_{3,2}(n)_{(4)} \\
     &\oplus &&\oplus\\
     &B_{3,0}(n)_{(2^3)}\ar@{~>}[rdd]\ar@{~>}[ruu]
     &\oplus&B_{3,2}(n)_{(3,1)}\\
     &&&\oplus\\
     &&B_{3,1}(n)_{(2,1^3)}\ar@{~>}[ruu]\ar@{~>}[r] &B^{(1)}_{3,2}(n)_{(2^2)} \\
     &&&\oplus\\
     &&&B^{(2)}_{3,2}(n)_{(2^2)}\ar@{~>}[r]&B_{3,3}(n)_{(1^3)}\ar@{~>}[r] &B_{3,4}(n)_{(2)}\\
     }
   $$
   where $B_{3,2}^{(i)}(n)_{(2^2)}$ is the irreducible component of $B_{3,2}(n)_{(2^2)}$ generated by $$\centre{
\begingroup%
  \makeatletter%
  \providecommand\color[2][]{%
    \errmessage{(Inkscape) Color is used for the text in Inkscape, but the package 'color.sty' is not loaded}%
    \renewcommand\color[2][]{}%
  }%
  \providecommand\transparent[1]{%
    \errmessage{(Inkscape) Transparency is used (non-zero) for the text in Inkscape, but the package 'transparent.sty' is not loaded}%
    \renewcommand\transparent[1]{}%
  }%
  \providecommand\rotatebox[2]{#2}%
  \newcommand*\fsize{\dimexpr\f@size pt\relax}%
  \newcommand*\lineheight[1]{\fontsize{\fsize}{#1\fsize}\selectfont}%
  \ifx\svgwidth\undefined%
    \setlength{\unitlength}{39.00118268bp}%
    \ifx\svgscale\undefined%
      \relax%
    \else%
      \setlength{\unitlength}{\unitlength * \real{\svgscale}}%
    \fi%
  \else%
    \setlength{\unitlength}{\svgwidth}%
  \fi%
  \global\let\svgwidth\undefined%
  \global\let\svgscale\undefined%
  \makeatother%
  \begin{picture}(1,0.96024559)%
    \lineheight{1}%
    \setlength\tabcolsep{0pt}%
    \put(0,0){\includegraphics[width=\unitlength,page=1]{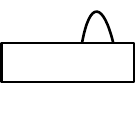}}%
    \put(0.29540421,0.46180856){\makebox(0,0)[lt]{\lineheight{1.45000005}\smash{\begin{tabular}[t]{l}$c_{(2,2)}$\end{tabular}}}}%
    \put(0,0){\includegraphics[width=\unitlength,page=2]{B321.pdf}}%
    \put(0.05880912,0.00903919){\makebox(0,0)[lt]{\lineheight{1.45000005}\smash{\begin{tabular}[t]{l}$v_1$\end{tabular}}}}%
    \put(0.29197851,0.01414043){\makebox(0,0)[lt]{\lineheight{1.45000005}\smash{\begin{tabular}[t]{l}$v_1$\end{tabular}}}}%
    \put(0.53139592,0.01448437){\makebox(0,0)[lt]{\lineheight{1.45000005}\smash{\begin{tabular}[t]{l}$v_2$\end{tabular}}}}%
    \put(0.76456531,0.0195856){\makebox(0,0)[lt]{\lineheight{1.45000005}\smash{\begin{tabular}[t]{l}$v_2$\end{tabular}}}}%
    \put(0,0){\includegraphics[width=\unitlength,page=3]{B321.pdf}}%
  \end{picture}%
\endgroup%
}+16\centre{
\begingroup%
  \makeatletter%
  \providecommand\color[2][]{%
    \errmessage{(Inkscape) Color is used for the text in Inkscape, but the package 'color.sty' is not loaded}%
    \renewcommand\color[2][]{}%
  }%
  \providecommand\transparent[1]{%
    \errmessage{(Inkscape) Transparency is used (non-zero) for the text in Inkscape, but the package 'transparent.sty' is not loaded}%
    \renewcommand\transparent[1]{}%
  }%
  \providecommand\rotatebox[2]{#2}%
  \newcommand*\fsize{\dimexpr\f@size pt\relax}%
  \newcommand*\lineheight[1]{\fontsize{\fsize}{#1\fsize}\selectfont}%
  \ifx\svgwidth\undefined%
    \setlength{\unitlength}{34.24794879bp}%
    \ifx\svgscale\undefined%
      \relax%
    \else%
      \setlength{\unitlength}{\unitlength * \real{\svgscale}}%
    \fi%
  \else%
    \setlength{\unitlength}{\svgwidth}%
  \fi%
  \global\let\svgwidth\undefined%
  \global\let\svgscale\undefined%
  \makeatother%
  \begin{picture}(1,0.56385483)%
    \lineheight{1}%
    \setlength\tabcolsep{0pt}%
    \put(0,0){\includegraphics[width=\unitlength,page=1]{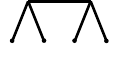}}%
    \put(-0.00441975,0.01029389){\makebox(0,0)[lt]{\lineheight{1.45000005}\smash{\begin{tabular}[t]{l}$v_1$\end{tabular}}}}%
    \put(0.27128854,0.01123905){\makebox(0,0)[lt]{\lineheight{1.45000005}\smash{\begin{tabular}[t]{l}$v_2$\end{tabular}}}}%
    \put(0.54660327,0.01084532){\makebox(0,0)[lt]{\lineheight{1.45000005}\smash{\begin{tabular}[t]{l}$v_1$\end{tabular}}}}%
    \put(0.82081489,0.01029373){\makebox(0,0)[lt]{\lineheight{1.45000005}\smash{\begin{tabular}[t]{l}$v_2$\end{tabular}}}}%
  \end{picture}%
\endgroup%
}\quad (i=1),$$ $$\centre{}\quad (i=2),$$
   respectively.
   Note that for $n=3$, $B_3(3)$ includes all of the above irreducible subrepresentations but $B_{3,1}(3)_{(2,1^3)}=0$ and there are all of the wavy arrows but the three wavy arrows that are directed to or coming from $B_{3,1}(3)_{(2,1^3)}$.
   For $n=2$, we have $$B_3(2)=(B_{3,0}(2)_{(6)} \oplus B_{3,0}(2)_{(4,2)})\oplus B_{3,2}(2)\oplus B_{3,4}(2).$$
   For $n=1$, we have $$B_3(1)=B_{3,0}(1)_{(6)} \oplus B_{3,2}(1)_{(4)} \oplus B_{3,4}(1)_{(2)}.$$
   For $n=1,2$, there are no wavy arrows because $B_{3,1}(n)=B_{3,3}(n)=0$.
  \end{remark}

  By Proposition \ref{propd3} and Remark \ref{compd3}, we have the surjectivity and the nontriviality of the bracket map for $n\geq 3$.
  Thus, by Theorem \ref{propindecomp}, one can obtain the following theorem, which improves Theorem \ref{propindecomp} for $d=3$.
  \begin{theorem}\label{decomp3}
    We have an indecomposable decomposition
    $$A_3(n)=A_3 P(n)\oplus A_3 Q(n)$$
    of $\Aut(F_n)$-modules for $n\geq 3$.
  \end{theorem}

  For $\lambda\vdash 4$, let $R_{\lambda}=\centre{
\begingroup%
  \makeatletter%
  \providecommand\color[2][]{%
    \errmessage{(Inkscape) Color is used for the text in Inkscape, but the package 'color.sty' is not loaded}%
    \renewcommand\color[2][]{}%
  }%
  \providecommand\transparent[1]{%
    \errmessage{(Inkscape) Transparency is used (non-zero) for the text in Inkscape, but the package 'transparent.sty' is not loaded}%
    \renewcommand\transparent[1]{}%
  }%
  \providecommand\rotatebox[2]{#2}%
  \newcommand*\fsize{\dimexpr\f@size pt\relax}%
  \newcommand*\lineheight[1]{\fontsize{\fsize}{#1\fsize}\selectfont}%
  \ifx\svgwidth\undefined%
    \setlength{\unitlength}{39.00118268bp}%
    \ifx\svgscale\undefined%
      \relax%
    \else%
      \setlength{\unitlength}{\unitlength * \real{\svgscale}}%
    \fi%
  \else%
    \setlength{\unitlength}{\svgwidth}%
  \fi%
  \global\let\svgwidth\undefined%
  \global\let\svgscale\undefined%
  \makeatother%
  \begin{picture}(1,0.91200334)%
    \lineheight{1}%
    \setlength\tabcolsep{0pt}%
    \put(0,0){\includegraphics[width=\unitlength,page=1]{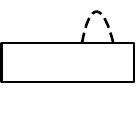}}%
    \put(0.33928657,0.37242658){\makebox(0,0)[lt]{\lineheight{1.45000005}\smash{\begin{tabular}[t]{l}$c_{\lambda}$\end{tabular}}}}%
    \put(0,0){\includegraphics[width=\unitlength,page=2]{Rlambda.pdf}}%
  \end{picture}%
\endgroup%
}\in A_3(4)$, where $c_{\lambda}$ is the Young symmetrizer.
  Let $S=R_{(2,2)}+ 16 \centre{
\begingroup%
  \makeatletter%
  \providecommand\color[2][]{%
    \errmessage{(Inkscape) Color is used for the text in Inkscape, but the package 'color.sty' is not loaded}%
    \renewcommand\color[2][]{}%
  }%
  \providecommand\transparent[1]{%
    \errmessage{(Inkscape) Transparency is used (non-zero) for the text in Inkscape, but the package 'transparent.sty' is not loaded}%
    \renewcommand\transparent[1]{}%
  }%
  \providecommand\rotatebox[2]{#2}%
  \newcommand*\fsize{\dimexpr\f@size pt\relax}%
  \newcommand*\lineheight[1]{\fontsize{\fsize}{#1\fsize}\selectfont}%
  \ifx\svgwidth\undefined%
    \setlength{\unitlength}{32.970128bp}%
    \ifx\svgscale\undefined%
      \relax%
    \else%
      \setlength{\unitlength}{\unitlength * \real{\svgscale}}%
    \fi%
  \else%
    \setlength{\unitlength}{\svgwidth}%
  \fi%
  \global\let\svgwidth\undefined%
  \global\let\svgscale\undefined%
  \makeatother%
  \begin{picture}(1,0.53288262)%
    \lineheight{1}%
    \setlength\tabcolsep{0pt}%
    \put(0,0){\includegraphics[width=\unitlength,page=1]{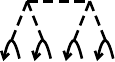}}%
  \end{picture}%
\endgroup%
}\in A_3(4)$ and  $T=\centre{
\begingroup%
  \makeatletter%
  \providecommand\color[2][]{%
    \errmessage{(Inkscape) Color is used for the text in Inkscape, but the package 'color.sty' is not loaded}%
    \renewcommand\color[2][]{}%
  }%
  \providecommand\transparent[1]{%
    \errmessage{(Inkscape) Transparency is used (non-zero) for the text in Inkscape, but the package 'transparent.sty' is not loaded}%
    \renewcommand\transparent[1]{}%
  }%
  \providecommand\rotatebox[2]{#2}%
  \newcommand*\fsize{\dimexpr\f@size pt\relax}%
  \newcommand*\lineheight[1]{\fontsize{\fsize}{#1\fsize}\selectfont}%
  \ifx\svgwidth\undefined%
    \setlength{\unitlength}{39.7264635bp}%
    \ifx\svgscale\undefined%
      \relax%
    \else%
      \setlength{\unitlength}{\unitlength * \real{\svgscale}}%
    \fi%
  \else%
    \setlength{\unitlength}{\svgwidth}%
  \fi%
  \global\let\svgwidth\undefined%
  \global\let\svgscale\undefined%
  \makeatother%
  \begin{picture}(1,0.59284079)%
    \lineheight{1}%
    \setlength\tabcolsep{0pt}%
    \put(0,0){\includegraphics[width=\unitlength,page=1]{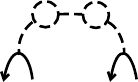}}%
  \end{picture}%
\endgroup%
}\in A_3(2)$.

  For $n=2$, we can check that $A_{3,2}(2)$ is semisimple as $\Aut(F_2)$-modules, that is,
  $$A_{3,2}(2)=A_3 R_{(4)}(2)\oplus A_3 R_{(3,1)}(2)\oplus A_3 S(2)\oplus A_3 U(2)\oplus A_3 T(2),$$
  where $U=\centre{
\begingroup%
  \makeatletter%
  \providecommand\color[2][]{%
    \errmessage{(Inkscape) Color is used for the text in Inkscape, but the package 'color.sty' is not loaded}%
    \renewcommand\color[2][]{}%
  }%
  \providecommand\transparent[1]{%
    \errmessage{(Inkscape) Transparency is used (non-zero) for the text in Inkscape, but the package 'transparent.sty' is not loaded}%
    \renewcommand\transparent[1]{}%
  }%
  \providecommand\rotatebox[2]{#2}%
  \newcommand*\fsize{\dimexpr\f@size pt\relax}%
  \newcommand*\lineheight[1]{\fontsize{\fsize}{#1\fsize}\selectfont}%
  \ifx\svgwidth\undefined%
    \setlength{\unitlength}{42.73034447bp}%
    \ifx\svgscale\undefined%
      \relax%
    \else%
      \setlength{\unitlength}{\unitlength * \real{\svgscale}}%
    \fi%
  \else%
    \setlength{\unitlength}{\svgwidth}%
  \fi%
  \global\let\svgwidth\undefined%
  \global\let\svgscale\undefined%
  \makeatother%
  \begin{picture}(1,0.82173075)%
    \lineheight{1}%
    \setlength\tabcolsep{0pt}%
    \put(0,0){\includegraphics[width=\unitlength,page=1]{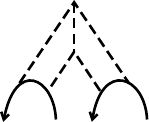}}%
  \end{picture}%
\endgroup%
}-\frac{1}{8}\centre{
\begingroup%
  \makeatletter%
  \providecommand\color[2][]{%
    \errmessage{(Inkscape) Color is used for the text in Inkscape, but the package 'color.sty' is not loaded}%
    \renewcommand\color[2][]{}%
  }%
  \providecommand\transparent[1]{%
    \errmessage{(Inkscape) Transparency is used (non-zero) for the text in Inkscape, but the package 'transparent.sty' is not loaded}%
    \renewcommand\transparent[1]{}%
  }%
  \providecommand\rotatebox[2]{#2}%
  \newcommand*\fsize{\dimexpr\f@size pt\relax}%
  \newcommand*\lineheight[1]{\fontsize{\fsize}{#1\fsize}\selectfont}%
  \ifx\svgwidth\undefined%
    \setlength{\unitlength}{42.82149062bp}%
    \ifx\svgscale\undefined%
      \relax%
    \else%
      \setlength{\unitlength}{\unitlength * \real{\svgscale}}%
    \fi%
  \else%
    \setlength{\unitlength}{\svgwidth}%
  \fi%
  \global\let\svgwidth\undefined%
  \global\let\svgscale\undefined%
  \makeatother%
  \begin{picture}(1,0.60253549)%
    \lineheight{1}%
    \setlength\tabcolsep{0pt}%
    \put(0,0){\includegraphics[width=\unitlength,page=1]{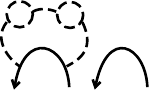}}%
  \end{picture}%
\endgroup%
}-\frac{1}{8}\centre{}-\frac{1}{8}\centre{
\begingroup%
  \makeatletter%
  \providecommand\color[2][]{%
    \errmessage{(Inkscape) Color is used for the text in Inkscape, but the package 'color.sty' is not loaded}%
    \renewcommand\color[2][]{}%
  }%
  \providecommand\transparent[1]{%
    \errmessage{(Inkscape) Transparency is used (non-zero) for the text in Inkscape, but the package 'transparent.sty' is not loaded}%
    \renewcommand\transparent[1]{}%
  }%
  \providecommand\rotatebox[2]{#2}%
  \newcommand*\fsize{\dimexpr\f@size pt\relax}%
  \newcommand*\lineheight[1]{\fontsize{\fsize}{#1\fsize}\selectfont}%
  \ifx\svgwidth\undefined%
    \setlength{\unitlength}{47.92795687bp}%
    \ifx\svgscale\undefined%
      \relax%
    \else%
      \setlength{\unitlength}{\unitlength * \real{\svgscale}}%
    \fi%
  \else%
    \setlength{\unitlength}{\svgwidth}%
  \fi%
  \global\let\svgwidth\undefined%
  \global\let\svgscale\undefined%
  \makeatother%
  \begin{picture}(1,0.53833857)%
    \lineheight{1}%
    \setlength\tabcolsep{0pt}%
    \put(0,0){\includegraphics[width=\unitlength,page=1]{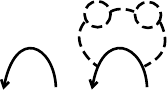}}%
  \end{picture}%
\endgroup%
}\in A_3(2).$
  We do not know whether or not the $\Aut(F_2)$-module $A_3(2)$ is semisimple.
  \begin{remark}\label{rad3}
    Since $A_{3,2}(2)$ is semisimple, we have $\Rad(A_{3,2}(2))=0$.
    On the other hand, we have $A_{3,3}(2)=A_{3,4}(2)\cong B_{3,4}(2)\neq 0$. Therefore, we have $\Rad(A_{3,2}(2))\neq A_{3,3}(2)$.
  \end{remark}

  For $n=1$, we have $\Aut(F_1)=\Z/2\Z$. We can easily check the following proposition.
  \begin{proposition}\label{A31}
    The $\Aut(F_1)$-action on $A_3(1)$ is trivial.
    Therefore, we have $A_3(1)=A_3 P(1)\oplus A_3 R_{(4)}(1) \oplus A_3 T(1).$
  \end{proposition}

 \subsection{The socle of $A_d(n)$ for small $d$}
  For an $\Aut(F_n)$-module $M$, let $\Soc(M)$ denote the \emph{socle} of $M$; that is,
  $$\Soc(M)=\sum\{K\subset M \mid K \text{ is simple}\}.$$

  Let us consider the cases for small $d$.
  Since $A_1(n)\cong \Sym^2(V_n)$ is simple, we have $$\Soc(A_1(n))=A_1(n)\quad (n\geq 1).$$

  By Theorem 6.9 of \cite{Mai1}, we have
  $$\Soc(A_2(n))=A_2 P(n)\oplus A_2 \ti{T}(n) \quad(n\geq 3, n=1),$$
  $$\Soc(A_2(n))=A_2(n)=A_2 P(n)\oplus A_2 W(n) \oplus A_2 \ti{T}(n) \quad(n=2),$$
  where $$\ti{T}=\scalebox{0.8}{$\centre{
\begingroup%
  \makeatletter%
  \providecommand\color[2][]{%
    \errmessage{(Inkscape) Color is used for the text in Inkscape, but the package 'color.sty' is not loaded}%
    \renewcommand\color[2][]{}%
  }%
  \providecommand\transparent[1]{%
    \errmessage{(Inkscape) Transparency is used (non-zero) for the text in Inkscape, but the package 'transparent.sty' is not loaded}%
    \renewcommand\transparent[1]{}%
  }%
  \providecommand\rotatebox[2]{#2}%
  \newcommand*\fsize{\dimexpr\f@size pt\relax}%
  \newcommand*\lineheight[1]{\fontsize{\fsize}{#1\fsize}\selectfont}%
  \ifx\svgwidth\undefined%
    \setlength{\unitlength}{24.72646255bp}%
    \ifx\svgscale\undefined%
      \relax%
    \else%
      \setlength{\unitlength}{\unitlength * \real{\svgscale}}%
    \fi%
  \else%
    \setlength{\unitlength}{\svgwidth}%
  \fi%
  \global\let\svgwidth\undefined%
  \global\let\svgscale\undefined%
  \makeatother%
  \begin{picture}(1,0.95248029)%
    \lineheight{1}%
    \setlength\tabcolsep{0pt}%
    \put(0,0){\includegraphics[width=\unitlength,page=1]{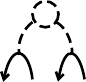}}%
  \end{picture}%
\endgroup%
}$},\; W=2\scalebox{0.8}{$\centre{
\begingroup%
  \makeatletter%
  \providecommand\color[2][]{%
    \errmessage{(Inkscape) Color is used for the text in Inkscape, but the package 'color.sty' is not loaded}%
    \renewcommand\color[2][]{}%
  }%
  \providecommand\transparent[1]{%
    \errmessage{(Inkscape) Transparency is used (non-zero) for the text in Inkscape, but the package 'transparent.sty' is not loaded}%
    \renewcommand\transparent[1]{}%
  }%
  \providecommand\rotatebox[2]{#2}%
  \newcommand*\fsize{\dimexpr\f@size pt\relax}%
  \newcommand*\lineheight[1]{\fontsize{\fsize}{#1\fsize}\selectfont}%
  \ifx\svgwidth\undefined%
    \setlength{\unitlength}{31.474021bp}%
    \ifx\svgscale\undefined%
      \relax%
    \else%
      \setlength{\unitlength}{\unitlength * \real{\svgscale}}%
    \fi%
  \else%
    \setlength{\unitlength}{\svgwidth}%
  \fi%
  \global\let\svgwidth\undefined%
  \global\let\svgscale\undefined%
  \makeatother%
  \begin{picture}(1,0.60753284)%
    \lineheight{1}%
    \setlength\tabcolsep{0pt}%
    \put(0,0){\includegraphics[width=\unitlength,page=1]{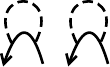}}%
  \end{picture}%
\endgroup%
}$}-\scalebox{0.8}{$\centre{
\begingroup%
  \makeatletter%
  \providecommand\color[2][]{%
    \errmessage{(Inkscape) Color is used for the text in Inkscape, but the package 'color.sty' is not loaded}%
    \renewcommand\color[2][]{}%
  }%
  \providecommand\transparent[1]{%
    \errmessage{(Inkscape) Transparency is used (non-zero) for the text in Inkscape, but the package 'transparent.sty' is not loaded}%
    \renewcommand\transparent[1]{}%
  }%
  \providecommand\rotatebox[2]{#2}%
  \newcommand*\fsize{\dimexpr\f@size pt\relax}%
  \newcommand*\lineheight[1]{\fontsize{\fsize}{#1\fsize}\selectfont}%
  \ifx\svgwidth\undefined%
    \setlength{\unitlength}{31.474021bp}%
    \ifx\svgscale\undefined%
      \relax%
    \else%
      \setlength{\unitlength}{\unitlength * \real{\svgscale}}%
    \fi%
  \else%
    \setlength{\unitlength}{\svgwidth}%
  \fi%
  \global\let\svgwidth\undefined%
  \global\let\svgscale\undefined%
  \makeatother%
  \begin{picture}(1,0.65427985)%
    \lineheight{1}%
    \setlength\tabcolsep{0pt}%
    \put(0,0){\includegraphics[width=\unitlength,page=1]{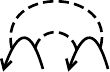}}%
  \end{picture}%
\endgroup%
}$}-\scalebox{0.8}{$\centre{
\begingroup%
  \makeatletter%
  \providecommand\color[2][]{%
    \errmessage{(Inkscape) Color is used for the text in Inkscape, but the package 'color.sty' is not loaded}%
    \renewcommand\color[2][]{}%
  }%
  \providecommand\transparent[1]{%
    \errmessage{(Inkscape) Transparency is used (non-zero) for the text in Inkscape, but the package 'transparent.sty' is not loaded}%
    \renewcommand\transparent[1]{}%
  }%
  \providecommand\rotatebox[2]{#2}%
  \newcommand*\fsize{\dimexpr\f@size pt\relax}%
  \newcommand*\lineheight[1]{\fontsize{\fsize}{#1\fsize}\selectfont}%
  \ifx\svgwidth\undefined%
    \setlength{\unitlength}{31.474021bp}%
    \ifx\svgscale\undefined%
      \relax%
    \else%
      \setlength{\unitlength}{\unitlength * \real{\svgscale}}%
    \fi%
  \else%
    \setlength{\unitlength}{\svgwidth}%
  \fi%
  \global\let\svgwidth\undefined%
  \global\let\svgscale\undefined%
  \makeatother%
  \begin{picture}(1,0.64458778)%
    \lineheight{1}%
    \setlength\tabcolsep{0pt}%
    \put(0,0){\includegraphics[width=\unitlength,page=1]{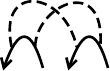}}%
  \end{picture}%
\endgroup%
}$}-\frac{1}{2}\scalebox{0.8}{$\centre{
\begingroup%
  \makeatletter%
  \providecommand\color[2][]{%
    \errmessage{(Inkscape) Color is used for the text in Inkscape, but the package 'color.sty' is not loaded}%
    \renewcommand\color[2][]{}%
  }%
  \providecommand\transparent[1]{%
    \errmessage{(Inkscape) Transparency is used (non-zero) for the text in Inkscape, but the package 'transparent.sty' is not loaded}%
    \renewcommand\transparent[1]{}%
  }%
  \providecommand\rotatebox[2]{#2}%
  \newcommand*\fsize{\dimexpr\f@size pt\relax}%
  \newcommand*\lineheight[1]{\fontsize{\fsize}{#1\fsize}\selectfont}%
  \ifx\svgwidth\undefined%
    \setlength{\unitlength}{31.88709807bp}%
    \ifx\svgscale\undefined%
      \relax%
    \else%
      \setlength{\unitlength}{\unitlength * \real{\svgscale}}%
    \fi%
  \else%
    \setlength{\unitlength}{\svgwidth}%
  \fi%
  \global\let\svgwidth\undefined%
  \global\let\svgscale\undefined%
  \makeatother%
  \begin{picture}(1,0.72812646)%
    \lineheight{1}%
    \setlength\tabcolsep{0pt}%
    \put(0,0){\includegraphics[width=\unitlength,page=1]{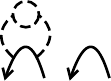}}%
  \end{picture}%
\endgroup%
}$}-\frac{1}{2}\scalebox{0.8}{$\centre{
\begingroup%
  \makeatletter%
  \providecommand\color[2][]{%
    \errmessage{(Inkscape) Color is used for the text in Inkscape, but the package 'color.sty' is not loaded}%
    \renewcommand\color[2][]{}%
  }%
  \providecommand\transparent[1]{%
    \errmessage{(Inkscape) Transparency is used (non-zero) for the text in Inkscape, but the package 'transparent.sty' is not loaded}%
    \renewcommand\transparent[1]{}%
  }%
  \providecommand\rotatebox[2]{#2}%
  \newcommand*\fsize{\dimexpr\f@size pt\relax}%
  \newcommand*\lineheight[1]{\fontsize{\fsize}{#1\fsize}\selectfont}%
  \ifx\svgwidth\undefined%
    \setlength{\unitlength}{33.12910267bp}%
    \ifx\svgscale\undefined%
      \relax%
    \else%
      \setlength{\unitlength}{\unitlength * \real{\svgscale}}%
    \fi%
  \else%
    \setlength{\unitlength}{\svgwidth}%
  \fi%
  \global\let\svgwidth\undefined%
  \global\let\svgscale\undefined%
  \makeatother%
  \begin{picture}(1,0.70082912)%
    \lineheight{1}%
    \setlength\tabcolsep{0pt}%
    \put(0,0){\includegraphics[width=\unitlength,page=1]{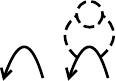}}%
  \end{picture}%
\endgroup%
}$} \in A_2(2).$$
  Note that $A_2 \ti{T}(n)=A_{2,2}(n)$

  By Proposition \ref{A31}, we have $\Soc(A_3(1))=A_3(1)$.
  \begin{proposition}\label{soc}
    For $n\geq 3$, we have
    $$\Soc(A_3(n))=A_3 P(n)\oplus A_3 R_{(4)}(n)\oplus A_3 R_{(3,1)}(n)\oplus A_3 S(n)\oplus A_3 T(n).$$
  \end{proposition}
  \begin{proof}
    A simple $\Aut(F_n)$-submodule $K\subset A_3(n)$ corresponds to an irreducible component of $B_3(n)$ via the PBW map.
    Therefore, by Remark \ref{compd3}, we have
    $$\Soc(A_3(n))\subset A_3 P(n)\oplus A_3 R_{(4)}(n)\oplus A_3 R_{(3,1)}(n)\oplus A_3 S(n)\oplus A_3 T(n).$$
    Moreover, we can check that
    \begin{gather*}
       A_3 P(n)\cong V_{(6)},\quad A_3 R_{(4)}(n)\cong V_{(4)},\quad
       A_3 R_{(3,1)}(n)\cong V_{(3,1)},\\
       A_3 S(n)\cong V_{(2,2)},\quad A_3 T(n)\cong V_{(2)}.
    \end{gather*}
    Hence, we have $$\Soc(A_3(n))\supset A_3 P(n)\oplus A_3 R_{(4)}(n)\oplus A_3 R_{(3,1)}(n)\oplus A_3 S(n)\oplus A_3 T(n)$$
    and the proof is complete.
  \end{proof}

 \subsection{The indecomposable decomposition of $A_4(n)$}
  Here, we consider the indecomposable decomposition of $A_4(n)$.

  Similarly, in degree $4$, we have $\GL(n;\Z)$-module homomorphisms
  \begin{gather*}
   \begin{split}
     \rho_1:B_{4,0}(n)_{(6,2)}&\rightarrow \Hom(\opegr^1(\IA(n)), B_{4,1}(n)_{(5,1^2)}),\\
     \rho_2:B_{4,0}(n)_{(4^2)}&\rightarrow \Hom(\opegr^1(\IA(n)), B_{4,1}(n)_{(3^2,1)}),\\
     \rho_3:B_{4,0}(n)_{(4,2^2)}&\rightarrow \Hom(\opegr^1(\IA(n)), B_{4,1}(n)_{(5,1^2)}),\\
     \rho_4:B_{4,0}(n)_{(4,2^2)}&\rightarrow \Hom(\opegr^1(\IA(n)), B_{4,1}(n)_{(4,1^3)}),\\
     \rho_5:B_{4,0}(n)_{(4,2^2)}&\rightarrow \Hom(\opegr^1(\IA(n)), B_{4,1}(n)_{(3^2,1)}),\\
     \rho_6:B_{4,0}(n)_{(4,2^2)}&\rightarrow \Hom(\opegr^1(\IA(n)), B_{4,1}(n)_{(3,2,1^2)}),\\
     \rho_7:B_{4,0}(n)_{(2^4)}&\rightarrow \Hom(\opegr^1(\IA(n)), B_{4,1}(n)_{(3,2,1^2)}),\\
     \rho_8:B_{4,0}(n)_{(2^4)}&\rightarrow \Hom(\opegr^1(\IA(n)), B_{4,1}(n)_{(2^2,1^3)}).
   \end{split}
  \end{gather*}

  \begin{proposition}\label{propd4}
    The $\GL(n;\Z)$-module homomorphisms $\rho_1, \rho_2,\rho_3$ and $\rho_5$ are injective for $n\geq 3$, $\rho_4,\rho_6$ and $\rho_7$ for $n\geq 4$ and $\rho_8$ for $n\geq 5$.
  \end{proposition}
  \begin{proof}
    As in the proof of Proposition \ref{propd3} in degree $3$, we will check that $\rho_1$ is injective for $n\geq 3$, $\rho_7$ for $n\geq 4$ and $\rho_8$ for $n\geq 5$. The others can be obtained in a similar way.

    For $n\geq 3$, we have
    $$[u,K_{3,1,2}]=14 w\neq 0\in B_{4,1}(n)_{(5,1^2)},$$
    where
    $$u=\frac{1}{1728}\centre{
\begingroup%
  \makeatletter%
  \providecommand\color[2][]{%
    \errmessage{(Inkscape) Color is used for the text in Inkscape, but the package 'color.sty' is not loaded}%
    \renewcommand\color[2][]{}%
  }%
  \providecommand\transparent[1]{%
    \errmessage{(Inkscape) Transparency is used (non-zero) for the text in Inkscape, but the package 'transparent.sty' is not loaded}%
    \renewcommand\transparent[1]{}%
  }%
  \providecommand\rotatebox[2]{#2}%
  \newcommand*\fsize{\dimexpr\f@size pt\relax}%
  \newcommand*\lineheight[1]{\fontsize{\fsize}{#1\fsize}\selectfont}%
  \ifx\svgwidth\undefined%
    \setlength{\unitlength}{74.25118142bp}%
    \ifx\svgscale\undefined%
      \relax%
    \else%
      \setlength{\unitlength}{\unitlength * \real{\svgscale}}%
    \fi%
  \else%
    \setlength{\unitlength}{\svgwidth}%
  \fi%
  \global\let\svgwidth\undefined%
  \global\let\svgscale\undefined%
  \makeatother%
  \begin{picture}(1,0.46353505)%
    \lineheight{1}%
    \setlength\tabcolsep{0pt}%
    \put(0,0){\includegraphics[width=\unitlength,page=1]{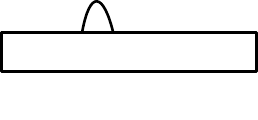}}%
    \put(0.38451146,0.23284819){\makebox(0,0)[lt]{\lineheight{1.45000005}\smash{\begin{tabular}[t]{l}$c_{(6,2)}$\end{tabular}}}}%
    \put(0,0){\includegraphics[width=\unitlength,page=2]{rhou.pdf}}%
    \put(0.02898502,0.00474793){\makebox(0,0)[lt]{\lineheight{1.45000005}\smash{\begin{tabular}[t]{l}$v_1$\end{tabular}}}}%
    \put(0.15336466,0.0064749){\makebox(0,0)[lt]{\lineheight{1.45000005}\smash{\begin{tabular}[t]{l}$v_1$\end{tabular}}}}%
    \put(0.27912107,0.00665556){\makebox(0,0)[lt]{\lineheight{1.45000005}\smash{\begin{tabular}[t]{l}$v_1$\end{tabular}}}}%
    \put(0.39863754,0.00489793){\makebox(0,0)[lt]{\lineheight{1.45000005}\smash{\begin{tabular}[t]{l}$v_1$\end{tabular}}}}%
    \put(0,0){\includegraphics[width=\unitlength,page=3]{rhou.pdf}}%
    \put(0.51673905,0.00585473){\makebox(0,0)[lt]{\lineheight{1.45000005}\smash{\begin{tabular}[t]{l}$v_1$\end{tabular}}}}%
    \put(0.63921359,0.0075817){\makebox(0,0)[lt]{\lineheight{1.45000005}\smash{\begin{tabular}[t]{l}$v_1$\end{tabular}}}}%
    \put(0,0){\includegraphics[width=\unitlength,page=4]{rhou.pdf}}%
    \put(0.76926017,0.00585466){\makebox(0,0)[lt]{\lineheight{1.45000005}\smash{\begin{tabular}[t]{l}$v_2$\end{tabular}}}}%
    \put(0.8917347,0.0075817){\makebox(0,0)[lt]{\lineheight{1.45000005}\smash{\begin{tabular}[t]{l}$v_2$\end{tabular}}}}%
  \end{picture}%
\endgroup%
}= \centre{
\begingroup%
  \makeatletter%
  \providecommand\color[2][]{%
    \errmessage{(Inkscape) Color is used for the text in Inkscape, but the package 'color.sty' is not loaded}%
    \renewcommand\color[2][]{}%
  }%
  \providecommand\transparent[1]{%
    \errmessage{(Inkscape) Transparency is used (non-zero) for the text in Inkscape, but the package 'transparent.sty' is not loaded}%
    \renewcommand\transparent[1]{}%
  }%
  \providecommand\rotatebox[2]{#2}%
  \newcommand*\fsize{\dimexpr\f@size pt\relax}%
  \newcommand*\lineheight[1]{\fontsize{\fsize}{#1\fsize}\selectfont}%
  \ifx\svgwidth\undefined%
    \setlength{\unitlength}{86.17435255bp}%
    \ifx\svgscale\undefined%
      \relax%
    \else%
      \setlength{\unitlength}{\unitlength * \real{\svgscale}}%
    \fi%
  \else%
    \setlength{\unitlength}{\svgwidth}%
  \fi%
  \global\let\svgwidth\undefined%
  \global\let\svgscale\undefined%
  \makeatother%
  \begin{picture}(1,0.18277808)%
    \lineheight{1}%
    \setlength\tabcolsep{0pt}%
    \put(0,0){\includegraphics[width=\unitlength,page=1]{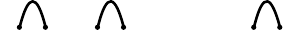}}%
    \put(-0.00175652,0.00422363){\makebox(0,0)[lt]{\lineheight{1.45000005}\smash{\begin{tabular}[t]{l}$v_1$\end{tabular}}}}%
    \put(0.14127659,0.0064574){\makebox(0,0)[lt]{\lineheight{1.45000005}\smash{\begin{tabular}[t]{l}$v_1$\end{tabular}}}}%
    \put(0.92878714,0.00554643){\makebox(0,0)[lt]{\lineheight{1.45000005}\smash{\begin{tabular}[t]{l}$v_2$\end{tabular}}}}%
    \put(0.27450554,0.00453982){\makebox(0,0)[lt]{\lineheight{1.45000005}\smash{\begin{tabular}[t]{l}$v_1$\end{tabular}}}}%
    \put(0.79819249,0.00489901){\makebox(0,0)[lt]{\lineheight{1.45000005}\smash{\begin{tabular}[t]{l}$v_2$\end{tabular}}}}%
    \put(0.40333958,0.004091){\makebox(0,0)[lt]{\lineheight{1.45000005}\smash{\begin{tabular}[t]{l}$v_1$\end{tabular}}}}%
    \put(0,0){\includegraphics[width=\unitlength,page=2]{rho1u4.pdf}}%
    \put(0.54458713,0.00504745){\makebox(0,0)[lt]{\lineheight{1.45000005}\smash{\begin{tabular}[t]{l}$v_1$\end{tabular}}}}%
    \put(0.66739569,0.00640112){\makebox(0,0)[lt]{\lineheight{1.45000005}\smash{\begin{tabular}[t]{l}$v_1$\end{tabular}}}}%
  \end{picture}%
\endgroup%
}-\centre{
\begingroup%
  \makeatletter%
  \providecommand\color[2][]{%
    \errmessage{(Inkscape) Color is used for the text in Inkscape, but the package 'color.sty' is not loaded}%
    \renewcommand\color[2][]{}%
  }%
  \providecommand\transparent[1]{%
    \errmessage{(Inkscape) Transparency is used (non-zero) for the text in Inkscape, but the package 'transparent.sty' is not loaded}%
    \renewcommand\transparent[1]{}%
  }%
  \providecommand\rotatebox[2]{#2}%
  \newcommand*\fsize{\dimexpr\f@size pt\relax}%
  \newcommand*\lineheight[1]{\fontsize{\fsize}{#1\fsize}\selectfont}%
  \ifx\svgwidth\undefined%
    \setlength{\unitlength}{86.17435255bp}%
    \ifx\svgscale\undefined%
      \relax%
    \else%
      \setlength{\unitlength}{\unitlength * \real{\svgscale}}%
    \fi%
  \else%
    \setlength{\unitlength}{\svgwidth}%
  \fi%
  \global\let\svgwidth\undefined%
  \global\let\svgscale\undefined%
  \makeatother%
  \begin{picture}(1,0.18277808)%
    \lineheight{1}%
    \setlength\tabcolsep{0pt}%
    \put(0,0){\includegraphics[width=\unitlength,page=1]{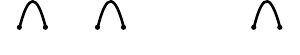}}%
    \put(-0.00175652,0.00422363){\makebox(0,0)[lt]{\lineheight{1.45000005}\smash{\begin{tabular}[t]{l}$v_1$\end{tabular}}}}%
    \put(0.14127659,0.0064574){\makebox(0,0)[lt]{\lineheight{1.45000005}\smash{\begin{tabular}[t]{l}$v_1$\end{tabular}}}}%
    \put(0.92878714,0.00554643){\makebox(0,0)[lt]{\lineheight{1.45000005}\smash{\begin{tabular}[t]{l}$v_2$\end{tabular}}}}%
    \put(0.27450554,0.00453982){\makebox(0,0)[lt]{\lineheight{1.45000005}\smash{\begin{tabular}[t]{l}$v_1$\end{tabular}}}}%
    \put(0.79819249,0.00489901){\makebox(0,0)[lt]{\lineheight{1.45000005}\smash{\begin{tabular}[t]{l}$v_1$\end{tabular}}}}%
    \put(0.40333958,0.004091){\makebox(0,0)[lt]{\lineheight{1.45000005}\smash{\begin{tabular}[t]{l}$v_1$\end{tabular}}}}%
    \put(0,0){\includegraphics[width=\unitlength,page=2]{rho1u42.pdf}}%
    \put(0.54458713,0.00504745){\makebox(0,0)[lt]{\lineheight{1.45000005}\smash{\begin{tabular}[t]{l}$v_1$\end{tabular}}}}%
    \put(0.66739569,0.00640112){\makebox(0,0)[lt]{\lineheight{1.45000005}\smash{\begin{tabular}[t]{l}$v_2$\end{tabular}}}}%
  \end{picture}%
\endgroup%
}\in B_{4,0}(n)_{(6,2)}$$
    and
    $$w=\frac{1}{336}\centre{
\begingroup%
  \makeatletter%
  \providecommand\color[2][]{%
    \errmessage{(Inkscape) Color is used for the text in Inkscape, but the package 'color.sty' is not loaded}%
    \renewcommand\color[2][]{}%
  }%
  \providecommand\transparent[1]{%
    \errmessage{(Inkscape) Transparency is used (non-zero) for the text in Inkscape, but the package 'transparent.sty' is not loaded}%
    \renewcommand\transparent[1]{}%
  }%
  \providecommand\rotatebox[2]{#2}%
  \newcommand*\fsize{\dimexpr\f@size pt\relax}%
  \newcommand*\lineheight[1]{\fontsize{\fsize}{#1\fsize}\selectfont}%
  \ifx\svgwidth\undefined%
    \setlength{\unitlength}{72.00117926bp}%
    \ifx\svgscale\undefined%
      \relax%
    \else%
      \setlength{\unitlength}{\unitlength * \real{\svgscale}}%
    \fi%
  \else%
    \setlength{\unitlength}{\svgwidth}%
  \fi%
  \global\let\svgwidth\undefined%
  \global\let\svgscale\undefined%
  \makeatother%
  \begin{picture}(1,0.51089883)%
    \lineheight{1}%
    \setlength\tabcolsep{0pt}%
    \put(0,0){\includegraphics[width=\unitlength,page=1]{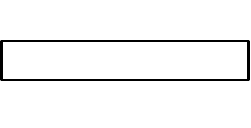}}%
    \put(0.35783744,0.24028311){\makebox(0,0)[lt]{\lineheight{1.45000005}\smash{\begin{tabular}[t]{l}$c_{(5,1^2)}$\end{tabular}}}}%
    \put(0,0){\includegraphics[width=\unitlength,page=2]{rhov.pdf}}%
    \put(0.27988667,0.0048963){\makebox(0,0)[lt]{\lineheight{1.45000005}\smash{\begin{tabular}[t]{l}$v_1$\end{tabular}}}}%
    \put(0.40815311,0.00667724){\makebox(0,0)[lt]{\lineheight{1.45000005}\smash{\begin{tabular}[t]{l}$v_1$\end{tabular}}}}%
    \put(0.53783936,0.00686354){\makebox(0,0)[lt]{\lineheight{1.45000005}\smash{\begin{tabular}[t]{l}$v_1$\end{tabular}}}}%
    \put(0.69068513,0.00505099){\makebox(0,0)[lt]{\lineheight{1.45000005}\smash{\begin{tabular}[t]{l}$v_2$\end{tabular}}}}%
    \put(0,0){\includegraphics[width=\unitlength,page=3]{rhov.pdf}}%
    \put(0.84489027,0.00603769){\makebox(0,0)[lt]{\lineheight{1.45000005}\smash{\begin{tabular}[t]{l}$v_3$\end{tabular}}}}%
    \put(0,0){\includegraphics[width=\unitlength,page=4]{rhov.pdf}}%
    \put(0.02989079,0.0048963){\makebox(0,0)[lt]{\lineheight{1.45000005}\smash{\begin{tabular}[t]{l}$v_1$\end{tabular}}}}%
    \put(0.15815722,0.00667716){\makebox(0,0)[lt]{\lineheight{1.45000005}\smash{\begin{tabular}[t]{l}$v_1$\end{tabular}}}}%
    \put(0,0){\includegraphics[width=\unitlength,page=5]{rhov.pdf}}%
  \end{picture}%
\endgroup%
}=
    \centre{
\begingroup%
  \makeatletter%
  \providecommand\color[2][]{%
    \errmessage{(Inkscape) Color is used for the text in Inkscape, but the package 'color.sty' is not loaded}%
    \renewcommand\color[2][]{}%
  }%
  \providecommand\transparent[1]{%
    \errmessage{(Inkscape) Transparency is used (non-zero) for the text in Inkscape, but the package 'transparent.sty' is not loaded}%
    \renewcommand\transparent[1]{}%
  }%
  \providecommand\rotatebox[2]{#2}%
  \newcommand*\fsize{\dimexpr\f@size pt\relax}%
  \newcommand*\lineheight[1]{\fontsize{\fsize}{#1\fsize}\selectfont}%
  \ifx\svgwidth\undefined%
    \setlength{\unitlength}{67.15510934bp}%
    \ifx\svgscale\undefined%
      \relax%
    \else%
      \setlength{\unitlength}{\unitlength * \real{\svgscale}}%
    \fi%
  \else%
    \setlength{\unitlength}{\svgwidth}%
  \fi%
  \global\let\svgwidth\undefined%
  \global\let\svgscale\undefined%
  \makeatother%
  \begin{picture}(1,0.29942078)%
    \lineheight{1}%
    \setlength\tabcolsep{0pt}%
    \put(0,0){\includegraphics[width=\unitlength,page=1]{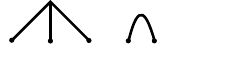}}%
    \put(0.63100605,0.01074403){\makebox(0,0)[lt]{\lineheight{1.45000005}\smash{\begin{tabular}[t]{l}$v_1$\end{tabular}}}}%
    \put(0.33000446,0.00594713){\makebox(0,0)[lt]{\lineheight{1.45000005}\smash{\begin{tabular}[t]{l}$v_3$\end{tabular}}}}%
    \put(-0.00225399,0.00649886){\makebox(0,0)[lt]{\lineheight{1.45000005}\smash{\begin{tabular}[t]{l}$v_1$\end{tabular}}}}%
    \put(0.49618747,0.01020103){\makebox(0,0)[lt]{\lineheight{1.45000005}\smash{\begin{tabular}[t]{l}$v_1$\end{tabular}}}}%
    \put(0.16437732,0.00524963){\makebox(0,0)[lt]{\lineheight{1.45000005}\smash{\begin{tabular}[t]{l}$v_2$\end{tabular}}}}%
    \put(0,0){\includegraphics[width=\unitlength,page=2]{rho1v4.pdf}}%
    \put(0.90861868,0.01074403){\makebox(0,0)[lt]{\lineheight{1.45000005}\smash{\begin{tabular}[t]{l}$v_1$\end{tabular}}}}%
    \put(0.77061706,0.01020099){\makebox(0,0)[lt]{\lineheight{1.45000005}\smash{\begin{tabular}[t]{l}$v_1$\end{tabular}}}}%
  \end{picture}%
\endgroup%
}\neq 0\in B_{4,1}(n)_{(5,1^2)}.$$
    Thus, we have $\rho_1\neq 0$ for $n\geq 3$.
    Since $B_{4,0}(n)_{(6,2)}$ is irreducible, $\rho_1$ is injective.

    For $n\geq 4$, we have
    $$[x,K_{1,4,3}]=-48 y\neq 0 \in B_{4,1}(n)_{(3,2,1^2)},$$
    where
    $$x=\centre{
\begingroup%
  \makeatletter%
  \providecommand\color[2][]{%
    \errmessage{(Inkscape) Color is used for the text in Inkscape, but the package 'color.sty' is not loaded}%
    \renewcommand\color[2][]{}%
  }%
  \providecommand\transparent[1]{%
    \errmessage{(Inkscape) Transparency is used (non-zero) for the text in Inkscape, but the package 'transparent.sty' is not loaded}%
    \renewcommand\transparent[1]{}%
  }%
  \providecommand\rotatebox[2]{#2}%
  \newcommand*\fsize{\dimexpr\f@size pt\relax}%
  \newcommand*\lineheight[1]{\fontsize{\fsize}{#1\fsize}\selectfont}%
  \ifx\svgwidth\undefined%
    \setlength{\unitlength}{74.25118142bp}%
    \ifx\svgscale\undefined%
      \relax%
    \else%
      \setlength{\unitlength}{\unitlength * \real{\svgscale}}%
    \fi%
  \else%
    \setlength{\unitlength}{\svgwidth}%
  \fi%
  \global\let\svgwidth\undefined%
  \global\let\svgscale\undefined%
  \makeatother%
  \begin{picture}(1,0.46353505)%
    \lineheight{1}%
    \setlength\tabcolsep{0pt}%
    \put(0,0){\includegraphics[width=\unitlength,page=1]{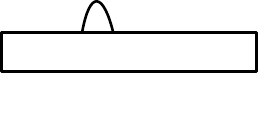}}%
    \put(0.38451146,0.23284819){\makebox(0,0)[lt]{\lineheight{1.45000005}\smash{\begin{tabular}[t]{l}$c_{(2^4)}$\end{tabular}}}}%
    \put(0,0){\includegraphics[width=\unitlength,page=2]{rho7x.pdf}}%
    \put(0.02898502,0.00474793){\makebox(0,0)[lt]{\lineheight{1.45000005}\smash{\begin{tabular}[t]{l}$v_1$\end{tabular}}}}%
    \put(0.15336466,0.0064749){\makebox(0,0)[lt]{\lineheight{1.45000005}\smash{\begin{tabular}[t]{l}$v_1$\end{tabular}}}}%
    \put(0.27912107,0.00665556){\makebox(0,0)[lt]{\lineheight{1.45000005}\smash{\begin{tabular}[t]{l}$v_2$\end{tabular}}}}%
    \put(0.39863754,0.00489793){\makebox(0,0)[lt]{\lineheight{1.45000005}\smash{\begin{tabular}[t]{l}$v_2$\end{tabular}}}}%
    \put(0,0){\includegraphics[width=\unitlength,page=3]{rho7x.pdf}}%
    \put(0.51673905,0.00585473){\makebox(0,0)[lt]{\lineheight{1.45000005}\smash{\begin{tabular}[t]{l}$v_3$\end{tabular}}}}%
    \put(0.63921359,0.0075817){\makebox(0,0)[lt]{\lineheight{1.45000005}\smash{\begin{tabular}[t]{l}$v_3$\end{tabular}}}}%
    \put(0,0){\includegraphics[width=\unitlength,page=4]{rho7x.pdf}}%
    \put(0.76926017,0.00585466){\makebox(0,0)[lt]{\lineheight{1.45000005}\smash{\begin{tabular}[t]{l}$v_4$\end{tabular}}}}%
    \put(0.8917347,0.0075817){\makebox(0,0)[lt]{\lineheight{1.45000005}\smash{\begin{tabular}[t]{l}$v_4$\end{tabular}}}}%
  \end{picture}%
\endgroup%
}\in B_{4,0}(n)_{(2^4)}$$
    and
    $$y=\centre{
\begingroup%
  \makeatletter%
  \providecommand\color[2][]{%
    \errmessage{(Inkscape) Color is used for the text in Inkscape, but the package 'color.sty' is not loaded}%
    \renewcommand\color[2][]{}%
  }%
  \providecommand\transparent[1]{%
    \errmessage{(Inkscape) Transparency is used (non-zero) for the text in Inkscape, but the package 'transparent.sty' is not loaded}%
    \renewcommand\transparent[1]{}%
  }%
  \providecommand\rotatebox[2]{#2}%
  \newcommand*\fsize{\dimexpr\f@size pt\relax}%
  \newcommand*\lineheight[1]{\fontsize{\fsize}{#1\fsize}\selectfont}%
  \ifx\svgwidth\undefined%
    \setlength{\unitlength}{72.00117926bp}%
    \ifx\svgscale\undefined%
      \relax%
    \else%
      \setlength{\unitlength}{\unitlength * \real{\svgscale}}%
    \fi%
  \else%
    \setlength{\unitlength}{\svgwidth}%
  \fi%
  \global\let\svgwidth\undefined%
  \global\let\svgscale\undefined%
  \makeatother%
  \begin{picture}(1,0.5608343)%
    \lineheight{1}%
    \setlength\tabcolsep{0pt}%
    \put(0,0){\includegraphics[width=\unitlength,page=1]{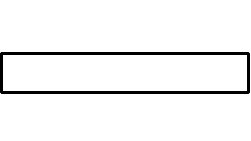}}%
    \put(0.35783743,0.24028307){\makebox(0,0)[lt]{\lineheight{1.45000005}\smash{\begin{tabular}[t]{l}$c_{(3,2,1^2)}$\end{tabular}}}}%
    \put(0,0){\includegraphics[width=\unitlength,page=2]{rho7y.pdf}}%
    \put(0.27988668,0.0048963){\makebox(0,0)[lt]{\lineheight{1.45000005}\smash{\begin{tabular}[t]{l}$v_1$\end{tabular}}}}%
    \put(0.4081531,0.00667724){\makebox(0,0)[lt]{\lineheight{1.45000005}\smash{\begin{tabular}[t]{l}$v_2$\end{tabular}}}}%
    \put(0.53783934,0.00686354){\makebox(0,0)[lt]{\lineheight{1.45000005}\smash{\begin{tabular}[t]{l}$v_2$\end{tabular}}}}%
    \put(0.69068511,0.00505099){\makebox(0,0)[lt]{\lineheight{1.45000005}\smash{\begin{tabular}[t]{l}$v_3$\end{tabular}}}}%
    \put(0,0){\includegraphics[width=\unitlength,page=3]{rho7y.pdf}}%
    \put(0.8448903,0.00603769){\makebox(0,0)[lt]{\lineheight{1.45000005}\smash{\begin{tabular}[t]{l}$v_4$\end{tabular}}}}%
    \put(0,0){\includegraphics[width=\unitlength,page=4]{rho7y.pdf}}%
    \put(0.02989079,0.0048963){\makebox(0,0)[lt]{\lineheight{1.45000005}\smash{\begin{tabular}[t]{l}$v_1$\end{tabular}}}}%
    \put(0.15815722,0.00667716){\makebox(0,0)[lt]{\lineheight{1.45000005}\smash{\begin{tabular}[t]{l}$v_1$\end{tabular}}}}%
    \put(0,0){\includegraphics[width=\unitlength,page=5]{rho7y.pdf}}%
  \end{picture}%
\endgroup%
}\in B_{4,1}(n)_{(3,2,1^2)}.$$
    Thus, $\rho_7$ is injective for $n\geq 4$.

    For $n\geq 5$, we have
    $$[x, K_{5,4,3}]=-48y' -32z,$$
    where
    $$y'=\centre{
\begingroup%
  \makeatletter%
  \providecommand\color[2][]{%
    \errmessage{(Inkscape) Color is used for the text in Inkscape, but the package 'color.sty' is not loaded}%
    \renewcommand\color[2][]{}%
  }%
  \providecommand\transparent[1]{%
    \errmessage{(Inkscape) Transparency is used (non-zero) for the text in Inkscape, but the package 'transparent.sty' is not loaded}%
    \renewcommand\transparent[1]{}%
  }%
  \providecommand\rotatebox[2]{#2}%
  \newcommand*\fsize{\dimexpr\f@size pt\relax}%
  \newcommand*\lineheight[1]{\fontsize{\fsize}{#1\fsize}\selectfont}%
  \ifx\svgwidth\undefined%
    \setlength{\unitlength}{72.00117926bp}%
    \ifx\svgscale\undefined%
      \relax%
    \else%
      \setlength{\unitlength}{\unitlength * \real{\svgscale}}%
    \fi%
  \else%
    \setlength{\unitlength}{\svgwidth}%
  \fi%
  \global\let\svgwidth\undefined%
  \global\let\svgscale\undefined%
  \makeatother%
  \begin{picture}(1,0.5608343)%
    \lineheight{1}%
    \setlength\tabcolsep{0pt}%
    \put(0,0){\includegraphics[width=\unitlength,page=1]{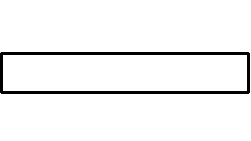}}%
    \put(0.35783743,0.24028307){\makebox(0,0)[lt]{\lineheight{1.45000005}\smash{\begin{tabular}[t]{l}$c_{(3,2,1^2)}$\end{tabular}}}}%
    \put(0,0){\includegraphics[width=\unitlength,page=2]{rho8y.pdf}}%
    \put(0.27988668,0.0048963){\makebox(0,0)[lt]{\lineheight{1.45000005}\smash{\begin{tabular}[t]{l}$v_5$\end{tabular}}}}%
    \put(0.4081531,0.00667724){\makebox(0,0)[lt]{\lineheight{1.45000005}\smash{\begin{tabular}[t]{l}$v_2$\end{tabular}}}}%
    \put(0.53783934,0.00686354){\makebox(0,0)[lt]{\lineheight{1.45000005}\smash{\begin{tabular}[t]{l}$v_2$\end{tabular}}}}%
    \put(0.69068511,0.00505099){\makebox(0,0)[lt]{\lineheight{1.45000005}\smash{\begin{tabular}[t]{l}$v_3$\end{tabular}}}}%
    \put(0,0){\includegraphics[width=\unitlength,page=3]{rho8y.pdf}}%
    \put(0.8448903,0.00603769){\makebox(0,0)[lt]{\lineheight{1.45000005}\smash{\begin{tabular}[t]{l}$v_4$\end{tabular}}}}%
    \put(0,0){\includegraphics[width=\unitlength,page=4]{rho8y.pdf}}%
    \put(0.02989079,0.0048963){\makebox(0,0)[lt]{\lineheight{1.45000005}\smash{\begin{tabular}[t]{l}$v_1$\end{tabular}}}}%
    \put(0.15815722,0.00667716){\makebox(0,0)[lt]{\lineheight{1.45000005}\smash{\begin{tabular}[t]{l}$v_1$\end{tabular}}}}%
    \put(0,0){\includegraphics[width=\unitlength,page=5]{rho8y.pdf}}%
  \end{picture}%
\endgroup%
}\in B_{4,1}(n)_{(3,2,1^2)}$$
    and
    $$z=\centre{
\begingroup%
  \makeatletter%
  \providecommand\color[2][]{%
    \errmessage{(Inkscape) Color is used for the text in Inkscape, but the package 'color.sty' is not loaded}%
    \renewcommand\color[2][]{}%
  }%
  \providecommand\transparent[1]{%
    \errmessage{(Inkscape) Transparency is used (non-zero) for the text in Inkscape, but the package 'transparent.sty' is not loaded}%
    \renewcommand\transparent[1]{}%
  }%
  \providecommand\rotatebox[2]{#2}%
  \newcommand*\fsize{\dimexpr\f@size pt\relax}%
  \newcommand*\lineheight[1]{\fontsize{\fsize}{#1\fsize}\selectfont}%
  \ifx\svgwidth\undefined%
    \setlength{\unitlength}{72.00117926bp}%
    \ifx\svgscale\undefined%
      \relax%
    \else%
      \setlength{\unitlength}{\unitlength * \real{\svgscale}}%
    \fi%
  \else%
    \setlength{\unitlength}{\svgwidth}%
  \fi%
  \global\let\svgwidth\undefined%
  \global\let\svgscale\undefined%
  \makeatother%
  \begin{picture}(1,0.51089881)%
    \lineheight{1}%
    \setlength\tabcolsep{0pt}%
    \put(0,0){\includegraphics[width=\unitlength,page=1]{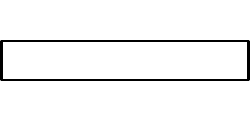}}%
    \put(0.35783743,0.24028307){\makebox(0,0)[lt]{\lineheight{1.45000005}\smash{\begin{tabular}[t]{l}$c_{(2^2,1^3)}$\end{tabular}}}}%
    \put(0,0){\includegraphics[width=\unitlength,page=2]{rho8z.pdf}}%
    \put(0.27988668,0.0048963){\makebox(0,0)[lt]{\lineheight{1.45000005}\smash{\begin{tabular}[t]{l}$v_2$\end{tabular}}}}%
    \put(0.4081531,0.00667724){\makebox(0,0)[lt]{\lineheight{1.45000005}\smash{\begin{tabular}[t]{l}$v_2$\end{tabular}}}}%
    \put(0.53783934,0.00686354){\makebox(0,0)[lt]{\lineheight{1.45000005}\smash{\begin{tabular}[t]{l}$v_3$\end{tabular}}}}%
    \put(0.69068511,0.00505099){\makebox(0,0)[lt]{\lineheight{1.45000005}\smash{\begin{tabular}[t]{l}$v_4$\end{tabular}}}}%
    \put(0,0){\includegraphics[width=\unitlength,page=3]{rho8z.pdf}}%
    \put(0.8448903,0.00603769){\makebox(0,0)[lt]{\lineheight{1.45000005}\smash{\begin{tabular}[t]{l}$v_5$\end{tabular}}}}%
    \put(0,0){\includegraphics[width=\unitlength,page=4]{rho8z.pdf}}%
    \put(0.02989079,0.0048963){\makebox(0,0)[lt]{\lineheight{1.45000005}\smash{\begin{tabular}[t]{l}$v_1$\end{tabular}}}}%
    \put(0.15815722,0.00667716){\makebox(0,0)[lt]{\lineheight{1.45000005}\smash{\begin{tabular}[t]{l}$v_1$\end{tabular}}}}%
    \put(0,0){\includegraphics[width=\unitlength,page=5]{rho8z.pdf}}%
  \end{picture}%
\endgroup%
}\neq 0\in B_{4,1}(n)_{(2^2,1^3)}.$$
    Therefore, we have $$\rho_8(x)(K_{5,4,3})=-32z\neq 0\in B_{4,1}(n)_{(2^2,1^3)}$$ and thus $\rho_8$ is injective for $n\geq 5$.
  \end{proof}

  By using Theorem \ref{propindecomp} and Proposition \ref{propd4} carefully, one can obtain the following theorem, which improves Theorem \ref{propindecomp} for $d=4$.
  \begin{theorem}\label{decomp4}
    We have an indecomposable decomposition
    $$A_4(n)=A_4 P(n)\oplus A_4 Q(n)$$
    of $\Aut(F_n)$-modules for $n\geq 7$.
  \end{theorem}
  We expect that Theorem \ref{decomp4} holds for $n\geq 3$.

\section{The $\Out(F_n)$-module structure of $A_d(n)$}
 In \cite{Mai1}, we observed that the $\Aut(F_n)$-action on $A_d(n)$ induces an action of $\Out(F_n)$ on $A_d(n)$.
 In this section, we obtain some results for $A_d(n)$ as $\Out(F_n)$-modules, which is induced by the results in Section \ref{s8}.

 Since the $\Aut(F_n)$-action on $A_d(n)$ factors through $\Out(F_n)$, any submodule of $A_d(n)$ as $\Aut(F_n)$-modules is a submodule of $A_d(n)$ as $\Out(F_n)$-modules, and vice versa.
 By Theorem \ref{rad}, we obtain the radical filtration of $A_d(n)$ as $\Out(F_n)$-modules.
 \begin{theorem}
   Let $n\geq 2d$.
   Then, the filtration of $A_d(n)$ by the number of trivalent vertices coincides with the radical filtration of $A_d(n)$ as $\Out(F_n)$-modules.
 \end{theorem}

 By Theorem \ref{propindecomp}, we obtain an indecomposable decomposition of $A_d(n)$ as $\Out(F_n)$-modules.
 \begin{theorem}
   Let $d\geq 2$.
   We have a direct decomposition $$A_d(n)=A_d P(n)\oplus A_d Q(n)$$ of $\Out(F_n)$-modules, which is indecomposable for $n\geq 2d$.
 \end{theorem}
 Theorems \ref{decomp3}, \ref{decomp4} also hold as $\Out(F_n)$-modules.
 Other results for $A_d(n)$ as $\Aut(F_n)$-modules such as Proposition \ref{soc} also hold.


\section{Indecomposable decomposition of the functor $A_d$}\label{s9}
 In this section, we obtain an indecomposable decomposition of the functor $A_d$ by using results in Section \ref{s8}.

 By Theorem \ref{decompositionofAd}, we obtain the following direct decomposition of the functor $A_d$.
 \begin{theorem}\label{directdecomp}
   We have a direct decomposition
   $$A_d=A_d P\oplus A_d Q$$
   in the functor category $\fVect^{\F^{\op}}$.
 \end{theorem}
 For $d=1$, we have $A_1 Q=0$ and the functor $A_1=A_1 P$ is simple.
 For $d=2$, we obtained this direct decomposition in Theorem 6.5 of \cite{Mai1}.
 Moreover, we proved that this direct decomposition of the functor $A_2$ is indecomposable (see Theorem 6.14 of \cite{Mai1}).

 By Theorem \ref{propindecomp}, we obtain the indecomposability of the direct decomposition of the functor $A_d$.
 \begin{proposition}\label{propindecompfunc}
  Let $d\geq 2$. The decomposition
  $$A_d=A_d P\oplus A_d Q$$
  of the functor $A_d$ is indecomposable in the functor category $\fVect^{\F^{\op}}$.
 \end{proposition}
 \begin{proof}
  Suppose that we have a decomposition
  $$A_d Q=G\oplus G'\in \fVect^{\F^{\op}}.$$
  Then we have $A_d Q(2d)=G(2d)\oplus G'(2d)$ as $\Aut(F_{2d})$-modules.
  By Theorem \ref{propindecomp}, the $\Aut(F_{2d})$-module $A_d Q(2d)$ is indecomposable. Therefore, we can assume that $G'(2d)=0$ and $A_d Q(2d)=G(2d)$.
  Since the subfunctor $A_d Q$ is generated by $Q\in A_d Q(2d)$, we have $A_d Q=G$.
  Hence, the subfunctor $A_d Q$ is also indecomposable.
  By Lemma \ref{l81}, $A_d P(2d)$ is also indecomposable. Therefore, by the similar argument, the subfunctor $A_d P$ is also indecomposable.
 \end{proof}


\appendix
\section{Presentation of the category $\A^{L}$}\label{sA}
 In this section, we construct a category $\prA$ and a full functor $F: \prA\rightarrow\A^{L}$ to study a presentation of the category $\A^{L}$, which we construct in Section \ref{ss42}.

 \subsection{The category $\prA$}\label{ssA1}
  In this section, we construct a category $\prA$, which has a generating set and some relations of the category $\A^{L}$.

  In a linear symmetric strict monoidal category $\catC$,
  let $H$ be a Hopf algebra and $L$ a Lie algebra.
  Define the \emph{adjoint action}
  $ad_H:H\otimes H\rightarrow H$ by
  $$ad_{H}= \mu^{[3]}(\id_{H^{\otimes 2}}\otimes S)(\id_{H}\otimes P_{H,H})(\Delta\otimes\id_{H}).$$
  We call a morphism $c:I\rightarrow L^{\otimes 2}$ a \emph{symmetric invariant $2$-tensor} if $c$ satisfies
  \begin{gather*}\label{symm2tensor}
    P_{L,L}c=c
  \end{gather*}
  and
  \begin{gather*}\label{invariant2tensor}
    ([\cdot,\cdot]\otimes \id_L)(\id_L\otimes c)=(\id_L\otimes [\cdot,\cdot])(c\otimes \id_L).
  \end{gather*}

  Define $\prA$ to be the category which is as a linear symmetric strict monoidal category, generated by
  \begin{itemize}
    \item a cocommutative Hopf algebra $(H,\mu,\eta,\Delta,\epsilon,S)$
    \item a Lie algebra with a symmetric invariant 2-tensor $(L,[\cdot,\cdot],c)$
    \item morphisms $i: L\rightarrow H$ and  $ad_{L}: H\otimes L \rightarrow L$
  \end{itemize}
  with the following 9 relations:
  \begin{enumerate}[($\prA$.1)]
    \item  $i\; [\cdot,\cdot]=-\mu(i\otimes i)+\mu P_{H,H}(i\otimes i),$
     \label{a1}
    \item  $\Delta i=i\otimes\eta+\eta\otimes i,$
     \label{a2}
    \item  $\epsilon i=0,$
     \label{a3}
    \item  $ad_{L}(\mu\otimes\id_L)=ad_{L}(\id_H\otimes ad_{L}),$
     \label{a4}
    \item  $ad_{L}(\eta\otimes\id_{L})=\id_{L},$
     \label{a5}
    \item  $(ad_{L}\otimes ad_{L}) (\id_{H}\otimes P_{H,L}\otimes \id_{L}) (\Delta\otimes c)=c\epsilon,$
     \label{a6}
    \item  $ad_{L}(\id_{H}\otimes[\cdot,\cdot])
    =[\cdot,\cdot](ad_{L}\otimes ad_{L}) (\id_{H}\otimes P_{H,L}\otimes\id_{L})(\Delta\otimes \id_{L^{\otimes 2}}),$
     \label{a7}
    \item  $i\; ad_{L}=ad_{H}\; i,$
     \label{a8}
    \item  $ad_{L}(i\otimes \id_{L})=-[\cdot,\cdot].$
     \label{a9}
  \end{enumerate}
  \begin{lemma}\label{lA11}
   In the category $\prA$, the following relations hold.
   \begin{enumerate}
     \item   $S i=- i.$
      \label{la1}
     \item   $ad_{H}(i\otimes i)=-i\; [\cdot,\cdot].$
      \label{la2}
   \end{enumerate}
  \end{lemma}
  \begin{proof}
   By $(\prA.\ref{a2})$ and $(\prA.\ref{a3})$ of the category $\prA$ and relations of Hopf algebras, we have
   $$i+S i =\mu(i\otimes S \eta)+\mu(\eta\otimes S i)=\mu(\id_H\otimes S)\Delta i=\eta \epsilon i=0.$$
   Thus, we have (\ref{la1}).
   By $(\prA.\ref{a8}), (\prA.\ref{a9})$, we have (\ref{la2}) as follows:
   $$ad_H(i\otimes i)=i\; ad_L(i\; \id_L)=-i\;[\cdot,\cdot].$$
  \end{proof}

  We review the definition of a Casimir Hopf algebra. Let $\catC$ be a linear symmetric strict monoidal category and $H$ be a cocommutative Hopf algebra in $\catC$. A \emph{Casimir 2-tensor} for $H$ is a morphism $c:I\rightarrow H^{\otimes 2}$ which is primitive, symmetric and invariant:
  \begin{gather}\label{casimir1}
       (\Delta \otimes \id_H)c= c_{13}+c_{23},
  \end{gather}
  \begin{gather}\label{casimir2}
       P_{H,H}c=c,
  \end{gather}
  \begin{gather}\label{casimir3}
       (ad_H\otimes ad_H)(\id_H\otimes P_{H,H}\otimes\id_H)(\Delta\otimes c)=c\epsilon,
  \end{gather}
  where $c_{13}:=(\id\otimes\eta\otimes\id)c$ and $c_{23}:=\eta\otimes c$.
  By a \emph{Casimir Hopf algebra}, we mean a cocommutative Hopf algebra $H$ equipped with a Casimir 2-tensor.

  \begin{lemma}\label{lA12}
     $(H,\mu,\eta,\Delta,\epsilon,S,\ti{c}:=(i\otimes i) c)$ is a Casimir Hopf algebra in $\prA$.
  \end{lemma}
  \begin{proof}
    Since $H$ is a cocommutative Hopf algebra in $\prA$, it suffices to check that $\ti{c}$ is a Casimir 2-tensor.
    By ($\prA$. \ref{a2}), we have (\ref{casimir1}) because
    $$(\Delta \otimes \id_H)\ti{c}= ((i\otimes\eta+\eta\otimes i)\otimes i)c =\ti{c}_{13}+\ti{c}_{23}.$$
    By the symmetricity of $c$, we have (\ref{casimir2}) because
    $$P_{H,H}\ti{c}=P_{H,H}(i\otimes i)c=(i\otimes i)P_{L,L}c=(i\otimes i)c=\ti{c}.$$
    By ($\prA$. \ref{a6}) and ($\prA$. \ref{a8}), we have (\ref{casimir3}) because
    \begin{gather*}
      \begin{split}
        &(ad_H\otimes ad_H)(\id_H\otimes P_{H,H}\otimes\id_H)(\Delta\otimes \ti{c})\\
        &=(ad_H\otimes ad_H)(\id_H\otimes P_{H,H}\otimes\id_H)(\Delta\otimes (i\otimes i))(\id_H\otimes c)\\
        &=(i\otimes i)(ad_L\otimes ad_L)(\id_H\otimes P_{H,L}\otimes\id_L)(\Delta\otimes c)\\
        &=(i\otimes i)c\epsilon\\
        &=\ti{c}\epsilon.
      \end{split}
    \end{gather*}
  \end{proof}

  The category $\A$ has a Casimir Hopf algebra $(H,c)=(1,\mu,\eta,\Delta,\epsilon, S, c)$, where $c=\centre{}.$
  Moreover, Theorem 5.11 in \cite{HM_k} implies that as a linear symmetric strict monoidal category, the category $\A$ is free on the Casimir Hopf algebra $(H,c)$.
  Therefore, we have a unique linear symmetric monoidal functor $F_{(H,\ti{c})}:\A \rightarrow \prA$.

 \subsection{Structure of the category $\A^{L}$}\label{ssA2}
  In Section \ref{ss425}, we observed that the category $\A^{L}$ has a cocommutative Hopf algebra $(H,\mu,\eta,\Delta,\epsilon,S)$ and morphisms
  $$[\cdot,\cdot]:L\otimes L\rightarrow L,\quad c_{L}:I\rightarrow L\otimes L,\quad i:L\rightarrow H,\quad ad_{L}:H\otimes L\rightarrow L.$$

  \begin{lemma}\label{lA21}
   In the category $\A^{L}$, $(L, [\cdot,\cdot], c_{L})$ is a Lie algebra with a symmetric invariant $2$-tensor.
  \end{lemma}
  \begin{proof}
   By the AS and IHX relations, it follows that $(L,[\cdot,\cdot])$ is a Lie algebra.
   Since we have
   $$P_{L,L}c_{L}=\centre{
\begingroup%
  \makeatletter%
  \providecommand\color[2][]{%
    \errmessage{(Inkscape) Color is used for the text in Inkscape, but the package 'color.sty' is not loaded}%
    \renewcommand\color[2][]{}%
  }%
  \providecommand\transparent[1]{%
    \errmessage{(Inkscape) Transparency is used (non-zero) for the text in Inkscape, but the package 'transparent.sty' is not loaded}%
    \renewcommand\transparent[1]{}%
  }%
  \providecommand\rotatebox[2]{#2}%
  \newcommand*\fsize{\dimexpr\f@size pt\relax}%
  \newcommand*\lineheight[1]{\fontsize{\fsize}{#1\fsize}\selectfont}%
  \ifx\svgwidth\undefined%
    \setlength{\unitlength}{30.75117737bp}%
    \ifx\svgscale\undefined%
      \relax%
    \else%
      \setlength{\unitlength}{\unitlength * \real{\svgscale}}%
    \fi%
  \else%
    \setlength{\unitlength}{\svgwidth}%
  \fi%
  \global\let\svgwidth\undefined%
  \global\let\svgscale\undefined%
  \makeatother%
  \begin{picture}(1,1)%
    \lineheight{1}%
    \setlength\tabcolsep{0pt}%
    \put(0,0){\includegraphics[width=\unitlength,page=1]{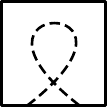}}%
  \end{picture}%
\endgroup%
}=\centre{}=c_{L}$$
   and
   $$([\cdot,\cdot]\otimes \id_L)(\id_L\otimes c_{L}) =\centre{
\begingroup%
  \makeatletter%
  \providecommand\color[2][]{%
    \errmessage{(Inkscape) Color is used for the text in Inkscape, but the package 'color.sty' is not loaded}%
    \renewcommand\color[2][]{}%
  }%
  \providecommand\transparent[1]{%
    \errmessage{(Inkscape) Transparency is used (non-zero) for the text in Inkscape, but the package 'transparent.sty' is not loaded}%
    \renewcommand\transparent[1]{}%
  }%
  \providecommand\rotatebox[2]{#2}%
  \newcommand*\fsize{\dimexpr\f@size pt\relax}%
  \newcommand*\lineheight[1]{\fontsize{\fsize}{#1\fsize}\selectfont}%
  \ifx\svgwidth\undefined%
    \setlength{\unitlength}{30.75117737bp}%
    \ifx\svgscale\undefined%
      \relax%
    \else%
      \setlength{\unitlength}{\unitlength * \real{\svgscale}}%
    \fi%
  \else%
    \setlength{\unitlength}{\svgwidth}%
  \fi%
  \global\let\svgwidth\undefined%
  \global\let\svgscale\undefined%
  \makeatother%
  \begin{picture}(1,1)%
    \lineheight{1}%
    \setlength\tabcolsep{0pt}%
    \put(0,0){\includegraphics[width=\unitlength,page=1]{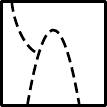}}%
  \end{picture}%
\endgroup%
}=\centre{
\begingroup%
  \makeatletter%
  \providecommand\color[2][]{%
    \errmessage{(Inkscape) Color is used for the text in Inkscape, but the package 'color.sty' is not loaded}%
    \renewcommand\color[2][]{}%
  }%
  \providecommand\transparent[1]{%
    \errmessage{(Inkscape) Transparency is used (non-zero) for the text in Inkscape, but the package 'transparent.sty' is not loaded}%
    \renewcommand\transparent[1]{}%
  }%
  \providecommand\rotatebox[2]{#2}%
  \newcommand*\fsize{\dimexpr\f@size pt\relax}%
  \newcommand*\lineheight[1]{\fontsize{\fsize}{#1\fsize}\selectfont}%
  \ifx\svgwidth\undefined%
    \setlength{\unitlength}{30.75117737bp}%
    \ifx\svgscale\undefined%
      \relax%
    \else%
      \setlength{\unitlength}{\unitlength * \real{\svgscale}}%
    \fi%
  \else%
    \setlength{\unitlength}{\svgwidth}%
  \fi%
  \global\let\svgwidth\undefined%
  \global\let\svgscale\undefined%
  \makeatother%
  \begin{picture}(1,1)%
    \lineheight{1}%
    \setlength\tabcolsep{0pt}%
    \put(0,0){\includegraphics[width=\unitlength,page=1]{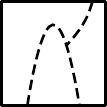}}%
  \end{picture}%
\endgroup%
}=(\id_L\otimes [\cdot,\cdot])(c_{L}\otimes \id_L),$$
   it follows that $c_{L}$ is a symmetric invariant $2$-tensor.
  \end{proof}
  \begin{remark}\label{casimirlie}
    The full subcategory of $\A^{L}$ with the free monoid generated by $L$ as the set of objects is isomorphic to the PROP $LIE^{c}$ for Casimir Lie algebras (see \cite{Hinich} for details).

    For each $m\geq 1,n\in\N$, the degree $0$ part $\A^{L}_0(L^{\otimes m}, H^{\otimes n})$ of the hom-set $\A^{L}(L^{\otimes m}, H^{\otimes n})$ has an $\Aut(F_n)$-module structure which is defined in a way similar to that of $A_d(n)$.
    For general $m,n$, the $\Aut(F_n)$-action on $\A^{L}_0(L^{\otimes m}, H^{\otimes n})$ does not factors through the outer automorphism group $\Out(F_n)$.
  \end{remark}

  \begin{proposition}\label{pA21}
   There exists a unique linear symmetric monoidal functor $F:\prA \rightarrow\A^{L}$ which maps $(L,[\cdot,\cdot],c_{L},i,ad_L)$ in $\prA$ to $(L,[\cdot,\cdot],c,i,ad_L)$ in $\A^{L}$ and which makes the following diagram commutative
    \begin{gather*}
     \xymatrix{
          \A\ar[rd]_{\text{inclu.}}\ar[rr]^{F_{(H,\ti{c})}}
          &&
           \prA\ar[ld]^{F}
          \\
           &
           \A^{L}
           &.
        }
    \end{gather*}
  \end{proposition}
  \begin{proof}
   We can check that morphisms of $\A^{L}$ satisfy the relations $(\prA.1),\cdots ,(\prA.9)$ by diagrammatic computation.
   Since $\prA$ is the linear symmetric strict monoidal category generated by $H$, $L$ and morphisms $i, ad_L$ with relations $(\prA.1),\cdots ,(\prA.9)$, we can construct a unique linear symmetric monoidal functor $F:\prA \rightarrow\A^{L}$ which maps $(H,L,c,i,ad_L)$ in $\prA$ to $(H,L,c_{L},i,ad_L)$ in $\A^{L}$.
  \end{proof}

 \subsection{The full functor $F: \prA \rightarrow \A^{L}$}\label{ssA3}
  We prove that the functor $F$ in Proposition \ref{pA21} is full.
  \begin{lemma}\label{lA31}
   A morphism in $\A^{L}$ can be written as a linear sum of the following diagrams:
   \begin{equation}\label{catDmor}
     \centre{
\begingroup%
  \makeatletter%
  \providecommand\color[2][]{%
    \errmessage{(Inkscape) Color is used for the text in Inkscape, but the package 'color.sty' is not loaded}%
    \renewcommand\color[2][]{}%
  }%
  \providecommand\transparent[1]{%
    \errmessage{(Inkscape) Transparency is used (non-zero) for the text in Inkscape, but the package 'transparent.sty' is not loaded}%
    \renewcommand\transparent[1]{}%
  }%
  \providecommand\rotatebox[2]{#2}%
  \newcommand*\fsize{\dimexpr\f@size pt\relax}%
  \newcommand*\lineheight[1]{\fontsize{\fsize}{#1\fsize}\selectfont}%
  \ifx\svgwidth\undefined%
    \setlength{\unitlength}{165.75118126bp}%
    \ifx\svgscale\undefined%
      \relax%
    \else%
      \setlength{\unitlength}{\unitlength * \real{\svgscale}}%
    \fi%
  \else%
    \setlength{\unitlength}{\svgwidth}%
  \fi%
  \global\let\svgwidth\undefined%
  \global\let\svgscale\undefined%
  \makeatother%
  \begin{picture}(1,0.88687506)%
    \lineheight{1}%
    \setlength\tabcolsep{0pt}%
    \put(0,0){\includegraphics[width=\unitlength,page=1]{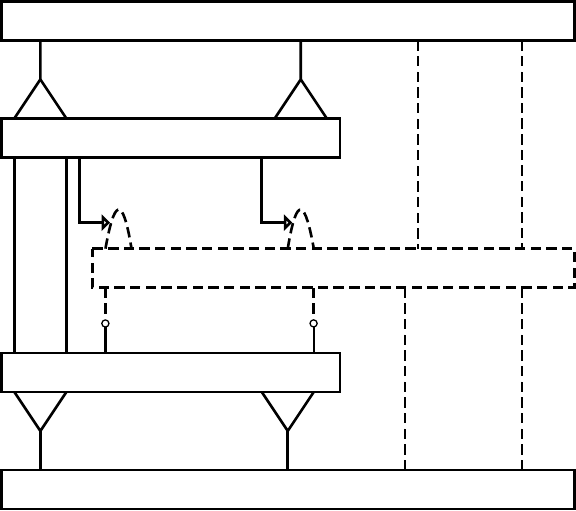}}%
    \put(0.10700486,0.8414417){\makebox(0,0)[lt]{\lineheight{1.45000005}\smash{\begin{tabular}[t]{l}\tiny{a diagram generated by $P_{H,L},P_{L,H}$}\end{tabular}}}}%
    \put(0,0){\includegraphics[width=\unitlength,page=2]{catD.pdf}}%
    \put(0.02335832,0.63895085){\makebox(0,0)[lt]{\lineheight{1.45000005}\smash{\begin{tabular}[t]{l}\tiny{a diagram generated by $P_{H,H}$}\end{tabular}}}}%
    \put(0.17998365,0.40884384){\makebox(0,0)[lt]{\lineheight{1.45000005}\smash{\begin{tabular}[t]{l}\tiny{a diagram generated by $c_{L},c^{\ast},P_{L,L},[\cdot,\cdot]$}\end{tabular}}}}%
    \put(0.02418545,0.22919437){\makebox(0,0)[lt]{\lineheight{1.45000005}\smash{\begin{tabular}[t]{l}\tiny{a diagram generated by $P_{H,H}$}\end{tabular}}}}%
    \put(0.10760468,0.0269861){\makebox(0,0)[lt]{\lineheight{1.45000005}\smash{\begin{tabular}[t]{l}\tiny{a diagram generated by $P_{H,L},P_{L,H}$}\end{tabular}}}}%
    \put(0.15050804,0.31493777){\makebox(0,0)[lt]{\lineheight{1.45000005}\smash{\begin{tabular}[t]{l}$i$\end{tabular}}}}%
    \put(0.55746821,0.31361898){\makebox(0,0)[lt]{\lineheight{1.45000005}\smash{\begin{tabular}[t]{l}$i$\end{tabular}}}}%
    \put(0,0){\includegraphics[width=\unitlength,page=3]{catD.pdf}}%
  \end{picture}%
\endgroup%
},
   \end{equation}
   where $\centre{}$ denotes $S$ or $\id_H$ and  $c^{\ast}=\centre{
\begingroup%
  \makeatletter%
  \providecommand\color[2][]{%
    \errmessage{(Inkscape) Color is used for the text in Inkscape, but the package 'color.sty' is not loaded}%
    \renewcommand\color[2][]{}%
  }%
  \providecommand\transparent[1]{%
    \errmessage{(Inkscape) Transparency is used (non-zero) for the text in Inkscape, but the package 'transparent.sty' is not loaded}%
    \renewcommand\transparent[1]{}%
  }%
  \providecommand\rotatebox[2]{#2}%
  \newcommand*\fsize{\dimexpr\f@size pt\relax}%
  \newcommand*\lineheight[1]{\fontsize{\fsize}{#1\fsize}\selectfont}%
  \ifx\svgwidth\undefined%
    \setlength{\unitlength}{30.75117737bp}%
    \ifx\svgscale\undefined%
      \relax%
    \else%
      \setlength{\unitlength}{\unitlength * \real{\svgscale}}%
    \fi%
  \else%
    \setlength{\unitlength}{\svgwidth}%
  \fi%
  \global\let\svgwidth\undefined%
  \global\let\svgscale\undefined%
  \makeatother%
  \begin{picture}(1,1)%
    \lineheight{1}%
    \setlength\tabcolsep{0pt}%
    \put(0,0){\includegraphics[width=\unitlength,page=1]{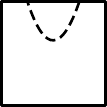}}%
  \end{picture}%
\endgroup%
}$.
  \end{lemma}
  Note that $c^{\ast}$ is not a morphism in $\A^{L}$ but just a diagram.
  \begin{proof}
   By using symmetries $P_{H,L}, P_{L,H}$, we can deform any diagram $f\in \A^{L}$ into a morphism in $\A^{L}(H^{\otimes m}\otimes L^{\otimes n}, H^{\otimes m'}\otimes L^{\otimes n'})$, so it suffices to consider a diagram $f$ in $\A^{L}(H^{\otimes m}\otimes L^{\otimes n}, H^{\otimes m'}\otimes L^{\otimes n'})$.

   We can decompose $f$ as follows:
   $f=f'\circ((P\circ\Delta^{[c_1,\cdots, c_m]})\otimes\id_{L^{\otimes n}}),$
   where $P$ is a tensor product of copies of $P_{H,H}$ and $\id_H$,
   $c_1,\cdots,c_m\geq 0$, and
   $f'$ is a diagram such that each handle has only one solid or dashed line. We can assume that handles of $U_m$ which include a dashed line are arranged right-hand side of $U_m$.

   By pulling univalent vertices that are attached to the solid lines toward the upper right-hand side of $U_m$, we can decompose $f'$ as $\centre{
\begingroup%
  \makeatletter%
  \providecommand\color[2][]{%
    \errmessage{(Inkscape) Color is used for the text in Inkscape, but the package 'color.sty' is not loaded}%
    \renewcommand\color[2][]{}%
  }%
  \providecommand\transparent[1]{%
    \errmessage{(Inkscape) Transparency is used (non-zero) for the text in Inkscape, but the package 'transparent.sty' is not loaded}%
    \renewcommand\transparent[1]{}%
  }%
  \providecommand\rotatebox[2]{#2}%
  \newcommand*\fsize{\dimexpr\f@size pt\relax}%
  \newcommand*\lineheight[1]{\fontsize{\fsize}{#1\fsize}\selectfont}%
  \ifx\svgwidth\undefined%
    \setlength{\unitlength}{94.50058901bp}%
    \ifx\svgscale\undefined%
      \relax%
    \else%
      \setlength{\unitlength}{\unitlength * \real{\svgscale}}%
    \fi%
  \else%
    \setlength{\unitlength}{\svgwidth}%
  \fi%
  \global\let\svgwidth\undefined%
  \global\let\svgscale\undefined%
  \makeatother%
  \begin{picture}(1,0.95237505)%
    \lineheight{1}%
    \setlength\tabcolsep{0pt}%
    \put(0,0){\includegraphics[width=\unitlength,page=1]{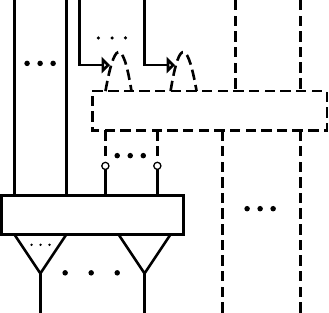}}%
    \put(0.19694119,0.26024667){\makebox(0,0)[lt]{\lineheight{1.45000005}\smash{\begin{tabular}[t]{l}$P_{H,H}$\end{tabular}}}}%
    \put(0.2639803,0.42937627){\makebox(0,0)[lt]{\lineheight{1.45000005}\smash{\begin{tabular}[t]{l}$i$\end{tabular}}}}%
    \put(0.49963907,0.432581){\makebox(0,0)[lt]{\lineheight{1.45000005}\smash{\begin{tabular}[t]{l}$i$\end{tabular}}}}%
    \put(0,0){\includegraphics[width=\unitlength,page=2]{catDs.pdf}}%
    \put(0.31696796,0.59568494){\makebox(0,0)[lt]{\lineheight{1.45000005}\smash{\begin{tabular}[t]{l}\tiny{uni-trivalent graph}\end{tabular}}}}%
  \end{picture}%
\endgroup%
}$ (see Lemma 5.16 \cite{HM_k}).

   Furthermore, any uni-trivalent graph can be obtained from morphisms $c_{L}$, $P_{L,L}$, $[\cdot,\cdot]$, $\id_L\in\A^{L}$ and $c^{\ast}$ by the tensor product and the composition, so the proof is complete.
  \end{proof}

  \begin{proposition}\label{pA31}
   The linear symmetric monoidal functor $F:\prA \rightarrow \A^{L}$ in Proposition \ref{pA21} is full.
  \end{proposition}
  \begin{proof}
   It suffices to show that morphisms of $\A^{L}$ are generated by $\mu,\eta,\Delta,\epsilon,S$, $[\cdot,\cdot],c_{L}$, $i$, $ad_{L}$ and symmetries.
   By Lemma \ref{lA31}, we need to prove that we can eliminate $c^{\ast}$ from the diagram \eqref{catDmor} by using the above morphisms in $\A^{L}$.

   By the definition of the category $\A^{L}$, for any $c^{\ast}$ in \eqref{catDmor}, if exists, either of the endpoints of $c^{\ast}$ is finally attached to one of the lower dashed lines.
   Therefore, there is $c_{L}$ between $c^{\ast}$ and the lower dashed line.
   If there are more than one such $c_{L}$, then we choose one such that there are the least trivalent vertices between $c^{\ast}$ and itself.
   By the AS relation, we have only to consider the case where the neighborhood of the $c_{L}$ and the $c^{\ast}$ is either $$\centre{
\begingroup%
  \makeatletter%
  \providecommand\color[2][]{%
    \errmessage{(Inkscape) Color is used for the text in Inkscape, but the package 'color.sty' is not loaded}%
    \renewcommand\color[2][]{}%
  }%
  \providecommand\transparent[1]{%
    \errmessage{(Inkscape) Transparency is used (non-zero) for the text in Inkscape, but the package 'transparent.sty' is not loaded}%
    \renewcommand\transparent[1]{}%
  }%
  \providecommand\rotatebox[2]{#2}%
  \newcommand*\fsize{\dimexpr\f@size pt\relax}%
  \newcommand*\lineheight[1]{\fontsize{\fsize}{#1\fsize}\selectfont}%
  \ifx\svgwidth\undefined%
    \setlength{\unitlength}{68.94303061bp}%
    \ifx\svgscale\undefined%
      \relax%
    \else%
      \setlength{\unitlength}{\unitlength * \real{\svgscale}}%
    \fi%
  \else%
    \setlength{\unitlength}{\svgwidth}%
  \fi%
  \global\let\svgwidth\undefined%
  \global\let\svgscale\undefined%
  \makeatother%
  \begin{picture}(1,0.7108776)%
    \lineheight{1}%
    \setlength\tabcolsep{0pt}%
    \put(0,0){\includegraphics[width=\unitlength,page=1]{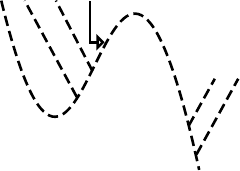}}%
    \put(0.45787489,0.51533796){\makebox(0,0)[lt]{\lineheight{1.45000005}\smash{\begin{tabular}[t]{l}\tiny$ad_L$\end{tabular}}}}%
    \put(0,0){\includegraphics[width=\unitlength,page=2]{cupc.pdf}}%
  \end{picture}%
\endgroup%
}=\centre{
\begingroup%
  \makeatletter%
  \providecommand\color[2][]{%
    \errmessage{(Inkscape) Color is used for the text in Inkscape, but the package 'color.sty' is not loaded}%
    \renewcommand\color[2][]{}%
  }%
  \providecommand\transparent[1]{%
    \errmessage{(Inkscape) Transparency is used (non-zero) for the text in Inkscape, but the package 'transparent.sty' is not loaded}%
    \renewcommand\transparent[1]{}%
  }%
  \providecommand\rotatebox[2]{#2}%
  \newcommand*\fsize{\dimexpr\f@size pt\relax}%
  \newcommand*\lineheight[1]{\fontsize{\fsize}{#1\fsize}\selectfont}%
  \ifx\svgwidth\undefined%
    \setlength{\unitlength}{47.14634764bp}%
    \ifx\svgscale\undefined%
      \relax%
    \else%
      \setlength{\unitlength}{\unitlength * \real{\svgscale}}%
    \fi%
  \else%
    \setlength{\unitlength}{\svgwidth}%
  \fi%
  \global\let\svgwidth\undefined%
  \global\let\svgscale\undefined%
  \makeatother%
  \begin{picture}(1,1.04012472)%
    \lineheight{1}%
    \setlength\tabcolsep{0pt}%
    \put(0,0){\includegraphics[width=\unitlength,page=1]{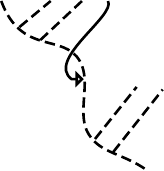}}%
    \put(0.18712362,0.44530336){\makebox(0,0)[lt]{\lineheight{1.45000005}\smash{\begin{tabular}[t]{l}\tiny$ad_L$\end{tabular}}}}%
    \put(0,0){\includegraphics[width=\unitlength,page=2]{cupceli.pdf}}%
    \put(0.61380049,0.84083934){\makebox(0,0)[lt]{\lineheight{1.45000005}\smash{\begin{tabular}[t]{l}\tiny$S$\end{tabular}}}}%
  \end{picture}%
\endgroup%
} \quad \text{or}\quad \centre{
\begingroup%
  \makeatletter%
  \providecommand\color[2][]{%
    \errmessage{(Inkscape) Color is used for the text in Inkscape, but the package 'color.sty' is not loaded}%
    \renewcommand\color[2][]{}%
  }%
  \providecommand\transparent[1]{%
    \errmessage{(Inkscape) Transparency is used (non-zero) for the text in Inkscape, but the package 'transparent.sty' is not loaded}%
    \renewcommand\transparent[1]{}%
  }%
  \providecommand\rotatebox[2]{#2}%
  \newcommand*\fsize{\dimexpr\f@size pt\relax}%
  \newcommand*\lineheight[1]{\fontsize{\fsize}{#1\fsize}\selectfont}%
  \ifx\svgwidth\undefined%
    \setlength{\unitlength}{71.702628bp}%
    \ifx\svgscale\undefined%
      \relax%
    \else%
      \setlength{\unitlength}{\unitlength * \real{\svgscale}}%
    \fi%
  \else%
    \setlength{\unitlength}{\svgwidth}%
  \fi%
  \global\let\svgwidth\undefined%
  \global\let\svgscale\undefined%
  \makeatother%
  \begin{picture}(1,0.65225591)%
    \lineheight{1}%
    \setlength\tabcolsep{0pt}%
    \put(0,0){\includegraphics[width=\unitlength,page=1]{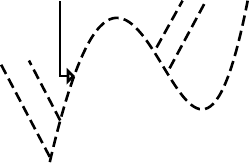}}%
    \put(0.32128469,0.32821074){\makebox(0,0)[lt]{\lineheight{1.45000005}\smash{\begin{tabular}[t]{l}\tiny$ad_L$\end{tabular}}}}%
    \put(0,0){\includegraphics[width=\unitlength,page=2]{ccup.pdf}}%
  \end{picture}%
\endgroup%
}=\centre{
\begingroup%
  \makeatletter%
  \providecommand\color[2][]{%
    \errmessage{(Inkscape) Color is used for the text in Inkscape, but the package 'color.sty' is not loaded}%
    \renewcommand\color[2][]{}%
  }%
  \providecommand\transparent[1]{%
    \errmessage{(Inkscape) Transparency is used (non-zero) for the text in Inkscape, but the package 'transparent.sty' is not loaded}%
    \renewcommand\transparent[1]{}%
  }%
  \providecommand\rotatebox[2]{#2}%
  \newcommand*\fsize{\dimexpr\f@size pt\relax}%
  \newcommand*\lineheight[1]{\fontsize{\fsize}{#1\fsize}\selectfont}%
  \ifx\svgwidth\undefined%
    \setlength{\unitlength}{61.50000019bp}%
    \ifx\svgscale\undefined%
      \relax%
    \else%
      \setlength{\unitlength}{\unitlength * \real{\svgscale}}%
    \fi%
  \else%
    \setlength{\unitlength}{\svgwidth}%
  \fi%
  \global\let\svgwidth\undefined%
  \global\let\svgscale\undefined%
  \makeatother%
  \begin{picture}(1,0.76475806)%
    \lineheight{1}%
    \setlength\tabcolsep{0pt}%
    \put(0,0){\includegraphics[width=\unitlength,page=1]{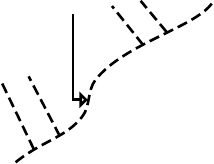}}%
    \put(0.43418841,0.27625652){\makebox(0,0)[lt]{\lineheight{1.45000005}\smash{\begin{tabular}[t]{l}\tiny$ad_L$\end{tabular}}}}%
    \put(0,0){\includegraphics[width=\unitlength,page=2]{ccupeli.pdf}}%
  \end{picture}%
\endgroup%
}.$$
   Hence, we can eliminate $c^{\ast}$ from \eqref{catDmor} and the proof is complete.
  \end{proof}



\end{document}